\documentclass[reqno]{amsart}
\usepackage{amssymb, latexsym}
\usepackage{amsmath}
\usepackage{amsthm}
\usepackage{graphicx}
\usepackage{titletoc}
\usepackage{pdfsync}
\usepackage{cancel}
\usepackage{color, xcolor}
\usepackage{framed}

\usepackage{graphicx}
\usepackage{epstopdf}
\usepackage{latexsym}
\usepackage{pdfsync}
\usepackage{cancel}

\usepackage{multirow}
\usepackage[font=bf, aboveskip=15pt]{caption}
\usepackage[toc, page]{appendix}
\usepackage[colorlinks=true]{hyperref}

\hypersetup{dvipdfmx,
       colorlinks,
       linkcolor=blue,
       filecolor=red,
       urlcolor=blue,
       citecolor=red,
       pdftitle={Concentration},
       pdfauthor={ Jun YANG},
       pdfsubject={   PDEs         },
       pdfkeywords={   Conjectures},
       pdfproducer={ps2pdf}}
\makeatletter
\def\@currentlabel{2.1}\label{e:dispaa}
\def\@currentlabel{2.21}\label{e:dispau}
\def\@currentlabel{2.22}\label{e:dispav}
\def\@currentlabel{2.23}\label{e:dispaw}
\def\@currentlabel{2.24}\label{e:dispax}
\def\theequation{\thesection.\@arabic\c@equation}

\makeatother
\def\theequation{\thesection.\@arabic\c@equation}
\makeatother

\renewcommand{\theequation}{\thesection.\arabic{equation}}
\newtheorem{theorem}{Theorem}[section]

\newtheorem{proposition}[theorem]{Proposition}
\newtheorem{remark}[theorem]{Remark}
\newtheorem{lemma}[theorem]{Lemma}

\allowdisplaybreaks

\UseRawInputEncoding

\begin{document}
\title[Ambrosetti-Malchiodi-Ni conjecture: clustering concentration layers]{On Ambrosetti-Malchiodi-Ni conjecture on two-dimensional smooth bounded domains:
 clustering concentration layers}

\author{Suting Wei$^\ddag$}
\thanks{$\ddag\, $ Department of Mathematics, South China
Agricultural University, Guangzhou, 510642, P. R. China.
Email: stwei@scau.edu.cn}

\author{Jun  Yang$^\S$}
\thanks{$\S\, $
School of Mathematics and Information Science,
Guangzhou University,
Guangzhou 510006, P. R. China.
Email: jyang2019@gzhu.edu.cn}

\thanks{Corresponding author: Jun Yang, jyang2019@gzhu.edu.cn}

\begin{abstract}
We consider the problem
$$
 \varepsilon^2 {\mathrm {div}}\big( \nabla_{{\mathfrak a}(y)} u\big)- V(y)u+u^p\, =\, 0, \quad  u>0  \quad\mbox{in }\Omega,
 \qquad
 \nabla_{{\mathfrak a}(y)} u\cdot \nu\, =\, 0\quad\mbox{on }   \partial \Omega,
$$
where $\Omega$ is a bounded domain in $\mathbb R^2$ with smooth boundary, the exponent $p$ is greater than $1$,  $\varepsilon>0$ is a small parameter,
$V$ is a uniformly positive smooth potential on $\bar{\Omega}$, and $\nu$ denotes the outward normal of $\partial \Omega$.
For two positive smooth functions ${\mathfrak a}_1(y), {\mathfrak a}_2(y)$ on $\bar\Omega$, the operator $\nabla_{{\mathfrak a}(y)}$ is given by
$$
\nabla_{{\mathfrak a}(y)} u=\Bigg({\mathfrak a}_1(y)\frac{\partial u}{\partial y_1}, \,  {\mathfrak a}_2(y)\frac{\partial u}{\partial y_2}\Bigg).
$$

(1). Let $\Gamma\subset{\bar\Omega}$ be a smooth curve intersecting orthogonally with
 $\partial \Omega$ at exactly two points and dividing $\Omega$ into two parts.
 Moreover,  $\Gamma$ is  a {\it non-degenerate geodesic} embedded in the Riemannian manifold ${\mathbb R}^2$
with metric $V^{2\sigma}(y)\big[{\mathfrak a}_2(y){\mathrm d}y_1^2+{\mathfrak a}_1(y){\mathrm d}y_2^2\big]$, where $\sigma=\frac {p+1}{p-1}-\frac 12$.
By assuming some additional constraints on the functions ${\mathfrak a}(y)$, $V(y)$ and the curves $\Gamma$, $\partial\Omega$,
we prove that there exists a sequence of $\varepsilon$ such that the problem has solutions $u_\varepsilon$ with clustering concentration layers directed along $\Gamma$, exponentially small in $\varepsilon$ at any positive distance from it.

(2). If ${\tilde\Gamma}$ is  a simple closed smooth curve in $\Omega$ (not touching the boundary $\partial\Omega$),
which is also a {\it non-degenerate geodesic} embedded in the Riemannian manifold ${\mathbb R}^2$
with metric $V^{2\sigma}(y)\big[{\mathfrak a}_2(y){\mathrm d}y_1^2+{\mathfrak a}_1(y){\mathrm d}y_2^2\big]$, then
a similar result of concentrated solutions is still true.



 \end{abstract}

\date{}\maketitle

 \section{Introduction}\label{section1}

 We consider the following problem for the existence of solutions with concentration phenomena
 \begin{align}
 \label{originalproblem}
 \varepsilon^2 {\mathrm {div}}\big( \nabla_{{\mathfrak a}(y)} u\big)- V(y)u+u^p\, =\, 0, \quad u>0\quad\mbox{in } \Omega,
\qquad\
\nabla_{{\mathfrak a}(y)} u\cdot \nu\, =\, 0\quad\mbox{on } \partial \Omega,
 \end{align}
where $\Omega$ is a bounded domain in $\mathbb R^{\mathbf d}$ with smooth boundary, $\varepsilon>0$ is a small parameter, $V$ is a uniformly positive, smooth
potential on $\bar{\Omega}$, and $\nu$ denotes the outward normal of $\partial \Omega$, the exponent $p>1$.
For ${\mathfrak a}(y)=\big({\mathfrak a}_1(y), \cdots, {\mathfrak a}_{\mathbf{d}}(y)\big)$,
the operator is defined in the form
$$
\nabla_{{\mathfrak a}(y)} u=\big({\mathfrak a}_1(y)u_{y_1}, \cdots, {\mathfrak a}_{\mathbf{d}}(y)u_{y_{\mathbf{d}}}\big),
$$
where
${\mathfrak a}_1(y), \cdots, {\mathfrak a}_{\mathbf{d}}(y)$ are positive smooth functions on $\bar\Omega$.

\medskip
\medskip
\noindent $\bullet$
In the case ${\mathfrak a}\equiv 1$ and $V\equiv 1$, problem \eqref{originalproblem} takes the form
\begin{equation}\label{original equationconstant}
 \varepsilon^2\Delta u-u+u^p\, =\, 0, \quad u>0 \quad{\rm in}\quad \Omega,
 \ \quad \nabla  u\cdot \nu\, =\, 0\quad \;\;{\rm on} \;\; \partial \Omega,
 \end{equation}
which is known as the stationary equation of Koeller-Segel system in chemotaxis \cite{LinNiTak1988}.
It can also be viewed as a limiting stationary equation of Gierer-Meinhardt system in biological pattern formation \cite{GierMein1972}.

 In the pioneering papers \cite{LinNiTak1988, NiTaka1991, NiTaka1993}, under the condition that $p$ is a subcritical Sobolev exponent, C.-S. Lin, W.-M. Ni and I. Takagi established,  for $\varepsilon$ sufficiently small, the existence of a least-energy
 solution $U_\varepsilon$ of (\ref{original equationconstant}) with only one local maximum point locating at the most curved point of $\partial\Omega$.
Such a solution is called a spike-layer, which has concentration phenomena at interior or boundary points.
For the existence of interior spikes, we refer the reader to the articles
 \cite{BatesFusc, DancerYan19992, dPFW1, GrosPistWei2000, GuiWei1999, Wei1998} and the references therein.
 On the other hand, boundary spikes related to the mean curvature
 of $\partial \Omega$ can be found in
 \cite{BatesDanShi, DancerYan19991, dPFW, GuiWeiWint2000, Li1998, Wei1997, WeiWinter1998},
 and the references therein. The coexistence of interior and boundary spikes was due to C. Gui and J. Wei \cite{GuiWei2000}.
 A good review of the subject up to 2004 can be found in \cite{Ni2}.

\medskip
There is a conjecture on higher-dimensional concentration by W.-M. Ni \cite{Ni1} (see also \cite{Ni2}):
\\
\noindent {\bf Conjecture 1.} {\em
For any integer $1\leq \mathbf{k}\leq {\mathbf d}-1$, there exists $p_\mathbf{k}\in (1, \infty)$ such that for all $1<p<p_\mathbf{k}$,
problem (\ref{original equationconstant}) has a solution $U_\varepsilon$ which concentrates on a $\mathbf{k}$-dimensional subset of $\bar{\Omega}$,
provided that $\varepsilon$ is sufficiently small.
\qed
}
\\
We here mention some results for the existence of higher dimensional boundary concentration phenomena in the papers \cite{MahMal2007, Malchi2004, Malchi2005, MalMont2002, MalMont2004}. The papers \cite{JWeiYang2007, JWeiYang2008} set up the existence of concentration on an interior line,
which connects the boundary $\partial\Omega$ and is non-degenerate in the sense of variation of arc-length. There are also some other results \cite{AoMussoWeiJDE, AoMussoWeiSIAM, dAprile2011, dAprilePist20101, dAprilePist20102} to exhibit concentration phenomena on interior line segments connecting the boundary of $\Omega$.
For higher dimensional extension, the reader can refer to \cite{AoYang,  DancerYan2007, delPKowParWei, Guoyang, LiPengYan2007}.
The reader can also refer to the survey paper by J. Wei \cite{Wei2008}.

\medskip
\medskip
\noindent $\bullet$
In the case ${\mathfrak a}\equiv 1$ and $V{\not\equiv}$ constant, problem (\ref{originalproblem})  on the whole space corresponds to the following problem
 \begin{align}
 \label{originalproblem-space}
 \varepsilon^2 \Delta u- V(y)u+u^p\, =\, 0, \quad u>0\qquad \mbox{in }  {\mathbb R}^{\mathbf d},
 \end{align}
where $\varepsilon>0$ is a small parameter, the exponent $p>1$,
and $V$ is a smooth function with
 \begin{equation*}
\displaystyle {\inf_{y\in {\mathbb R}^{\mathbf d}}} V(y)>0.
\end{equation*}
Started by \cite{FloerWein1986}, solutions exhibiting concentration around one or more points of space under various assumptions on the
potential and the nonlinearity were given by many authors
\cite{{AmbBadCing1997}, AmbMalSecc2001, CingLazzo, delPFelmer1996, delPFelmer1997, delPFelmer1998, delPFelmer2002, FelmTorr2002, wangx}.
On the other hand, radially symmetric solutions with concentration on sphere of radius $r_0$ can be constructed in \cite{AMNconjecture1},
whenever $V$ is a radial function and $r_0 >0$ is a non-degenerate critical point of the function
\begin{equation*}
M(r)\, =\, r^{{\mathbf d}-1}V^{\sigma}(r), \quad {\rm where}\quad \sigma\, =\, \frac {p+1}{p-1}-\frac 12.
\end{equation*}
Based on heuristic arguments, in 2003,  A. Ambrosetti, A. Malchiodi and W.-M. Ni  raised the following conjecture (p.465, \cite{AMNconjecture1}):
\\
\noindent {\bf Conjecture 2.} {\em
Let ${\mathcal K}$ be a  non-degenerate $\mathbf{k}$-dimensional stationary manifold  of the following functional
\[
\int_{\mathcal K} V^{\frac{p+1}{p-1} -\frac{1}{2} ({\mathbf d}-\mathbf{k})},
\]
where the exponent $p$ is subcritical w.r.t. ${\mathbb R}^{\mathbf{d}-\mathbf{k}}$.
Then there exists a solution to (\ref{originalproblem-space}) concentrating near ${\mathcal K}$, at least along a subsequence $\varepsilon_j \to 0$.
\qed
}
\\
The validity of this conjecture was confirmed in two-dimensional general case with concentration on stationary and non-degenerate curves for $\varepsilon$ satisfying a gap condition due to the resonance character of the problem, see \cite{delPKowWei2007}.
More results can be found in \cite{bartschpeng-1, bartschpeng-2, WangWeiYang, MahMalMont2009}.

\medskip
Let us go back to the problem on smooth bounded domains with homogeneous Neumann boundary condition,
\begin{align}
 \label{originalproblem-domain}
 \varepsilon^2  \Delta u - V(y)u+u^p\, =\, 0, \quad u>0\quad\mbox{in } \Omega,
\qquad\
\nabla u\cdot \nu\, =\, 0\quad\mbox{on } \partial \Omega,
 \end{align}
where $\Omega$ is a bounded domain in $\mathbb R^{\mathbf d}$ with smooth boundary, $\varepsilon>0$ is a small parameter, $V$ is a uniformly positive, smooth
 potential on $\bar{\Omega}$, and $\nu$ denotes the outward normal of $\partial \Omega$  and the exponent $p>1$.
If $\Omega$ is a unit ball $B_1(0)$, the existence of radial solutions to (\ref{originalproblem-domain}) was shown in \cite{AMNconjecture2},
where  the concentration lying on spheres in $B_1(0)$ will approach the boundary with speed $O(\varepsilon|\log\varepsilon|)$ as $\varepsilon\rightarrow 0$.
For the interior concentration phenomena connecting the boundary $\partial\Omega$,
there is a conjecture by A. Ambrosetti, A. Malchiodi and W.-M. Ni in 2004 (p. 327, \cite{AMNconjecture2}), which can be stated as:
\\
\noindent {\bf Conjecture 3.} {\em
Let ${\mathcal K}$ be a $\mathbf{k}$-dimensional manifold intersecting $\partial\Omega$ perpendicularly, which is also stationary and non-degenerate with respect to the following functional
\[
\int_{\mathcal K}  V^{\frac{p+1}{p-1} -\frac{1}{2} ({\mathbf d}-\mathbf{k})},
\]
where the exponent $p$ is subcritical w.r.t. ${\mathbb R}^{\mathbf{d}-\mathbf{k}}$.
Then there exists a solution to (\ref{originalproblem-domain}) concentrating near ${\mathcal K}$, at least along a subsequence $\varepsilon_j \to 0$.
\qed
}
\\
In \cite{weixuyang}, S. Wei, B. Xu and J. Yang considered (\ref{originalproblem-domain}) and provided an affirmative answer to the mentioned conjecture only in the case: ${\mathbf d}=2$ and $\mathbf{k}=1$ for the existence of solutions with single concentration layer connecting the boundary $\partial\Omega$.

\bigskip
\noindent $\bullet$
In present paper, we will consider a little bit more general case, i.e., problem  \eqref{originalproblem}, and investigate first the existence of clustering phenomena of multiple concentration layers, which connect the boundary $\partial\Omega$.
As the descriptions in ¡°Concluding remarks¡± of \cite{delPKowWei2010}, the difficulties arise from the multiple resonance phenomena, see also Remark \ref{remark14}.
On the other hand, the much more complicated situation is the balance between neighbouring layers plus the interaction among the interior concentration layers, the boundary of $\Omega$ and the competition between  $\mathcal{{\mathfrak a}}$ and $V$.
It is natural to introduce the new ingredient (\ref{boundaryadmissibility}) to handle this delicate thing. Whence, in the present paper we focus on the two dimensional case of problem \eqref{originalproblem}. We make the following assumptions:
\begin{itemize}
  \item[ \bf(A1).]
{\em
 Let $\Omega$ be a smooth and bounded domain in $\Bbb R^2, $\;
 $\Gamma$ be a curve intersecting $\partial \Omega$ at exactly two points, saying  $P_1, P_2$,  and, at these points $\Gamma\bot \partial \Omega$.
 In the small neighborhoods of $P_1, P_2$, the boundary $\partial \Omega$ are two curves, say $\mathcal{C}_1$ and $\mathcal{C}_2$, which can be represented
 by the graphs of two functions respectively:
 $$
 {y}_2=\varphi_1({y}_1)\quad\mbox{with } (0, \varphi_1(0))=P_1,
 $$
 $$
 {y}_2=\varphi_2({y}_1)\quad\mbox{with } (0, \varphi_2(0))=P_2.
 $$
Without loss of generality, we can assume that $\Gamma$ has length $1$, and then denote $k_1$, $k_2$ the signed curvatures of ${\mathcal C}_1$ and ${\mathcal C}_2$ respectively, also $k$ the curvature of $\Gamma$.
}
  \item[ \bf(A2).]
{\em
 $\Gamma$ separates the domain $\Omega$ into two disjoint components $\Omega_1$ and $\Omega_2$.
}
  \item[ \bf(A3).]
{\em
The functions ${\mathfrak a}_1(y)$ and ${\mathfrak a}_2(y)$ satisfy the condition
\begin{align}\label{a1=a2}
{\mathfrak a}_1={\mathfrak a}_2\quad \mbox{at the points } P_1 \mbox{ and } P_2.
\end{align}
The curve $\Gamma$ is a non-degenerate geodesic embedded in the Riemannian manifold ${\mathbb R}^2$ with the following metric
\begin{equation}\label{sigma}
 V^{2\sigma}(y) \big[ {\mathfrak a}_2(y){\mathrm d}y_1^2 +{\mathfrak a}_1(y) {\mathrm d}y_2^2 \big],
 \quad\mbox{with}\ \sigma\, =\, \frac{p+1}{p-1}-\frac{1}{2}.
\end{equation}
This will be clarified in the next section (see (\ref{stationary}) and (\ref{nondegeneracy})).
\qed
}
\end{itemize}

 \medskip
Let $w$ denote the unique positive solution of the problem
 \begin{equation}
 \label{wsolution}
 w{''}-w+w^{p}\, =\, 0, \quad w>0\quad\mbox{in }{\mathbb R}, \quad w{'}(0)\, =\, 0, \quad w(\pm \infty)\, =\, 0.
 \end{equation}
We can formulate the first result.

\begin{theorem}
\label{theorem 1.1}
Let ${\mathbf d}=2$, $p>1$ and recall the assumptions in {\bf (A1)}-{\bf (A3)} as well as the modified Fermi coordinates $(t, \theta)$ in \eqref{Fermicoordinates-modified}.
Moreover, we assume that
\begin{equation}\label{taupositivity}
\begin{split}
\tau_2(\theta)
\, \equiv\, \mathcal{H}_2'(\theta)-\mathcal{H}_3(\theta)
+2\frac{|\beta'(\theta)|^2}{\beta^2(\theta)}\mathcal{H}_1(\theta)
-\frac{\beta''(\theta)}{\beta(\theta)}\mathcal{H}_1(\theta)
-\frac{\beta'(\theta)}{\beta(\theta)}\mathcal{H}_1'(\theta)
\, >\, 0,
\end{split}
\end{equation}
and also the validity of the admissibility conditions
\begin{align}\label{boundaryadmissibility}
\frac{{\mathfrak b}_2}{{\mathfrak b}_1}+\frac{\beta'(0)}{\beta(0)}=0,
\qquad
\frac{{\mathfrak b}_7}{{\mathfrak b}_6}+\frac{\beta' (1)}{\beta(1)}=0,
\end{align}
where the function $\beta>0$ is defined in (\ref{alpha-beta}),
the functions $\mathcal{H}_1$, $\mathcal{H}_2$ and $\mathcal{H}_3$ are given in \eqref{mathcalH1}-\eqref{mathcalH3},
and the constants ${\mathfrak b}_1$, ${\mathfrak b}_2$, ${\mathfrak b}_6$, ${\mathfrak b}_7$ are given in \eqref{b1b2} and \eqref{b6b7}.
Then for each $N$, there exists a sequence of $\varepsilon$, say $\{\varepsilon_l \}$, such that problem \eqref{originalproblem} has a positive solution $u_{\varepsilon_l}$ with exactly N concentration layers at mutual distances $O(\varepsilon_l|\ln\varepsilon_l|)$.
 In addition, the center of mass for N concentration layers collapses to $\Gamma$ at speed $O(\varepsilon_l^{1+\mu})$ for some small positive constant $\mu$.
More precisely, $u_{\varepsilon_l}$ has the form
 \begin{equation}
 \label{taketheform}
 u_{\varepsilon_l}( y_1, y_2)\, \approx\, \sum_{j=1}^N V(0, \theta)^{\frac {1}{p-1}}
 \, w\left(
 \sqrt{\frac{V(0, \theta)  \Big({\mathfrak a}_1(0, \theta)|n_1(\theta)|^2+{\mathfrak a}_2(0, \theta)|n_2(\theta)|^2\Big)}
 {|{\mathfrak a}_1(0, \theta)|^2+|{\mathfrak a}_2(0, \theta)|^2} }
 \frac {\, t-\varepsilon_l\, f_j(\theta)\, } {\varepsilon_l}
 \right),
 \end{equation}
 where $n(\theta)=(n_1(\theta), n_2(\theta))$ is the unit normal to $\Gamma$.
The functions $f_j$'s satisfy
\begin{align}
\|f_j\|_{\infty}\leq\, C|\ln\varepsilon_l|^2,
  \qquad
  \sum_{j=1}^Nf_j=O(\varepsilon_l^\mu)\quad\mbox{with } \mu>0,
\label{fproperties1}
\end{align}
\begin{align}
\min\limits_{1\leq j\leq N-1}(f_{j+1}-f_j)
\approx 2 |\ln\varepsilon_l|\,\sqrt{
\frac {|{\mathfrak a}_1(0, \theta)|^2+|{\mathfrak a}_2(0, \theta)|^2}
{V(0, \theta)  \Big({\mathfrak a}_1(0, \theta)|n_1(\theta)|^2+{\mathfrak a}_2(0, \theta)|n_2(\theta)|^2\Big)}
}
\, ,
\label{fproperties2}
\end{align}
and solve the Jacobi-Toda system, for $j=1, \cdots, N$,
\begin{align}
&\varepsilon_l^2\, \varsigma
\Big[\,
\mathcal{H}_1f''_j
+\mathcal{H}_1'f'_j
+ \big(\mathcal{H}_2'-\mathcal{H}_3\big)f_j
\, \Big]
 -e^{-\beta(f_j-f_{j-1})}+e^{-\beta(f_{j+1}-f_j)}\approx0\quad\mbox{in } {(0, 1)},
\label{fJacobiToda}
\end{align}
with boundary conditions
\begin{align}
{\mathfrak b}_6f_j'(1)-{\mathfrak b}_7f_j(1)\approx 0,
\qquad
{\mathfrak b}_1f_j'(0)-{\mathfrak b}_2f_j(0)\approx 0,
\qquad\forall\, j=1, \cdots, N,
\label{fboundary}
\end{align}
where the function $\varsigma>0$ is defined in (\ref{varsigma}) with the conventions
$
f_0=-\infty, \ f_{N+1}=\infty.
$
\qed
\end{theorem}

\medskip
Here is the second result for the existence of interior clustering concentration layers, which do not touch the boundary $\partial\Omega$.

\begin{theorem}
\label{theorem 1.2}
Let ${\mathbf d}=2$ and $p>1$.
Suppose that $({\hat t}, {\hat\theta})$ are the modified Fermi coordinates given in \eqref{Fermicoordinates-modified-tilde}.
Assume that ${\hat\Gamma}$ is a simple closed smooth curve with unit length in $\Omega$, which is also a non-degenerate geodesic embedded in the Riemannian manifold ${\mathbb R}^2$ with the following metric
\begin{equation}\label{sigma22222222}
 V^{2\sigma}(y) \big[ {\mathfrak a}_2(y){\mathrm d}y_1^2 +{\mathfrak a}_1(y) {\mathrm d}y_2^2 \big],
 \quad\mbox{with}\ \sigma\, =\, \frac{p+1}{p-1}-\frac{1}{2},
\end{equation}
see (\ref{stationary4}) and (\ref{nondegeneracy4}).
Moreover, we assume that
\begin{equation}\label{taupositivity2}
\begin{split}
{\hat\tau}_2({\hat\theta})
\, \equiv\, \widehat{\mathcal H}_2'({\hat\theta})-\widehat{\mathcal H}_3({\hat\theta})
+2\frac{|\beta'({\hat\theta})|^2}{\beta^2({\hat\theta})}\widehat{\mathcal H}_1({\hat\theta})
-\frac{\beta''({\hat\theta})}{\beta({\hat\theta})}\widehat{\mathcal H}_1({\hat\theta})
-\frac{\beta'({\hat\theta})}{\beta({\hat\theta})}\widehat{\mathcal H}_1'({\hat\theta})
\, >\, 0,
\end{split}
\end{equation}
where the functions $\beta>0$, $\widehat{\mathcal H}_1$, $\widehat{\mathcal H}_2$ and $\widehat{\mathcal H}_3$
are given in \eqref{alpha-beta}, \eqref{mathcalH1tilde}-\eqref{mathcalH3tilde}.
Then for each $N$, there exists a sequence of $\varepsilon$, say $\{{\hat\varepsilon}_l \}$, such that problem \eqref{originalproblem} has a positive solution $u_{{\hat\varepsilon}_l}$ with exactly N concentration layers at mutual distances $O({\hat\varepsilon}_l|\ln{\hat\varepsilon}_l|)$.
 In addition, the center of mass for N concentration layers collapses to ${\hat\Gamma}$ at speed $O({\hat\varepsilon}_l^{1+\mu})$ for some small positive constant $\mu$.
More precisely, $u_{{\hat\varepsilon}_l}$ has the form
\begin{equation}
 \label{taketheform2}
 u_{{\hat\varepsilon}_l}( y_1, y_2)\, \approx\, \sum_{j=1}^N V(0, {\hat\theta})^{\frac {1}{p-1}}\,
 w\left(\sqrt{\frac{V(0, {\hat\theta})  \Big({\mathfrak a}_1(0, {\hat\theta})|{\hat n}_1({\hat\theta})|^2+{\mathfrak a}_2(0, {\hat\theta})|{\hat n}_2({\hat\theta})|^2\Big)}
 {|{\mathfrak a}_1(0, {\hat\theta})|^2+|{\mathfrak a}_2(0, {\hat\theta})|^2} }
 \frac {\, {\hat t} -{\hat\varepsilon}_l\, {\hat f}_j({\hat\theta})\, } {{\hat\varepsilon}_l}\right),
\end{equation}
where ${\hat n}({\hat\theta})=({\hat n}_1({\hat\theta}), {\hat n}_2({\hat\theta}))$ is the unit normal to ${\hat\Gamma}$.
The functions ${\hat f}_j$'s satisfy
\begin{equation}
\|{\hat f}_j\|_{\infty}\leq\, C|\ln{\hat\varepsilon}_l|^2,
\qquad
  \sum_{j=1}^N {\hat f}_j=O({\hat\varepsilon}_l^\mu)\quad\mbox{with } \mu>0,
\end{equation}
\begin{equation}
\min\limits_{1\leq j\leq N-1}({\hat f}_{j+1}-{\hat f}_j)
\approx 2|\ln{\hat\varepsilon}_l|\,
\sqrt{
\frac{|{\mathfrak a}_1(0, {\hat\theta})|^2+|{\mathfrak a}_2(0, {\hat\theta})|^2}
{V(0, {\hat\theta})  \Big({\mathfrak a}_1(0, {\hat\theta})|{\hat n}_1({\hat\theta})|^2+{\mathfrak a}_2(0, {\hat\theta})|{\hat n}_2({\hat\theta})|^2\Big)}
 }
 \, ,
\end{equation}
and solve the Jacobi-Toda system, for $j=1, \cdots, N$,
\begin{align}
&{\hat\varepsilon}_l^2\, \varsigma
\Big[\,
\widehat{\mathcal H}_1 {\hat f}''_j
+\widehat{\mathcal H}_1' {\hat f}'_j
+ \big(\widehat{\mathcal H}_2'-\widehat{\mathcal H}_3\big){\hat f}_j
\, \Big]
-e^{-\beta({\hat f}_j-{\hat f}_{j-1})}
+e^{-\beta({\hat f}_{j+1}-{\hat f}_j)}\approx0\quad\mbox{in } {(0, 1)},
\end{align}
with boundary conditions
$$
{\hat f}_j'(0)={\hat f}_j'(1),
\qquad
{\hat f}_j(0)={\hat f}_j(1),
\qquad\forall\, j=1, \cdots, N,
$$
where $\varsigma>0$ is defined in (\ref{varsigma}) with the conventions
$
{\hat f}_0=-\infty, \ {\hat f}_{N+1}=\infty.
$
\qed
\end{theorem}

Here are some words for further discussions.
Since the solutions have exponential decaying as $y$ leaves away the curve ${\tilde\Gamma}$, the proof of Theorem \ref{theorem 1.2}
is much more simpler than that of Theorem \ref{theorem 1.1}.
Whence, in the present paper, we will only provide the details to show the validity of Theorem \ref{theorem 1.1}.
Based on the same reason, a same result as in Theorem \ref{theorem 1.2} also holds for the first equation of \eqref{originalproblem}
in the whole space $\mathbb R^2$ under the condition $u(y)\rightarrow 0$ as $|y|\rightarrow\infty$.
Whence, Theorems \ref{theorem 1.1} and \ref{theorem 1.2} for the existence of cluster of multiple concentration layers
can be concerned as the extensions of the results in \cite{weixuyang} and \cite{delPKowWei2007},
where solutions with single concentration layer were constructed for partial confirmation of the two dimensional cases of Conjectures 2 and 3.

\medskip
However, in addition to the interaction between neighbouring layers in the cluster of multiple concentration layers,
the new ingredient is the role of the term ${\mathfrak a}(y)=\big({\mathfrak a}_1(y), {\mathfrak a}_2(y)\big)$.
This is the reason that we shall set up the new local  coordinates,
see \eqref{Fermicoordinates2} together with \eqref{Fermicoordinates-modified} and also \eqref{Fermicoordinates-modified-tilde}.
In the procedure of variational calculus, the deformations of the curves $\Gamma$ and ${\tilde\Gamma}$ are no longer directed along their normal directions,
see \eqref{deformation1}-\eqref{weightedlength1} and \eqref{deformation2}-\eqref{weightedlength2}.
The term ${\mathfrak a}$ will play an effect in the variational properties of the curves,
see the notions of non-degenerate stationary curves in Section \ref{Stationary and non-degenerate curves}.

\begin{remark}\label{remark12}
{\it
The Toda system was used first in \cite{delPKowWei2008}  to construct the clustered interfaces
for Allen-Cahn model in a two dimensional bounded domain.
Later, M. del Pino, M. Kowalczyk, J. Wei and J. Yang \cite{delPKowWeiYang} used the Jacobi-Toda
system in the construction of clustered phase transition layers for Allen-Cahn model on general Riemannian
manifolds.
The reader can refer to  \cite{delPKowWei2010, delPKowWei2013, JWeiYang2008, JWeiYang2010, weiyang20201, weiyang20202, YangYang2013}for more results.

For a ${\mathbf d}$-dimensional smooth compact Riemannian manifold $(\tilde{\mathcal M}, {\mathfrak g})$, M. del Pino, M. Kowalczyk, J. Wei and J. Yang \cite{delPKowWeiYang} considered the singularly perturbed Allen-Cahn equation
$$
\varepsilon^2\Delta_{ {\mathfrak g}} {u}\, +\, (1 - {u}^2)u \, =\, 0\quad \mbox{in } \tilde{\mathcal M},
$$
where $\varepsilon$ is a small parameter.
We let in what follows ${\mathcal K}$  be a minimal $({\mathbf d}-1)$-dimensional
embedded submanifold of $\tilde{\mathcal M}$,  which divides $\tilde{\mathcal M}$ into two open
components $\tilde{\mathcal M}_\pm$.
(The latter condition is not needed in some cases.)
Assume that
${\mathcal K}$ is non-degenerate in the sense that it does not support
non-trivial  Jacobi fields, and that
\begin{align}\label{posivity}
|\mathcal{A}_{{\mathcal K}}|^2+\mbox{Ric}_{\mathfrak g}(\nu_{{\mathcal K}}, \nu_{{\mathcal K}})>0
\quad \mbox{along } {\mathcal K}.
\end{align}
Then for each integer $N\geq 2$, they established the existence of a sequence $\varepsilon = \varepsilon_j\to 0$,
and solutions $u_{\varepsilon}$ with $N$-transition layers near ${\mathcal K}$, with mutual distance $O(\varepsilon|\ln \varepsilon|)$.

\medskip
As the above geometric language, we consider ${\mathbb R}^2$ as a manifold with the metric
$$
{\tilde{\mathfrak g}}=V^{2\sigma}(y) \big[ {\mathfrak a}_2(y){\mathrm d}y_1^2 +{\mathfrak a}_1(y) {\mathrm d}y_2^2 \big],
$$
with $\sigma$ in (\ref{sigma}) and $\Gamma$ as its submanifold with boundary.
In the manifold $({\mathbb R}^2, {\tilde{\mathfrak g}})$, $\Gamma$ is a non-degenerate  geodesic with endpoints on $\partial\Omega$.
In other words, in order to construct the clustering phase transition layers connecting the boundary $\partial\Omega$ in Theorem \ref{theorem 1.1},
we need the  condition (\ref{taupositivity}), which is similar  as (\ref{posivity}) in \cite{delPKowWeiYang}.
The reader can refer to Section \ref{section6.2}.
\qed
}
\end{remark}

\begin{remark}\label{remark13}
{\it
At $P_1$ and $P_2$ (the intersection points of $\Gamma$ and $\partial\Omega$),
the conditions in (\ref{boundaryadmissibility}) set up relations between the terms ${\mathfrak a}$, $V$
and the geometric properties of the curves $\partial\Omega$ and $\Gamma$.
For example, by recalling the unit normal $n(\theta)=(n_1(\theta), n_2(\theta))$  to $\Gamma$ and also the curvatures $k_1$ and $k_2$ of $\partial\Omega$ at $P_1$ and $P_2$,
the first one in \eqref{boundaryadmissibility} can be exactly expressed in the following form
\begin{align}
\frac{\sqrt{2}}{{\mathfrak a}_1(0, 0)}\Big[\partial_t{\mathfrak a}_1(0, 0)- \partial_t{\mathfrak a}_2(0, 0) \Big]n_1(0)n_2(0)
\, +\,
\sqrt{2}\Big({\tilde{\mathfrak a}'}_1(0)|n_1(0)|^2+{\tilde{\mathfrak a}'}_2(0) |n_2(0)|^2\Big)
\, +\,
k_1
=\frac{\beta'(0)}{\beta(0)},
\end{align}
with the conventions
$$
{\tilde{\mathfrak a}}_1(\theta)=\frac{{\mathfrak a}_1(0, \theta)}{\sqrt{|{\mathfrak a}_1(0, \theta)|^2+|{\mathfrak a}_2(0, \theta)|^2}\, },
\qquad
{\tilde{\mathfrak a}}_2(\theta)=\frac{{\mathfrak a}_2(0, \theta)}{\sqrt{|{\mathfrak a}_1(0, \theta)|^2+|{\mathfrak a}_2(0, \theta)|^2}\, },
$$

$$
\beta(\theta)
=\sqrt{\frac{ V(0, \theta)\big(a_1(0, \theta)|n_1(\theta)|^2+a_2(0, \theta)|n_2(\theta)|^2\big)}
{|a_1(0, \theta)|^2+|a_2(0, \theta)|^2}},
$$
where we have used \eqref{tildea1a2}, \eqref{alpha-beta}, \eqref{b1=mathfrakw0}, \eqref{b2=mathfrakw1}.
These conditions  will be used to decompose the interaction of neighbouring layers on the boundary $\partial\Omega$, see Remark \ref{remark61}.
On the other hand, we still have to deal with the delicate boundary terms for the reduced equations in Section \ref{section6.2}.
\qed
}
\end{remark}

\begin{remark}\label{remark14}
We construct solutions with multiple clustering concentration layers only for a sequence of $\varepsilon$ due to the coexistence of two types of resonances,
see also the fourth open question in "Concluding Remarks" of \cite{delPKowWei2010}.
The first one is due to the instability of the profile function $w$,
see Proposition \ref{proposition7point1}, in which we can impose the following gap condition for $\varepsilon$
\begin{align}\label{gapconditionofve}
|\lambda_*-j^2\varepsilon^2|\geq {\tilde c}\, \varepsilon, \quad \forall\, j \in \mathbb{N},
\end{align}
where ${\tilde c}$ is a given small positive constant.
In the above, $\lambda_{*}$ is a positive constant given by
 \begin{align}
\label{definenumber}
\lambda_{*}\, =\, \frac {\lambda_0 \ell^2}{\pi^2},
\end{align}
where $\lambda_0$ and $\ell$ are the positive constants given in \eqref{lambda0} and \eqref{ell}.
More details about this resonance phenomena were described in \cite{delPKowWei2007}.
The other one comes from the Jacobi-Toda system which was concerned in \cite{delPKowWeiYang}.
In this case, we shall choose a sequence of $\varepsilon$ from those satisfying (\ref{gapconditionofve}), see Proposition \ref{proposition7point2}.
\qed
\end{remark}

\medskip
By the rescaling
\begin{equation}\label{rescaling}
  y\, =\, \varepsilon {\tilde y}
\end{equation}
in ${\mathbb R}^2$, problem (\ref{originalproblem}) will be rewritten as
\begin{align}
\label{problemafterscaling}
{\mathrm {div}}\big( \nabla_{{\mathfrak a}(\varepsilon \tilde y)} u\big)- V(\varepsilon \tilde y)u+u^p\, =\, 0\quad{\rm in }\ \Omega_\varepsilon,
\qquad
\nabla_{{\mathfrak a}(\varepsilon \tilde y)} u \cdot\nu_{\varepsilon}\, =\, 0\quad\mbox{on }\ \partial \Omega_\varepsilon,
\end{align}
where $\Omega_\varepsilon=\Omega/\varepsilon$ and $\Gamma_\varepsilon=\Gamma/\varepsilon$, $\nu_\varepsilon$ is the unit outer normal of $\partial \Omega_\varepsilon$.
The remaining part of this paper is devoted to the proof of Theorem \ref{theorem 1.1}, which will be organized as follows:
\begin{itemize}
\item[1.] In Section \ref{section2}, we will set up a coordinate system in a neighborhood of $\Gamma$. Next we write down the local form of (\ref{problemafterscaling}),
especially the differential operators the differential operators ${\rm div}\big(\nabla_{\mathfrak a(y)}u\big)$ and $\nabla_{{\mathfrak a}(y)}u\cdot\nu$.
This local coordinate system also help us set up the stationary and non-degeneracy conditions for the curve $\Gamma$,
see (\ref{stationary}) and (\ref{nondegeneracy}).

\item[2.] We will set up an outline of the proof in Section \ref{thegluingprocedure}, which involves the gluing procedure from \cite {delPKowWei2007}, so that we can transform \eqref{problemafterscaling} into a projected form,
see \eqref{system-1}-\eqref{system-4}.

\item[3.] In Section \ref{section4}, we are devoted to the constructing of a local approximate solution in such a way that it solves the nonlinear problem  locally up to order $O(\varepsilon^2)$.

\item[4.]
To get a real solution, the well-known infinite dimensional reduction method \cite{delPKowWei2007}
will be needed in Sections \ref{section5}-\ref{sectionsolvingreducedequation}.
In fact, the reduced problem involves a Toda-Jacobi system and  inherits the resonance phenomena, which will be handled by complicated Fourier analysis.
\end{itemize}

\section{Geometric preliminaries}
\label{section2}
\setcounter{equation}{0}
In this section, we will set up  a coordinate system in a neighborhood of $\Gamma$.
This system is similar to the modified Fermi coordinates in \cite{weixuyang}.
However, some adaptions should be introduced due to the existence of the term ${\mathfrak a}$ in \eqref{originalproblem},
which make the geometric computations much more complicated.
The differential operators in \eqref{originalproblem} will be then derived in the local coordinates.
The notion of a stationary and non-degenerate curve $\Gamma$  will be also derived in the last part of this section.

\subsection{Modified Fermi coordinates}\label{Fermi}\

Recall the assumptions  {\bf (A1)}-{\bf (A3)} in Section \ref{section1} and notation therein.
For basic notions of curves, such as the signed curvature, the reader can refer to the book by do Carmo \cite{docarmo}.
\\[2mm]
\noindent {\bf Step  1.} Let the natural parameterization of the curve $\Gamma$ be as follows.
$$
\gamma_0:[0, 1]\rightarrow \Gamma\subset \bar\Omega\subset {\mathbb R}^2.
$$
For some small positive number $\sigma_0$, one can make a smooth extension and define the mapping
$$
\gamma=(\gamma_1, \gamma_2):(-\sigma_0, 1+\sigma_0)\rightarrow {\mathbb R}^2,
$$
such that
$$
\gamma(\tilde {\theta})\, =\, \gamma_0({\tilde\theta}), \quad \forall\, {\tilde\theta}\in [0, 1].
$$
There holds the Frenet formula
\begin{align}\label{Frenet}
\gamma''\, =\, kn \quad\textrm{and}\quad n'\, =\, -k\gamma',
\end{align}
where $k$, $n=(n_1, n_2)$ are the curvature and the normal of $\gamma$.
The relations
\begin{align}
|\gamma_1'({\tilde\theta})|^2\, +\, |\gamma_2'({\tilde\theta})|^2=1,
\qquad
|n_1({\tilde\theta})|^2+|n_2({\tilde\theta})|^2=1,
\qquad
\gamma_1'({\tilde\theta})n_1({\tilde\theta})+ \gamma_2'({\tilde\theta})n_2({\tilde\theta})=0,
\label{perpendicular}
\end{align}
will give that
\begin{equation}\label{relationofga&n}
 \big(\gamma_1'({\tilde\theta}), \gamma_2'({\tilde\theta})\big)\, =\, \big(-n_2({\tilde\theta}), n_1({\tilde\theta})\big),
\end{equation}
and
\begin{align}\label{relationofga&n2}
-\gamma_1'({\tilde\theta})n_2({\tilde\theta})+\gamma_2'({\tilde\theta})n_1({\tilde\theta})\, =\, 1.
\end{align}

\medskip
Choosing $\delta_0>0$ very small, and setting
$$
{\mathfrak S}_1\, \equiv\, (-\delta_0, \delta_0)\times (-\sigma_0, 1+\sigma_0),
$$
we construct the following mapping
\begin{align}
{\hat{\mathbb H}}: {\mathfrak S}_1\rightarrow {\hat{\mathbb H}}({\mathfrak S}_1)\, \equiv\, {\hat\Omega}_{\delta_0, \sigma_0}
\quad
\mbox{with}
\quad
{\hat{\mathbb H}}(\tilde {t}, \tilde {\theta})\, =\, \gamma({\tilde\theta})+{\tilde t}\, n({\tilde\theta}).
\label{Fermicoordinates0}
\end{align}
Note that ${\hat{\mathbb H}}$ is a diffeomorphism (locally) and ${\hat{\mathbb H}}(0, {\tilde\theta})=\gamma({\tilde\theta})$.
By this, we will write  the functions ${\mathfrak a}_1$, ${\mathfrak a}_2$ in the forms ${\mathfrak a}_1({\tilde t}, {\tilde\theta})$ and ${\mathfrak a}_2({\tilde t}, {\tilde\theta})$
and then set
\begin{align}\label{tildea1a2}
{\tilde{\mathfrak a}}_1({\tilde\theta})=\frac{{\mathfrak a}_1(0, {\tilde\theta})}{\sqrt{|{\mathfrak a}_1(0, {\tilde\theta})|^2+|{\mathfrak a}_2(0, {\tilde\theta})|^2}\, },
\qquad
{\tilde{\mathfrak a}}_2({\tilde\theta})=\frac{{\mathfrak a}_2(0, {\tilde\theta})}{\sqrt{|{\mathfrak a}_1(0, {\tilde\theta})|^2+|{\mathfrak a}_2(0, {\tilde\theta})|^2}\, }.
\end{align}
Note that
\begin{align}\label{tildea1=a2}
{\tilde{\mathfrak a}}_1(0)={\tilde{\mathfrak a}}_2(0)={\tilde{\mathfrak a}}_1(1)={\tilde{\mathfrak a}}_2(1)=\frac{1}{\sqrt 2},
\end{align}
due to the assumptions in \eqref{a1=a2}.
After that, we construct another mapping
\begin{align}
{\mathbb H}: {\mathfrak S}_1\rightarrow {\mathbb H}({\mathfrak S}_1)\, \equiv\, \Omega_{\delta_0, \sigma_0}
\quad
\mbox{with}
\quad
\mathbb{H}(\tilde {t}, \tilde {\theta})\, =\, \gamma({\tilde\theta})
\, +\,
{\tilde t}\, \big({\tilde{\mathfrak a}}_1({\tilde\theta})n_1({\tilde\theta}), \, {\tilde{\mathfrak a}}_2( {\tilde\theta})n_2({\tilde\theta})\big),
\label{Fermicoordinates2}
\end{align}
in such a way that it is a local diffeomorphism  and ${\mathbb H}(0, {\tilde\theta})=\gamma({\tilde\theta})$.
This is due to the fact that ${\mathfrak a}_1$ and ${\mathfrak a}_2$ are positive functions.

\medskip
\noindent {\bf Step 2}.
Recall ${\mathcal C}_1$, ${\mathcal C}_2$ given in the assumptions {\bf (A1)}-{\bf (A3)} in Section \ref{section1} and then  denote the preimages
$$
{\tilde{\mathcal C}}_1\, \equiv\, {\mathbb H}^{-1}\left(\mathcal C_1\right)
\quad\textrm{and}\quad
{\tilde{\mathcal C}}_2\, \equiv\, {\mathbb H}^{-1}\left(\mathcal C_2\right),
$$
which are two smooth curves in $({\tilde t}, {\tilde\theta})$ coordinates of \eqref{Fermicoordinates2} and can be parameterized respectively by $\left({\tilde t}, {\tilde \varphi}_1({\tilde t})\right)$ and  $ \left({\tilde t}, {\tilde \varphi}_2({\tilde t})\right)$
for some smooth functions ${\tilde \varphi}_1({\tilde t})\, {\rm and}\; {\tilde \varphi}_2({\tilde t})$ with the properties
\begin{align}
{\tilde \varphi}_1(0)=0,
\qquad
{\tilde \varphi}_2(0)=1.
\label{fact0}
\end{align}
We define a  mapping
$$
\tilde{{\mathbb H}}:\, {\mathfrak S}_1\rightarrow {\mathfrak S}_2\, \equiv\, \tilde{{\mathbb H}}({\mathfrak S}_1)\subset \Bbb R^2,
$$
such that
$$   t={\tilde t}, \quad \theta=\frac {{\tilde\theta}-{\tilde \varphi}_1({\tilde t})}{{\tilde \varphi}_2({\tilde t})-{\tilde \varphi}_1({\tilde t})}. $$
\noindent This transformation will straighten up the curves ${\tilde{\mathcal C}}_1$ and ${\tilde{\mathcal C}}_2$.
It is obvious that
\begin{align*}
{\tilde {\mathbb H}}^{-1}(0, \theta)\, =\, (0, {\theta}), \ \theta\in [0, 1].
\end{align*}
Moreover, we have

\begin{lemma}
\label{lemma2.1}
There hold
\begin{align}\label{fact2}
\tilde\varphi_1'(0)=0,
\quad
\tilde\varphi_2'(0)=0,
\quad
\tilde{\varphi}_1''(0)={\tilde k}_1,
\quad
\tilde{\varphi}_2''(0)={\tilde k}_2,
\end{align}
where
\begin{equation}
  {\tilde k}_1
  =
  \frac
  {\big(|{\tilde{\mathfrak a}}_1(0)n_1(0)|^2
  \, +\,
  |{\tilde{\mathfrak a}}_2(0)n_2(0)|^2\big)^{3/2}}
  {{\tilde{\mathfrak a}}_1(0) |n_1(0)|^2
  \, +\,
  {\tilde{\mathfrak a}}_2(0) |n_2(0)|^2}
  \, k_1
  =\frac{1}{2}k_1,
    \label{tildek1}
\end{equation}
\begin{equation}
 {\tilde k}_2 =
  \frac
  {\big(|{\tilde{\mathfrak a}}_1(1)n_1(1)|^2
  \, +\,
  |{\tilde{\mathfrak a}}_2(1)n_2(1)|^2\big)^{3/2}}
  {{\tilde{\mathfrak a}}_1(1) |n_1(1)|^2
  \, +\,
  {\tilde{\mathfrak a}}_2(1) |n_2(1)|^2}
  \, k_2
  =\frac{1}{2}k_2.
  \label{tildek2}
\end{equation}
\end{lemma}
\proof
In fact, the curves can be expressed in the following forms
\begin{align*}
{\mathcal C}_1:&\ {\mathbb H}\big({\tilde t}, {\tilde \varphi}_1({\tilde t})\big)
\, =\, \gamma\big({\tilde \varphi}_1({\tilde t})\big)
\, +\,
{\tilde t}\,
\Big({\tilde{\mathfrak a}}_1\big( {\tilde \varphi}_1({\tilde t})\big) n_1\big({\tilde \varphi}_1({\tilde t})\big),
\,
{\tilde{\mathfrak a}}_2\big( {\tilde \varphi}_1({\tilde t})\big) n_2\big({\tilde \varphi}_1({\tilde t})\big)\Big),
\\[2mm]
{\mathcal C}_2:&\ {\mathbb H}\big({\tilde t}, {\tilde \varphi}_2({\tilde t})\big)
\, =\,\gamma\big({\tilde \varphi}_2({\tilde t})\big)
\, +\,
{\tilde t}
\, \Big({\tilde{\mathfrak a}}_1\big( {\tilde \varphi}_2({\tilde t})\big) n_1\big({\tilde \varphi}_2({\tilde t})\big),
\,
{\tilde{\mathfrak a}}_2\big( {\tilde \varphi}_2({\tilde t})\big) n_2\big({\tilde \varphi}_2({\tilde t})\big)\Big),
\\[2mm]
\Gamma:&\ {\mathbb H}({\tilde\theta})\, =\, \gamma({\tilde\theta}).
\end{align*}
It follows that the tangent vectors of ${\mathcal C}_1$ at $P_1$ can be written as
\begin{align*}
\frac {{\mathrm d} {\mathcal C}_1}{{\mathrm d}{\tilde t}}\Big|_{{\tilde t}=0}
\, =\,
&\frac {\partial \gamma}{\partial {\tilde\theta}}\Big|_{{\tilde\theta}\, =\, {{\tilde\varphi}_1(0)}}
\cdot \frac {{\mathrm d} {{\tilde\varphi}_1}}{{\mathrm d}  {\tilde t}}\Big|_{{\tilde t}=0}
\, +\,
\Big({\tilde{\mathfrak a}}_1\big( {\tilde \varphi}_1({\tilde t})\big) n_1\big({\tilde \varphi}_1({\tilde t})\big),
\,
{\tilde{\mathfrak a}}_2\big( {\tilde \varphi}_1({\tilde t})\big) n_2\big({\tilde \varphi}_1({\tilde t})\big)\Big)\Big|_{{\tilde t}=0}
\\[2mm]
\, =\, &
\gamma'(0){\tilde\varphi}_1'(0)
\, +\,
\big({\tilde{\mathfrak a}}_1(0) n_1(0), \, {\tilde{\mathfrak a}}_2(0) n_2(0)\big)
,
\end{align*}
and the tangent vector of $\Gamma$ at $P_1$ is $\gamma'(0)$.
According to the condition: $\Gamma\bot \partial \Omega$ at $P_1$,
we have that
$$
\left\langle\frac {{\mathrm d}{\mathcal C}_1}{{\mathrm d}{\tilde t}}\Big|_{{\tilde t}=0}, \, \, \gamma'(0)  \right\rangle\, =\, 0.
$$
By \eqref{perpendicular} and \eqref{tildea1=a2}, we have
$$
{\tilde{\mathfrak a}}_1(0) n_1(0)\gamma_1'(0)\, +\, {\tilde{\mathfrak a}}_2(0) n_2(0)\gamma_2'(0)=0,
$$
and then drive from the above to get
$$ \tilde{\varphi}_1'(0)=0.$$
Similarly, we can show ${\tilde \varphi}_2'(0)=0$.

\medskip
The curve ${\mathcal C}_1$ can be expressed in the following form
\begin{equation*}
\begin{split}
  {\mathcal C}_1: {\mathbb H}\big({\tilde t}, {\tilde \varphi}_1({\tilde t})\big)
\,
&=\Big(\gamma_1\big({\tilde \varphi}_1({\tilde t})\big)+{\tilde t} {\tilde{\mathfrak a}}_1\big( {\tilde \varphi}_1({\tilde t})\big)n_1\big({\tilde \varphi}_1({\tilde t})\big),
\, \,
\gamma_2\big({\tilde \varphi}_1({\tilde t})\big)+{\tilde t} {\tilde{\mathfrak a}}_2\big( {\tilde \varphi}_1({\tilde t})\big) n_2\big({\tilde \varphi}_1({\tilde t})\big)\Big)
\\[2mm]
&\equiv\, \big( y_1({\tilde t})\, , \, y_2({\tilde t})\big).
\end{split}
\end{equation*}
The calculations
\begin{align*}
y_1'({\tilde t})
=&\,
\gamma_1'({\tilde\varphi}_1)\cdot \frac {{\mathrm d}{{\tilde\varphi}_1}}{{\mathrm d} {\tilde t}}
+{\tilde{\mathfrak a}}_1({\tilde\varphi}_1) n_1({\tilde\varphi}_1)
+{\tilde t}\cdot{\tilde{\mathfrak a}}_1({\tilde\varphi}_1)\cdot n_1'({\tilde\varphi}_1)\cdot \frac {{\mathrm d} {{\tilde\varphi}_1}}{{\mathrm d} {\tilde t}}
+{\tilde t}\cdot{\tilde{\mathfrak a}}_1'({\tilde\varphi}_1)\cdot n_1({\tilde\varphi}_1)\cdot \frac {{\mathrm d} {{\tilde\varphi}_1}}{{\mathrm d} {\tilde t}},
\\[2mm]
y_2'({\tilde t})
=&\,
\gamma_2'({\tilde\varphi}_1)\cdot \frac {{\mathrm d}{{\tilde\varphi}_1}}{{\mathrm d} {\tilde t}}
+{\tilde{\mathfrak a}}_2({\tilde\varphi}_1) n_2({\tilde\varphi}_1)
+{\tilde t}\cdot{\tilde{\mathfrak a}}_2({\tilde\varphi}_1)\cdot n_2'({\tilde\varphi}_1)\cdot \frac {{\mathrm d} {{\tilde\varphi}_1}}{{\mathrm d} {\tilde t}}
+{\tilde t}\cdot{\tilde{\mathfrak a}}_2'({\tilde\varphi}_2)\cdot n_2({\tilde\varphi}_1)\cdot \frac {{\mathrm d} {{\tilde\varphi}_1}}{{\mathrm d} {\tilde t}},
\end{align*}
and
\begin{align*}
y_1''({\tilde t})=&\,
\gamma_1''({\tilde\varphi}_1)\cdot \Big(\frac {{\mathrm d} {{\tilde\varphi}_1}}{{\mathrm d} {\tilde t}}\Big)^2
+\gamma_1'({\tilde\varphi}_1)\cdot \frac {{\mathrm d}^2 {{\tilde\varphi}_1}}{{\mathrm d} {\tilde t}^2}
+2\, {\tilde{\mathfrak a}}_1'({\tilde\varphi}_1)\cdot n_1({\tilde\varphi}_1)\cdot\frac {{\mathrm d} {{\tilde\varphi}_1}}{{\mathrm d} {\tilde t}}
\\[2mm]
&
+2\, {\tilde{\mathfrak a}}_1({\tilde\varphi}_1)\cdot n_1'({\tilde\varphi}_1)\cdot\frac {{\mathrm d} {{\tilde\varphi}_1}}{{\mathrm d} {\tilde t}}
+2{\tilde t}\, {\tilde{\mathfrak a}}_1'({\tilde\varphi}_1) n_1'({\tilde\varphi}_1)\cdot\Big(\frac {{\mathrm d} {{\tilde\varphi}_1}}{{\mathrm d} {\tilde t}}\Big)^2
+{\tilde t}\, {\tilde{\mathfrak a}}_1({\tilde\varphi}_1) n_1''({\tilde\varphi}_1)\cdot\Big(\frac {{\mathrm d} {{\tilde\varphi}_1}}{{\mathrm d} {\tilde t}}\Big)^2
\\[2mm]
&
+{\tilde t}\, {\tilde{\mathfrak a}}_1({\tilde\varphi}_1) n_1'({\tilde\varphi}_1)\cdot\frac {{\mathrm d}^2 {{\tilde\varphi}_1}}{{\mathrm d} {\tilde t}^2}
+{\tilde t}\, {\tilde{\mathfrak a}}_1''({\tilde\varphi}_1) n_1({\tilde\varphi}_1)\cdot\Big(\frac {{\mathrm d} {{\tilde\varphi}_1}}{{\mathrm d} {\tilde t}}\Big)^2
+{\tilde t}\, {\tilde{\mathfrak a}}_1'({\tilde\varphi}_1) n_1({\tilde\varphi}_1)\cdot\frac {{\mathrm d}^2 {{\tilde\varphi}_1}}{{\mathrm d} {\tilde t}^2},
 \end{align*}
\begin{align*}
y_2''({\tilde t})=&\,
\gamma_2''({\tilde\varphi}_1)\cdot \Big(\frac {{\mathrm d} {{\tilde\varphi}_1}}{{\mathrm d} {\tilde t}}\Big)^2
+\gamma_2'({\tilde\varphi}_1)\cdot \frac {{\mathrm d}^2 {{\tilde\varphi}_1}}{{\mathrm d} {\tilde t}^2}
+2\, {\tilde{\mathfrak a}}_2'({\tilde\varphi}_1)\cdot n_2({\tilde\varphi}_1)\cdot\frac {{\mathrm d} {{\tilde\varphi}_1}}{{\mathrm d} {\tilde t}}
\\[2mm]
&
+2\, {\tilde{\mathfrak a}}_2({\tilde\varphi}_1)\cdot n_2'({\tilde\varphi}_1)\cdot\frac {{\mathrm d} {{\tilde\varphi}_1}}{{\mathrm d} {\tilde t}}
+2{\tilde t}\, {\tilde{\mathfrak a}}_2'({\tilde\varphi}_1) n_2'({\tilde\varphi}_1)\cdot\Big(\frac {{\mathrm d} {{\tilde\varphi}_1}}{{\mathrm d} {\tilde t}}\Big)^2
+{\tilde t}\, {\tilde{\mathfrak a}}_2({\tilde\varphi}_1) n_2''({\tilde\varphi}_1)\cdot\Big(\frac {{\mathrm d} {{\tilde\varphi}_1}}{{\mathrm d} {\tilde t}}\Big)^2
\\[2mm]
&
+{\tilde t}\, {\tilde{\mathfrak a}}_2({\tilde\varphi}_1) n_2'({\tilde\varphi}_1)\cdot\frac {{\mathrm d}^2 {{\tilde\varphi}_1}}{{\mathrm d} {\tilde t}^2}
+{\tilde t}\, {\tilde{\mathfrak a}}_2''({\tilde\varphi}_1) n_2({\tilde\varphi}_1)\cdot\Big(\frac {{\mathrm d} {{\tilde\varphi}_1}}{{\mathrm d} {\tilde t}}\Big)^2
+{\tilde t}\, {\tilde{\mathfrak a}}_2'({\tilde\varphi}_1) n_2({\tilde\varphi}_1)\cdot\frac {{\mathrm d}^2 {{\tilde\varphi}_1}}{{\mathrm d} {\tilde t}^2},
 \end{align*}
imply that
 \begin{align*}
  |y_1'(0)|^2+|y_2'(0)|^2
  \, =\,
  |{\tilde{\mathfrak a}}_1(0)n_1(0)|^2
  \, +\,
  |{\tilde{\mathfrak a}}_2(0)n_2(0)|^2,
 \end{align*}
and
\begin{align*}
y_1'(0)\, y_2''(0)-y_1''(0)\, y_2'(0)
=&\, \big({\tilde{\mathfrak a}}_1(0)n_1(0)\, \gamma_2'(0)
\,-\, {\tilde{\mathfrak a}}_2(0)n_2(0)\, \gamma_1'(0)\big)\, {\tilde\varphi}_1''(0)
\\[2mm]
\, =&\,
\big({\tilde{\mathfrak a}}_1(0) |n_1(0)|^2
\, +\,
{\tilde{\mathfrak a}}_2(0) |n_2(0)|^2\big)\, {\tilde\varphi}_1''(0).
\end{align*}
Therefore, the signed curvature of the curve $\mathcal{C}_1$ at the point $P_1$ is
\begin{equation*}
k_1=\frac{\, y_1'(0)y_2''(0)-y_1''(0)y_2'(0)\, }{\big[\, |y_1'(0)|^2+|y_2'(0)|^2\, \big]^{\frac{3}{2}}}
  =
  \frac{{\tilde{\mathfrak a}}_1(0) |n_1(0)|^2
  \, +\,
  {\tilde{\mathfrak a}}_2(0) |n_2(0)|^2}
  {\big(|{\tilde{\mathfrak a}}_1(0)n_1(0)|^2
  \, +\,
  |{\tilde{\mathfrak a}}_2(0)n_2(0)|^2\big)^{3/2}}
  \, \tilde{\varphi}_1''(0)
  =2\tilde{\varphi}_1''(0).
\end{equation*}
Similarly, we can show
\begin{equation*}
k_2=
\frac{{\tilde{\mathfrak a}}_1(1) |n_1(1)|^2
\, +\,
{\tilde{\mathfrak a}}_2(1) |n_2(1)|^2}
{\big(|{\tilde{\mathfrak a}}_1(1)n_1(1)|^2
 \, +\,
|{\tilde{\mathfrak a}}_2(1)n_2(1)|^2\big)^{3/2}} \, \tilde{\varphi}_2''(0)
=2\tilde{\varphi}_2''(0).
\end{equation*}

\medskip
\noindent {\bf Step 3}.
We define the {\bf modified Fermi coordinates}
\begin{equation}
\label{Fermicoordinates-modified}
(y_1, y_2)\, =\, F(t, \theta)
\, =\,
{\mathbb H}\circ{\tilde {\mathbb H}}^{-1}(t, \theta)\, :\, (-\delta_0, \delta_0)\times(-\sigma_0, 1+\sigma_0)\rightarrow{\mathbb R}^2
\end{equation}
for  given small positive constants $\sigma_0$ and $\delta_0$.
More precisely,  we write $F(t, \theta)=\big(F_1(t, \theta), \, F_2(t, \theta)\big)$ with
\begin{equation}\label{F1}
F_i(t, \theta)
\, =\, \gamma_i\big(\Theta(t, \theta)\big)
\, +\,
t\, {\tilde{\mathfrak a}}_i\big(\Theta(t, \theta)\big) n_i\big(\Theta(t, \theta)\big),
\quad
i=1,2,
\end{equation}
where
\begin{equation}\label{Theta}
\Theta(t, \theta)
\, \equiv\, \big({\tilde \varphi}_2(t)-{\tilde \varphi}_1(t)\big)\theta
\, +\, {\tilde \varphi}_1(t).
\end{equation}
From \eqref{fact0}-\eqref{tildek2}, we have
\begin{align}
\Theta(0, \theta)\, =\, \theta,
\qquad
\Theta_t(0, \theta)\, =\, 0,
\qquad
\Theta_\theta(0, \theta)\, =\, 1,
\qquad
\Theta_{\theta t}(0, \theta)\, =\, 0,
\label{Thetaderivative1}
\\[2mm]
\Theta_{tt}(0, \theta)\, =\, \big({\tilde k}_2-{\tilde k}_1\big)\theta +{\tilde k}_1,
\qquad
\Theta_{tt\theta}(0, \theta)\, =\, {\tilde k}_2-{\tilde k}_1.
\label{Thetaderivative2}
\end{align}
These quantities will play an important role in the further settings.

\medskip
We now derive the asymptotical behaviors of the coordinates.
For given $i=1$ or $i=2$,
consider the derivative of first order
\begin{align}
\begin{aligned}
\frac{{\partial} F_i}{{\partial} t}
&\, =\,
\gamma_i'(\Theta)\cdot\Theta_t
\, +\, {\tilde{\mathfrak a}}_i(\Theta) n_i(\Theta)
\, +\, t\, {\tilde{\mathfrak a}}_i'(\Theta) n_i(\Theta)\cdot {\Theta}_t
\, +\, t\, {\tilde{\mathfrak a}}_i(\Theta) n_i'(\Theta)\cdot {\Theta}_t,
\end{aligned}
\label{partialFipartialt}
\end{align}
and also the derivative of second order
\begin{equation}
\begin{split}
\frac{{\partial}^2 F_i}{{\partial} t^2}
&\, =\, \gamma_i''(\Theta)\cdot |\Theta_t|^2
\, +\, \gamma_i'(\Theta)\cdot \Theta_{tt}
\, +\, 2{\tilde{\mathfrak a}}_i'(\Theta) n_i(\Theta)\cdot \Theta_{t}
\, +\, 2{\tilde{\mathfrak a}}_i(\Theta) n_i'(\Theta)\cdot \Theta_{t}
\\[2mm]
&\quad
\, +\, 2t\, {\tilde{\mathfrak a}}_i'(\Theta) n_i'(\Theta)\cdot|\Theta_t|^2
\, +\, t\, {\tilde{\mathfrak a}}_i''(\Theta) n_i(\Theta)\cdot|\Theta_t|^2
\, +\, t\, {\tilde{\mathfrak a}}_i'(\Theta) n_i(\Theta)\cdot\Theta_{tt}
\\[2mm]
&\quad
\, +\, t\, {\tilde{\mathfrak a}}_i(\Theta) n_i''(\Theta)\cdot|\Theta_t|^2
\, +\, t\, {\tilde{\mathfrak a}}_i(\Theta) n_i'(\Theta)\cdot\Theta_{tt}.
\end{split}
\label{partialFipartialtt}
\end{equation}
These imply that
\begin{equation}
\frac{{\partial} F_i}{{\partial} t}(0, \theta)\, =\, {\tilde{\mathfrak a}}_i(\theta)n_i(\theta),
\end{equation}
and
\begin{equation}
 \frac{{\partial}^2F_i}{{\partial} t^2}(0, \theta)
\, =\, \gamma_i'(\theta)\cdot \Theta_{tt}(0, \theta)
\, =\, \big({\tilde k}_2-{\tilde k}_1\big)\theta +{\tilde k}_1
\, \equiv\, q_i(\theta).
\label{q_i}
\end{equation}
Hence
\begin{equation}
 \Theta_{tt}(0, \theta)\, =\, q_1 \gamma_1'+q_2 \gamma_2'\, =\, -q_1n_2+q_2n_1.
 \label{2.22}
\end{equation}
Here comes the derivative of third order
\begin{align}\label{m_i}
\frac{{\partial}^3F_i}{{\partial} t^3}(0, \theta)
&\, =\,
\bigg[ \gamma_i'''(\Theta)\cdot (\Theta_t)^3
+3\gamma_i''(\Theta)\cdot \Theta_t\cdot \Theta_{tt}
+\gamma_i'(\Theta)\cdot \Theta_{ttt}
\nonumber\\[2mm]
&\qquad
+3{\tilde{\mathfrak a}}_i''(\Theta) n_i(\Theta)\cdot |\Theta_t|^2
+8{\tilde{\mathfrak a}}_i'(\Theta) n_i'(\Theta)\cdot |\Theta_t|^2
+3{\tilde{\mathfrak a}}_i'(\Theta) n_i(\Theta)\cdot \Theta_{tt}
\nonumber\\[2mm]
&\qquad
+3{\tilde{\mathfrak a}}_i(\Theta) n_i''(\Theta)\cdot |\Theta_t|^2
+3{\tilde{\mathfrak a}}_i(\Theta) n_i'(\Theta)\cdot \Theta_{tt}
\bigg]\Big|_{(t, \theta)=(0, \theta)}
\nonumber\\[2mm]
&\, =\, \gamma_i'(\theta)\cdot \Theta_{ttt}(0, \theta)
+3{\tilde{\mathfrak a}}_i'(\theta) n_i(\theta)\cdot \Theta_{tt}(0, \theta)
+3{\tilde{\mathfrak a}}_i(\theta) n_i'(\theta)\cdot \Theta_{tt}(0, \theta)
\nonumber\\[2mm]
&\, \equiv\, m_i(\theta).
\end{align}
These results will be collected in the following way.
\begin{lemma}\label{derivativeofF}
The mapping $F$ has the following properties:
\\

\noindent {\rm\textbf{(1).}}
$\, \, F(0, \theta)\, =\, \gamma(\theta)$,
\qquad
$\, \, \frac{{\partial} F}{{\partial} t}(0, \theta)\, =\, \big({\tilde{\mathfrak a}}_1(\theta)n_1(\theta), {\tilde{\mathfrak a}}_2(\theta)n_2(\theta)\big)$,
\\

\noindent {\rm\textbf{(2).}}
$\, \, \, \frac{{\partial}^2 F}{{\partial} t^2}(0, \theta)=q(\theta)
\quad\mbox{with}\quad
q(\theta)\, =\, \big(q_1(\theta), q_2(\theta)\big)\quad \bot\, n(\theta)
$,
\\

\noindent {\rm\textbf{(3).}}
$\, \, \frac{{\partial}^3F}{{\partial} t^3}(0, \theta)=m(\theta)
\quad\mbox{with}\quad
m(\theta)\, =\, \big(m_1(\theta), m_2(\theta)\big)
$.
\\

\noindent Here $m_i$'s and $q_i$'s are given in \eqref{q_i}-\eqref{m_i}.
\qed
\end{lemma}

As a conclusion, as $t$ is small enough, there holds the expansion
\begin{align}
F(t, \theta)
\, =\, &
\gamma(\theta)
+t\, \big({\tilde{\mathfrak a}}_1(\theta)n_1(\theta), {\tilde{\mathfrak a}}_2(\theta)n_2(\theta)\big)
+\frac {t^2}{2}q(\theta)
\nonumber\\[2mm]
&+\frac {t^3}{6}m(\theta)
+O(t^4),
\quad
\forall\, \theta\in [0, 1],  t\in (-\delta_0, \delta_0),
\label{expansionofFermi}
\end{align}
where $\delta_0>0$ is a small constant.
This gives us that
\begin{align}
\frac{{\partial} F}{{\partial} t}(t, \theta)
\, =\, &
\big({\tilde{\mathfrak a}}_1(\theta)n_1(\theta), {\tilde{\mathfrak a}}_2(\theta)n_2(\theta)\big)
+t q(\theta)
+\frac{t^2}{2}m(\theta)
+O(t^3), \label{F_t}
\\[2mm]
\frac{{\partial} F_1}{{\partial} \theta}(t, \theta)
\, =\, &\gamma_1'(\theta)
+t{\tilde{\mathfrak a}}_1'(\theta)n_1(\theta)
-tk{\tilde{\mathfrak a}}_1(\theta)\gamma_1'(\theta)
+\frac{t^2}{2}q_1'(\theta)+O(t^3)
\nonumber\\[2mm]
\, =\, &-n_2(\theta)
+t{\tilde{\mathfrak a}}_1'(\theta)n_1(\theta)
+tk{\tilde{\mathfrak a}}_1(\theta)n_2(\theta)
+\frac{t^2}{2}q_1'(\theta)+O(t^3), \label{F1_theta}
\\[2mm]
\frac{{\partial} F_2}{{\partial} \theta}(t, \theta)
\, =\, &
\gamma_2'(\theta)
+t{\tilde{\mathfrak a}}_1'(\theta)n_2(\theta)
-tk{\tilde{\mathfrak a}}_1(\theta)\gamma_2'(\theta)
+\frac{t^2}{2}q_2'(\theta)+O(t^3)
\nonumber\\[2mm]
\, =\, &
n_1(\theta)
+t{\tilde{\mathfrak a}}_1'(\theta)n_2(\theta)
-tk{\tilde{\mathfrak a}}_1(\theta)n_1(\theta)
+\frac{t^2}{2}q_2'(\theta)+O(t^3), \label{F2_theta}
\end{align}
where we have used \eqref{relationofga&n} and the Frenet formula \eqref{Frenet}.
Moreover, there hold
\begin{align}\label{qi'}
q_i'(\theta)
\, =\, \gamma_i''(\theta)\cdot\Theta_{tt}(0, \theta)
\, +\, \gamma_i'(\theta)\Theta_{tt\theta}(0, \theta),
\end{align}
and specially,
\begin{align*}
q'(0)\, =\, \gamma{''}(0){\tilde k}_1\, +\, \gamma'(0)({\tilde k}_2-{\tilde k}_1),
\qquad
q'(1)\, =\, \gamma{''}(1){\tilde k}_2\, +\, \gamma'(1)({\tilde k}_2-{\tilde k}_1).
\end{align*}

\subsection{The metric}
In the local coordinates $(t, \theta)$ in \eqref{Fermicoordinates-modified},
here are the preparing computations for the metric matrix:
\begin{align}
\frac{{\partial} F_i}{{\partial} t}\frac{{\partial} F_i}{{\partial} t}
&\, =\, ({\tilde{\mathfrak a}}_i)^2|n_i|^2
+2t\, {\tilde{\mathfrak a}}_i n_i q_i
+t^2\, {\tilde{\mathfrak a}}_i n_i m_i
+t^2|q_i|^2
+O(t^3),
\label{Fit-Fit}
\end{align}
\begin{align}\label{Fit-Fitheta}
\frac{{\partial} F_1}{{\partial} t} \frac{{\partial} F_1}{{\partial} \theta}
\, =\, &
-{\tilde{\mathfrak a}}_1 n_1 n_2
+t {\tilde{\mathfrak a}}_1{\tilde{\mathfrak a}'}_1 |n_1|^2
+tk {\tilde{\mathfrak a}}_1^2 n_1 n_2
-t q_1n_2
+\frac{t^2}{2}{\tilde{\mathfrak a}}_1 q_1'n_1
\nonumber\\[2mm]
&+t^2{\tilde{\mathfrak a}' }_1q_1n_1
+t^2k {\tilde{\mathfrak a}}_1  q_1 n_2
-\frac{t^2}{2}m_1 n_2+O(t^3),
\\[2mm]
\frac{{\partial} F_2}{{\partial} t} \frac{{\partial} F_2}{{\partial} \theta}
\, =\, &
{\tilde{\mathfrak a}}_2 n_2 n_1
+t {\tilde{\mathfrak a}}_2{\tilde{\mathfrak a}'}_2 |n_2|^2
-tk {\tilde{\mathfrak a}}_2^2 n_2n_1
+t q_2n_1
+\frac{t^2}{2}{\tilde{\mathfrak a}}_2 q_2'n_2
\nonumber\\[2mm]
&+t^2{\tilde{\mathfrak a}' }_2q_2n_2
-t^2k {\tilde{\mathfrak a}}_2  q_2n_1
+\frac{t^2}{2}m_2n_1+O(t^3),
\end{align}
\begin{align}
\frac{{\partial} F_1}{{\partial} \theta} \frac{{\partial} F_1}{{\partial} \theta}
\, =\, & |n_2|^2
-2t {\tilde{\mathfrak a}'}_1  n_1n_2
-2t k{\tilde{\mathfrak a}}_1 |n_2|^2
-t^2 q_1'n_2
+t^2({\tilde{\mathfrak a}'}_1)^2|n_1|^2
\nonumber\\[2mm]
&+2t^2k{\tilde{\mathfrak a}}_1 {\tilde{\mathfrak a}'}_1 n_1n_2
+t^2k^2 ({\tilde{\mathfrak a}}_1)^2 |n_2|^2
+O(t^3),
\label{Fitheta-Fitheta}
\end{align}
and
\begin{align}
\frac{{\partial} F_2}{{\partial} \theta} \frac{{\partial} F_2}{{\partial} \theta}
\, =\, & |n_1|^2
+2t {\tilde{\mathfrak a}'}_2 n_1 n_2
-2t k{\tilde{\mathfrak a}}_2 |n_1|^2
+t^2 q_2'n_1
+t^2({\tilde{\mathfrak a}'}_2)^2|n_2|^2
\nonumber\\[2mm]
&-2t^2k{\tilde{\mathfrak a}}_2 {\tilde{\mathfrak a}'}_2 n_1n_2
+t^2k^2 ({\tilde{\mathfrak a}}_2)^2 |n_1|^2
+O(t^3).
\label{F2theta-F2theta}
\end{align}

\medskip
The elements of the metric matrix are:
\begin{align}
g_{11}\, =\, &
\frac{{\partial} F_1}{{\partial} t}\frac{{\partial} F_1}{{\partial} t}
\, +\,
\frac{{\partial} F_2}{{\partial} t}\frac{{\partial} F_2}{{\partial} t}
\nonumber\\[2mm]
\, =\, &\big[{\tilde{\mathfrak a}}_1^2 |n_1|^2 +{\tilde{\mathfrak a}}_2^2 |n_2|^2 \big]
+2t\, \big[{\tilde{\mathfrak a}}_1 n_1 q_1 + {\tilde{\mathfrak a}}_2 n_2 q_2\big]
\nonumber\\[2mm]
&+t^2\, \big[ {\tilde{\mathfrak a}}_1 n_1 m_1 + {\tilde{\mathfrak a}}_2 n_2 m_2 \big]
+t^2\big[|q_1|^2+|q_2|^2\big]+O(t^3),
\nonumber\\[2mm]
g_{12}\, =\, &
\frac{{\partial} F_1}{{\partial} t}\frac{{\partial} F_1}{{\partial} \theta}
\, +\,
\frac{{\partial} F_2}{{\partial} t}\frac{{\partial} F_2}{{\partial} \theta}
\nonumber\\[2mm]
=&\big[-{\tilde{\mathfrak a}}_1n_1n_2+{\tilde{\mathfrak a}}_2n_1n_2\big]
+t\big[{\tilde{\mathfrak a}}_1{\tilde{\mathfrak a}'}_1|n_1|^2+{\tilde{\mathfrak a}}_2{\tilde{\mathfrak a}'}_2 |n_2|^2\big]
-tk\big[-{\tilde{\mathfrak a}}_1^2n_1n_2+{\tilde{\mathfrak a}}_2^2n_1n_2\big]
\nonumber\\[2mm]
&+t \big[-q_1n_2+q_2n_1\big]
+\frac{t^2}{2}\big[{\tilde{\mathfrak a}}_1q_1'n_1+{\tilde{\mathfrak a}}_2q_2'n_2\big]
+t^2\big[{\tilde{\mathfrak a}'}_1q_1n_1+{\tilde{\mathfrak a}'}_2q_2n_2\big]
\nonumber\\[2mm]
&-t^2k\big[-{\tilde{\mathfrak a}}_1q_1n_2+{\tilde{\mathfrak a}}_2q_2n_1\big]
+\frac{t^2}{2}\big[-m_1n_2+m_2n_1\big]+O(t^3),
\end{align}
and
\begin{align}
g_{22}\, =\, &
\frac{{\partial} F_1}{{\partial} \theta}\frac{{\partial} F_1}{{\partial} \theta}
\, +\,
\frac{{\partial} F_2}{{\partial} \theta}\frac{{\partial} F_2}{{\partial} \theta}
\nonumber\\[2mm]
=& 1+2t \big[-{\tilde{\mathfrak a}'}_1n_1 n_2+{\tilde{\mathfrak a}'}_2n_1 n_2\big]
-2tk\big[{\tilde{\mathfrak a}}_1|n_2|^2+{\tilde{\mathfrak a}}_2|n_1|^2\big]
+t^2\big[-q_1'n_2+q_2'n_1\big]
\nonumber\\[2mm]
&+t^2\big[({\tilde{\mathfrak a}'}_1)^2|n_1|^2+({\tilde{\mathfrak a}'}_2)^2|n_2|^2\big]
-2t^2k\big[-{\tilde{\mathfrak a}}_1{\tilde{\mathfrak a}'}_1n_1n_2+{\tilde{\mathfrak a}}_2{\tilde{\mathfrak a}'}_2n_1n_2\big]
\nonumber\\[2mm]
&+t^2k^2 \big[{\tilde{\mathfrak a}}_1^2 |n_2|^2+{\tilde{\mathfrak a}}_2^2 |n_1|^2\big]+O(t^3).
\end{align}

\medskip
So the  determinant  of the metric matrix is
$$g=g_{11}g_{22}-g_{12}g_{12}$$
where
\begin{align*}
g_{11}g_{22}
=&\big[{\tilde{\mathfrak a}}_1^2 |n_1|^2 +{\tilde{\mathfrak a}}_2^2 |n_2|^2 \big]
+2t\big[{\tilde{\mathfrak a}}_1 n_1 q_1 + {\tilde{\mathfrak a}}_2 n_2 q_2\big]
+2t\big[{\tilde{\mathfrak a}}_1^2 |n_1|^2 +{\tilde{\mathfrak a}}_2^2 |n_2|^2 \big]\big[-{\tilde{\mathfrak a}'}_1n_1 n_2+{\tilde{\mathfrak a}'}_2n_1 n_2\big]
\\[2mm]
&-2tk\big[{\tilde{\mathfrak a}}_1|n_2|^2+{\tilde{\mathfrak a}}_2|n_1|^2\big]\big[{\tilde{\mathfrak a}}_1^2 |n_1|^2 +{\tilde{\mathfrak a}}_2^2 |n_2|^2 \big]
+t^2\, \big[ {\tilde{\mathfrak a}}_1 n_1 m_1 + {\tilde{\mathfrak a}}_2 n_2 m_2 \big]
+t^2\big[|q_1|^2+|q_2|^2\big]
\\[2mm]
&+4t^2\big[{\tilde{\mathfrak a}}_1 n_1 q_1 + {\tilde{\mathfrak a}}_2 n_2 q_2\big]\big[-{\tilde{\mathfrak a}'}_1+{\tilde{\mathfrak a}'}_2\big]n_1n_2
+t^2\big[{\tilde{\mathfrak a}}_1^2 |n_1|^2 +{\tilde{\mathfrak a}}_2^2 |n_2|^2 \big]\big[-q_1'n_2+q_2'n_1\big]
\\[2mm]
&+t^2\big[{\tilde{\mathfrak a}}_1^2 |n_1|^2 +{\tilde{\mathfrak a}}_2^2 |n_2|^2 \big]\big[({\tilde{\mathfrak a}'}_1)^2|n_1|^2+({\tilde{\mathfrak a}'}_2)^2|n_2|^2\big]
-4t^2k\big[{\tilde{\mathfrak a}}_1|n_2|^2+{\tilde{\mathfrak a}}_2|n_1|^2\big]\big[{\tilde{\mathfrak a}}_1 n_1 q_1 + {\tilde{\mathfrak a}}_2 n_2 q_2\big]
\nonumber\\[2mm]
&-2t^2k\big[{\tilde{\mathfrak a}}_1^2 |n_1|^2 +{\tilde{\mathfrak a}}_2^2 |n_2|^2 \big]\big[-{\tilde{\mathfrak a}}_1{\tilde{\mathfrak a}'}_1n_1n_2+{\tilde{\mathfrak a}}_2{\tilde{\mathfrak a}'}_2n_1n_2\big]
\nonumber\\[2mm]
&+t^2k^2\big[{\tilde{\mathfrak a}}_1^2 |n_1|^2 +{\tilde{\mathfrak a}}_2^2 |n_2|^2 \big] \big[{\tilde{\mathfrak a}}_1^2 |n_2|^2+{\tilde{\mathfrak a}}_2^2 |n_1|^2\big]+O(t^3),
\end{align*}
and
\begin{align*}
 g_{12}g_{12}
=&\big[-{\tilde{\mathfrak a}}_1 n_1 n_2+{\tilde{\mathfrak a}}_2 n_2n_1\big]^2
+2t\big[-{\tilde{\mathfrak a}}_1 n_1 n_2+{\tilde{\mathfrak a}}_2 n_2n_1\big]\big[{\tilde{\mathfrak a}}_1{\tilde{\mathfrak a}'}_1|n_1|^2+{\tilde{\mathfrak a}}_2{\tilde{\mathfrak a}'}_2 |n_2|^2\big]
\\[2mm]
&+2t \big[-{\tilde{\mathfrak a}}_1 n_1 n_2+{\tilde{\mathfrak a}}_2 n_2n_1\big]
\big[-q_1n_2+q_2n_1\big]
-2tk\big[-{\tilde{\mathfrak a}}_1 n_1 n_2+{\tilde{\mathfrak a}}_2 n_2n_1\big]
\big[-{\tilde{\mathfrak a}}_1^2 n_1n_2+{\tilde{\mathfrak a}}_2^2 n_2n_1\big]
\\[2mm]
&+t^2 \big[{\tilde{\mathfrak a}}_1{\tilde{\mathfrak a}'}_1|n_1|^2+{\tilde{\mathfrak a}}_2{\tilde{\mathfrak a}'}_2 |n_2|^2\big]^2
+2t^2 \big[{\tilde{\mathfrak a}}_1{\tilde{\mathfrak a}'}_1|n_1|^2+{\tilde{\mathfrak a}}_2{\tilde{\mathfrak a}'}_2 |n_2|^2\big]\big[-q_1n_2+q_2n_1\big]
+t^2 \big[-q_1n_2+q_2n_1\big]^2
\\[2mm]
&+t^2\big[{\tilde{\mathfrak a}}_1q_1'n_1+{\tilde{\mathfrak a}}_2q_2'n_2\big]\big[-{\tilde{\mathfrak a}}_1 n_1n_2+{\tilde{\mathfrak a}}_2 n_2n_1\big]
+2t^2\big[{\tilde{\mathfrak a}'}_1q_1n_1+{\tilde{\mathfrak a}'}_2q_2n_2\big]\big[-{\tilde{\mathfrak a}}_1 n_1 n_2+{\tilde{\mathfrak a}}_2 n_2n_1\big]
\\[2mm]
&+t^2\big[-m_1n_2+m_2n_1\big]\big[-{\tilde{\mathfrak a}}_1 n_1 n_2+{\tilde{\mathfrak a}}_2 n_2n_1\big]
-2t^2k \big[{\tilde{\mathfrak a}}_1{\tilde{\mathfrak a}'}_1|n_1|^2+{\tilde{\mathfrak a}}_2{\tilde{\mathfrak a}'}_2 |n_2|^2\big]\big[-{\tilde{\mathfrak a}}_1^2 n_1n_2+{\tilde{\mathfrak a}}_2^2 n_2n_1\big]
\\[2mm]
&-2t^2k\big[-{\tilde{\mathfrak a}}_1^2 n_1n_2+{\tilde{\mathfrak a}}_2^2 n_2n_1\big]
\big[-q_1n_2+q_2n_1\big]
-2t^2k \big[-{\tilde{\mathfrak a}}_1q_1n_2+{\tilde{\mathfrak a}}_2q_2n_1\big]\big[-{\tilde{\mathfrak a}}_1 n_1 n_2+{\tilde{\mathfrak a}}_2 n_2n_1\big]
\\[2mm]
&+t^2k^2\big[-{\tilde{\mathfrak a}}_1^2 n_1n_2+{\tilde{\mathfrak a}}_2^2 n_2n_1\big]^2+O(t^3).
\end{align*}

\medskip
We now make a rearrangement of the terms in $g$ and then consider the following terms.

\noindent $\clubsuit$ Term 1:
\begin{align*}
 &\big[{\tilde{\mathfrak a}}_1^2 |n_1|^2 +{\tilde{\mathfrak a}}_2^2 |n_2|^2 \big]
 - \big[-{\tilde{\mathfrak a}}_1 n_1 n_2+{\tilde{\mathfrak a}}_2 n_2n_1\big]^2
=\big[{\tilde{\mathfrak a}}_1 |n_1|^2+{\tilde{\mathfrak a}}_2|n_2|^2\big]^2.
\end{align*}

\noindent $\clubsuit$ Term 2:
\begin{align*}
 &2t\, \big[{\tilde{\mathfrak a}}_1 n_1 q_1 + {\tilde{\mathfrak a}}_2 n_2 q_2\big]
 - 2t \big[-{\tilde{\mathfrak a}}_1 n_1 n_2+{\tilde{\mathfrak a}}_2 n_1 n_2\big]\big[-q_1n_2+q_2n_1\big]
 \\[2mm]
&=2t \big[ {\tilde{\mathfrak a}}_1 n_1 q_1 \big( 1- |n_2|^2 \big)
+{\tilde{\mathfrak a}}_2 n_2q_2\big( 1- |n_1|^2 \big)
+{\tilde{\mathfrak a}}_2q_1 n_1|n_2|^2
+{\tilde{\mathfrak a}}_1q_2 |n_1|^2n_2 \big]
\\[2mm]
&=2t \big[ {\tilde{\mathfrak a}}_1 n_1 q_1 |n_1|^2
+{\tilde{\mathfrak a}}_2 n_2q_2|n_2|^2
+{\tilde{\mathfrak a}}_2q_1 n_1|n_2|^2
+{\tilde{\mathfrak a}}_1q_2 |n_1|^2n_2  \big]
\\[2mm]
&=2t \big[ {\tilde{\mathfrak a}}_1 |n_1|^2 \big( q_1n_1+q_2n_2\big)
-{\tilde{\mathfrak a}}_2 |n_2|^2\big( -q_2n_2- q_1n_1 \big) \big]
=0,
\end{align*}
due to the fact that $q\bot n$ given in Lemma \ref{derivativeofF}.

\noindent $\clubsuit$ Term 3:
\begin{align*}
 &2t \big[-{\tilde{\mathfrak a}'}_1n_1n_2+{\tilde{\mathfrak a}'}_2n_1n_2\big]\big[{\tilde{\mathfrak a}}_1^2 |n_1|^2 +{\tilde{\mathfrak a}}_2^2 |n_2|^2 \big]
 - 2t\big[-{\tilde{\mathfrak a}}_1 n_1n_2+{\tilde{\mathfrak a}}_2 n_1n_2\big]\big[{\tilde{\mathfrak a}}_1{\tilde{\mathfrak a}'}_1|n_1|^2+{\tilde{\mathfrak a}}_2{\tilde{\mathfrak a}'}_2 |n_2|^2\big]
 \\[2mm]
&=2t \big[
{\tilde{\mathfrak a}'}_2{\tilde{\mathfrak a}}_1^2 n_2|n_1|^3
-{\tilde{\mathfrak a}'}_1{\tilde{\mathfrak a}}_2^2 n_1|n_2|^3
-{\tilde{\mathfrak a}}_2 {\tilde{\mathfrak a}}_1{\tilde{\mathfrak a}'}_1n_2|n_1|^3
+{\tilde{\mathfrak a}}_1{\tilde{\mathfrak a}}_2{\tilde{\mathfrak a}'}_2 n_1 |n_2|^3 \big]
\\[2mm]
&=2t \big(
{\tilde{\mathfrak a}'}_2{\tilde{\mathfrak a}}_1 -{\tilde{\mathfrak a}'}_1{\tilde{\mathfrak a}}_2\big) \big( {\tilde{\mathfrak a}}_1 |n_1|^2+{\tilde{\mathfrak a}}_2|n_2|^2\big) n_1n_2.
\end{align*}

\noindent $\clubsuit$ Term 4:
\begin{align*}
&-2tk\big[{\tilde{\mathfrak a}}_1|n_2|^2+{\tilde{\mathfrak a}}_2|n_1|^2\big]\big[{\tilde{\mathfrak a}}_1^2 |n_1|^2 +{\tilde{\mathfrak a}}_2^2 |n_2|^2 \big]
+2tk\big[-{\tilde{\mathfrak a}}_1 n_1n_2+{\tilde{\mathfrak a}}_2 n_1n_2\big]\big[-{\tilde{\mathfrak a}}_1^2 n_1n_2+{\tilde{\mathfrak a}}_2^2 n_1n_2\big]
\\[2mm]
&=-2tk \Big[
{\tilde{\mathfrak a}}_2{\tilde{\mathfrak a}}_1^2 |n_1|^2|n_1|^2
+{\tilde{\mathfrak a}}_1{\tilde{\mathfrak a}}_2^2 |n_2|^2|n_2|^2
+{\tilde{\mathfrak a}}_2{\tilde{\mathfrak a}}_1^2 |n_1|^2|n_2|^2
+{\tilde{\mathfrak a}}_1{\tilde{\mathfrak a}}_2^2 |n_1|^2|n_2|^2  \Big]
\\[2mm]
&=-2tk{\tilde{\mathfrak a}}_1{\tilde{\mathfrak a}}_2\big[{\tilde{\mathfrak a}}_1|n_1|^2+{\tilde{\mathfrak a}}_2|n_2|^2\big].
\end{align*}

\noindent $\clubsuit$ Term 5:
\begin{align*}
&t^2\, \big[ {\tilde{\mathfrak a}}_1 n_1 m_1 + {\tilde{\mathfrak a}}_2 n_2 m_2 \big]
-t^2\big[-m_1n_2+m_2n_1\big]\big[-{\tilde{\mathfrak a}}_1 n_1 n_2+{\tilde{\mathfrak a}}_2 n_2n_1\big]
\\[2mm]
&=t^2\, \big[ {\tilde{\mathfrak a}}_1 n_1 m_1  \big( 1- |n_2|^2\big)
+{\tilde{\mathfrak a}}_2 n_2 m_2 \big( 1- |n_1|^2\big)
+{\tilde{\mathfrak a}}_1m_2n_1^2n_2
+{\tilde{\mathfrak a}}_2m_1n_1n_2^2 \big]
\\[2mm]
&=t^2\big(  m_1n_1+ m_2n_2\big)\big[ {\tilde{\mathfrak a}}_1 |n_1|^2
+{\tilde{\mathfrak a}}_2 |n_2|^2 \big].
\end{align*}

\noindent $\clubsuit$ Term 6:
\begin{align*}
t^2\big[|q_1|^2+|q_2|^2\big]-t^2 \big[-q_1n_2+q_2n_1\big]^2
&=t^2\big[|q_1|^2+|q_2|^2\big]
-t^2 \big(q_1^2|n_2|^2
-2q_2q_1n_1n_2+q_2^2|n_1|^2 \big)
\\[2mm]
&=t^2\big[|q_1|^2 \big( 1- |n_2|^2\big)
+|q_2|^2 \big( 1-|n_1|^2 \big)
+2q_2q_1n_1n_2  \big]
\\[2mm]
&=t^2\big[q_1n_1+q_2n_2\big]^2
=0,
\end{align*}
due to the fact that $q\bot n$ given in Lemma \ref{derivativeofF}.

\noindent $\clubsuit$ Term 7:
\begin{align*}
&4t^2 \big[-{\tilde{\mathfrak a}'}_1n_1n_2+{\tilde{\mathfrak a}'}_2n_1 n_2\big]\big[{\tilde{\mathfrak a}}_1 n_1 q_1 + {\tilde{\mathfrak a}}_2 n_2 q_2\big]
-2t^2\big[{\tilde{\mathfrak a}'}_1q_1n_1+{\tilde{\mathfrak a}'}_2q_2n_2\big]\big[-{\tilde{\mathfrak a}}_1 n_1 n_2+{\tilde{\mathfrak a}}_2 n_1 n_2\big]
\\[2mm]
&-2t^2 \big[{\tilde{\mathfrak a}}_1{\tilde{\mathfrak a}'}_1|n_1|^2+{\tilde{\mathfrak a}}_2{\tilde{\mathfrak a}'}_2 |n_2|^2\big]\big[-q_1n_2+q_2n_1\big]
\\[2mm]
&=t^2 \big[
2{\tilde{\mathfrak a}'}_2{\tilde{\mathfrak a}}_1  n_1n_2 \big( q_1n_1+q_2n_2 \big)
-2{\tilde{\mathfrak a}'}_1{\tilde{\mathfrak a}}_2 n_1n_2 \big(  q_2n_2+q_1n_1 \big)
\\[2mm]
&\quad
+2{\tilde{\mathfrak a}'}_2q_1n_2 {\tilde{\mathfrak a}}_1|n_1|^2
+2{\tilde{\mathfrak a}'}_2q_1n_2 {\tilde{\mathfrak a}}_2|n_2|^2
-2{\tilde{\mathfrak a}'}_1q_2n_1 {\tilde{\mathfrak a}}_2|n_2|^2
-2{\tilde{\mathfrak a}'}_1q_2n_1 {\tilde{\mathfrak a}}_1|n_1|^2  \big]
\\[2mm]
&= 2t^2\big[{\tilde{\mathfrak a}}_2'q_1n_2-{\tilde{\mathfrak a}}_1'q_2n_1\big]\big[{\tilde{\mathfrak a}}_1 |n_1|^2+{\tilde{\mathfrak a}}_2|n_2|^2\big],
\end{align*}
where we have used  the fact that $q\bot n$ given in Lemma \ref{derivativeofF}.

\noindent $\clubsuit$ Term 8:
\begin{align*}
&t^2\big[{\tilde{\mathfrak a}}_1^2 |n_1|^2
+{\tilde{\mathfrak a}}_2^2 |n_2|^2 \big]\big[-q_1'n_2+q_2'n_1\big]
-t^2\big[{\tilde{\mathfrak a}}_1q_1'n_1+{\tilde{\mathfrak a}}_2q_2'n_2\big]
\big[-{\tilde{\mathfrak a}}_1 n_1n_2+{\tilde{\mathfrak a}}_2 n_2n_1\big]
\\[2mm]
&= t^2 \big[
-{\tilde{\mathfrak a}}_2^2 |n_2|^3q_1'
+{\tilde{\mathfrak a}}_1^2 |n_1|^3q_2'
+{\tilde{\mathfrak a}}_1{\tilde{\mathfrak a}}_2q_2'n_1n_2^2
-{\tilde{\mathfrak a}}_2{\tilde{\mathfrak a}}_1q_1'n_1^2 n_2\big]
\\[2mm]
&=t^2 \big[{\tilde{\mathfrak a}}_1q_2'n_1-{\tilde{\mathfrak a}}_2q_1'n_2\big]\big[{\tilde{\mathfrak a}}_1 |n_1|^2+{\tilde{\mathfrak a}}_2|n_2|^2\big].
\end{align*}

\noindent $\clubsuit$ Term 9:
\begin{align*}
&t^2\big[{\tilde{\mathfrak a}}_1^2 |n_1|^2 +{\tilde{\mathfrak a}}_2^2 |n_2|^2 \big]\big[({\tilde{\mathfrak a}'}_1)^2|n_1|^2+({\tilde{\mathfrak a}'}_2)^2|n_2|^2\big]
-t^2 \big[{\tilde{\mathfrak a}}_1{\tilde{\mathfrak a}'}_1|n_1|^2+{\tilde{\mathfrak a}}_2{\tilde{\mathfrak a}'}_2 |n_2|^2\big]^2
\\[2mm]
&=t^2 \big[
{\tilde{\mathfrak a}}_2^2 ({\tilde{\mathfrak a}'}_1)^2|n_1|^2|n_2|^2
+{\tilde{\mathfrak a}}_1^2 ({\tilde{\mathfrak a}'}_2)^2|n_1|^2|n_2|^2
-2{\tilde{\mathfrak a}}_1{\tilde{\mathfrak a}'}_1{\tilde{\mathfrak a}}_2{\tilde{\mathfrak a}'}_2 |n_1|^2|n_2|^2 \big]
\\[2mm]
&=t^2\big[{\tilde{\mathfrak a}}_1{\tilde{\mathfrak a}}_2'-{\tilde{\mathfrak a}}_1'{\tilde{\mathfrak a}}_2\big]^2|n_1|^2|n_2|^2.
\end{align*}

\noindent $\clubsuit$ Term 10:
\begin{align*}
&-4t^2k\big[{\tilde{\mathfrak a}}_1|n_2|^2+{\tilde{\mathfrak a}}_2|n_1|^2\big]\big[{\tilde{\mathfrak a}}_1 n_1 q_1 + {\tilde{\mathfrak a}}_2 n_2 q_2\big]
+2t^2k\big[-{\tilde{\mathfrak a}}_1^2 n_1n_2+{\tilde{\mathfrak a}}_2^2 n_2n_1\big]\big[-q_1n_2+q_2n_1\big]
\\[2mm]
&+2t^2k \big[-{\tilde{\mathfrak a}}_1q_1n_2+{\tilde{\mathfrak a}}_2q_2n_1 \big]
\big[-{\tilde{\mathfrak a}}_1 n_1n_2+{\tilde{\mathfrak a}}_2 n_2n_1\big]
\\[2mm]
&=-t^2k \Big[
2{\tilde{\mathfrak a}}_2{\tilde{\mathfrak a}}_1 q_1 n_1 \big(|n_1|^2+|n_2|^2\big)
+2{\tilde{\mathfrak a}}_1{\tilde{\mathfrak a}}_2q_2 n_2\big(|n_2|^2+|n_1|^2\big)
\\[2mm]
&\qquad \qquad+2{\tilde{\mathfrak a}}_2q_1n_1 \big( {\tilde{\mathfrak a}}_1 |n_1|^2+{\tilde{\mathfrak a}}_2|n_2|^2\big)
+2{\tilde{\mathfrak a}}_1 q_2 n_2 \big( {\tilde{\mathfrak a}}_2 |n_2|^2 +{\tilde{\mathfrak a}}_1|n_1|^2\big) \Big]
\\[2mm]
&=-t^2k \Big[
2{\tilde{\mathfrak a}}_2{\tilde{\mathfrak a}}_1 \big( q_1 n_1 +q_2 n_2 \big)
+2 \big( {\tilde{\mathfrak a}}_2q_1n_1+{\tilde{\mathfrak a}}_1 q_2 n_2 \big)
\big( {\tilde{\mathfrak a}}_1 |n_1|^2+{\tilde{\mathfrak a}}_2|n_2|^2\big)  \Big]
\\[2mm]
&=-2t^2k\big( {\tilde{\mathfrak a}}_2q_1n_1+{\tilde{\mathfrak a}}_1 q_2 n_2\big)   \big( {\tilde{\mathfrak a}}_1 |n_1|^2+{\tilde{\mathfrak a}}_2|n_2|^2\big),
\end{align*}
where we have used  the fact that $q\bot n$ given in Lemma \ref{derivativeofF}.

\noindent $\clubsuit$ Term 11:
\begin{align*}
&-2t^2k\big[{\tilde{\mathfrak a}}_1^2 |n_1|^2 +{\tilde{\mathfrak a}}_2^2 |n_2|^2 \big]\big[-{\tilde{\mathfrak a}}_1{\tilde{\mathfrak a}'}_1n_1n_2+{\tilde{\mathfrak a}}_2{\tilde{\mathfrak a}'}_2n_2n_1\big]
+2t^2k\big[{\tilde{\mathfrak a}}_1{\tilde{\mathfrak a}'}_1|n_1|^2+{\tilde{\mathfrak a}}_2{\tilde{\mathfrak a}'}_2 |n_2|^2\big]\big[-{\tilde{\mathfrak a}}_1^2 n_1n_2+{\tilde{\mathfrak a}}_2^2 n_2n_1\big]
\\[2mm]
&=-2t^2k \big[
-{\tilde{\mathfrak a}}_2^2 {\tilde{\mathfrak a}}_1{\tilde{\mathfrak a}'}_1n_1|n_2|^3
+{\tilde{\mathfrak a}}_1^2 {\tilde{\mathfrak a}}_2{\tilde{\mathfrak a}'}_2|n_1|^3n_2
+{\tilde{\mathfrak a}}_2{\tilde{\mathfrak a}'}_2 {\tilde{\mathfrak a}}_1^2 n_1|n_2|^3
-{\tilde{\mathfrak a}}_1{\tilde{\mathfrak a}'}_1{\tilde{\mathfrak a}}_2^2|n_1|^3 n_2   \big]
\\[2mm]
&= 2t^2k \big[{\tilde{\mathfrak a}}_2^2 {\tilde{\mathfrak a}}_1{\tilde{\mathfrak a}'}_1-{\tilde{\mathfrak a}}_1^2 {\tilde{\mathfrak a}}_2{\tilde{\mathfrak a}'}_2\big]n_1n_2.
\end{align*}

\noindent $\clubsuit$ Term 12:
\begin{align*}
&t^2k^2\big[{\tilde{\mathfrak a}}_1^2 |n_1|^2 +{\tilde{\mathfrak a}}_2^2 |n_2|^2 \big] \big[{\tilde{\mathfrak a}}_1^2 |n_2|^2+{\tilde{\mathfrak a}}_2^2 |n_1|^2\big]
-t^2k^2\big[-{\tilde{\mathfrak a}}_1^2 n_1n_2+{\tilde{\mathfrak a}}_2^2 n_2n_1\big]^2
\\[2mm]
&=t^2k^2 \big[ {\tilde{\mathfrak a}}_1^2 {\tilde{\mathfrak a}}_2^2 |n_1|^2|n_1|^2
+{\tilde{\mathfrak a}}_1^2 {\tilde{\mathfrak a}}_2^2 |n_2|^2|n_2|^2
+2{\tilde{\mathfrak a}}_1^2 {\tilde{\mathfrak a}}_2^2 n_1^2n_2^2\big]
\\[2mm]
&=t^2k^2{\tilde{\mathfrak a}}_1^2{\tilde{\mathfrak a}}_2^2 \big( n_1^2+n_2^2 \big)^2 =t^2k^2{\tilde{\mathfrak a}}_1^2{\tilde{\mathfrak a}}_2^2.
\end{align*}

\medskip
Therefore, we obtain that
\begin{align}\label{mathfrakh1}
g\, =\, &\mbox{det}(g_{ij})=g_{11}g_{22}-g_{12}g_{12}
\nonumber\\[2mm]
=&\big[{\tilde{\mathfrak a}}_1 |n_1|^2+{\tilde{\mathfrak a}}_2|n_2|^2\big]^2
+
t \Big\{ 2 \big({\tilde{\mathfrak a}'}_2{\tilde{\mathfrak a}}_1 -{\tilde{\mathfrak a}'}_1{\tilde{\mathfrak a}}_2\big)\big[{\tilde{\mathfrak a}}_1 |n_1|^2+{\tilde{\mathfrak a}}_2|n_2|^2\big] n_1n_2
-2k{\tilde{\mathfrak a}}_1{\tilde{\mathfrak a}}_2\big[{\tilde{\mathfrak a}}_1n_1^2+{\tilde{\mathfrak a}}_2n_2^2\big]
\Big\}
\nonumber\\[2mm]
&+
t^2\Big\{ 2\big[{\tilde{\mathfrak a}}_2'q_1n_2-{\tilde{\mathfrak a}}_1'q_2n_1\big]\big[{\tilde{\mathfrak a}}_1 |n_1|^2+{\tilde{\mathfrak a}}_2|n_2|^2\big]
+ \big[{\tilde{\mathfrak a}}_1q_2'n_1-{\tilde{\mathfrak a}}_2q_1'n_2\big]\big[{\tilde{\mathfrak a}}_1 |n_1|^2+{\tilde{\mathfrak a}}_2|n_2|^2\big]
\nonumber\\[2mm]
&\qquad
+\big(m_1n_1+m_2n_2\big)\big[ {\tilde{\mathfrak a}}_1 |n_1|^2 +{\tilde{\mathfrak a}}_2 |n_2|^2\big]
+\big[{\tilde{\mathfrak a}}_1{\tilde{\mathfrak a}}_2'-{\tilde{\mathfrak a}}_1'{\tilde{\mathfrak a}}_2\big]^2|n_1|^2|n_2|^2
\nonumber\\[2mm]
&\qquad
+2 k \big[{\tilde{\mathfrak a}}_2^2 {\tilde{\mathfrak a}}_1{\tilde{\mathfrak a}'}_1-{\tilde{\mathfrak a}}_1^2 {\tilde{\mathfrak a}}_2{\tilde{\mathfrak a}'}_2\big]n_1n_2
-2 k\big( {\tilde{\mathfrak a}}_2q_1n_1+{\tilde{\mathfrak a}}_1 q_2n_2\big)\big[{\tilde{\mathfrak a}}_1 |n_1|^2+{\tilde{\mathfrak a}}_2|n_2|^2\big]
+ k^2{\tilde{\mathfrak a}}_1^2{\tilde{\mathfrak a}}_2^2
\Big\}
+O(t^3)
\nonumber\\[2mm]
\, \equiv\, &{\mathfrak h}_1(\theta)
\, +\, t{\mathfrak h}_2(\theta)
\, +\, t^2{\mathfrak h}_3(\theta)+O(t^3).
\end{align}

\medskip
We now compute the inverse of the metric matrix. By the formula
$$
\big[1+at+bt^2+O(t^3)\big]^{-1}\, =\, 1-at+(a^2-b)t^2+O(t^3),
$$
we have
\begin{align}
\frac{1}{g}\, =\, &\frac{1}{{\mathfrak h}_1}
-\frac{{\mathfrak h}_2}{{\mathfrak h}_1^2}t
+\Big[\frac{{\mathfrak h}_2^2}{{\mathfrak h}_1^3}
-\frac{{\mathfrak h}_3}{{\mathfrak h}_1^2}\Big]t^2
+O(t^3)
\nonumber\\[2mm]
\, \equiv\, &\frac{1}{{\mathfrak h}_1}
+{\mathfrak g}_1t
+{\mathfrak g}_2t^2
+O(t^3).
\label{g-inverse}
\end{align}
Note that if the function has the following asymptotic expansion
$$f(s)\, =\, 1+as+bs^2+O(s^3), \quad f(0)\, =\, 1,
$$
then for $s$ close to zero
$$
\sqrt{f(s)}\, =\, 1+\frac{a}{2}s+\frac12\left(b-\frac14a^2\right)s^2+O(s^3).
$$
We can get the following formulas
\begin{equation}
 \sqrt{g}\, =\,
\sqrt{{\mathfrak h}_1}
+\frac{1}{2}\frac{{\mathfrak h}_2}{\sqrt{{\mathfrak h}_1}} t
+\Big[\frac{1}{2}\frac{{\mathfrak h}_3}{\sqrt{{\mathfrak h}_1}}t^2
-\frac{1}{8}\frac{{\mathfrak h}_2^2}{(\sqrt{{\mathfrak h}_1})^3}\Big]
+O(t^3), \label{g-inverse1}
\end{equation}
and also
\begin{align}
\frac{1}{\sqrt{g}}&\, =\,
\frac{1}{\sqrt{{\mathfrak h}_1}}
-\frac{1}{2}\frac{{\mathfrak h}_2}{(\sqrt{{\mathfrak h}_1})^3} t
+\Big[\frac{3}{8}\frac{{\mathfrak h}_2^2}{(\sqrt{{\mathfrak h}_1})^5}
-\frac{1}{2}\frac{{\mathfrak h}_3}{(\sqrt{{\mathfrak h}_1})^3}\Big]t^2+O(t^3)
\nonumber\\[2mm]
&\, \equiv\, \frac{1}{\sqrt{{\mathfrak h}_1}}
+{\mathfrak r}_1t
+{\mathfrak r}_2t^2+O(t^3).
\label{g-inverse2}
\end{align}
Whence
\begin{align}
g^{12}\, =\, &-\frac{g_{12}}{ g}
\nonumber\\
\, =\, &\frac{1}{{\mathfrak h}_1}\big[{\tilde{\mathfrak a}}_1-{\tilde{\mathfrak a}}_2\big]n_1n_2
\nonumber\\[2mm]
&+t\Big\{-\frac{1}{{\mathfrak h}_1}\big[{\tilde{\mathfrak a}}_1{\tilde{\mathfrak a}'}_1|n_1|^2+{\tilde{\mathfrak a}}_2{\tilde{\mathfrak a}'}_2 |n_2|^2\big]
+k\frac{1}{{\mathfrak h}_1}\big[-{\tilde{\mathfrak a}}_1^2+{\tilde{\mathfrak a}}_2^2\big]n_1n_2
\nonumber\\[2mm]
&\qquad-\frac{1}{{\mathfrak h}_1}\big[-q_1n_2+q_2n_1\big]
+{\mathfrak g}_1\big[{\tilde{\mathfrak a}}_1-{\tilde{\mathfrak a}}_2\big]n_1n_2\Big\}
\nonumber\\[2mm]
&
+t^2\Big\{-\frac{1}{2}\frac{1}{{\mathfrak h}_1}\big[{\tilde{\mathfrak a}}_1q_1'n_1+{\tilde{\mathfrak a}}_2q_2'n_2\big]
-\frac{1}{{\mathfrak h}_1}\big[{\tilde{\mathfrak a}'}_1q_1n_1+{\tilde{\mathfrak a}'}_2q_2n_2\big]
+k\frac{1}{{\mathfrak h}_1}\big[-{\tilde{\mathfrak a}}_1q_1n_2+{\tilde{\mathfrak a}}_2q_2n_1\big]
\nonumber\\[2mm]
&\qquad+\frac{1}{2}\frac{1}{{\mathfrak h}_1}\big[-m_1n_2+m_2n_1\big]
+{\mathfrak g}_1\big[{\tilde{\mathfrak a}}_1{\tilde{\mathfrak a}'}_1|n_1|^2+{\tilde{\mathfrak a}}_2{\tilde{\mathfrak a}'}_2 |n_2|^2\big]
\nonumber\\[2mm]
&\qquad+k{\mathfrak g}_1\big[-{\tilde{\mathfrak a}}_1^2+{\tilde{\mathfrak a}}_2^2\big]n_1n_2
-{\mathfrak g}_1\big[-q_1n_2+q_2n_1\big]
+{\mathfrak g}_2\big[{\tilde{\mathfrak a}}_1-{\tilde{\mathfrak a}}_2\big]n_1n_2\Big\}+O(t^3)
\nonumber\\[2mm]
\, \equiv\, &{\mathfrak g}_3+{\mathfrak g}_4t+{\mathfrak g}_5t^2+O(t^3),
\label{g12inverse}
\end{align}
and
\begin{align}
g^{22}\, =\, \frac{g_{11}}{g}
\, =\, &\frac{1}{{\mathfrak h}_1}\big[{\tilde{\mathfrak a}}_1^2 |n_1|^2 +{\tilde{\mathfrak a}}_2^2 |n_2|^2 \big]
+t\Big\{2\frac{1}{{\mathfrak h}_1}\big[{\tilde{\mathfrak a}}_1 n_1 q_1 + {\tilde{\mathfrak a}}_2 n_2 q_2\big]
+{\mathfrak g}_1\big[{\tilde{\mathfrak a}}_1^2 |n_1|^2 +{\tilde{\mathfrak a}}_2^2 |n_2|^2\big]\Big\}
\nonumber\\[2mm]
&
+t^2\Big\{\frac{1}{{\mathfrak h}_1}\big[ {\tilde{\mathfrak a}}_1 n_1 m_1+{\tilde{\mathfrak a}}_2 n_2 m_2 \big]
+\frac{1}{{\mathfrak h}_1}\big[|q_1|^2+|q_2|^2\big]
\nonumber\\[2mm]
&\qquad \quad+2{\mathfrak g}_1\big[{\tilde{\mathfrak a}}_1 n_1 q_1 + {\tilde{\mathfrak a}}_2 n_2 q_2\big]
+{\mathfrak g}_2\big[{\tilde{\mathfrak a}}_1^2 |n_1|^2 +{\tilde{\mathfrak a}}_2^2 |n_2|^2\big]\Big\}+O(t^3)
\nonumber\\[2mm]
\, \equiv\, &\, {\mathfrak g}_6+{\mathfrak g}_7t+{\mathfrak g}_8t^2+O(t^3).
\label{g22inverse}
\end{align}
Similar asymptotic expression holds for the term $g^{11}=g_{22}/g$.

\subsection{Local forms of the differential operators in \eqref{originalproblem}}\
\label{section2.3}

In this section, we are devoted to presenting the expressions of the differential operators ${\rm div}\big(\nabla_{\mathfrak a(y)}u\big)$
and $\nabla_{{\mathfrak a}(y)}u\cdot\nu$ in problem \eqref{originalproblem}.

\medskip
\noindent{\bf{Part 1: the operator ${\rm div}\big(\nabla_{\mathfrak a(y)}u\big)$}}

We recall the relation of \eqref{Fermicoordinates-modified}
and then have
$$
{\mathrm d}y_1=\frac{\partial F_1}{\partial t}{\mathrm d}t\, +\, \frac{\partial F_1}{\partial \theta}{\mathrm d}\theta,
\qquad
{\mathrm d}y_2=\frac{\partial F_2}{\partial t}{\mathrm d}t\, +\, \frac{\partial F_2}{\partial \theta}{\mathrm d}\theta.
$$
This implies that
\begin{align}
{\mathrm d}t=\frac{1}{\, \frac{\partial F_1}{\partial t}\frac{\partial F_2}{\partial \theta}\,-\, \frac{\partial F_2}{\partial t}\frac{\partial F_1}{\partial \theta}\, }
\left[\frac{\partial F_2}{\partial \theta}{\mathrm d}y_1\,-\, \frac{\partial F_1}{\partial \theta}{\mathrm d}y_2\right],
\label{dt}
\end{align}
\begin{align}
{\mathrm d}\theta=\frac{1}{\, \frac{\partial F_1}{\partial t}\frac{\partial F_2}{\partial \theta}\,-\, \frac{\partial F_2}{\partial t}\frac{\partial F_1}{\partial \theta}\, }
\left[\frac{\partial F_1}{\partial t}{\mathrm d}y_2\,-\, \frac{\partial F_2}{\partial t}{\mathrm d}y_1\right].
\label{dtheta}
\end{align}
On the other hand, there hold
\begin{align}
\frac{\partial u}{\partial t}=\frac{\partial u}{\partial y_1}\frac{\partial F_1}{\partial t}\, +\, \frac{\partial u}{\partial y_2}\frac{\partial F_2}{\partial t},
\qquad
\frac{\partial u}{\partial \theta}=\frac{\partial u}{\partial y_1}\frac{\partial F_1}{\partial \theta}\, +\, \frac{\partial u}{\partial y_2}\frac{\partial F_2}{\partial \theta},
\end{align}
which give that
\begin{align}
\frac{\partial u}{\partial y_1}
=\frac{1}{\, \frac{\partial F_1}{\partial t}\frac{\partial F_2}{\partial \theta}\,-\, \frac{\partial F_2}{\partial t}\frac{\partial F_1}{\partial \theta}\, }
\left[ \frac{\partial u}{\partial t}\frac{\partial F_2}{\partial \theta}\,-\, \frac{\partial u}{\partial \theta}\frac{\partial F_2}{\partial t} \right],
\label{partialuy1}
\end{align}
\begin{align}
\frac{\partial u}{\partial y_2}
=\frac{1}{\, \frac{\partial F_1}{\partial t}\frac{\partial F_2}{\partial \theta}\,-\, \frac{\partial F_2}{\partial t}\frac{\partial F_1}{\partial \theta}\, }
\left[ \frac{\partial u}{\partial \theta}\frac{\partial F_1}{\partial t}\,-\, \frac{\partial u}{\partial t}\frac{\partial F_1}{\partial \theta} \right].
\label{partialuy2}
\end{align}

\medskip
We now compute
\begin{align}\label{udiri}
\nabla_{{\mathfrak a}(y)}u
\, =\, &{\mathfrak a}_1\frac{\partial u}{\partial y_1}\frac{\partial }{\partial y_1}\, +\, {\mathfrak a}_2\frac{\partial u}{\partial y_2}\frac{\partial }{\partial y_2}
\nonumber\\[2mm]
\, =\, &{\mathfrak a}_1\frac{\partial u}{\partial y_1}\left[ \frac{\partial t}{\partial y_1}\frac{\partial }{\partial t}\, +\, \frac{\partial \theta}{\partial y_1}\frac{\partial }{\partial \theta}\right]
\, +\,
{\mathfrak a}_2\frac{\partial u}{\partial y_2}\left[ \frac{\partial t}{\partial y_2}\frac{\partial }{\partial t}\, +\, \frac{\partial \theta}{\partial y_2}\frac{\partial }{\partial \theta}\right]
\nonumber\\[2mm]
\, =\, &\left[ {\mathfrak a}_1\frac{\partial t}{\partial y_1}\frac{\partial u}{\partial y_1}\, +\, {\mathfrak a}_2\frac{\partial t}{\partial y_2}\frac{\partial u}{\partial y_2}\right]
\frac{\partial }{\partial t}
\, +\,
\left[ {\mathfrak a}_1\frac{\partial \theta}{\partial y_1}\frac{\partial u}{\partial y_1}\, +\, {\mathfrak a}_2\frac{\partial \theta}{\partial y_2}\frac{\partial u}{\partial y_2}\right]
\frac{\partial }{\partial \theta}.
\end{align}
By substituting \eqref{dt}, \eqref{dtheta}, \eqref{partialuy1} and \eqref{partialuy2} in \eqref{udiri}, we obtain
\begin{align}
\nabla_{{\mathfrak a}(y)}u
\, =\, &\frac{1}{g}\left[{\mathfrak a}_1\frac{\partial F_2}{\partial \theta}\frac{\partial F_2}{\partial \theta}\, +\, {\mathfrak a}_2\frac{\partial F_1}{\partial \theta}\frac{\partial F_1}{\partial \theta}\right]\frac{\partial u}{\partial t}\frac{\partial }{\partial t}
\,-\,
\frac{1}{g}\left[{\mathfrak a}_1\frac{\partial F_2}{\partial \theta}\frac{\partial F_2}{\partial t}\, +\, {\mathfrak a}_2\frac{\partial F_1}{\partial \theta}\frac{\partial F_1}{\partial t}\right]\frac{\partial u}{\partial \theta}\frac{\partial }{\partial t}
\nonumber\\[2mm]
&\,-\,
\frac{1}{g}\left[{\mathfrak a}_1\frac{\partial F_2}{\partial t}\frac{\partial F_2}{\partial \theta}\, +\, {\mathfrak a}_2\frac{\partial F_1}{\partial t}\frac{\partial F_1}{\partial \theta}\right]\frac{\partial u}{\partial t}\frac{\partial }{\partial \theta}
\, +\,
\frac{1}{g}\left[{\mathfrak a}_1\frac{\partial F_2}{\partial t}\frac{\partial F_2}{\partial t}\, +\, {\mathfrak a}_2\frac{\partial F_1}{\partial t}\frac{\partial F_1}{\partial t}\right]\frac{\partial u}{\partial \theta}\frac{\partial }{\partial \theta}.
\end{align}
Furthermore,  by setting
\begin{align}
{\tilde g}^{11}\, =\, {\mathfrak a}_1\frac{\partial F_2}{\partial \theta}\frac{\partial F_2}{\partial \theta}\, +\, {\mathfrak a}_2\frac{\partial F_1}{\partial \theta}\frac{\partial F_1}{\partial \theta},
\\
{\tilde g}^{21}\, =\,{\tilde g}^{12}\, =\, {\mathfrak a}_1\frac{\partial F_2}{\partial \theta}\frac{\partial F_2}{\partial t}\, +\, {\mathfrak a}_2\frac{\partial F_1}{\partial \theta}\frac{\partial F_1}{\partial t},
\\
{\tilde g}^{22}\, =\, {\mathfrak a}_1\frac{\partial F_2}{\partial t}\frac{\partial F_2}{\partial t}
\, +\, {\mathfrak a}_2\frac{\partial F_1}{\partial t}\frac{\partial F_1}{\partial t},
\end{align}
 the definition of ${\rm div}$ operator will give that
\begin{align}
{\rm div}\big(\nabla_{{\mathfrak a}(y)}  u\big)
\, =\, &\frac{1}{\sqrt{g}}\frac{\partial}{\partial t}
\left[
\frac{1}{\sqrt{g}}{\tilde g}^{11}\frac{\partial u}{\partial t}
\right]
\,-\,
\frac{1}{\sqrt{g}}\frac{\partial}{\partial t}
\left[
\frac{1}{\sqrt{g}}{\tilde g}^{12}\frac{\partial u}{\partial \theta}
\right]
\nonumber\\[2mm]
&\,-\,
\frac{1}{\sqrt{g}}\frac{\partial}{\partial \theta}
\left[
\frac{1}{\sqrt{g}}{\tilde g}^{21}\frac{\partial u}{\partial t}
\right]
\, +\,
\frac{1}{\sqrt{g}}\frac{\partial}{\partial \theta}
\left[
\frac{1}{\sqrt{g}}{\tilde g}^{22}\frac{\partial u}{\partial \theta}
\right].
\label{divoperator}
\end{align}

\medskip
Here are the computations of all coefficients in \eqref{divoperator}.
By using the Taylor expansion
\begin{equation}\label{expressionofai}
 {\mathfrak a}_i(t, \theta)= {\mathfrak a}_i(0, \theta)+t \partial_{t} {\mathfrak a}_i(0, \theta)
 +\frac{t^2}{2}\partial_{tt} {\mathfrak a}_i (0, \theta)+O(t^3), ~~\forall\, \theta\in [0, 1], ~t\in (-\delta_0, \delta_0),
\end{equation}
and recalling \eqref{Fit-Fit}-\eqref{Fitheta-Fitheta},
we obtain
\begin{align}
{\tilde g}^{11}\, =\, &{\mathfrak a}_1\frac{\partial F_2}{\partial \theta}\frac{\partial F_2}{\partial \theta}\, +\, {\mathfrak a}_2\frac{\partial F_1}{\partial \theta}\frac{\partial F_1}{\partial \theta}
\nonumber\\[2mm]
\, =\, &{\mathfrak a}_1\Big[|n_1|^2
+2t {\tilde{\mathfrak a}'}_2 n_1 n_2
-2t k{\tilde{\mathfrak a}}_2 |n_1|^2
+t^2 q_2'n_1
+t^2({\tilde{\mathfrak a}'}_2)^2|n_2|^2
-2t^2k{\tilde{\mathfrak a}}_2 {\tilde{\mathfrak a}'}_2 n_2n_1
+t^2k^2 {\tilde{\mathfrak a}}_2^2 |n_1|^2\Big]
\nonumber\\[2mm]
&+{\mathfrak a}_2\Big[ |n_2|^2
-2t {\tilde{\mathfrak a}'}_1 n_2 n_1
-2t k{\tilde{\mathfrak a}}_1 |n_2|^2
-t^2 q_1'n_2
+t^2({\tilde{\mathfrak a}'}_1)^2|n_1|^2
+2t^2k{\tilde{\mathfrak a}}_1 {\tilde{\mathfrak a}'}_1 n_1n_2
+t^2k^2 {\tilde{\mathfrak a}}_1^2 |n_2|^2\Big]+O(t^3)
\nonumber\\[2mm]
\, =\, &
\Big[{\mathfrak a}_1(0, \theta) |n_1|^2 +{\mathfrak a}_2(0, \theta) |n_2|^2\Big]
\nonumber\\[2mm]
&+t\Big\{
2\big[{\mathfrak a}_1(0, \theta) {\tilde{\mathfrak a}'}_2 n_1 n_2 -{\mathfrak a}_2(0, \theta){\tilde{\mathfrak a}'}_1 n_2 n_1\big]
-2k\big[{\mathfrak a}_1(0, \theta) {\tilde{\mathfrak a}}_2|n_1|^2 +{\mathfrak a}_2(0, \theta){\tilde{\mathfrak a}}_1 |n_2|^2\big]
\nonumber\\[2mm]
&\qquad
+\big[\partial_{t} {\mathfrak a}_1(0, \theta) |n_1|^2 +\partial_{t} {\mathfrak a}_2(0, \theta) |n_2|^2\big]
\Big\}
\nonumber\\[2mm]
&+t^2\Big\{
\big[{\mathfrak a}_1(0, \theta) n_1q_2' -{\mathfrak a}_2(0, \theta) n_2q_1'\big]
+ \big[{\mathfrak a}_1(0, \theta) ({\tilde{\mathfrak a}'}_2)^2|n_2|^2+{\mathfrak a}_2(0, \theta) ({\tilde{\mathfrak a}'}_1)^2|n_1|^2\big]
\nonumber\\[2mm]
&\qquad
-2 k\big[{\mathfrak a}_1(0, \theta){\tilde{\mathfrak a}}_2 {\tilde{\mathfrak a}'}_2 n_2n_1-{\mathfrak a}_2(0, \theta){\tilde{\mathfrak a}}_1 {\tilde{\mathfrak a}'}_1n_1n_2\big]
+ k^2\big[{\mathfrak a}_1(0, \theta) {\tilde{\mathfrak a}}_2|n_1|^2 +{\mathfrak a}_2(0, \theta) {\tilde{\mathfrak a}}_1|n_2|^2\big]
\nonumber\\[2mm]
&\qquad
+2\big[\partial_{t} {\mathfrak a}_1(0, \theta) {\tilde{\mathfrak a}'}_2n_1n_2   -\partial_{t} {\mathfrak a}_2(0, \theta){\tilde{\mathfrak a}'}_1n_2 n_1\big]
-2 k\big[\partial_{t} {\mathfrak a}_1(0, \theta){\tilde{\mathfrak a}}_2 |n_1|^2  +\partial_{t} {\mathfrak a}_2(0, \theta){\tilde{\mathfrak a}}_1 |n_2|^2\big]
\nonumber\\[2mm]
&\qquad
+\frac{1}{2} \big[\partial_{tt} {\mathfrak a}_1(0, \theta) |n_1|^2  +\partial_{tt} {\mathfrak a}_2(0, \theta) |n_2|^2\big]
\Big\}
+ O(t^3)
\nonumber\\[2mm]
\, \equiv\, &{\mathfrak f}_0
+t{\mathfrak f}_1
+t^2{\mathfrak f}_2+ O(t^3).
\label{m22}
\end{align}
Similarly, there hold
\begin{align}
 {\tilde g}^{12}\, =\, &{\mathfrak a}_1\frac{\partial F_2}{\partial \theta}\frac{\partial F_2}{\partial t}\, +\, {\mathfrak a}_2\frac{\partial F_1}{\partial \theta}\frac{\partial F_1}{\partial t}
\nonumber\\[2mm]
\, =\, &{\mathfrak a}_1\Big[{\tilde{\mathfrak a}}_2 n_2n_1
+t {\tilde{\mathfrak a}}_2{\tilde{\mathfrak a}'}_2|n_2|^2
-tk {\tilde{\mathfrak a}}_2^2 n_2n_1
+t q_2n_1
+\frac{t^2}{2}{\tilde{\mathfrak a}}_2q_2'n_2
\Big]
\nonumber\\[2mm]
&+{\mathfrak a}_2\Big[-{\tilde{\mathfrak a}}_1 n_1 n_2
+t {\tilde{\mathfrak a}}_1{\tilde{\mathfrak a}'}_1 |n_1|^2
+tk {\tilde{\mathfrak a}}_1^2 n_1n_2
-t q_1n_2
\Big]+O(t^2)
\nonumber\\[2mm]
\, =\, &
t\Big\{
\big[{\mathfrak a}_1(0, \theta){\tilde{\mathfrak a}}_2{\tilde{\mathfrak a}'}_2|n_2|^2+{\mathfrak a}_2(0, \theta){\tilde{\mathfrak a}}_1{\tilde{\mathfrak a}'}_1|n_1|^2\big]
-k\big[{\mathfrak a}_1(0, \theta){\tilde{\mathfrak a}}_2^2 n_1n_2 -{\mathfrak a}_2(0, \theta){\tilde{\mathfrak a}}_1^2 n_2n_1 \big]
\nonumber\\[2mm]
&\quad
+ \big[{\mathfrak a}_1(0, \theta)n_1q_2-{\mathfrak a}_2(0, \theta)n_2q_1 \big]
+ \big[\partial_{t} {\mathfrak a}_1(0, \theta){\tilde{\mathfrak a}}_2n_1n_2 -\partial_{t} {\mathfrak a}_2(0, \theta){\tilde{\mathfrak a}}_1n_2n_1 \big]
\Big\}
+O(t^2)
\nonumber\\
\, \equiv\, &t{\mathfrak l}_1
+O(t^2),
\label{m12}
\end{align}
and
\begin{align}
{\tilde g}^{22}\, =\, &{\mathfrak a}_1\frac{\partial F_2}{\partial t}\frac{\partial F_2}{\partial t}
\, +\, {\mathfrak a}_2\frac{\partial F_1}{\partial t}\frac{\partial F_1}{\partial t}
\nonumber\\[2mm]
\, =\, &{\mathfrak a}_1\Big[{\tilde{\mathfrak a}}_2^2 |n_2|^2
+2t\, {\tilde{\mathfrak a}}_2 n_2 q_2
+t^2\, {\tilde{\mathfrak a}}_2 n_2 m_2
+t^2|q_2|^2\Big]
\nonumber\\[2mm]
&+{\mathfrak a}_2\Big[{\tilde{\mathfrak a}}_1^2 |n_1|^2
+2t\, {\tilde{\mathfrak a}}_1n_1q_1
+t^2\, {\tilde{\mathfrak a}}_1n_1m_1
+t^2|q_1|^2\Big]+O(t^3)
\nonumber\\[2mm]
\, =\, &\Big[{\mathfrak a}_1(0, \theta){\tilde{\mathfrak a}}_2^2 |n_2|^2+{\mathfrak a}_2(0, \theta){\tilde{\mathfrak a}}_1^2 |n_1|^2\Big]
+t\Big\{
2\, {\mathfrak a}_1(0, \theta){\tilde{\mathfrak a}}_2n_2q_2
+2\, {\mathfrak a}_2(0, \theta){\tilde{\mathfrak a}}_1n_1q_1
\nonumber\\[2mm]
&+\partial_{t} {\mathfrak a}_1 (0, \theta){\tilde{\mathfrak a}}_2^2 |n_2|^2
+\partial_{t} {\mathfrak a}_2 (0, \theta){\tilde{\mathfrak a}}_1^2 |n_1|^2
\Big\}
+ O(t^3)
\nonumber\\[2mm]
\, \equiv\, &{\mathfrak w}_0
+t{\mathfrak w}_1
+ O(t^2).
\label{m11}
\end{align}

By recalling \eqref{g-inverse}-\eqref{g-inverse2} and \eqref{m22}-\eqref{m11}, we can obtain that
\begin{align*}
 \frac{{\tilde g}^{11}}{g}
\, =\, &\Big[\frac{1}{{\mathfrak h}_1}
+{\mathfrak g}_1t
+{\mathfrak g}_2t^2
+O(t^3)\Big]\times
\Big[{\mathfrak f}_0
+t{\mathfrak f}_1
+t^2{\mathfrak f}_2+O(t^3)\Big]
\\[2mm]
\, =\, & \frac{1}{{\mathfrak h}_1}{\mathfrak f}_0
+t\frac{1}{{\mathfrak h}_1}{\mathfrak f}_1
+t{\mathfrak g}_1{\mathfrak f}_0
+t^2\frac{1}{{\mathfrak h}_1}{\mathfrak f}_2
+t^2{\mathfrak g}_1{\mathfrak f}_1
+t^2{\mathfrak g}_2{\mathfrak f}_0
+O(t^3),
\\[2mm]
\frac{{\tilde g}^{11}}{\sqrt{g}}
\, =\, &\Big[\frac{1}{\sqrt{{\mathfrak h}_1}}
+{\mathfrak r}_1t
+{\mathfrak r}_2t^2+O(t^3)\Big]
\times\Big[{\mathfrak f}_0
+t{\mathfrak f}_1
+t^2{\mathfrak f}_2+O(t^3)\Big]
\\
\, =\, &\frac{1}{\sqrt{{\mathfrak h}_1}}{\mathfrak f}_0
+t\frac{1}{\sqrt{{\mathfrak h}_1}}{\mathfrak f}_1
+t{\mathfrak r}_1{\mathfrak f}_0
+t^2\frac{1}{\sqrt{{\mathfrak h}_1}}{\mathfrak f}_2
+t^2{\mathfrak r}_1{\mathfrak f}_1
+t^2{\mathfrak r}_2{\mathfrak f}_0
+O(t^3),
\\[2mm]
 \frac{{\tilde g}^{22}}{g}
\, =\, &\Big[\frac{1}{{\mathfrak h}_1}
+{\mathfrak g}_1t
+{\mathfrak g}_2t^2+O(t^3)\Big]
\times \Big[{\mathfrak w}_0
+t{\mathfrak w}_1
+ O(t^2)\Big]
\\[2mm]
\, =\, &\frac{1}{{\mathfrak h}_1}{\mathfrak w}_0
+t\frac{1}{{\mathfrak h}_1}{\mathfrak w}_1
+t{\mathfrak g}_1{\mathfrak w}_0
+O(t^2),
\\[2mm]
\frac{{\tilde g}^{22}}{\sqrt{g}}
\, =\, &\Big[\frac{1}{\sqrt{{\mathfrak h}_1}}
+{\mathfrak r}_1t
+{\mathfrak r}_2t^2+O(t^3)\Big]
\times\big[{\mathfrak w}_0
+t{\mathfrak w}_1
+t^2{\mathfrak w}_2+O(t^3)\big]
\\
\, =\, &\frac{1}{\sqrt{{\mathfrak h}_1}}{\mathfrak w}_0
+t\frac{1}{\sqrt{{\mathfrak h}_1}}{\mathfrak w}_1
+t{\mathfrak r}_1{\mathfrak w}_0
+O(t^2),
\\[2mm]
\frac{{\tilde g}^{12}}{g}
\, =\, &\Big[\frac{1}{{\mathfrak h}_1}
+{\mathfrak g}_1t
+{\mathfrak g}_2t^2
+O(t^3)\Big]
\times \big[t{\mathfrak l}_1
+O(t^2)\big]
\, =\, t\frac{1}{{\mathfrak h}_1}{\mathfrak l}_1
+ O(t^2),
\\[2mm]
\frac{{\tilde g}^{12}}{\sqrt{g}}
\, =\, &\Big[\frac{1}{\sqrt{{\mathfrak h}_1}}
+{\mathfrak r}_1t
+{\mathfrak r}_2t^2
+O(t^3)\Big]
\times \Big[t{\mathfrak l}_1
+O(t^2)\Big]
\, =\, t\frac{1}{\sqrt{{\mathfrak h}_1}}{\mathfrak l}_1
+ O(t^2).
\end{align*}
Then, there hold
\begin{align*}
\partial_t \Big[\frac{1}{\sqrt{g}} {\tilde g}^{12}  \Big]
=\frac{1}{\sqrt{{\mathfrak h}_1}}{\mathfrak l}_1
+ O(t),
\qquad
 \partial_{\theta} \Big[\frac{1}{\sqrt{g}} {\tilde g}^{12} \Big]
=t \Big[\frac{1}{\sqrt{{\mathfrak h}_1}}{\mathfrak l}_1\Big]'
+ O(t^2),
\end{align*}
and
\begin{align*}
\partial_t \Big[\frac{1}{\sqrt{g}} {\tilde g}^{11} \Big]
\, =\, &\frac{1}{\sqrt{{\mathfrak h}_1}}{\mathfrak f}_1
+{\mathfrak r}_1{\mathfrak f}_0
+2t\frac{1}{\sqrt{{\mathfrak h}_1}}{\mathfrak f}_2
+2t{\mathfrak r}_1{\mathfrak f}_1
+2t{\mathfrak r}_2{\mathfrak f}_0
+O(t^2),
\\[2mm]
\partial_{\theta} \Big[\frac{1}{\sqrt{g}} {\tilde g}^{22}  \Big]
\, =\, &\Big[ \frac{1}{\sqrt{{\mathfrak h}_1}}{\mathfrak w}_0\Big]'
+t \Big[\frac{1}{\sqrt{{\mathfrak h}_1}}{\mathfrak w}_1\Big]'
+t \big[{\mathfrak r}_1{\mathfrak w}_0\big]'
+ O(t^2).
\end{align*}
Combining the expression of $(\sqrt{g})^{-1}$ as in \eqref{g-inverse2}, we can obtain that

\begin{align*}
\frac{1}{\sqrt{g}}\partial_t \Big[\frac{1}{\sqrt{g}} {\tilde g}^{12}  \Big]
\, =\, &\Big[\frac{1}{\sqrt{{\mathfrak h}_1}}
+{\mathfrak r}_1t
+{\mathfrak r}_2t^2
+O(t^3)\Big]
\times
\Big[\frac{1}{\sqrt{{\mathfrak h}_1}}{\mathfrak l}_1 + O(t) \Big]
\, =\, \frac{1}{{\mathfrak h}_1}{\mathfrak l}_{1}
+ O(t),
\end{align*}
\begin{align*}
\frac{1}{\sqrt{g}} \partial_{\theta} \Big[\frac{1}{\sqrt{g}} {\tilde g}^{12} \Big]
\, =\, &\Big[\frac{1}{\sqrt{{\mathfrak h}_1}}
+{\mathfrak r}_1t
+{\mathfrak r}_2t^2
+O(t^3)\Big]
\times
\Big[t\, \Big(\frac{1}{\sqrt{{\mathfrak h}_1}}{\mathfrak l}_1\Big)'  + O(t^2) \Big]
\, =\, t\, \frac{1}{\sqrt{{\mathfrak h}_1}}\, \Big[\frac{1}{\sqrt{{\mathfrak h}_1}}{\mathfrak l}_1\Big]'+ O(t^2),
\end{align*}
and
\begin{align*}
\frac{1}{\sqrt{g}}\partial_t\Big[\frac{1}{\sqrt{g}} {\tilde g}^{11}  \Big]
\, =\, &\frac{1}{{\mathfrak h}}{\mathfrak f}_1
+{\mathfrak r}_1\frac{1}{\sqrt{{\mathfrak h}_1}}{\mathfrak f}_0
+2t\frac{1}{{\mathfrak h}_1}{\mathfrak f}_2
+3t\frac{1}{\sqrt{{\mathfrak h}_1}} {\mathfrak r}_1 {\mathfrak f}_1
+2t \frac{1}{\sqrt{{\mathfrak h}_1}} {\mathfrak r}_2 {\mathfrak f}_0
+t{\mathfrak r}_1{\mathfrak r}_1{\mathfrak f}_0
+O(t^2),
\end{align*}

\begin{align*}
\frac{1}{\sqrt{g}} \partial_{\theta} \Big[\frac{1}{\sqrt{g}} {\tilde g}^{22} \Big]
\, =\, &\frac{1}{\sqrt{{\mathfrak h}_1}} \Big[ \frac{1}{\sqrt{{\mathfrak h}_1}}{\mathfrak w}_0\Big]'
+t\frac{1}{\sqrt{{\mathfrak h}_1}} \Big[\frac{1}{\sqrt{{\mathfrak h}_1}}{\mathfrak w}_1\Big]'
+t\frac{1}{\sqrt{{\mathfrak h}_1}} \big[{\mathfrak r}_1{\mathfrak w}_0\big]'
+t{\mathfrak r}_1 \Big[ \frac{1}{\sqrt{{\mathfrak h}_1}}{\mathfrak w}_0\Big]'
+ O(t^2).
\end{align*}

\medskip
\medskip
\noindent{\bf{Notation 1: }}
{\em
By collecting all the computations in the above, we set the following conventions.
\begin{equation}\label{h1}
h_1(\theta)=\frac{{\mathfrak f}_0(\theta)}{{\mathfrak h}_1(\theta)}
=\frac{|{\mathfrak a}_1(0, \theta)|^2+|{\mathfrak a}_2(0, \theta)|^2} {{\mathfrak a}_1(0, \theta)|n_1(\theta)|^2+{\mathfrak a}_2(0, \theta)|n_2(\theta)|^2},
\end{equation}

\begin{equation}\label{h2}
 h_2(\theta)=\frac{{\mathfrak w}_0(\theta)}{{\mathfrak h}_1(\theta)}
 =\frac{{\mathfrak a}_1(0, \theta){\mathfrak a}_2(0, \theta)} {{\mathfrak a}_1(0, \theta)|n_1(\theta)|^2+{\mathfrak a}_2(0, \theta)|n_2(\theta)|^2},
\end{equation}

\begin{equation}\label{h3h4}
 h_3(\theta)=\frac{1}{{\mathfrak h}_1}{\mathfrak f}_1
+{\mathfrak r}_1\frac{1}{\sqrt{{\mathfrak h}_1}}{\mathfrak f}_0
=\frac{1}{{\mathfrak h}_1}{\mathfrak f}_1
-\frac{1}{2}\frac{{\mathfrak h}_2}{{\mathfrak h}_1^{2}}{\mathfrak f}_0,
\qquad
h_4(\theta)=\frac{1}{\sqrt{{\mathfrak h}_1}} \Big[ \frac{1}{\sqrt{{\mathfrak h}_1}}{\mathfrak w}_0\Big]'
-\frac{1}{{\mathfrak h}_1}{\mathfrak l}_{1},
\end{equation}

\begin{align}
h_5(\theta)
\, =\, &
2\frac{1}{{\mathfrak h}_1}{\mathfrak f}_2
+3{\mathfrak r}_1\frac{1}{\sqrt{{\mathfrak h}_1}}{\mathfrak f}_1
+2{\mathfrak r}_2\frac{1}{\sqrt{{\mathfrak h}_1}}{\mathfrak f}_0
+{\mathfrak r}_1{\mathfrak r}_1{\mathfrak f}_0
-\frac{1}{\sqrt{{\mathfrak h}_1}} \Big[\frac{1}{\sqrt{{\mathfrak h}_1}}{\mathfrak l}_1\Big]'
\nonumber\\[2mm]
\, =\, &2\frac{1}{{\mathfrak h}_1}{\mathfrak f}_2
-\frac{3}{2}\frac{{\mathfrak h}_2}{{\mathfrak h}_1^{2}}{\mathfrak f}_1
+\frac{{\mathfrak h}_2^2}{{\mathfrak h}_1^{3}}{\mathfrak f}_0
-\frac{{\mathfrak h}_3}{{\mathfrak h}_1^{2}}{\mathfrak f}_0
-\frac{1}{\sqrt{{\mathfrak h}_1}} \Big[\frac{1}{\sqrt{{\mathfrak h}_1}}{\mathfrak l}_1\Big]',
\label{h5}
\end{align}

\begin{equation}\label{h6h7}
 h_6(\theta)=-2\frac{1}{{\mathfrak h}_1}{\mathfrak l}_1,
\qquad
h_7(\theta)\, =\,
\frac{1}{{\mathfrak h}_1}{\mathfrak f}_2
+{\mathfrak g}_1{\mathfrak f}_1
+{\mathfrak g}_2{\mathfrak f}_0
\, =\,
\frac{1}{{\mathfrak h}_1}{\mathfrak f}_2
-\frac{{\mathfrak h}_2}{{\mathfrak h}_1^2}{\mathfrak f}_1
+\Big[\frac{{\mathfrak h}_2^2}{{\mathfrak h}_1^3}
-\frac{{\mathfrak h}_3}{{\mathfrak h}_1^2}\Big]{\mathfrak f}_0,
\end{equation}
\begin{equation}
h_8(\theta)\, =\,
\frac{1}{{\mathfrak h}_1}{\mathfrak f}_1
+{\mathfrak g}_1{\mathfrak f}_0
\, =\,
\frac{1}{{\mathfrak h}_1}{\mathfrak f}_1
-\frac{{\mathfrak h}_2}{{\mathfrak h}_1^2}{\mathfrak f}_0.
\label{h8}
\end{equation}
\qed
}

Hence, the term ${\rm div}\big(\nabla_{{\mathfrak a}(y)}  u\big) $ in \eqref{divoperator} has the following form in the modified Fermi coordinate system
\begin{align}
{\rm div}\big(\nabla_{{\mathfrak a}(y)}  u\big)
\, =\, &\frac{{\tilde g}^{11}}{g} u_{tt}
\, +\, \frac{1}{\sqrt{g}}\partial_t \Big[\frac{1}{\sqrt{g}} {\tilde g}^{11}  \Big] u_{t}
\,-\, 2\frac{{\tilde g}^{12}}{g}  u_{\theta t}
\,-\, \frac{1}{\sqrt{g}} \partial_t \Big[\frac{1}{\sqrt{g}} {\tilde g}^{12}  \Big] u_{\theta}
\nonumber\\[2mm]
&\,-\, \frac{1}{\sqrt{g}} \partial_{\theta} \Big[\frac{1}{\sqrt{g}} {\tilde g}^{12}  \Big] u_{t}
\, +\, \frac{{\tilde g}^{22}}{g} u_{\theta\theta}
\, +\, \frac{1}{\sqrt{g}} \partial_{\theta} \Big[\frac{1}{\sqrt{g}} {\tilde g}^{22}  \Big]u_{\theta}
\nonumber\\[2mm]
\, =\, &h_1(\theta) u_{tt}+h_2(\theta)u_{\theta\theta}
+h_3(\theta)u_{t}+h_4(\theta) u_\theta+\bar{B}_1(u)+\bar{B}_0(u),
\label{laplacelocal}
\end{align}
where
\begin{equation}
\bar{B}_1(u)\, =\, h_5(\theta)tu_{t}+h_6(\theta) tu_{t\theta}
+h_7(\theta)t^2 u_{tt}+h_8(\theta)tu_{tt},
\label{B1bar}
\end{equation}
and
\begin{equation}
\bar{B}_0(u)
\, =\, h_9(t, \theta)t u_{\theta\theta}
+h_{10}(t, \theta)t^2 u_{\theta\theta}
+h_{11}(t, \theta)t^2u_{\theta t}
+h_{12}(t, \theta)t u_{\theta}.
\label{B0bar}
\end{equation}
Here, $h_9, \cdots, h_{12}$ are smooth functions.

\bigskip
\noindent{\bf{Part 2: the operator $\nabla_{{\mathfrak a}(y)}u\cdot\nu$}}

We finally show the local expression of $\nabla_{{\mathfrak a}(y)}u\cdot\nu$ in (\ref{originalproblem}).
Suppose that, in the local coordinates $(t, \theta)$ of \eqref{Fermicoordinates-modified}, the unit outer normal of $\partial\Omega$ is expressed in the form
$$
\nu\, =\, \sigma_1\frac{{\partial} F}{{\partial} t}+\sigma_2\frac{{\partial} F}{{\partial} \theta}.
$$
If $\theta=0$ or $\theta=1$, the expression of $F(t, \theta)$ in \eqref{Fermicoordinates-modified} gives the curves $\mathcal{C}_1$ or $\mathcal{C}_2$.
We then have
$$
\langle {{\partial} F}/{{\partial} t}, \, \nu\rangle=0\quad \mbox{at } \theta=0, 1.
$$
For the convenience of notation, in the following lines of this part,
we will always take $\theta=0$ or $\theta=1$ without any further announcement.
Hence
$$
\sigma_2\neq 0, \quad \frac{{\partial} F}{{\partial} \theta}\neq 0, \quad
{\rm and}\quad
\sigma_1g_{11}+\sigma_2g_{12}\, =\, 0.$$
On the other hand, $\langle \nu, \nu \rangle=1$, that is
$$
\left<\sigma_1\frac{{\partial} F}{{\partial} t}+\sigma_2\frac{{\partial} F}{{\partial} \theta}, \, \sigma_1\frac{{\partial} F}{{\partial} t}+\sigma_2\frac{{\partial} F}{{\partial} \theta}\right>\, =\, 1,
$$
which implies that
$$
\sigma_1^2g_{11}+\sigma_2^2g_{22}+2\sigma_1\sigma_2g_{12}\, =\, 1.
$$
Combining above two equations, one can get
$$\sigma_1\, =\, \pm\frac{g^{12}}{\sqrt{g^{22}}},
\qquad
\sigma_2\, =\, \pm\sqrt{g^{22}}.$$
By choosing the sign $"+"$ and using \eqref{g12inverse}-\eqref{g22inverse}, it is easy to check that
\begin{align}
\sigma_1
\, =\, &\Big[\, {\mathfrak g}_3+{\mathfrak g}_4t+{\mathfrak g}_5t^2+O(t^3)\, \Big]
\times\Big\{\, \frac{1}{{\mathfrak g}_6}-\frac{{\mathfrak g}_7}{{\mathfrak g}_6^2}t
+\Big[\frac{1}{2}\frac{{\mathfrak g}_7^2}{{\mathfrak g}_6^{3}}-\frac{{\mathfrak g}_8}{{\mathfrak g}_6^{2}}\Big]t^2+O(t^3)\, \Big\}
\nonumber\\[2mm]
\, =\, &\frac{{\mathfrak g}_3}{{\mathfrak g}_6}
+t\, \Big[\, \frac{{\mathfrak g}_4}{{\mathfrak g}_6}-\frac{{\mathfrak g}_7{\mathfrak g}_3}{{\mathfrak g}_6^2}\, \Big]
+t^2\, \Big[\, \frac{{\mathfrak g}_5}{{\mathfrak g}_6}
-\frac{{\mathfrak g}_7{\mathfrak g}_4}{{\mathfrak g}_6^2}
+{\mathfrak g}_3{\mathfrak g}_5\, \Big]+O(t^3)
\nonumber\\[2mm]
\, \equiv\, &{\mathfrak y}_1(\theta)+t{\mathfrak y}_2(\theta)+t^2{\mathfrak y}_3(\theta)+O(t^3),
\label{sigma1}
\end{align}
and
\begin{align}\label{sigma2}
\sigma_2\, =\, \sqrt{g^{22}}
\, =\, &\sqrt{{\mathfrak g}_6}
+\frac{1}{2}\frac{{\mathfrak g}_7}{\sqrt{{\mathfrak g}_6}}t
+\frac{1}{2}\Big[\frac{{\mathfrak g}_8}{\sqrt{{\mathfrak g}_6}}
-\frac{1}{4}\frac{{\mathfrak g}_7^2}{(\sqrt{{\mathfrak g}_6})^3}\Big]t^2+O(t^3)
\nonumber\\[2mm]
\, \equiv\, &{\mathfrak y}_4(\theta)+t{\mathfrak y}_5(\theta)+t^2{\mathfrak y}_6(\theta)+O(t^3).
\end{align}

\medskip
In the modified Fermi coordinates $(t, \theta)$ in (\ref{Fermicoordinates-modified}), the normal derivative
$\nabla_{{\mathfrak a}(y)}u\cdot\nu$  has a local form as follows
\begin{align*}
 \nabla_{{\mathfrak a}(y)}u\cdot\nu
\, =\, &\Big({\mathfrak a}_1(y)\frac{\partial u}{\partial y_1}, {\mathfrak a}_2(y)\frac{\partial u}{\partial y_2}  \Big)\Big( \sigma_1 \frac{\partial F_1}{\partial t}+\sigma_2 \frac{\partial F_1}{\partial \theta}, \sigma_1 \frac{\partial F_2}{\partial t}+\sigma_2 \frac{\partial F_2}{\partial \theta}\Big)
 \\[2mm]
\, =\, &{\mathfrak a}_1(y)\frac{\partial u}{\partial y_1} \Big( \sigma_1 \frac{\partial F_1}{\partial t}+\sigma_2 \frac{\partial F_1}{\partial \theta} \Big)
\, +\, {\mathfrak a}_2(y)\frac{\partial u}{\partial y_2} \Big(\sigma_1 \frac{\partial F_2}{\partial t}+\sigma_2 \frac{\partial F_2}{\partial \theta}\Big)
\\[2mm]
\, =\, & {\mathfrak a}_1(y)\frac{1}{\sqrt{g}}
\left[ \frac{\partial u}{\partial t}\frac{\partial F_2}{\partial \theta}
\,-\, \frac{\partial u}{\partial \theta}\frac{\partial F_2}{\partial t} \right]
\Big( \sigma_1 \frac{\partial F_1}{\partial t}+\sigma_2 \frac{\partial F_1}{\partial \theta} \Big)
\\[2mm]
&\, +\, {\mathfrak a}_2(y)\frac{1}{\sqrt{g}}
\left[ \frac{\partial u}{\partial \theta}\frac{\partial F_1}{\partial t}
\,-\, \frac{\partial u}{\partial t}\frac{\partial F_1}{\partial \theta} \right]
\Big(\sigma_1 \frac{\partial F_2}{\partial t}+\sigma_2 \frac{\partial F_2}{\partial \theta}\Big)
\\[2mm]
\, =\, & \frac{\sigma_1}{\sqrt{g}}
\left[{\mathfrak a}_1(y)\frac{\partial F_1}{\partial t}\frac{\partial F_2}{\partial \theta}
-{\mathfrak a}_2(y)\frac{\partial F_1}{\partial \theta}\frac{\partial F_2}{\partial t}\right]\frac{\partial u}{\partial t}
+\frac{\sigma_2 }{\sqrt{g}}
\left[{\mathfrak a}_1(y)\frac{\partial F_1}{\partial \theta}\frac{\partial F_2}{\partial \theta}
-{\mathfrak a}_2(y)\frac{\partial F_1}{\partial \theta}\frac{\partial F_2}{\partial \theta} \right]\frac{\partial u}{\partial t}
\\[2mm]
&+\frac{\sigma_1}{\sqrt{g}}
\left[ {\mathfrak a}_2(y)\frac{\partial F_1}{\partial t} \frac{\partial F_2}{\partial t}
-{\mathfrak a}_1(y)\frac{\partial F_1}{\partial t} \frac{\partial F_2}{\partial t}\right] \frac{\partial u}{\partial \theta}
+\frac{\sigma_2 }{\sqrt{g}}
\left[{\mathfrak a}_2(y) \frac{\partial F_1}{\partial t} \frac{\partial F_2}{\partial \theta}
-{\mathfrak a}_1(y) \frac{\partial F_1}{\partial \theta}\frac{\partial F_2}{\partial t}\right] \frac{\partial u}{\partial \theta}.
\end{align*}
 According to the expressions of $\frac{\partial F_i}{\partial t}$ and $\frac{\partial F_i}{\partial \theta}$ as in \eqref{F_t}-\eqref{F2_theta}, it is easy to derive that
\begin{align}
\frac{\partial F_1}{\partial t} \frac{\partial F_2}{\partial t}
\, =\, &{\tilde{\mathfrak a}}_1{\tilde{\mathfrak a}}_2n_1n_2
+t \big( {\tilde{\mathfrak a}}_2q_1n_2+{\tilde{\mathfrak a}}_1q_2n_1 \big)
+ t^2q_1q_2+\frac{t^2}{2}\big( {\tilde{\mathfrak a}}_2m_1n_2+{\tilde{\mathfrak a}}_1m_2n_1 \big) +O(t^3), \label{FtFt}
\\[2mm]
\frac{\partial F_1}{\partial \theta}\frac{\partial F_2}{\partial \theta}
\, =\, &-n_1n_2
\, +\, t \big(- {\tilde{\mathfrak a}' }_2|n_2|^2+ {\tilde{\mathfrak a}' }_1|n_1|^2\big)
\, +\, tk\big({\tilde{\mathfrak a}}_1+{\tilde{\mathfrak a}}_2 \big) n_1n_2
\, +\, \frac{t^2}{2}\big(-q_2'n_2+ q_1'n_1\big)
\nonumber\\[2mm]
&\, +\, t^2{\tilde{\mathfrak a}' }_1{\tilde{\mathfrak a}' }_2n_1n_2
\,-\, t^2k \big( {\tilde{\mathfrak a}' }_1{\tilde{\mathfrak a}}_2 |n_1|^2- {\tilde{\mathfrak a}}_1 {\tilde{\mathfrak a}' }_2|n_2|^2 \big)
\,-\, t^2k^2{\tilde{\mathfrak a}}_1 {\tilde{\mathfrak a}}_2 n_1n_2
\, +\, O(t^3), \label{FteFte}
\end{align}
and
\begin{align}
\frac{\partial F_1}{\partial t} \frac{\partial F_2}{\partial \theta}
\, =\, &{\tilde{\mathfrak a}}_1|n_1|^2
+t{\tilde{\mathfrak a}}_1{\tilde{\mathfrak a}' }_2n_1n_2
\,-\, tk{\tilde{\mathfrak a}}_1{\tilde{\mathfrak a}}_2|n_1|^2
\, +\, t q_1n_1
\, +\, t^2q_1{\tilde{\mathfrak a}' }_2n_2
\nonumber\\[2mm]
&\,-\, t^2k {\tilde{\mathfrak a}}_2q_1n_1
\, +\, \frac{t^2}{2}{\tilde{\mathfrak a}}_1n_1q_2'\, +\, \frac{t^2}{2}m_1n_1+O(t^3),
\label{FtFte}
\\
\frac{\partial F_1}{\partial \theta}\frac{\partial F_2}{\partial t}
\, =\, &-{\tilde{\mathfrak a}}_2|n_2|^2
\, +\, t{\tilde{\mathfrak a}' }_1{\tilde{\mathfrak a}}_2n_1n_2
\, +\, tk{\tilde{\mathfrak a}}_1{\tilde{\mathfrak a}}_2|n_2|^2
\,-\, t q_2n_2
\, +\, t^2{\tilde{\mathfrak a}' }_1n_1 q_2
\nonumber\\[2mm]
&\, +\, t^2k{\tilde{\mathfrak a}}_1 q_2n_2
\, +\, \frac{t^2}{2}{\tilde{\mathfrak a}}_2n_2q_1'
\,-\, \frac{t^2}{2}m_2n_2\, +\, O(t^3).
\label{FteFt}
\end{align}

By using the Taylor expansion of ${\mathfrak a}_i(y)$ as in \eqref{expressionofai} and \eqref{FtFt}-\eqref{FteFt},  it is easy to derive that
\begin{align}
&{\mathfrak a}_1(y)\frac{\partial F_1}{\partial t}\frac{\partial F_2}{\partial \theta}
\,-\, {\mathfrak a}_2(y)\frac{\partial F_1}{\partial \theta}\frac{\partial F_2}{\partial t}
\nonumber\\[2mm]
&\, =\, {\mathfrak a}_1(y)\big[{\tilde{\mathfrak a}}_1|n_1|^2
+t{\tilde{\mathfrak a}}_1{\tilde{\mathfrak a}' }_2n_1n_2
\,-\, tk{\tilde{\mathfrak a}}_1{\tilde{\mathfrak a}}_2|n_1|^2
\, +\, t q_1n_1
+O(t^2)\big]
\nonumber\\[2mm]
&\qquad\,-\, {\mathfrak a}_2(y)\big[-{\tilde{\mathfrak a}}_2|n_2|^2
\, +\, t{\tilde{\mathfrak a}' }_1{\tilde{\mathfrak a}}_2n_1n_2
\, +\, tk{\tilde{\mathfrak a}}_1{\tilde{\mathfrak a}}_2|n_2|^2
\,-\, t q_2n_2
\, +\, O(t^2)\big]
\nonumber\\[2mm]
&\, =\, \big[{\mathfrak a}_1(0, \theta){\tilde{\mathfrak a}}_1|n_1|^2\, +\, {\mathfrak a}_2(0, \theta){\tilde{\mathfrak a}}_2|n_2|^2 \big]
\nonumber\\[2mm]
&\qquad+t\Big\{\big[ {\mathfrak a}_1(0, \theta){\tilde{\mathfrak a}}_1{\tilde{\mathfrak a}'}_2n_1n_2-{\mathfrak a}_2(0, \theta){\tilde{\mathfrak a}'}_1{\tilde{\mathfrak a}}_2n_1n_2 \big]
\nonumber\\[2mm]
&\qquad \qquad
- k\big[ {\mathfrak a}_1(0, \theta){\tilde{\mathfrak a}}_1{\tilde{\mathfrak a}}_2|n_1|^2-{\mathfrak a}_2(0, \theta){\tilde{\mathfrak a}}_1{\tilde{\mathfrak a}}_2|n_2|^2 \big]
+\big[{\mathfrak a}_1(0, \theta)q_1n_1+{\mathfrak a}_2(0, \theta)q_2n_2 \big]
\nonumber\\[2mm]
&\qquad \qquad
+\big[\partial_t{\mathfrak a}_1(0, \theta){\tilde{\mathfrak a}}_1|n_1|^2+\partial_t{\mathfrak a}_2(0, \theta){\tilde{\mathfrak a}}_2|n_2|^2\big]\Big\}
\, +\, O(t^2)
\nonumber\\[2mm]
&\, \equiv\, \mathfrak{p}_1(\theta)+ t\mathfrak{p}_2(\theta)\, +\, O(t^2),
\label{p1p2}
\end{align}
\begin{align}
&{\mathfrak a}_1(y)\frac{\partial F_1}{\partial \theta}\frac{\partial F_2}{\partial \theta}\,-\, {\mathfrak a}_2(y)\frac{\partial F_1}{\partial \theta}\frac{\partial F_2}{\partial \theta}
\nonumber\\[2mm]
&\, =\, \Big[{\mathfrak a}_1(y)-{\mathfrak a}_2(y)\Big]\times\Big[-n_1n_2
+t \big({\tilde{\mathfrak a}' }_1|n_1|^2- {\tilde{\mathfrak a}' }_2|n_2|^2\big)
+tk\big({\tilde{\mathfrak a}}_1+{\tilde{\mathfrak a}}_2 \big) n_1n_2
+\frac{t^2}{2}\big(q_1'n_1-q_2'n_2\big)
\nonumber\\[2mm]
&\qquad \qquad \qquad \qquad \qquad +t^2{\tilde{\mathfrak a}' }_1{\tilde{\mathfrak a}' }_2n_1n_2
-t^2k \big( {\tilde{\mathfrak a}' }_1{\tilde{\mathfrak a}}_2 |n_1|^2- {\tilde{\mathfrak a}}_1 {\tilde{\mathfrak a}' }_2|n_2|^2 \big)
-t^2k^2{\tilde{\mathfrak a}}_1 {\tilde{\mathfrak a}}_2 n_1n_2
+O(t^3)\Big]
\nonumber\\[2mm]
&\, =\, -\big[{\mathfrak a}_1(0, \theta)- {\mathfrak a}_2(0, \theta) \big]n_1n_2
\nonumber\\[2mm]
&\quad\, \, +t\Big\{\big[{\mathfrak a}_1(0, \theta)- {\mathfrak a}_2(0, \theta) \big] \big(- {\tilde{\mathfrak a}' }_2|n_2|^2+ {\tilde{\mathfrak a}' }_1|n_1|^2\big)
\nonumber\\[2mm]
&\qquad\qquad+k\big[{\mathfrak a}_1(0, \theta)- {\mathfrak a}_2(0, \theta) \big]\big({\tilde{\mathfrak a}}_1+{\tilde{\mathfrak a}}_2 \big) n_1n_2
-\big[\partial_t{\mathfrak a}_1(0, \theta)- \partial_t{\mathfrak a}_2(0, \theta) \big]n_1n_2\Big\}
\nonumber\\[2mm]
&\quad\, \, +t^2\Big\{\frac{1}{2}\big[{\mathfrak a}_1(0, \theta)- {\mathfrak a}_2(0, \theta) \big]\big(-q_2'n_2+ q_1'n_1\big)
+\big[{\mathfrak a}_1(0, \theta)- {\mathfrak a}_2(0, \theta) \big]{\tilde{\mathfrak a}' }_1{\tilde{\mathfrak a}' }_2n_1n_2
\nonumber\\[2mm]
&\qquad \qquad -k\big[{\mathfrak a}_1(0, \theta)- {\mathfrak a}_2(0, \theta) \big] \big( {\tilde{\mathfrak a}' }_1{\tilde{\mathfrak a}}_2 |n_1|^2- {\tilde{\mathfrak a}}_1 {\tilde{\mathfrak a}' }_2|n_2|^2 \big)
-k^2\big[{\mathfrak a}_1(0, \theta)- {\mathfrak a}_2(0, \theta) \big] {\tilde{\mathfrak a}}_1 {\tilde{\mathfrak a}}_2 n_1n_2
\nonumber\\[2mm]
&\qquad \qquad
+\big[\partial_t{\mathfrak a}_1(0, \theta)- \partial_t{\mathfrak a}_2(0, \theta) \big] \big(- {\tilde{\mathfrak a}' }_2|n_2|^2+ {\tilde{\mathfrak a}' }_1|n_1|^2\big)
\nonumber\\[2mm]
&\qquad \qquad+k\big[\partial_t{\mathfrak a}_1(0, \theta)- \partial_t{\mathfrak a}_2(0, \theta) \big]\big({\tilde{\mathfrak a}}_1+{\tilde{\mathfrak a}}_2 \big) n_1n_2\Big\}
+O(t^3)
\nonumber\\[2mm]
&\, \equiv\, \mathfrak{p}_3(\theta)+t \mathfrak{p}_4(\theta)\, +\, t^2 \mathfrak{p}_5(\theta)+O(t^3),
\label{p3p4}
\end{align}
\begin{align}
{\mathfrak a}_2(y)\frac{\partial F_1}{\partial t} \frac{\partial F_2}{\partial t}
\,-\, {\mathfrak a}_1(y)\frac{\partial F_1}{\partial t} \frac{\partial F_2}{\partial t}
\, =\, &\big[{\mathfrak a}_1(0, \theta)- {\mathfrak a}_2(0, \theta) \big]{\tilde{\mathfrak a}}_1{\tilde{\mathfrak a}}_2n_1n_2
\nonumber\\[2mm]
&+t \Big\{\big[{\mathfrak a}_1(0, \theta)- {\mathfrak a}_2(0, \theta) \big]\big( {\tilde{\mathfrak a}}_2q_1n_2+{\tilde{\mathfrak a}}_1q_2n_1 \big)
\nonumber\\[2mm]
&\qquad+\big[\partial_t{\mathfrak a}_1(0, \theta)- \partial_t{\mathfrak a}_2(0, \theta) \big]{\tilde{\mathfrak a}}_1{\tilde{\mathfrak a}}_2n_1n_2 \Big\}
+O(t^2)
\nonumber\\[2mm]
\, \equiv\, &\mathfrak{p}_6(\theta)\, +\, t\mathfrak{p}_7(\theta)+O(t^2),
\label{p6p7}
\end{align}
and
\begin{align}
&{\mathfrak a}_2(y) \frac{\partial F_1}{\partial t} \frac{\partial F_2}{\partial \theta}
\,-\, {\mathfrak a}_1(y) \frac{\partial F_1}{\partial \theta}\frac{\partial F_2}{\partial t}
\nonumber\\[2mm]
&\, =\, {\mathfrak a}_2(y) \big[{\tilde{\mathfrak a}}_1|n_1|^2
+t{\tilde{\mathfrak a}}_1{\tilde{\mathfrak a}' }_2n_1n_2
\,-\, tk{\tilde{\mathfrak a}}_1{\tilde{\mathfrak a}}_2|n_1|^2
\, +\, t q_1n_1
\, +\, t^2{\tilde{\mathfrak a}' }_2n_2q_1
\nonumber\\[2mm]
&\qquad \qquad\,-\, t^2k {\tilde{\mathfrak a}}_2q_1n_1
\, +\, \frac{t^2}{2}{\tilde{\mathfrak a}}_1n_1q_2'
\, +\, \frac{t^2}{2}m_1n_1+O(t^3)\big]
\nonumber\\
&\quad\,-\, {\mathfrak a}_1(y)\big[-{\tilde{\mathfrak a}}_2|n_2|^2
\, +\, t{\tilde{\mathfrak a}' }_1{\tilde{\mathfrak a}}_2n_1n_2
\, +\, tk{\tilde{\mathfrak a}}_1{\tilde{\mathfrak a}}_2|n_2|^2
\,-\, t q_2n_2
\, +\, t^2{\tilde{\mathfrak a}' }_1n_1 q_2
\nonumber\\[2mm]
&\qquad \qquad \quad\, +\, t^2k{\tilde{\mathfrak a}}_1 q_2n_2
\, +\, \frac{t^2}{2}{\tilde{\mathfrak a}}_2n_2q_1'
\,-\, \frac{t^2}{2}m_2n_2\, +\, O(t^3)\big]
\nonumber\\[2mm]
&\, =\, \big[{\mathfrak a}_2(0, \theta){\tilde{\mathfrak a}}_1|n_1|^2+{\mathfrak a}_1(0, \theta){\tilde{\mathfrak a}}_2|n_2|^2 \big]
\nonumber\\[2mm]
&\qquad+t\Big\{\big[\partial_t{\mathfrak a}_2(0, \theta){\tilde{\mathfrak a}}_1|n_1|^2 +\partial_t{\mathfrak a}_1(0, \theta){\tilde{\mathfrak a}}_2|n_2|^2\big]
+\big[{\mathfrak a}_2(0, \theta){\tilde{\mathfrak a}}_1{\tilde{\mathfrak a}' }_2
-{\mathfrak a}_1(0, \theta){\tilde{\mathfrak a}' }_1{\tilde{\mathfrak a}}_2\big]n_1n_2
\nonumber\\[2mm]
&\qquad \qquad\,-\, k\big[{\mathfrak a}_2(0, \theta){\tilde{\mathfrak a}}_1{\tilde{\mathfrak a}}_2|n_1|^2+{\mathfrak a}_1(0, \theta){\tilde{\mathfrak a}}_1{\tilde{\mathfrak a}}_2|n_2|^2\big]
+\big[{\mathfrak a}_2(0, \theta)q_1n_1+{\mathfrak a}_1(0, \theta)q_2n_2  \big]\Big\}
\nonumber\\[2mm]
&\qquad+t^2\Big\{\big[\partial_t{\mathfrak a}_2(0, \theta){\tilde{\mathfrak a}}_1{\tilde{\mathfrak a}' }_2
-\partial_t{\mathfrak a}_1(0, \theta){\tilde{\mathfrak a}' }_1{\tilde{\mathfrak a}}_2\big]n_1n_2
\nonumber\\[2mm]
&\qquad \qquad\,-\, k\big[ \partial_t{\mathfrak a}_2(0, \theta){\tilde{\mathfrak a}}_1{\tilde{\mathfrak a}}_2|n_1|^2+\partial_t{\mathfrak a}_1(0, \theta){\tilde{\mathfrak a}}_1{\tilde{\mathfrak a}}_2|n_2|^2 \big]
\, +\, \big[\partial_t{\mathfrak a}_2(0, \theta)q_1n_1+\partial_t{\mathfrak a}_1(0, \theta)q_2n_2 \big]
\nonumber\\[2mm]
&\qquad \qquad+\big[{\mathfrak a}_2(0, \theta){\tilde{\mathfrak a}' }_2n_2q_1-{\mathfrak a}_1(0, \theta){\tilde{\mathfrak a}' }_1n_1q_2\big]
\,-\, k\big[{\mathfrak a}_2(0, \theta){\tilde{\mathfrak a}}_2q_1n_1-{\mathfrak a}_1(0, \theta){\tilde{\mathfrak a}}_2q_2n_2  \big]
\nonumber\\[2mm]
&\qquad \qquad\, +\, \frac{1}{2} \big[{\mathfrak a}_2(0, \theta){\tilde{\mathfrak a}}_1n_1q_2'-{\mathfrak a}_1(0, \theta){\tilde{\mathfrak a}}_2n_2q_1'  \big]
\, +\, \frac{1}{2} \big[{\mathfrak a}_2(0, \theta)m_1n_1+{\mathfrak a}_1(0, \theta)m_2n_2 \big]\Big\}
+O(t^3)
\nonumber\\[2mm]
&\, \equiv\, \mathfrak{p}_8(\theta)\, +\, t\mathfrak{p}_9(\theta)
\, +\, t^2\mathfrak{p}_{10}(\theta)+O(t^3).
\label{p8p9p10}
\end{align}

\medskip
On the other hand, by recalling \eqref{g-inverse2}, \eqref{sigma1} and \eqref{sigma2}, we can obtain that
\begin{align*}
\frac{1}{\sqrt{g}} \sigma_1
\, =\, & \Big[\frac{1}{\sqrt{{\mathfrak h}_1}}+{\mathfrak r}_1t+{\mathfrak r}_3t^2+O(t^3) \Big]\times\big[{\mathfrak y}_1(\theta)+t{\mathfrak y}_2(\theta)+t^2{\mathfrak y}_3(\theta)+O(t^3)\big]
\\[2mm]
\, =\, &\frac{1}{\sqrt{{\mathfrak h}_1}}{\mathfrak y}_1(\theta)
+t{\mathfrak r}_1{\mathfrak y}_1(\theta)
+t\frac{1}{\sqrt{{\mathfrak h}_1}}{\mathfrak y}_2(\theta)
+t^2\Big[{\mathfrak r}_1{\mathfrak y}_2(\theta)
+\frac{1}{\sqrt{{\mathfrak h}_1}}{\mathfrak y}_3(\theta)
+{\mathfrak c}_3{\mathfrak y}_1(\theta)
\Big]+O(t^3)
\\[2mm]
\, \equiv\, &\frac{1}{\sqrt{{\mathfrak h}_1}}{\mathfrak y}_1(\theta)
+t{\mathfrak r}_1{\mathfrak y}_1(\theta)
+t\frac{1}{\sqrt{{\mathfrak h}_1}}{\mathfrak y}_2(\theta)
+t^2{\mathfrak k}_1+O(t^3),
\end{align*}
\begin{align*}
 \frac{1}{\sqrt{g}} \sigma_2
 \, =\, & \Big[\frac{1}{\sqrt{{\mathfrak h}_1}}+{\mathfrak r}_1t+{\mathfrak r}_3t^2+O(t^3) \big]\times \big[{\mathfrak y}_4(\theta)+t{\mathfrak y}_5(\theta)+t^2{\mathfrak y}_6(\theta)+O(t^3)\Big]
\\[2mm]
\, =\, &\frac{1}{\sqrt{{\mathfrak h}_1}}{\mathfrak y}_4(\theta)
+t{\mathfrak r}_1{\mathfrak y}_4(\theta)
+t\frac{1}{\sqrt{{\mathfrak h}_1}}{\mathfrak y}_5(\theta)
+t^2\Big[{\mathfrak r}_3{\mathfrak y}_4(\theta)
+{\mathfrak y}_5(\theta){\mathfrak r}_1
+\frac{1}{\sqrt{{\mathfrak h}_1}}{\mathfrak y}_6(\theta)
\Big]+O(t^3)
\\[2mm]
\, \equiv\, &\frac{1}{\sqrt{{\mathfrak h}_1}}{\mathfrak y}_4(\theta)
+t{\mathfrak r}_1{\mathfrak y}_4(\theta)
+t\frac{1}{\sqrt{{\mathfrak h}_1}}{\mathfrak y}_5(\theta)
+t^2{\mathfrak k}_2+O(t^3).
\end{align*}

Therefore, we obtain that
\begin{align*}
&\frac{1}{\sqrt{g}} \sigma_1
\left[{\mathfrak a}_1(y)\frac{\partial F_1}{\partial t}\frac{\partial F_2}{\partial\theta}
\,-\, {\mathfrak a}_2(y)\frac{\partial F_1}{\partial \theta}\frac{\partial F_2}{\partial t}\right]u_t
\, +\, \frac{1}{\sqrt{g}}\sigma_2
\left[{\mathfrak a}_1(y)\frac{\partial F_1}{\partial \theta}\frac{\partial F_2}{\partial \theta}
\,-\, {\mathfrak a}_2(y)\frac{\partial F_1}{\partial \theta}\frac{\partial F_2}{\partial \theta} \right]u_t
\\[2mm]
\, =\, &\Big[\frac{1}{\sqrt{{\mathfrak h}_1}}{\mathfrak y}_1(\theta)
+t{\mathfrak r}_1{\mathfrak y}_1(\theta)
+t\frac{1}{\sqrt{{\mathfrak h}_1}}{\mathfrak y}_2(\theta)
+t^2{\mathfrak k}_1+O(t^3)\Big]\times
\Big[\mathfrak{p}_1(\theta)+ t\mathfrak{p}_2(\theta)\, +\, O(t^2)\Big]u_t
\\[2mm]
&\, +\, \Big[\frac{1}{\sqrt{{\mathfrak h}_1}}{\mathfrak y}_4(\theta)
+t{\mathfrak r}_1{\mathfrak y}_4(\theta)
+t\frac{1}{\sqrt{{\mathfrak h}_1}}{\mathfrak y}_5(\theta)
+t^2{\mathfrak k}_2+O(t^3)\Big]\times
\Big[\mathfrak{p}_3(\theta)+t \mathfrak{p}_4(\theta)+ t^2 \mathfrak{p}_5(\theta)+O(t^3)\Big]u_t
\\[2mm]
\, \equiv\, &\frac{1}{\sqrt{{\mathfrak h}_1}}{\mathfrak y}_1(\theta)\mathfrak{p}_1(\theta)u_t
+\frac{1}{\sqrt{{\mathfrak h}_1}}{\mathfrak y}_4(\theta)\mathfrak{p}_3(\theta)u_t
\\[2mm]
&+\Big[{\mathfrak r}_1{\mathfrak y}_1(\theta)\mathfrak{p}_1(\theta)
+\frac{1}{\sqrt{{\mathfrak h}_1}}{\mathfrak y}_2(\theta)\mathfrak{p}_1(\theta)
+\frac{1}{\sqrt{{\mathfrak h}_1}}{\mathfrak y}_1(\theta)\mathfrak{p}_2(\theta)
\\[2mm]
&\quad+\frac{1}{\sqrt{{\mathfrak h}_1}}{\mathfrak y}_4(\theta)\mathfrak{p}_4(\theta)
+{\mathfrak r}_1{\mathfrak y}_4(\theta)\mathfrak{p}_3(\theta)
+\frac{1}{\sqrt{{\mathfrak h}_1}}{\mathfrak y}_5(\theta)\mathfrak{p}_3(\theta)\Big]tu_t
+b_1(\theta)t^2u_t+O(t^3),
\end{align*}
and
\begin{align*}
&\frac{1}{\sqrt{g}} \sigma_1
\left[ {\mathfrak a}_2(y)\frac{\partial F_1}{\partial t} \frac{\partial F_2}{\partial t}
\,-\, {\mathfrak a}_1(y)\frac{\partial F_1}{\partial t} \frac{\partial F_2}{\partial t}\right]u_{\theta}
\, +\, \frac{1}{\sqrt{g}}\sigma_2
\left[{\mathfrak a}_2(y) \frac{\partial F_1}{\partial t} \frac{\partial F_2}{\partial \theta}
\,-\, {\mathfrak a}_1(y) \frac{\partial F_1}{\partial \theta}\frac{\partial F_2}{\partial t}\right] u_{\theta}
\\[2mm]
\, =\, &\Big[\frac{1}{\sqrt{{\mathfrak h}_1}}{\mathfrak y}_1(\theta)
+t{\mathfrak r}_1{\mathfrak y}_1(\theta)
+t\frac{1}{\sqrt{{\mathfrak h}_1}}{\mathfrak y}_2(\theta)
+t^2{\mathfrak k}_1+O(t^3)\Big]\times
\Big[\mathfrak{p}_6(\theta)\, +\, t\mathfrak{p}_7(\theta)+O(t^2)\Big]u_{\theta}
\\[2mm]
&\, +\, \frac{1}{\sqrt{{\mathfrak h}_1}}{\mathfrak y}_4(\theta)
+t{\mathfrak r}_1{\mathfrak y}_4(\theta)
+t\frac{1}{\sqrt{{\mathfrak h}_1}}{\mathfrak y}_5(\theta)
+t^2{\mathfrak k}_2+O(t^3)\Big]\times
\Big[\mathfrak{p}_8(\theta)+t\mathfrak{p}_9(\theta)
+t^2\mathfrak{p}_{10}(\theta)+O(t^3)\Big]u_{\theta}
\\[2mm]
\, \equiv\, &\frac{1}{\sqrt{{\mathfrak h}_1}}{\mathfrak y}_1(\theta)\mathfrak{p}_6(\theta)u_{\theta}
+\frac{1}{\sqrt{{\mathfrak h}_1}}{\mathfrak y}_4(\theta)\mathfrak{p}_8(\theta)u_{\theta}
+b_3(\theta)tu_{\theta}
+b_2(\theta)t^2u_{\theta}+O(t^3).
\end{align*}

\medskip
\medskip
Whence, the terms in $\nabla_{{\mathfrak a}(y)}u\cdot\nu$ will be rearranged in the form
\begin{align*}
\nabla_{{\mathfrak a}(y)}u\cdot\nu
\, =\, &\frac{1}{\sqrt{{\mathfrak h}_1}}{\mathfrak y}_1(\theta)\mathfrak{p}_1(\theta)u_t
+\frac{1}{\sqrt{{\mathfrak h}_1}}{\mathfrak y}_4(\theta)\mathfrak{p}_3(\theta)u_t
+\frac{1}{\sqrt{{\mathfrak h}_1}}{\mathfrak y}_1(\theta)\mathfrak{p}_6(\theta)u_{\theta}
+\frac{1}{\sqrt{{\mathfrak h}_1}}{\mathfrak y}_4(\theta)\mathfrak{p}_8(\theta)u_{\theta}
\\[2mm]
&+\Big[{\mathfrak r}_1{\mathfrak y}_1(\theta)\mathfrak{p}_1(\theta)
+\frac{1}{\sqrt{{\mathfrak h}_1}}{\mathfrak y}_2(\theta)\mathfrak{p}_1(\theta)
+\frac{1}{\sqrt{{\mathfrak h}_1}}{\mathfrak y}_1(\theta)\mathfrak{p}_2(\theta)
\\[2mm]
&\qquad+\frac{1}{\sqrt{{\mathfrak h}_1}}{\mathfrak y}_4(\theta)\mathfrak{p}_4(\theta)
+{\mathfrak r}_1{\mathfrak y}_4(\theta)\mathfrak{p}_3(\theta)
+\frac{1}{\sqrt{{\mathfrak h}_1}}{\mathfrak y}_5(\theta)\mathfrak{p}_3(\theta)\Big]tu_t
\\[2mm]
&+b_3(\theta)tu_{\theta}\, +\, b_1(\theta)t^2u_t+b_2(\theta)t^2u_{\theta}+O(t^3).
\end{align*}
More precisely, we shall derive the boundary operator in the following way.
According to \eqref{tildea1=a2}, we can obtain that
\begin{align*}
{\mathfrak y}_1(0)\, =\, {\mathfrak y}_1(1)\, =\, 0, \quad
{\mathfrak p}_3(0)\, =\, {\mathfrak p}_3(1)\, =\, 0.
\end{align*}
For $\theta\, =\, 0$, the boundary operator becomes
\begin{equation}
\begin{split}
&
{\mathfrak b}_1u_{\theta}
\, +\,
{\mathfrak b}_2t u_t
+{\mathfrak b}_3t^2 u_t
+{\mathfrak b}_4tu_{\theta}
+{\mathfrak b}_5t^2 u_{\theta}+{\bar D}_0^0(u),
\end{split}
\label{boundaryoriginal0}
\end{equation}
where
\begin{align}
{\mathfrak b}_1=\frac{1}{\sqrt{{\mathfrak h}_1(0)}}{\mathfrak y}_4(0)\mathfrak{p}_8(0),
\qquad
{\mathfrak b}_2=\frac{1}{\sqrt{{\mathfrak h}_1(0)}}\big[{\mathfrak y}_4(0)\mathfrak{p}_4(0)+{\mathfrak y}_2(0)\mathfrak{p}_1(0)\big],
\label{b1b2}
\\[2mm]
{\mathfrak b}_3=b_1(0),
\qquad
{\mathfrak b}_4=b_3(0),
\qquad
{\mathfrak b}_5=b_2(0),
\qquad
{\bar D}_0^0(u)\, =\, \sigma_3(t)\, u_t\, +\, \sigma_4(t)u_\theta.
\end{align}
On the other hand, for $\theta=1$, it has the form
\begin{equation}
\begin{split}
{\mathfrak b}_6u_{\theta}
\, +\,
{\mathfrak b}_7 t u_t
+{\mathfrak b}_8t^2 u_t
+{\mathfrak b}_9tu_{\theta}
+{\mathfrak b}_{10}t^2 u_{\theta}+{\bar D}_0^1(u),
\end{split}
\label{boundaryoriginal1}
\end{equation}
with the notation
\begin{align}
{\mathfrak b}_6=\frac{1}{\sqrt{{\mathfrak h}_1(1)}}{\mathfrak y}_4(1)\mathfrak{p}_8(1),
\qquad
{\mathfrak b}_7=\frac{1}{\sqrt{{\mathfrak h}_1(1)}}\big[{\mathfrak y}_4(1)\mathfrak{p}_4(1)+{\mathfrak y}_2(1)\mathfrak{p}_1(1)\big],
\label{b6b7}
\\[2mm]
{\mathfrak b}_8=b_1(1),
\qquad
{\mathfrak b}_9=b_3(1),
\qquad
{\mathfrak b}_{10}=b_2(1),
\qquad
{\bar D}_0^1(u)\, =\, \sigma_5(t)\, u_t\, +\, \sigma_6(t)u_\theta.
\end{align}
In the above, the functions $\sigma_3, \cdots, \sigma_6$ are smooth functions of $t$ with the properties
$$
|\sigma_i(t)|\leq C|t|^3, \quad i=3, 4, 5, 6.
$$

\subsection{Stationary and non-degenerate curves}\label{Stationary and non-degenerate curves}\

In the following, for a simple smooth curve $\Gamma$ connecting the boundary $\partial\Omega$,
we will make precisely the notion of a non-degenerate  geodesic with respect to
the metric ${\mathrm d}s^2= V^{2\sigma}(y)\big[{\mathfrak a}_2(y){\mathrm d}{y}_1^2+{\mathfrak a}_1(y){\mathrm d}{y}_2^2\big]$.
Consider the deformation of $\Gamma$ in the form
\begin{equation}
\label{deformation1}
\Gamma_{{\mathfrak t}}(\theta)=\big(\Gamma_{{\mathfrak t}1}(\theta), \Gamma_{{\mathfrak t}2}(\theta)\big): \gamma\big(\Theta({\mathfrak t}(\theta), \theta)\big)
 \, +\, {\mathfrak t}(\theta)\, \Big({\tilde{\mathfrak a}}_1(\theta)n_1\big(\Theta({\mathfrak t}(\theta), \theta)\big),
 {\tilde{\mathfrak a}}_2(\theta)n_2\big(\Theta({\mathfrak t}(\theta), \theta)\big)\Big),
\end{equation}
where $\theta\in [0, 1]$ and ${\mathfrak t}$ is a smooth function of $\theta$ with small $L^{\infty}$-norm.
Note that  the end points of  $\Gamma_{{\mathfrak t}}$ stay on $\partial \Omega$. We denote ${\mathfrak a}_1(F(t, \theta)), {\mathfrak a}_2(F(t, \theta))$ and $ V(F(t, \theta))$
as ${\mathfrak a}_1(t, \theta)$, ${\mathfrak a}_2(t, \theta)$, $ V(t, \theta)$.
The weighted length of the curve $\Gamma_{{\mathfrak t}}$ is given by the functional of ${\mathfrak t}$
\begin{equation}
\begin{split}
{\mathcal J}({\mathfrak t})&\equiv
\int_0^1 V^{\sigma}\big(\Gamma_{{\mathfrak t}}(\theta)\big)
\sqrt{{\mathfrak a}_2\left(\Gamma_{{\mathfrak t}}(\theta)\right)\big|\Gamma'_{{\mathfrak t}1}(\theta)\big|^2
+ {\mathfrak a}_1\left(\Gamma_{{\mathfrak t}}(\theta)\right)\big|\Gamma'_{{\mathfrak t}2}(\theta)\big|^2}\, {\rm d}\theta
\\[2mm]
&=\int_0^1  V^{\sigma}\big(\Gamma_{{\mathfrak t}}(\theta)\big) \sqrt{{\mathcal W}\big( {\mathfrak t}(\theta)\big)}
\, {\rm d}\theta.
\end{split}
\label{weightedlength1}
\end{equation}
where we have denoted
\begin{align}
{\mathcal W}\big( {\mathfrak t}(\theta)\big)
=& \Big(\Theta_t({\mathfrak t}(\theta), \theta){\mathfrak t}'(\theta)+\Theta_{\theta}({\mathfrak t}(\theta), \theta)\Big)^2
\nonumber\\[2mm]
&\qquad\times\Big\{\, \big[1-k(\theta){\mathfrak t}(\theta){\tilde{\mathfrak a}}_1\big]^2{\mathfrak a}_2\big(\Gamma_{{\mathfrak t}}(\theta)\big)\big|\gamma_1'\big(\Theta({\mathfrak t}(\theta), \theta)\big)\big|^2
\nonumber\\[2mm]
&\qquad\qquad+\big[1-k(\theta){\mathfrak t}(\theta){\tilde{\mathfrak a}}_2\big]^2{\mathfrak a}_1\big(\Gamma_{{\mathfrak t}}(\theta)\big)\big|\gamma_2'\big(\Theta({\mathfrak t}(\theta), \theta)\big)\big|^2\, \Big\}
\nonumber\\[2mm]
&\, +\, |{\mathfrak t}'(\theta)|^2\bigg[{\mathfrak a}_2\big(\Gamma_{{\mathfrak t}}(\theta)\big)\big|{\tilde{\mathfrak a}}_1(\theta)\big|^2\big|n_1\big(\Theta({\mathfrak t}(\theta), \theta)\big)\big|^2
+{\mathfrak a}_1\big(\Gamma_{{\mathfrak t}}(\theta)\big)\big|{\tilde{\mathfrak a}}_2(\theta)\big|^2\big|n_2\big(\Theta({\mathfrak t}(\theta), \theta)\big)\big|^2\bigg]
\nonumber\\[2mm]
&+2{\mathfrak t}'(\theta){\mathfrak t}(\theta)\bigg[{\mathfrak a}_2\big(\Gamma_{{\mathfrak t}}(\theta)\big){\tilde{\mathfrak a}}_1(\theta){\tilde{\mathfrak a}}_1'(\theta)|n_1\big(\Theta({\mathfrak t}(\theta), \theta)\big)|^2
+{\mathfrak a}_1\big(\Gamma_{{\mathfrak t}}(\theta)\big){\tilde{\mathfrak a}}_2(\theta){\tilde{\mathfrak a}}_2'(\theta)|n_2\big(\Theta({\mathfrak t}(\theta), \theta)\big)|^2\bigg]
\nonumber\\[2mm]
&+|{\mathfrak t}(\theta)|^2\bigg[{\mathfrak a}_2\big(\Gamma_{{\mathfrak t}}(\theta)\big)\big|{\tilde{\mathfrak a}}_1'(\theta)\big|^2|n_1\big(\Theta({\mathfrak t}(\theta), \theta)\big)|^2
+{\mathfrak a}_1\big(\Gamma_{{\mathfrak t}}(\theta)\big)\big|{\tilde{\mathfrak a}}_2'(\theta)\big|^2|n_2\big(\Theta({\mathfrak t}(\theta), \theta)\big)|^2\bigg]
\nonumber\\[2mm]
&\, +\, 2{\mathfrak t}'(\theta)\, \Big(\Theta_t({\mathfrak t}(\theta), \theta){\mathfrak t}'(\theta)+\Theta_{\theta}({\mathfrak t}(\theta), \theta)\Big)
\nonumber\\[2mm]
&\qquad\times
 \Big\{\big[1-k(\theta){\mathfrak t}(\theta){\tilde{\mathfrak a}}_1\big]{\mathfrak a}_2\big(\Gamma_{{\mathfrak t}}(\theta)\big){\tilde{\mathfrak a}}_1(\theta)n_1\big(\Theta({\mathfrak t}(\theta), \theta)\big)\gamma_1'\big(\Theta({\mathfrak t}(\theta), \theta)\big)
\nonumber\\[2mm]
&\qquad \qquad+\big[1-k(\theta){\mathfrak t}(\theta){\tilde{\mathfrak a}}_2\big]{\mathfrak a}_1\big(\Gamma_{{\mathfrak t}}(\theta)\big){\tilde{\mathfrak a}}_2(\theta)n_2\big(\Theta({\mathfrak t}(\theta), \theta)\big)\gamma_2'\big(\Theta({\mathfrak t}(\theta), \theta)\big)\Big\}
\nonumber\\[2mm]
&\, +\, 2{\mathfrak t}(\theta)\Big(\Theta_t({\mathfrak t}(\theta), \theta){\mathfrak t}'(\theta)+\Theta_{\theta}({\mathfrak t}(\theta), \theta)\Big)
\nonumber\\[2mm]
&\qquad\times
\, \Big\{\big[1-k(\theta){\mathfrak t}(\theta){\tilde{\mathfrak a}}_1\big]{\mathfrak a}_2\big(\Gamma_{{\mathfrak t}}(\theta)\big){\tilde{\mathfrak a}}_1'(\theta)n_1\big(\Theta({\mathfrak t}(\theta), \theta)\big)\gamma_1'\big(\Theta({\mathfrak t}(\theta), \theta)\big)
\nonumber\\[2mm]
&\qquad \qquad
+\big[1-k(\theta){\mathfrak t}(\theta){\tilde{\mathfrak a}}_2\big]{\mathfrak a}_1\big(\Gamma_{{\mathfrak t}}(\theta)\big){\tilde{\mathfrak a}}_2'(\theta)n_2\big(\Theta({\mathfrak t}(\theta), \theta)\big)\gamma_2'\big(\Theta({\mathfrak t}(\theta), \theta)\big)\big]\Big\}.
\label{mathcalM}
\end{align}

\medskip
\noindent{\bf{Step 1. }}
The first variation of ${\mathcal J}$ at ${\mathfrak t}$ along the direction $h$ is given by
\begin{align}
{\mathcal J}'({\mathfrak t})[h]
=\, \frac {\mathrm d}{{\mathrm d}s}{\mathcal J}({\mathfrak t}+sh)\Big|_{s=0}
&=\, \int_0^1 \sigma V^{\sigma-1}\big(\Gamma_{{\mathfrak t}}(\theta)\big)\,  V_t\big(\Gamma_{{\mathfrak t}}(\theta)\big)\, h\,
\, \sqrt{{\mathcal W}\big( {\mathfrak t}(\theta)\big)}
\, {\rm d}\theta
\nonumber\\[2mm]
&\quad\ +\
\frac{1}{2}\int_0^1\frac{ V^{\sigma}\big(\Gamma_{{\mathfrak t}}(\theta)\big)}{\sqrt{\mathcal W\big( {\mathfrak t}(\theta)\big)}}
\frac {\mathrm d}{{\mathrm d}s}{\mathcal W}({\mathfrak t}+sh)\Big|_{s=0}
\, {\rm d}\theta,
\end{align}
where
\begin{align}
\frac {\mathrm d}{{\mathrm d}s}{\mathcal W}({\mathfrak t}+sh)\Big|_{s=0}
\, =\, &{\mathcal W}_1\big( {\mathfrak t}(\theta)\big)[h]
\ +\
{\mathcal W}_2\big( {\mathfrak t}(\theta)\big)[h]
\ +\
{\mathcal W}_3\big( {\mathfrak t}(\theta)\big)[h]
\nonumber\\[2mm]
&\, +\, {\mathcal W}_4\big( {\mathfrak t}(\theta)\big)[h]
\, +\, {\mathcal W}_5\big( {\mathfrak t}(\theta)\big)[h]
\, +\, {\mathcal W}_6\big( {\mathfrak t}(\theta)\big)[h],
\end{align}
in which ${\mathcal W}_i\big( {\mathfrak t}(\theta)\big)[h], (i=1, \cdots, 6)$ are given by
\begin{align}
{\mathcal W}_1\big( {\mathfrak t}(\theta)\big)[h]
&= \Big(\Theta_t({\mathfrak t}(\theta), \theta){\mathfrak t}'(\theta)+\Theta_{\theta}({\mathfrak t}(\theta), \theta)\Big)^2
\nonumber\\[2mm]
&\qquad\quad\times\Big\{\, 2\big[1-k(\theta){\mathfrak t}(\theta){\tilde{\mathfrak a}}_1\big]\, (-k{\tilde{\mathfrak a}}_1)\, {\mathfrak a}_2\big(\Gamma_{{\mathfrak t}}(\theta)\big)\big|\gamma_1'\big(\Theta({\mathfrak t}(\theta), \theta)\big)\big|^2
\nonumber\\[2mm]
&\qquad \qquad+2\big[1-k(\theta){\mathfrak t}(\theta){\tilde{\mathfrak a}}_2\big]\, (-k{\tilde{\mathfrak a}}_2)\, {\mathfrak a}_1\big(\Gamma_{{\mathfrak t}}(\theta)\big)\big|\gamma_2'\big(\Theta({\mathfrak t}(\theta), \theta)\big)\big|^2\, \Big\}h
\nonumber\\[2mm]
&\quad\, +\, 2\Big(\Theta_t({\mathfrak t}(\theta), \theta){\mathfrak t}'(\theta)+\Theta_{\theta}({\mathfrak t}(\theta), \theta)\Big)
\Big(
\Theta_{tt}({\mathfrak t}(\theta), \theta)\, {\mathfrak t}'(\theta)\, h
\, +\, \Theta_{t}({\mathfrak t}(\theta), \theta)\, h'
\, +\, \Theta_{\theta t}({\mathfrak t}(\theta), \theta)\, h
\Big)
\nonumber\\[2mm]
&\qquad\quad\times\Big\{\, \big[1-k(\theta){\mathfrak t}(\theta){\tilde{\mathfrak a}}_1\big]^2{\mathfrak a}_2\big(\Gamma_{{\mathfrak t}}(\theta)\big)\big|\gamma_1'\big(\Theta({\mathfrak t}(\theta), \theta)\big)\big|^2
\nonumber\\[2mm]
&\qquad \qquad \quad+\big[1-k(\theta){\mathfrak t}(\theta){\tilde{\mathfrak a}}_2\big]^2{\mathfrak a}_1\big(\Gamma_{{\mathfrak t}}(\theta)\big)\big|\gamma_2'\big(\Theta({\mathfrak t}(\theta), \theta)\big)\big|^2\, \Big\}
\nonumber\\[2mm]
&\quad+\Big(\Theta_t({\mathfrak t}(\theta), \theta){\mathfrak t}'(\theta)+\Theta_{\theta}({\mathfrak t}(\theta), \theta)\Big)^2
\nonumber\\[2mm]
&\qquad\quad\times\Big\{\, \big[1-k(\theta){\mathfrak t}(\theta){\tilde{\mathfrak a}}_1\big]^2\partial_t{\mathfrak a}_2\big(\Gamma_{{\mathfrak t}}(\theta)\big)\big|\gamma_1'\big(\Theta({\mathfrak t}(\theta), \theta)\big)\big|^2
\nonumber\\[2mm]
&\qquad\quad\quad\, +\, \big[1-k(\theta){\mathfrak t}(\theta){\tilde{\mathfrak a}}_1\big]^2
2{\mathfrak a}_2\big(\Gamma_{{\mathfrak t}}(\theta)\big)\gamma_1'\big(\Theta({\mathfrak t}(\theta), \theta)\big)\gamma_1''\big(\Theta({\mathfrak t}(\theta), \theta)\big)
\Theta_t({\mathfrak t}(\theta), \theta)
\nonumber\\[2mm]
&\qquad\quad\quad\, +\,
\big[1-k(\theta){\mathfrak t}(\theta){\tilde{\mathfrak a}}_2\big]^2\partial_t{\mathfrak a}_1\big(\Gamma_{{\mathfrak t}}(\theta)\big)\big|\gamma_2'\big(\Theta({\mathfrak t}(\theta), \theta)\big)\big|^2
\nonumber\\[2mm]
&\qquad\quad\quad\, +\,
\big[1-k(\theta){\mathfrak t}(\theta){\tilde{\mathfrak a}}_2\big]^22{\mathfrak a}_1\big(\Gamma_{{\mathfrak t}}(\theta)\big)\gamma_2'\big(\Theta({\mathfrak t}(\theta), \theta)\big)\gamma_2''\big(\Theta({\mathfrak t}(\theta), \theta)\big)
\Theta_t({\mathfrak t}(\theta), \theta)\Big\}h,
\label{mathcalM1}
\end{align}
\begin{align}
{\mathcal W}_2\big( {\mathfrak t}(\theta)\big)[h]
&=2{\mathfrak t}'(\theta)\, \bigg[{\mathfrak a}_2\big(\Gamma_{{\mathfrak t}}(\theta)\big)\big|{\tilde{\mathfrak a}}_1(\theta)\big|^2\cdot\big|n_1\big(\Theta({\mathfrak t}(\theta), \theta)\big)\big|^2
+{\mathfrak a}_1\big(\Gamma_{{\mathfrak t}}(\theta)\big)\big|{\tilde{\mathfrak a}}_2(\theta)\big|^2\cdot\big|n_2\big(\Theta({\mathfrak t}(\theta), \theta)\big)\big|^2\bigg]\, h'
\nonumber\\[2mm]
&\quad\, +\,
|{\mathfrak t}'(\theta)|^2\Big\{\partial_t{\mathfrak a}_2\big(\Gamma_{{\mathfrak t}}(\theta)\big)\big|{\tilde{\mathfrak a}}_1(\theta)\big|^2\cdot\big|n_1\big(\Theta({\mathfrak t}(\theta), \theta)\big)\big|^2
\, +\, \partial_t{\mathfrak a}_1\big(\Gamma_{{\mathfrak t}}(\theta)\big)\big|{\tilde{\mathfrak a}}_2(\theta)\big|^2\cdot\big|n_2\big(\Theta({\mathfrak t}(\theta), \theta)\big)\big|^2
\nonumber\\[2mm]
&\qquad\qquad\qquad\, +\, {\mathfrak a}_2\big(\Gamma_{{\mathfrak t}}(\theta)\big)\big|{\tilde{\mathfrak a}}_1(\theta)\big|^2\,2 n_1\big(\Theta({\mathfrak t}(\theta), \theta)\big)n_1'\big(\Theta({\mathfrak t}(\theta), \theta)\big)\Theta_{t}({\mathfrak t}(\theta), \theta)
\nonumber\\[2mm]
&\qquad\qquad\qquad\, +\, {\mathfrak a}_1\big(\Gamma_{{\mathfrak t}}(\theta)\big)\big|{\tilde{\mathfrak a}}_2(\theta)\big|^2\,2 n_2\big(\Theta({\mathfrak t}(\theta), \theta)\big)n_2'\big(\Theta({\mathfrak t}(\theta), \theta)\big)\Theta_{t}({\mathfrak t}(\theta), \theta)\Big\}h,
\label{mathcalM2}
\end{align}
\begin{align}
{\mathcal W}_3\big( {\mathfrak t}(\theta)\big)[h]
&=2\, {\mathfrak t}(\theta)\Big[{\mathfrak a}_2\big(\Gamma_{{\mathfrak t}}(\theta)\big){\tilde{\mathfrak a}}_1(\theta){\tilde{\mathfrak a}}_1'(\theta)|n_1\big(\Theta({\mathfrak t}(\theta), \theta)\big)|^2
+{\mathfrak a}_1\big(\Gamma_{{\mathfrak t}}(\theta)\big){\tilde{\mathfrak a}}_2(\theta){\tilde{\mathfrak a}}_2'(\theta)|n_2\big(\Theta({\mathfrak t}(\theta), \theta)\big)|^2\Big]h'
\nonumber\\[2mm]
&\quad\, +\, 2{\mathfrak t}'(\theta)\Big[{\mathfrak a}_2\big(\Gamma_{{\mathfrak t}}(\theta)\big){\tilde{\mathfrak a}}_1(\theta){\tilde{\mathfrak a}}_1'(\theta)|n_1\big(\Theta({\mathfrak t}(\theta), \theta)\big)|^2
+{\mathfrak a}_1\big(\Gamma_{{\mathfrak t}}(\theta)\big){\tilde{\mathfrak a}}_2(\theta){\tilde{\mathfrak a}}_2'(\theta)|n_2\big(\Theta({\mathfrak t}(\theta), \theta)\big)|^2\Big]h
\nonumber\\[2mm]
&\quad\, +\, 2{\mathfrak t}'(\theta){\mathfrak t}(\theta)
\Big\{\partial_t{\mathfrak a}_2\big(\Gamma_{{\mathfrak t}}(\theta)\big){\tilde{\mathfrak a}}_1(\theta){\tilde{\mathfrak a}}_1'(\theta)|n_1\big(\Theta({\mathfrak t}(\theta), \theta)\big)|^2
\, +\, \partial_t{\mathfrak a}_1\big(\Gamma_{{\mathfrak t}}(\theta)\big){\tilde{\mathfrak a}}_2(\theta){\tilde{\mathfrak a}}_2'(\theta)|n_2\big(\Theta({\mathfrak t}(\theta), \theta)\big)|^2
\nonumber\\[2mm]
&\qquad \qquad \qquad \quad\, +\, {\mathfrak a}_2\big(\Gamma_{{\mathfrak t}}(\theta)\big){\tilde{\mathfrak a}}_1(\theta){\tilde{\mathfrak a}}_1'(\theta)2n_1\big(\Theta({\mathfrak t}(\theta), \theta)n_1'\big(\Theta({\mathfrak t}(\theta), \theta)\Theta_t({\mathfrak t}(\theta))
\nonumber\\[2mm]
&\qquad \qquad \qquad \quad\, +\, {\mathfrak a}_1\big(\Gamma_{{\mathfrak t}}(\theta)\big){\tilde{\mathfrak a}}_2(\theta){\tilde{\mathfrak a}}_2'(\theta)2n_2\big(\Theta({\mathfrak t}(\theta), \theta)\big)n_2'\big(\Theta({\mathfrak t}(\theta), \theta)\big)\Theta_t({\mathfrak t}(\theta))\Big\}h,
\label{mathcalM3}
\end{align}
\begin{align}
{\mathcal W}_4\big( {\mathfrak t}(\theta)\big)[h]
&=2{\mathfrak t}(\theta) \Big[{\mathfrak a}_2\big(\Gamma_{{\mathfrak t}}(\theta)\big)\big|{\tilde{\mathfrak a}}_1'(\theta)\big|^2|n_1\big(\Theta({\mathfrak t}(\theta), \theta)\big)|^2
+{\mathfrak a}_1\big(\Gamma_{{\mathfrak t}}(\theta)\big)\big|{\tilde{\mathfrak a}}_2'(\theta)\big|^2|n_2\big(\Theta({\mathfrak t}(\theta), \theta)\big)|^2\Big]h
\nonumber\\[2mm]
&\quad\, +\, |{\mathfrak t}(\theta)|^2
\Big\{\partial_t{\mathfrak a}_2\big(\Gamma_{{\mathfrak t}}(\theta)\big)\big|{\tilde{\mathfrak a}}_1'(\theta)\big|^2|n_1\big(\Theta({\mathfrak t}(\theta), \theta)\big)|^2
\, +\, \partial_t{\mathfrak a}_1\big(\Gamma_{{\mathfrak t}}(\theta)\big)\big|{\tilde{\mathfrak a}}_2'(\theta)\big|^2|n_2\big(\Theta({\mathfrak t}(\theta), \theta)\big)|^2
\nonumber\\[2mm]
&\qquad \qquad \qquad \, +\, {\mathfrak a}_2\big(\Gamma_{{\mathfrak t}}(\theta)\big)\big|{\tilde{\mathfrak a}}_1'(\theta)\big|^22n_1\big(\Theta({\mathfrak t}(\theta), \theta)\big)n_1'\big(\Theta({\mathfrak t}(\theta), \theta)\big)\Theta_t({\mathfrak t}(\theta), \theta)
\nonumber\\[2mm]
&\qquad \qquad \qquad \, +\, {\mathfrak a}_1\big(\Gamma_{{\mathfrak t}}(\theta)\big)\big|{\tilde{\mathfrak a}}_2'(\theta)\big|^22n_2\big(\Theta({\mathfrak t}(\theta), \theta)\big)n_2'\big(\Theta({\mathfrak t}(\theta), \theta)\big)\Theta_t({\mathfrak t}(\theta), \theta)
\Big\}h,
\label{mathcalM4}
\end{align}
\begin{align}
{\mathcal W}_5\big( {\mathfrak t}(\theta)\big)[h]
&=2\, \Big(\Theta_t({\mathfrak t}(\theta), \theta){\mathfrak t}'(\theta)+\Theta_{\theta}({\mathfrak t}(\theta), \theta)\Big)
\nonumber\\[2mm]
&\qquad\times
\, \Big\{\big[1-k(\theta){\mathfrak t}(\theta){\tilde{\mathfrak a}}_1\big]{\mathfrak a}_2\big(\Gamma_{{\mathfrak t}}(\theta)\big){\tilde{\mathfrak a}}_1(\theta)n_1\big(\Theta({\mathfrak t}(\theta), \theta)\big)\gamma_1'\big(\Theta({\mathfrak t}(\theta), \theta)\big)
\nonumber\\[2mm]
&\qquad \qquad
+\big[1-k(\theta){\mathfrak t}(\theta){\tilde{\mathfrak a}}_2\big]{\mathfrak a}_1\big(\Gamma_{{\mathfrak t}}(\theta)\big){\tilde{\mathfrak a}}_2(\theta)n_2\big(\Theta({\mathfrak t}(\theta), \theta)\big)\gamma_2'\big(\Theta({\mathfrak t}(\theta), \theta)\big)\Big\}h'
\nonumber\\[2mm]
&\quad\, +\, 2{\mathfrak t}'(\theta)\, \Big(
\Theta_{tt}({\mathfrak t}(\theta), \theta)\, {\mathfrak t}'(\theta)\, h
\, +\, \Theta_{t}({\mathfrak t}(\theta), \theta)\, h'
\, +\, \Theta_{\theta t}({\mathfrak t}(\theta), \theta)\, h
\Big)
\nonumber\\[2mm]
&\qquad\times
\, \Big\{\big[1-k(\theta){\mathfrak t}(\theta){\tilde{\mathfrak a}}_1\big]{\mathfrak a}_2\big(\Gamma_{{\mathfrak t}}(\theta)\big){\tilde{\mathfrak a}}_1(\theta)n_1\big(\Theta({\mathfrak t}(\theta), \theta)\big)\gamma_1'\big(\Theta({\mathfrak t}(\theta), \theta)\big)
\nonumber\\[2mm]
&\qquad \qquad
+\big[1-k(\theta){\mathfrak t}(\theta){\tilde{\mathfrak a}}_2\big]{\mathfrak a}_1\big(\Gamma_{{\mathfrak t}}(\theta)\big){\tilde{\mathfrak a}}_2(\theta)n_2\big(\Theta({\mathfrak t}(\theta), \theta)\big)\gamma_2'\big(\Theta({\mathfrak t}(\theta), \theta)\big)\Big\}h
\nonumber\\[2mm]
&\quad\, +\, 2{\mathfrak t}'(\theta)\, \Big(\Theta_t({\mathfrak t}(\theta), \theta){\mathfrak t}'(\theta)+\Theta_{\theta}({\mathfrak t}(\theta), \theta)\Big)
\nonumber\\[2mm]
&\qquad\times
\, \Big\{\big[-k(\theta){\tilde{\mathfrak a}}_1\big]{\mathfrak a}_2\big(\Gamma_{{\mathfrak t}}(\theta)\big){\tilde{\mathfrak a}}_1(\theta)n_1\big(\Theta({\mathfrak t}(\theta), \theta)\big)\gamma_1'\big(\Theta({\mathfrak t}(\theta), \theta)\big)
\nonumber\\[2mm]
&\qquad\qquad +\big[1-k(\theta){\mathfrak t}(\theta){\tilde{\mathfrak a}}_1\big]\partial_t{\mathfrak a}_2\big(\Gamma_{{\mathfrak t}}(\theta)\big){\tilde{\mathfrak a}}_1(\theta)n_1\big(\Theta({\mathfrak t}(\theta), \theta)\big)\gamma_1'\big(\Theta({\mathfrak t}(\theta), \theta)\big)
\nonumber\\[2mm]
&\qquad\qquad -\big[1-k(\theta){\mathfrak t}(\theta){\tilde{\mathfrak a}}_1\big]{\mathfrak a}_2\big(\Gamma_{{\mathfrak t}}(\theta)\big){\tilde{\mathfrak a}}_1(\theta)n_1'\big(\Theta({\mathfrak t}(\theta), \theta)\big)\Theta_t({\mathfrak t}(\theta))n_2\big(\Theta({\mathfrak t}(\theta), \theta)\big)
\nonumber\\[2mm]
&\qquad\qquad -\big[1-k(\theta){\mathfrak t}(\theta){\tilde{\mathfrak a}}_1\big]{\mathfrak a}_2\big(\Gamma_{{\mathfrak t}}(\theta)\big){\tilde{\mathfrak a}}_1(\theta)n_1\big(\Theta({\mathfrak t}(\theta), \theta)\big)n_2'\big(\Theta({\mathfrak t}(\theta), \theta)\big)\Theta_t({\mathfrak t}(\theta), \theta)
\nonumber\\[2mm]
&\qquad \qquad
+\big[-k(\theta){\tilde{\mathfrak a}}_2\big]{\mathfrak a}_1\big(\Gamma_{{\mathfrak t}}(\theta)\big){\tilde{\mathfrak a}}_2(\theta)n_2\big(\Theta({\mathfrak t}(\theta), \theta)\big)\gamma_2'\big(\Theta({\mathfrak t}(\theta), \theta)\big)
\nonumber\\[2mm]
&\qquad \qquad
+\big[1-k(\theta){\mathfrak t}(\theta){\tilde{\mathfrak a}}_2\big]{\partial_t\mathfrak a}_1\big(\Gamma_{{\mathfrak t}}(\theta)\big){\tilde{\mathfrak a}}_2(\theta)n_2\big(\Theta({\mathfrak t}(\theta), \theta)\big)\gamma_2'\big(\Theta({\mathfrak t}(\theta), \theta)\big)
\nonumber\\[2mm]
&\qquad \qquad
+\big[1-k(\theta){\mathfrak t}(\theta){\tilde{\mathfrak a}}_2\big]{\mathfrak a}_1\big(\Gamma_{{\mathfrak t}}(\theta)\big){\tilde{\mathfrak a}}_2(\theta)n_2'\big(\Theta({\mathfrak t}(\theta), \theta)\big)n_1\big(\Theta({\mathfrak t}(\theta), \theta)\big)\Theta_t({\mathfrak t}(\theta), \theta)
\nonumber\\[2mm]
&\qquad \qquad
+\big[1-k(\theta){\mathfrak t}(\theta){\tilde{\mathfrak a}}_2\big]{\mathfrak a}_1\big(\Gamma_{{\mathfrak t}}(\theta)\big){\tilde{\mathfrak a}}_2(\theta)n_2\big(\Theta({\mathfrak t}(\theta), \theta)\big)n_2'\big(\Theta({\mathfrak t}(\theta), \theta)\big)\Theta({\mathfrak t}(\theta), \theta)
\Big\}h,
\label{mathcalM5}
\end{align}
and
\begin{align}
{\mathcal W}_6\big( {\mathfrak t}(\theta)\big)[h]
&=2\, \Big(\Theta_t({\mathfrak t}(\theta), \theta){\mathfrak t}'(\theta)+\Theta_{\theta}({\mathfrak t}(\theta), \theta)\Big)
\nonumber\\[2mm]
&\qquad\times
\, \Big\{\big[1-k(\theta){\mathfrak t}(\theta){\tilde{\mathfrak a}}_1\big]{\mathfrak a}_2\big(\Gamma_{{\mathfrak t}}(\theta)\big){\tilde{\mathfrak a}}_1'(\theta)n_1\big(\Theta({\mathfrak t}(\theta), \theta)\big)\gamma_1'\big(\Theta({\mathfrak t}(\theta), \theta)\big)
\nonumber\\[2mm]
&\qquad \qquad
+\big[1-k(\theta){\mathfrak t}(\theta){\tilde{\mathfrak a}}_2\big]{\mathfrak a}_1\big(\Gamma_{{\mathfrak t}}(\theta)\big){\tilde{\mathfrak a}}_2'(\theta)n_2\big(\Theta({\mathfrak t}(\theta), \theta)\big)\gamma_2'\big(\Theta({\mathfrak t}(\theta), \theta)\big)\Big\}h
\nonumber\\[2mm]
&+2{\mathfrak t}(\theta)\,
\Big(
\Theta_{tt}({\mathfrak t}(\theta), \theta)\, {\mathfrak t}'(\theta)\, h
\, +\, \Theta_{t}({\mathfrak t}(\theta), \theta)\, h'
\, +\, \Theta_{\theta t}({\mathfrak t}(\theta), \theta)\, h
\Big)
\nonumber\\[2mm]
&\qquad\times
\, \Big\{\big[1-k(\theta){\mathfrak t}(\theta){\tilde{\mathfrak a}}_1\big]{\mathfrak a}_2\big(\Gamma_{{\mathfrak t}}(\theta)\big){\tilde{\mathfrak a}}_1'(\theta)n_1\big(\Theta({\mathfrak t}(\theta), \theta)\big)\gamma_1'\big(\Theta({\mathfrak t}(\theta), \theta)\big)
\nonumber\\[2mm]
&\qquad \qquad
+\big[1-k(\theta){\mathfrak t}(\theta){\tilde{\mathfrak a}}_2\big]{\mathfrak a}_1\big(\Gamma_{{\mathfrak t}}(\theta)\big){\tilde{\mathfrak a}}_2'(\theta)n_2\big(\Theta({\mathfrak t}(\theta), \theta)\big)\gamma_2'\big(\Theta({\mathfrak t}(\theta), \theta)\big)\Big\}
\nonumber\\[2mm]
&+2{\mathfrak t}(\theta)\, \Big(\Theta_t({\mathfrak t}(\theta), \theta){\mathfrak t}'(\theta)+\Theta_{\theta}({\mathfrak t}(\theta), \theta)\Big)
\nonumber\\[2mm]
&\qquad\times
\, \Big\{\big[-k(\theta){\tilde{\mathfrak a}}_1\big]{\mathfrak a}_2\big(\Gamma_{{\mathfrak t}}(\theta)\big){\tilde{\mathfrak a}}_1'(\theta)n_1\big(\Theta({\mathfrak t}(\theta), \theta)\big)\gamma_1'\big(\Theta({\mathfrak t}(\theta), \theta)\big)
\nonumber\\[2mm]
&\qquad \qquad
+\big[1-k(\theta){\mathfrak t}(\theta){\tilde{\mathfrak a}}_1\big]{\partial_t\mathfrak a}_2\big(\Gamma_{{\mathfrak t}}(\theta)\big){\tilde{\mathfrak a}}_1'(\theta)n_1\big(\Theta({\mathfrak t}(\theta), \theta)\big)\gamma_1'\big(\Theta({\mathfrak t}(\theta), \theta)\big)
\nonumber\\[2mm]
&\qquad \qquad
-\big[1-k(\theta){\mathfrak t}(\theta){\tilde{\mathfrak a}}_1\big]{\mathfrak a}_2\big(\Gamma_{{\mathfrak t}}(\theta)\big){\tilde{\mathfrak a}}_1'(\theta)n_1'\big(\Theta({\mathfrak t}(\theta), \theta)\big)n_2\big(\Theta({\mathfrak t}(\theta), \theta)\big)\Theta_t({\mathfrak t}(\theta), \theta)
\nonumber\\[2mm]
&\qquad  \qquad
+\big[-k(\theta){\tilde{\mathfrak a}}_2\big]{\mathfrak a}_1\big(\Gamma_{{\mathfrak t}}(\theta)\big){\tilde{\mathfrak a}}_2'(\theta)n_2\big(\Theta({\mathfrak t}(\theta), \theta)\big)\gamma_2'\big(\Theta({\mathfrak t}(\theta), \theta)\big)
\nonumber\\[2mm]
&\qquad \qquad
+\big[1-k(\theta){\mathfrak t}(\theta){\tilde{\mathfrak a}}_2\big]\partial_t{\mathfrak a}_1\big(\Gamma_{{\mathfrak t}}(\theta)\big){\tilde{\mathfrak a}}_2'(\theta)n_2\big(\Theta({\mathfrak t}(\theta), \theta)\big)\gamma_2'\big(\Theta({\mathfrak t}(\theta), \theta)\big)
\nonumber\\[2mm]
&\qquad \qquad
+\big[1-k(\theta){\mathfrak t}(\theta){\tilde{\mathfrak a}}_2\big]{\mathfrak a}_1\big(\Gamma_{{\mathfrak t}}(\theta)\big){\tilde{\mathfrak a}}_2'(\theta)n_2\big(\Theta({\mathfrak t}(\theta), \theta)\big)n_1'\big(\Theta({\mathfrak t}(\theta), \theta)\big)\Theta_t({\mathfrak t}(\theta), \theta)\Big\}h.
\label{mathcalM6}
\end{align}

By substituting the relations in \eqref{Thetaderivative1}-\eqref{Thetaderivative2} into \eqref{mathcalM1}-\eqref{mathcalM6}, we get
\begin{align}
{\mathcal W}(0)\, =\, &{\mathfrak a}_1(0, \theta)\big|n_1\big|^2+{\mathfrak a}_2(0, \theta)\big|n_2\big|^2={\mathfrak f}_0(\theta),
\label{mathcalM(0)}
\\[2mm]
{\mathcal W}_1(0)[h]\, =\, &-2\, k \bigg[\, {\tilde{\mathfrak a}}_1\, {\mathfrak a}_2(0, \theta)\big|n_2\big|^2
+{\tilde{\mathfrak a}}_2\, {\mathfrak a}_1(0, \theta)\big|n_1\big|^2\, \bigg]h
\nonumber\\[2mm]
&+\bigg[\partial_t{\mathfrak a}_2(0, \theta)\big|n_2\big|^2
\, +\, \partial_t{\mathfrak a}_1(0, \theta)\big|n_1\big|^2\bigg]\, h,
\\[2mm]
{\mathcal W}_2(0)[h]\, =\, &0,
\qquad
{\mathcal W}_3(0)[h]=0,
\qquad
{\mathcal W}_4(0)[h]=0,
\label{mathcalM23(0)}
\\[2mm]
{\mathcal W}_5(0)[h]
\, =\, &2h'\big[-{\mathfrak a}_2(0, \theta){\tilde{\mathfrak a}}_1(\theta)
+{\mathfrak a}_1(0, \theta){\tilde{\mathfrak a}}_2(\theta)\big]n_1n_2\, =\, 0,
\\[2mm]
{\mathcal W}_6(0)[h]\, =\, &
2h\, \big[-{\mathfrak a}_2(0, \theta){\tilde{\mathfrak a}}_1'(\theta)
+{\mathfrak a}_1(0, \theta){\tilde{\mathfrak a}}_2'(\theta)\big]n_1n_2.
\label{mathcalM6(0)}
\end{align}
From the definitions of ${\mathfrak f}_1$ as in \eqref{m22}, we know that
\begin{equation*}
{\mathcal W}_1(0)[h]+{\mathcal W}_2(0)[h]+{\mathcal W}_3(0)[h]+{\mathcal W}_4(0)[h]+{\mathcal W}_5(0)[h]+{\mathcal W}_6(0)[h]=  {\mathfrak f}_1.
\end{equation*}
The above computations give that
\begin{align}
{\mathcal J}'(0)[h]
=&\int_0^1
\Bigg\{
\frac{\sigma V_t(0, \theta)}{ V^{1-\sigma}(0, \theta)}\, \sqrt{{\mathfrak f}_0}
+\frac{1}{2}\frac{ V^{\sigma}(0, \theta)}{\, \sqrt{{\mathfrak f}_0}} {\mathfrak f}_1\Bigg\}
\, h\, {\rm d}\theta.
\end{align}

The curve  $\Gamma$ is said to be {\bf stationary} with respect to the weighted length in \eqref{weightedlength1}
if the first variation of
${\mathcal J}$ at ${\mathfrak t}=0$ is equal to zero.
That is, for any smooth function $h(\theta)$ defined at $[0, 1]$ there holds
 \begin{equation*}
 \begin{split}
{\mathcal J}'(0)[h]=0.
 \end{split}
 \end{equation*}
This is equivalent to the relation
\begin{equation}
\label{stationary}
\frac{\sigma V_t(0, \theta)}{ V^{1-\sigma}(0, \theta)}\, \sqrt{{\mathfrak f}_0}
+\frac{1}{2}\frac{ {\mathfrak f}_1}{\, \sqrt{{\mathfrak f}_0}} V^{\sigma}(0, \theta)
=0,
\quad\forall\, \theta\in(0, 1),
\end{equation}
where ${\mathfrak f}_0$ and ${\mathfrak f}_1$ are given in \eqref{m22}.
Specially, if the parameter $\sigma$ in \eqref{weightedlength1} is $\frac{p+1}{p-1}-\frac{1}{2}$, then \eqref{stationary} has the form
\begin{equation}
\label{stationary1}
\Big(\frac{p+1}{p-1}-\frac{1}{2}\Big) V_t(0, \theta)\, \sqrt{{\mathfrak f}_0}
+\frac{1}{2}\frac{ V(0, \theta)}{\, \sqrt{{\mathfrak f}_0}} {\mathfrak f}_1
=0,
\quad\forall\, \theta\in(0, 1),
\end{equation}
i.e.,
\begin{align}\label{stationary2}
k\, =\, &\frac{{\mathfrak a}_1(0, \theta){\tilde{\mathfrak a}'}_2 n_1 n_2 -{\mathfrak a}_2(0, \theta){\tilde{\mathfrak a}'}_1 n_2 n_1}{{\mathfrak a}_1(0, \theta) {\tilde{\mathfrak a}}_2|n_1|^2 +{\mathfrak a}_2(0, \theta){\tilde{\mathfrak a}}_1 |n_2|^2}
+\frac{1}{2}\frac{\partial_{t} {\mathfrak a}_1(0, \theta) |n_1|^2 +\partial_{t} {\mathfrak a}_2(0, \theta) |n_2|^2}{{\mathfrak a}_1(0, \theta) {\tilde{\mathfrak a}}_2|n_1|^2 +{\mathfrak a}_2(0, \theta){\tilde{\mathfrak a}}_1 |n_2|^2}
\nonumber\\[2mm]
&+\sigma\frac{ V_t(0, \theta)}{ V(0, \theta)} \frac{{\mathfrak a}_1(0, \theta) |n_1|^2 +{\mathfrak a}_2(0, \theta) |n_2|^2}{{\mathfrak a}_1(0, \theta) {\tilde{\mathfrak a}}_2|n_1|^2 +{\mathfrak a}_2(0, \theta){\tilde{\mathfrak a}}_1 |n_2|^2},
\qquad  \forall\, \theta\in(0, 1).
\end{align}

\medskip
\noindent{\bf{Step 2. }}
We now consider the second variation of ${\mathcal J}$
\begin{align}
{\mathcal J}''(0)[h, f]
=\, &\frac {\mathrm d}{{\mathrm d}s}{\mathcal J}'(0+sf)[h]\Big|_{s=0}
\nonumber\\[2mm]
=\, &\int_0^1
\left[\,
\frac{\sigma V_{tt}(0, \theta)\, }{ V^{1-\sigma}\big(0, \theta\big)}
+
\frac{\sigma(\sigma-1)\big| V_t(0, \theta)\big|^2\, }{ V^{2-\sigma}(0, \theta)}
\, \right]
\, \sqrt{{\mathcal W}(0)}
\, hf\, {\rm d}\theta
\nonumber\\[2mm]
&\, +\,
\frac{\sigma}{2}\int_0^1\frac{ V_t(0, \theta) h\, }{ V^{1-\sigma}\big(0, \theta\big)}
\, \frac{1}{\sqrt{{\mathcal W}(0)}}
\frac {\mathrm d}{{\mathrm d}s}{\mathcal W}(sf)\Big|_{s=0}
\, {\rm d}\theta
\nonumber\\[2mm]
& \, +\,
 \frac{\sigma}{2}\int_0^1\frac{\, V_t(0, \theta)f\, }{ V^{1-\sigma}\big(0, \theta\big)}
\, \frac{1}{\sqrt{{\mathcal W}(0)}}
\frac {\mathrm d}{{\mathrm d}s}{\mathcal W}(sh)\Big|_{s=0}
 \, {\rm d}\theta
\nonumber\\[2mm]
& \,-\,
 \frac{1}{4}\int_0^1\frac{ V^{\sigma}\big(0, \theta\big)\, }{\ \big(\sqrt{{\mathcal W}(0)}\big)^3\ }
\frac {\mathrm d}{{\mathrm d}s}{\mathcal W}(sf)\Big|_{s=0}
\frac {\mathrm d}{{\mathrm d}s}{\mathcal W}(sh)\Big|_{s=0}
 \, {\rm d}\theta
\nonumber\\[2mm]
&\, +\,
\frac{1}{2}\int_0^1\frac{ V^{\sigma}\big(0, \theta\big)}{\sqrt{\mathcal W(0)}}
\frac {\mathrm d}{{\mathrm d}s}\big[ \sum_{i=1}^{6}{\mathcal W}_i\big(sf\big)[h]\big]\Big|_{s=0}
\, {\rm d}\theta.
\end{align}

From the definitions of ${\mathcal W}(f)$ and ${\mathfrak f}_1$  in \eqref{mathcalM}, \eqref{m22},  we can obtain that
\begin{align*}
\frac {\mathrm d}{{\mathrm d}s}{\mathcal W}(sf)\Big|_{s=0}
\, =\, &-2k \big[{\mathfrak a}_2(0, \theta){\tilde{\mathfrak a}}_1|n_2|^2+{\mathfrak a}_1(0, \theta){\tilde{\mathfrak a}}_2|n_1|^2\big]f
+\big[\partial_t{\mathfrak a}_2(0, \theta)|n_2|^2+\partial_t{\mathfrak a}_1(0, \theta)|n_1|^2\big]f
\\[2mm]
&+2\big[{\mathfrak a}_2(0, \theta){\tilde{\mathfrak a}}_1 n_1\gamma_1'
+{\mathfrak a}_1(0, \theta){\tilde{\mathfrak a}}_2 n_2\gamma_2' \big]f'
+2\big[{\mathfrak a}_2(0, \theta){\tilde{\mathfrak a}}_1' n_1\gamma_1'
+{\mathfrak a}_1(0, \theta){\tilde{\mathfrak a}}_2' n_2\gamma_2' \big]f
\\[2mm]
\, =\, &-2k \big[{\mathfrak a}_2(0, \theta){\tilde{\mathfrak a}}_1|n_2|^2+{\mathfrak a}_1(0, \theta){\tilde{\mathfrak a}}_2|n_1|^2\big]f
+ \big[\partial_t{\mathfrak a}_2(0, \theta)|n_2|^2+\partial_t{\mathfrak a}_1(0, \theta)|n_1|^2\big]f
\\[2mm]
&+2\big[-{\mathfrak a}_2(0, \theta){\tilde{\mathfrak a}}_1'
+{\mathfrak a}_1(0, \theta){\tilde{\mathfrak a}}_2'  \big]n_1n_2f
\\[2mm]
\, =\, &{\mathfrak f}_1f.
\end{align*}
Moreover,
\begin{align*}
 \frac {\mathrm d}{{\mathrm d}s}{\mathcal W}(sh)\Big|_{s=0}
={\mathfrak f}_1h,
\qquad
\frac {\mathrm d}{{\mathrm d}s}{\mathcal W}(sf)\Big|_{s=0}  \frac {\mathrm d}{{\mathrm d}s}{\mathcal W}(sh)\Big|_{s=0}
={\mathfrak f}_1^2hf.
\end{align*}
On the other hand, we use \eqref{mathcalM1}-\eqref{mathcalM6} together with the relations in \eqref{Thetaderivative1}-\eqref{Thetaderivative2} to derive the following
\begin{align}
&\frac {\mathrm d}{{\mathrm d}s} {\mathcal W}_1\big(sf\big)[h] \Big|_{s=0}
\nonumber\\[2mm]
&= 2\, k^2\big[\, {\mathfrak a}_2(0, \theta){\tilde{\mathfrak a}}_1|n_2|^2+{\mathfrak a}_1(0, \theta){\tilde{\mathfrak a}}_2|n_1|^2\, \big]\, fh
 -4\, k\big[\, \partial_t{\mathfrak a}_2(0, \theta){\tilde{\mathfrak a}}_1|n_2|^2+\partial_t{\mathfrak a}_1(0, \theta){\tilde{\mathfrak a}}_2|n_1|^2\, \big]\, fh
\nonumber\\[2mm]
&\quad\, +\, 2\big[\Theta_{tt}(0, \theta)f'h+\Theta_{tt}(0, \theta)fh'+\Theta_{\theta tt}(0, \theta)fh\big]
\big[{\mathfrak a}_2(0, \theta)|n_2|^2+{\mathfrak a}_1(0, \theta)|n_1|^2\big]
\nonumber\\[2mm]
&\quad\, +\, \big[\partial_{tt}{\mathfrak a}_2(0, \theta)|n_2|^2+\partial_{tt}{\mathfrak a}_1(0, \theta)|n_1|^2\big]fh
\, +\, 2\Theta_{tt}(0, \theta)\big[{\mathfrak a}_2(0, \theta)\gamma_1'\gamma_1''+{\mathfrak a}_1(0, \theta)\gamma_2'\gamma_2''\big]fh,
\label{mathcalM1deri}
\end{align}
\begin{align}
\frac {\mathrm d}{{\mathrm d}s} {\mathcal W}_2\big(sf\big)[h] \Big|_{s=0}
\, =\, 2f'\, h'\, \big[{\mathfrak a}_1(0, \theta){\tilde{\mathfrak a}}_2^2 |n_2|^2+{\mathfrak a}_2(0, \theta){\tilde{\mathfrak a}}_1^2 |n_1|^2\big]
\, =\, 2f'h'\, {\mathfrak w}_0,
\label{mathcalM2deri}
\end{align}
\begin{align}
\frac {\mathrm d}{{\mathrm d}s} {\mathcal W}_3\big(sf\big)[h] \Big|_{s=0}
\, =\, &2\, (h\, f'+h'f)
\, \big[{\mathfrak a}_2(0, \theta){\tilde{\mathfrak a}}_1{\tilde{\mathfrak a}}_1'|n_1|^2
 + {\mathfrak a}_1(0, \theta){\tilde{\mathfrak a}}_2{\tilde{\mathfrak a}}_2'|n_2|^2\big],
\label{mathcalM3deri}
\\[3mm]
\frac {\mathrm d}{{\mathrm d}s} {\mathcal W}_4\big(sf\big)[h] \Big|_{s=0}
\, =\, &2fh\big[{\mathfrak a}_2(0, \theta)|{\tilde{\mathfrak a}}_1'|^2|n_1|^2
 + {\mathfrak a}_1(0, \theta)|{\tilde{\mathfrak a}}_2'|^2|n_2|^2\big],
\label{mathcalM4deri}
\\[3mm]
\frac {\mathrm d}{{\mathrm d}s} {\mathcal W}_5\big(sf\big)[h] \Big|_{s=0}
\, =\, &-2k\, (h\, f'+h'f)
\, \big[{\mathfrak a}_2(0, \theta)|{\tilde{\mathfrak a}}_1|^2\gamma_1'n_1
+ {\mathfrak a}_1(0, \theta)|{\tilde{\mathfrak a}}_2|^2\gamma_2'n_2\big]
\nonumber\\[2mm]
&\, +\, 2(h\, f'+h'f)\, \big[\partial_t{\mathfrak a}_2(0, \theta){\tilde{\mathfrak a}}_1\gamma_1'n_1
\, +\, \partial_t{\mathfrak a}_1(0, \theta){\tilde{\mathfrak a}}_2\gamma_2'n_2\big]
\nonumber\\[2mm]
\, =\, &-2k\, (h\, f'+h'f)
\, \big[-{\mathfrak a}_2(0, \theta)|{\tilde{\mathfrak a}}_1|^2
+ {\mathfrak a}_1(0, \theta)|{\tilde{\mathfrak a}}_2|^2\big]n_1n_2
\nonumber\\[2mm]
&\, +\, 2(h\, f'+h'f)\, \big[-\partial_t{\mathfrak a}_2(0, \theta){\tilde{\mathfrak a}}_1
+ \partial_t{\mathfrak a}_1(0, \theta){\tilde{\mathfrak a}}_2\big]n_1n_2,
\label{mathcalM5deri}
\end{align}
and
\begin{align}
\frac {\mathrm d}{{\mathrm d}s} {\mathcal W}_6\big(sf\big)[h] \Big|_{s=0}
&=-4k\, \big[{\mathfrak a}_2(0, \theta){\tilde{\mathfrak a}}_1{\tilde{\mathfrak a}}_1'\gamma_1'n_1
+ {\mathfrak a}_1(0, \theta){\tilde{\mathfrak a}}_2{\tilde{\mathfrak a}}_2'\gamma_2'n_2\big]fh
\nonumber\\[2mm]
&\quad\, +\, 4\big[\partial_t{\mathfrak a}_2(0, \theta){\tilde{\mathfrak a}}_1'\gamma_1'n_1
+\partial_t{\mathfrak a}_1(0, \theta){\tilde{\mathfrak a}}_2'\gamma_2'n_2\big]fh
\nonumber\\
&=-4k\, \big[-{\mathfrak a}_2(0, \theta){\tilde{\mathfrak a}}_1{\tilde{\mathfrak a}}_1'
+ {\mathfrak a}_1(0, \theta){\tilde{\mathfrak a}}_2{\tilde{\mathfrak a}}_2'\big]n_1n_2fh
\nonumber\\[2mm]
&\quad\, +\, 4\big[-\partial_t{\mathfrak a}_2(0, \theta){\tilde{\mathfrak a}}_1'
+\partial_t{\mathfrak a}_1(0, \theta){\tilde{\mathfrak a}}_2'\big]n_1n_2 fh,
\label{mathcalM6deri}
\end{align}
where ${\mathfrak w}_0$ is given in \eqref{m11}.
Hence, we obtain
\begin{align}
&\frac {\mathrm d}{{\mathrm d}s}\big[ \sum_{i=1}^{6}{\mathcal W}_i\big(sf\big)[h] \big]\Big|_{s=0}
\nonumber\\[2mm]
&=\Big\{
-4k\big[-{\mathfrak a}_2(0, \theta){\tilde{\mathfrak a}}_1{\tilde{\mathfrak a}}_1'
+{\mathfrak a}_1(0, \theta){\tilde{\mathfrak a}}_2{\tilde{\mathfrak a}}_2'\big]n_1n_2
\, +\, 4\big[-\partial_t{\mathfrak a}_2(0, \theta){\tilde{\mathfrak a}}_1'
+\partial_t{\mathfrak a}_1(0, \theta){\tilde{\mathfrak a}}_2'\big]n_1n_2
\nonumber\\[2mm]
&\qquad+2\big[{\mathfrak a}_2(0, \theta)|{\tilde{\mathfrak a}}_1'|^2|n_1|^2
 + {\mathfrak a}_1(0, \theta)|{\tilde{\mathfrak a}}_2'|^2|n_2|^2\big]
 +2\, k^2\big[{\mathfrak a}_2(0, \theta){\tilde{\mathfrak a}}_1|n_2|^2+{\mathfrak a}_1(0, \theta){\tilde{\mathfrak a}}_2|n_1|^2\, \big]
 \nonumber\\[2mm]
&\qquad-4k\big[\, \partial_t{\mathfrak a}_2(0, \theta){\tilde{\mathfrak a}}_1|n_2|^2+\partial_t{\mathfrak a}_1(0, \theta){\tilde{\mathfrak a}}_2|n_1|^2\, \big]
\, +\, 2\Theta_{\theta tt}(0, \theta)\big[{\mathfrak a}_2(0, \theta)|n_2|^2+{\mathfrak a}_1(0, \theta)|n_1|^2\big]
\nonumber\\[2mm]
&\qquad+\big[\partial_{tt}{\mathfrak a}_2(0, \theta)|n_2|^2+\partial_{tt}{\mathfrak a}_1(0, \theta)|n_1|^2\big]
+2\Theta_{tt}(0, \theta)\big[{\mathfrak a}_2(0, \theta)\gamma_1'\gamma_1''+{\mathfrak a}_1(0, \theta)\gamma_2'\gamma_2''\big]
 \Big\}fh
\nonumber\\[2mm]
&\quad+\Big\{2\big[{\mathfrak a}_2(0, \theta)|n_2|^2+{\mathfrak a}_1(0, \theta)|n_1|^2\big]\Theta_{tt}(0, \theta)
+2k\big[{\mathfrak a}_2(0, \theta){\tilde{\mathfrak a}}_1{\tilde{\mathfrak a}}_1'|n_1|^2
+ {\mathfrak a}_1(0, \theta){\tilde{\mathfrak a}}_2{\tilde{\mathfrak a}}_2'|n_2|^2\big]
 \nonumber\\[2mm]
&\qquad \quad-2k\big[-{\mathfrak a}_2(0, \theta)|{\tilde{\mathfrak a}}_1'|^2
+{\mathfrak a}_1(0, \theta)|{\tilde{\mathfrak a}}_2'|^2\big]n_1n_2
 \nonumber\\[2mm]
&\qquad \quad+2\big[-\partial_t{\mathfrak a}_2(0, \theta){\tilde{\mathfrak a}}_1
+\partial_t{\mathfrak a}_1(0, \theta){\tilde{\mathfrak a}}_2\big]n_1n_2
\Big\}(f'h+fh')\, +\, 2f'h'\, {\mathfrak w}_0
 \nonumber\\[2mm]
&\equiv {\mathfrak f}_3 fh \, +\, {\mathfrak f}_4 (f' h +fh') \, +\, 2f'h'\, {\mathfrak w}_0.
\end{align}
The combinations of the above computations together with \eqref{mathcalM(0)}-\eqref{mathcalM6(0)} will give that
\begin{align}
{\mathcal J}''(0)[h, f]
=\, &\int_0^1
\left[\, \frac{\sigma V_{tt}(0, \theta)\, }{ V^{1-\sigma}\big(0, \theta\big)}
+\frac{\sigma(\sigma-1)\big| V_t(0, \theta)\big|^2\, }{ V^{2-\sigma}(0, \theta)}\, \right]\, \sqrt{{\mathfrak f}_0}\, hf\, {\rm d}\theta
\nonumber\\[2mm]
&\, +\,
\sigma\int_0^1\frac{ V_t(0, \theta) \, }{ V^{1-\sigma}\big(0, \theta\big)}\, \frac{1}{\sqrt{{\mathfrak f}_0}}{\mathfrak f}_1hf\, {\rm d}\theta
\,-\,
\frac{1}{4}\int_0^1\frac{ V^{\sigma}\big(0, \theta\big)}{\ \big(\sqrt{{\mathfrak f}_0}\, \big)^3\ }{\mathfrak f}_1^2hf \, {\rm d}\theta
\nonumber\\[2mm]
&\, +\,
\frac{1}{2}\int_0^1\frac{ V^{\sigma}\big(0, \theta\big)}{\sqrt{{\mathfrak f}_0}}
\big[ {\mathfrak f}_3 fh \, +\, {\mathfrak f}_4 (f' h +fh') \, +\, 2f'h'\, {\mathfrak w}_0 \big]\Big|_{s=0}
\, {\rm d}\theta.
\label{secondvariation}
\end{align}

\medskip
\noindent{\bf{Notation 2: }}\label{notation2}
By recalling ${\mathfrak w}_0, {\mathfrak l}_1$  and ${\mathfrak f}_0, {\mathfrak f}_2, {\mathfrak f}_3$
in \eqref{m11}, \eqref{m12} and \eqref{m22},
we  introduce following functions:

\begin{align}\label{mathcalH1}
{\mathcal H}_1(\theta)=&
\frac{ V^{\sigma}\big(0, \theta\big)}{\sqrt{{\mathfrak f}_0}}\, {\mathfrak w}_0,
\end{align}

\begin{align}\label{mathcalH2}
{\mathcal H}_2(\theta)
=&\frac{ V^{\sigma}\big(0, \theta\big)}{2\sqrt{{\mathfrak f}_0}}
\Big\{2\big[{\mathfrak a}_2(0, \theta)|n_2|^2+{\mathfrak a}_1(0, \theta)|n_1|^2\big]\Theta_{tt}(0, \theta)
+2\big[{\mathfrak a}_2(0, \theta){\tilde{\mathfrak a}}_1{\tilde{\mathfrak a}}_1'|n_1|^2
+{\mathfrak a}_1(0, \theta){\tilde{\mathfrak a}}_2{\tilde{\mathfrak a}}_2'|n_2|^2\big]
\nonumber\\[2mm]
&\qquad \qquad -2k\big[-{\mathfrak a}_2(0, \theta)|{\tilde{\mathfrak a}}_1'|^2
+ {\mathfrak a}_1(0, \theta)|{\tilde{\mathfrak a}}_2'|^2\big]n_1n_2
+2\big[-\partial_t{\mathfrak a}_2(0, \theta){\tilde{\mathfrak a}}_1
+\partial_t{\mathfrak a}_1(0, \theta){\tilde{\mathfrak a}}_2\big]n_1n_2
\Big\}
\nonumber\\
=&\frac{ V^{\sigma}\big(0, \theta\big)}{\sqrt{{\mathfrak f}_0}}\Big\{{\mathfrak l}_1
-\big[{\mathfrak a}_1(0, \theta)n_1q_2+{\mathfrak a}_2(0, \theta)n_2q_1\big]
+\big[{\mathfrak a}_2(0, \theta)|n_2|^2+{\mathfrak a}_1(0, \theta)|n_1|^2\big]\Theta_{tt}(0, \theta)\Big\}
\nonumber\\[2mm]
=&\frac{ V^{\sigma}\big(0, \theta\big)}{\sqrt{{\mathfrak f}_0}}{\mathfrak l}_1,
\end{align}
and
\begin{align}\label{mathcalH3}
{\mathcal H}_3(\theta)
=&\frac{ V^{\sigma}\big(0, \theta\big)}{2\sqrt{{\mathfrak f}_0}}
\Big\{-4k\, \big[{\mathfrak a}_1(0, \theta){\tilde{\mathfrak a}}_2{\tilde{\mathfrak a}}_2'-{\mathfrak a}_2(0, \theta){\tilde{\mathfrak a}}_1{\tilde{\mathfrak a}}_1'
\big]n_1n_2
+4\big[\partial_t{\mathfrak a}_1(0, \theta){\tilde{\mathfrak a}}_2'
-\partial_t{\mathfrak a}_2(0, \theta){\tilde{\mathfrak a}}_1' \big]n_1n_2
\nonumber\\[2mm]
&\qquad\qquad \quad+2\big[{\mathfrak a}_2(0, \theta)|{\tilde{\mathfrak a}}_1'|^2|n_1|^2
+{\mathfrak a}_1(0, \theta)|{\tilde{\mathfrak a}}_2'|^2|n_2|^2\big]
+2k^2\big[{\mathfrak a}_2(0, \theta){\tilde{\mathfrak a}}_1|n_2|^2+{\mathfrak a}_1(0, \theta){\tilde{\mathfrak a}}_2|n_1|^2\big]
 \nonumber\\[2mm]
&\qquad\qquad \quad-4k\big[\partial_t{\mathfrak a}_2(0, \theta){\tilde{\mathfrak a}}_1|n_2|^2+\partial_t{\mathfrak a}_1(0, \theta){\tilde{\mathfrak a}}_2|n_1|^2\big]
+\big[\partial_{tt}{\mathfrak a}_2(0, \theta)|n_2|^2+\partial_{tt}{\mathfrak a}_1(0, \theta)|n_1|^2\big]
\nonumber\\[2mm]
&\qquad\qquad \quad+2\Theta_{\theta tt}(0, \theta)\big[{\mathfrak a}_2(0, \theta)|n_2|^2+{\mathfrak a}_1(0, \theta)|n_1|^2\big]
\nonumber\\[2mm]
&\qquad\qquad \quad+2\Theta_{tt}(0, \theta)\big[{\mathfrak a}_2(0, \theta)\gamma_1'\gamma_1''+{\mathfrak a}_1(0, \theta)\gamma_2'\gamma_2''\big]
 \Big\}
\nonumber\\[2mm]
&+\left[\, \frac{\sigma V_{tt}(0, \theta)\, }{ V^{1-\sigma}\big(0, \theta\big)}
+\frac{\sigma(\sigma-1)\big| V_t(0, \theta)\big|^2\, }{ V^{2-\sigma}(0, \theta)}
\, \right]\sqrt{{\mathfrak f}_0}
\, +\, \sigma\frac{ V_t(0, \theta) \, }{ V^{1-\sigma}\big(0, \theta\big)}\frac{1}{\sqrt{{\mathfrak f}_0}}{\mathfrak f}_1
\,-\, \frac{1}{4}\frac{ V^{\sigma}\big(0, \theta\big)}{\ \big(\sqrt{{\mathfrak f}_0}\, \big)^3\ }{\mathfrak f}_1^2
\nonumber\\[2mm]
\, =\, &\frac{ V^{\sigma}\big(0, \theta\big)}{2\sqrt{{\mathfrak f}_0}}
\Big\{2{\mathfrak f}_2-\big[{\mathfrak a}_1(0, \theta) n_1q_2' -{\mathfrak a}_2(0, \theta) n_2q_1'\big]
+2\Theta_{\theta tt}(0, \theta)\big[{\mathfrak a}_2(0, \theta)|n_2|^2+{\mathfrak a}_1(0, \theta)|n_1|^2\big]
\nonumber\\[2mm]
&\qquad\qquad \quad+2\Theta_{tt}(0, \theta)\big[{\mathfrak a}_2(0, \theta)\gamma_1'\gamma_1''+{\mathfrak a}_1(0, \theta)\gamma_2'\gamma_2''\big]\Big\}
\nonumber\\[2mm]
&+\left[\,
\frac{\sigma V_{tt}(0, \theta)\, }{ V^{1-\sigma}\big(0, \theta\big)}
+\frac{\sigma(\sigma-1)\big| V_t(0, \theta)\big|^2\, }{ V^{2-\sigma}(0, \theta)}
\, \right]\, \sqrt{{\mathfrak f}_0}
\, +\, \sigma\frac{ V_t(0, \theta) \, }{ V^{1-\sigma}\big(0, \theta\big)}
\, \frac{1}{\sqrt{{\mathfrak f}_0}}{\mathfrak f}_1
\,-\, \frac{1}{4}\frac{ V^{\sigma}\big(0, \theta\big)}{\ \big(\sqrt{{\mathfrak f}_0}\, \big)^3\ }{\mathfrak f}_1^2
\nonumber\\[2mm]
=&\frac{ V^{\sigma}\big(0, \theta\big)}{\sqrt{{\mathfrak f}_0}}{\mathfrak f}_2
+\left[\,
\frac{\sigma V_{tt}(0, \theta)\, }{ V^{1-\sigma}\big(0, \theta\big)}
+\frac{\sigma(\sigma-1)\big| V_t(0, \theta)\big|^2\, }{ V^{2-\sigma}(0, \theta)}
\, \right]\, \sqrt{{\mathfrak f}_0}
\nonumber \\[2mm]
&\, +\, \sigma\frac{ V_t(0, \theta) \, }{ V^{1-\sigma}\big(0, \theta\big)}
\, \frac{1}{\sqrt{{\mathfrak f}_0}}{\mathfrak f}_1
\,-\, \frac{1}{4}\frac{ V^{\sigma}\big(0, \theta\big)}{\ \big(\sqrt{{\mathfrak f}_0}\, \big)^3\ }{\mathfrak f}_1^2.
\end{align}
In the above, we have used the facts \eqref{q_i} and \eqref{qi'}.


\medskip
The terms in \eqref{secondvariation} can be rearranged in the following way
\begin{align}
{\mathcal J}''(0)[h, f]
=\, &\int_0^1\Big( {\mathcal H}_1f' \, +\, {\mathcal H}_2 f \Big) h'\, {\rm d}\theta
\, +\,
\int_0^1\Big( {\mathcal H}_2f' \, +\, {\mathcal H}_3 f \Big) h\, {\rm d}\theta
\nonumber\\[2mm]
=&\Big( {\mathcal H}_1f' \, +\, {\mathcal H}_2 f \Big) h\Big|^1_{\theta=0}
\,-\,
\int_0^1\big[ {\mathcal H}_1f'' \, +\,  {\mathcal H}_1'f'    \, +\, \big({\mathcal H}_2'-{\mathcal H}_3\big) f  \big]h\, {\rm d}\theta.
\end{align}
Note that
\begin{align}
{\mathcal H}_1(1)f'(1) \, +\, {\mathcal H}_2(1) f(1)
\, =\, &\frac{ V^{\sigma}\big(0, 1\big)}{\sqrt{{\mathfrak f}_0(1)}}
\big[
{\mathfrak w}_0(1)f'(1)
\ +\
{\mathfrak l}_1(1)f(1)
\big]
\label{boundaryoffunctional1000}
\\[2mm]
\, =\, &\frac{ V^{\sigma}\big(0, 1\big)}{2\sqrt{{\mathfrak f}_0(1)}}
\big[
{\mathfrak b}_6f'(1)
\ -\
{\mathfrak b}_7f(1)
\big],
\label{boundaryoffunctional1}
\end{align}
and
\begin{align}
{\mathcal H}_1(0)f'(0) \, +\, {\mathcal H}_2(0) f(0)
\, =\, &\frac{ V^{\sigma}\big(0, 0\big)}{\sqrt{{\mathfrak f}_0(0)}}
\big[
{\mathfrak w}_0(0)f'(0)
\ +\
{\mathfrak l}_1(0)f(0)
\big]
\label{boundaryoffunctional2000}
\\[2mm]
\, =\, &\frac{ V^{\sigma}\big(0, 0\big)}{2\sqrt{{\mathfrak f}_0(0)}}
\big[
{\mathfrak b}_1f'(0)
\ -\
{\mathfrak b}_2f(0)
\big].
\label{boundaryoffunctional2}
\end{align}
For more details of the derivation of the equalities \eqref{boundaryoffunctional1} and \eqref{boundaryoffunctional2},
the reader can refer to the computations in Appendix \ref{appendixE}.

\medskip
For a stationary curve $\Gamma$, we say that it is {\bf non-degenerate} in the sense that if
 \begin{equation}
 \label{2-relation}
  {\mathcal J}{''}(0)[h, f]=0, \quad \forall\, h\in H^1(0, 1),
  \end{equation}
then $f\equiv 0$.
It is  equivalent to that the boundary problem
 \begin{equation}
 \label{nondegeneracy}
 \begin{split}
{\mathcal H}_1f'' \, +\,  {\mathcal H}_1'f'    \, +\, \big({\mathcal H}_2'-{\mathcal H}_3\big) f=0
\quad \mbox{in}\  (0, 1),
\\[2mm]
{\mathfrak b}_1f'(0)
\ -\
{\mathfrak b}_2f(0)=0,
\qquad
{\mathfrak b}_6f'(1)
\ -\
{\mathfrak b}_7f(1)=0,
\end{split}
\end{equation}
has only the trivial solution.

\medskip
Paralleling with the above arguments, we finally give the following Remark.
\begin{remark}
For a simple closed curve ${\hat\Gamma}$ in ${\mathbb R}^2$ with unit length,
we construct another {\bf modified Fermi coordinates}
\begin{align}
y\, =\, {\hat\gamma}({\hat\theta})
\, +\,
{\hat t}\, \big({\tilde{\mathfrak a}}_1({\hat\theta}){\hat n}_1({\hat\theta}), \, {\tilde{\mathfrak a}}_2( {\hat\theta}){\hat n}_2({\hat\theta})\big),
\label{Fermicoordinates-modified-tilde}
\end{align}
where ${\hat\gamma}({\hat\theta})=({\hat\gamma}_1({\hat\theta}), {\hat\gamma}_2({\hat\theta}))$ is a natural parametrization of ${\hat\Gamma}$ in ${\mathbb R}^2$,
${\hat n}({\hat\theta})=({\hat n}_1({\hat\theta}), {\hat n}_2({\hat\theta}))$ is the unit normal of ${\hat\Gamma}$,
and the functions ${\tilde{\mathfrak a}}_1$ and ${\tilde{\mathfrak a}}_2$ have similar expressions as in \eqref{tildea1a2} for ${\hat\theta}\in [0, 1]$.
Consider the deformation of ${\hat\Gamma}$ in the form
\begin{equation}
\label{deformation2}
{\hat\Gamma}_{{\mathfrak t}}({\hat\theta})=\big({\hat\Gamma}_{{\mathfrak t}1}({\hat\theta}), {\hat\Gamma}_{{\mathfrak t}2}({\hat\theta})\big)
\, :\, {\hat\gamma}({\hat\theta})
 \, +\,
 {\mathfrak t}({\hat\theta})\, \Big({\tilde{\mathfrak a}}_1({\hat\theta}){\hat n}_1({\hat\theta}), \, {\tilde{\mathfrak a}}_2({\hat\theta}){\hat n}_2({\hat\theta})\Big),
\end{equation}
where  ${\mathfrak t}$ is a smooth function of ${\hat\theta}$ with small $L^{\infty}$-norm,
and then the length functional
\begin{equation}
\begin{split}
{\hat {\mathcal J}}({\mathfrak t})&\equiv
\int_0^1 V^{\sigma}\big({\hat\Gamma}_{{\mathfrak t}}({\hat\theta})\big)
\sqrt{{\mathfrak a}_2\big({\hat\Gamma}_{{\mathfrak t}}({\hat\theta})\big)\big|{\hat\Gamma}'_{{\mathfrak t}1}({\hat\theta})\big|^2
+ {\mathfrak a}_1\big({\hat\Gamma}_{{\mathfrak t}}({\hat\theta})\big)\big|{\hat\Gamma}'_{{\mathfrak t}2}({\hat\theta})\big|^2}\, {\rm d}{\hat\theta}.
\end{split}
\label{weightedlength2}
\end{equation}
We can do the same variational calculations to the functional ${\hat {\mathcal J}}$ and then derive the following notions.
If the curvature ${\hat k}$ of ${\hat\Gamma}$ satisfies
\begin{align}
{\hat k}=&\frac{{\mathfrak a}_1(0, {\hat\theta}) {{\tilde{\mathfrak a}}'}_2 {\hat n}_1 {\hat n}_2 -{\mathfrak a}_2(0, {\hat\theta}){{\tilde{\mathfrak a}}'}_1 {\hat n}_2 {\hat n}_1}
{\, {\mathfrak a}_1(0, {\hat\theta}) {\tilde{\mathfrak a}}_2|{\hat n}_1|^2 +{\mathfrak a}_2(0, {\hat\theta}){\tilde{\mathfrak a}}_1 |{\hat n}_2|^2\, }
\, +\,
\frac{1}{2}\frac{\partial_{\hat t} {\mathfrak a}_1(0, {\hat\theta}) |{\hat n}_1|^2 +\partial_{\hat t} {\mathfrak a}_2(0, {\hat\theta}) |{\hat n}_2|^2}
{\, \big[{\mathfrak a}_1(0, {\hat\theta}) {\tilde{\mathfrak a}}_2|{\hat n}_1|^2 +{\mathfrak a}_2(0, {\hat\theta}){\tilde{\mathfrak a}}_1 |{\hat n}_2|^2\big]\, }
\nonumber\\[2mm]
&\, +\,
\sigma\frac{ V_{\hat t}(0, {\hat\theta})}{ V(0, {\hat\theta})}
\frac{{\mathfrak a}_1(0, {\hat\theta}) |{\hat n}_1|^2 +{\mathfrak a}_2(0, {\hat\theta}) |{\hat n}_2|^2}
{\, \big[{\mathfrak a}_1(0, {\hat\theta}) {\tilde{\mathfrak a}}_2|{\hat n}_1|^2 +{\mathfrak a}_2(0, {\hat\theta}){\tilde{\mathfrak a}}_1 |{\hat n}_2|^2\big]\, },
\quad\forall\, {\hat\theta}\in(0, 1),
\label{stationary4}
\end{align}
then the curve ${\hat\Gamma}$ is said to be {\bf stationary}.
Set the notation
\begin{align}
{\widehat{\mathcal H}}_1=\frac{ V^{\sigma}\big(0, {\hat\theta}\big)}{\sqrt{\hat{\mathfrak f}_0}}\, \big[{\mathfrak a}_1(0, {\hat\theta}){\tilde{\mathfrak a}}_2^2 |{\hat n}_2|^2+{\mathfrak a}_2(0, {\hat\theta}){\tilde{\mathfrak a}}_1^2 |{\hat n}_1|^2\big]
\label{mathcalH1tilde}
\end{align}
\begin{align}
{\widehat{\mathcal H}}_2
=&\frac{ V^{\sigma}\big(0, {\hat\theta}\big)}{2\sqrt{\hat{\mathfrak f}_0}}
\Big\{2\big[{\mathfrak a}_2(0, {\hat\theta}){\tilde{\mathfrak a}}_1{\tilde{\mathfrak a}}_1'|{\hat n}_1|^2
+{\mathfrak a}_1(0, {\hat\theta}){\tilde{\mathfrak a}}_2{\tilde{\mathfrak a}}_2'|{\hat n}_2|^2\big]
-2{\hat k}\big[{\mathfrak a}_1(0, {\hat\theta})|{\tilde{\mathfrak a}}_2'|^2
-{\mathfrak a}_2(0, {\hat\theta})|{\tilde{\mathfrak a}}_1'|^2
\big]{\hat n}_1{\hat n}_2
\nonumber\\[2mm]
&\qquad \qquad \quad
+2\big[-\partial_{\hat t}{\mathfrak a}_2(0, {\hat\theta}){\tilde{\mathfrak a}}_1
+\partial_{\hat t}{\mathfrak a}_1(0, {\hat\theta}){\tilde{\mathfrak a}}_2\big]{\hat n}_1{\hat n}_2
\Big\}
\label{mathcalH2tilde}
\end{align}
\begin{align}
{\widehat{\mathcal H}}_3
=&\frac{ V^{\sigma}\big(0, {\hat\theta}\big)}{2\sqrt{\hat{\mathfrak f}_0}}
\Big\{-4{\hat k}\, \big[{\mathfrak a}_1(0, {\hat\theta}){\tilde{\mathfrak a}}_2{\tilde{\mathfrak a}}_2'-{\mathfrak a}_2(0, {\hat\theta}){\tilde{\mathfrak a}}_1{\tilde{\mathfrak a}}_1'
\big]{\hat n}_1{\hat n}_2
+4\big[\partial_{\hat t}{\mathfrak a}_1(0, {\hat\theta}){\tilde{\mathfrak a}}_2'
-\partial_{\hat t}{\mathfrak a}_2(0, {\hat\theta}){\tilde{\mathfrak a}}_1' \big]{\hat n}_1{\hat n}_2
\nonumber\\[2mm]
&\qquad\qquad+2\big[{\mathfrak a}_2(0, {\hat\theta})|{\tilde{\mathfrak a}}_1'|^2|{\hat n}_1|^2
+{\mathfrak a}_1(0, {\hat\theta})|{\tilde{\mathfrak a}}_2'|^2|{\hat n}_2|^2\big]
+2{\hat k}^2\big[{\mathfrak a}_2\big(0, {\hat\theta}\big){\tilde{\mathfrak a}}_1|{\hat n}_2|^2+{\mathfrak a}_1\big(0, {\hat\theta}\big){\tilde{\mathfrak a}}_2|{\hat n}_1|^2\big]
 \nonumber\\[2mm]
&\qquad\qquad-4{\hat k}\big[\partial_{\hat t}{\mathfrak a}_2\big(0, {\hat\theta}\big){\tilde{\mathfrak a}}_1|{\hat n}_2|^2+\partial_{\hat t}{\mathfrak a}_1\big(0, {\hat\theta}\big){\tilde{\mathfrak a}}_2|{\hat n}_1|^2\big]
+\big[\partial_{{\hat t}{\hat t}}{\mathfrak a}_2\big(0, {\hat\theta}\big)|{\hat n}_2|^2
+\partial_{{\hat t}{\hat t}}{\mathfrak a}_1\big(0, {\hat\theta}\big)|{\hat n}_1|^2\big]\Big\}
\nonumber\\[2mm]
&+\left[\,
\frac{\sigma V_{{\hat t}{\hat t}}(0, {\hat\theta})\, }{ V^{1-\sigma}\big(0, {\hat\theta}\big)}
+
\frac{\sigma(\sigma-1)\big| V_{\hat t}(0, {\hat\theta})\big|^2\, }{ V^{2-\sigma}(0, {\hat\theta})}
\, \right]
\, \sqrt{\hat{\mathfrak f}_0}
+\sigma \frac{ V_{\hat t}(0, {\hat\theta}) \, }{ V^{1-\sigma}\big(0, {\hat\theta}\big)}
\, \frac{1}{\sqrt{\hat{\mathfrak f}_0}}\hat{\mathfrak f}_1
\,-\,
 \frac{1}{4}\frac{ V^{\sigma}\big(0, {\hat\theta}\big)}{\ \big(\sqrt{\hat{\mathfrak f}_0}\, \big)^3\ }\hat{\mathfrak f}_1^2,
\label{mathcalH3tilde}
\end{align}
where
$$
\hat{\mathfrak f}_0({\hat\theta})={\mathfrak a}_1(0, {\hat\theta}) |{\hat n}_1({\hat\theta})|^2 +{\mathfrak a}_2(0, {\hat\theta}) |{\hat n}_2({\hat\theta})|^2,
$$
\begin{align*}
\hat{\mathfrak f}_1({\hat\theta})
=&\,2\bigg[{\mathfrak a}_1(0, {\hat\theta}) {\tilde{\mathfrak a}'}_2({\hat\theta}) {\hat n}_1({\hat\theta}) {\hat n}_2({\hat\theta})
\,-\,
{\mathfrak a}_2(0, {\hat\theta}){\tilde{\mathfrak a}'}_1({\hat\theta}) {\hat n}_2({\hat\theta}) {\hat n}_1({\hat\theta})\bigg]
\nonumber\\[2mm]
&\,-\,
2{\hat k}({\hat\theta})\bigg[{\mathfrak a}_1(0, {\hat\theta}) {\tilde{\mathfrak a}}_2({\hat\theta})|{\hat n}_1({\hat\theta})|^2
\,+\,
{\mathfrak a}_2(0, {\hat\theta}){\tilde{\mathfrak a}}_1({\hat\theta}) |{\hat n}_2({\hat\theta})|^2\bigg]
\nonumber\\[2mm]
&
\,+\,\bigg[\partial_{\hat t} {\mathfrak a}_1(0, {\hat\theta}) |{\hat n}_1({\hat\theta})|^2 +\partial_{\hat t} {\mathfrak a}_2(0, {\hat\theta}) |{\hat n}_2({\hat\theta})|^2\bigg].
\end{align*}
For a simple closed curve ${\hat\Gamma}$ satisfying the stationary condition, if the boundary problem
 \begin{equation}
 \label{nondegeneracy4}
 \begin{split}
{\widehat{\mathcal H}}_1f''
\, +\,
 {\widehat{\mathcal H}}_1'f'
\, +\,
\big({\widehat{\mathcal H}}_2'-{\widehat{\mathcal H}}_3\big) f&=0
\quad \mbox{in}\  (0, 1),
\\[2mm]
f'(1)=f'(0),
\qquad
f(1)&=f(0),
\end{split}
\end{equation}
has only the trivial solution, we call ${\hat\Gamma}$ a {\bf non-degenerate} stationary curve.
\qed
\end{remark}

\section{Outline of the proof}\label{thegluingprocedure}\label{section3}
\setcounter{equation}{0}

Recall that $w$ is the solution to \eqref{wsolution}. In fact, $w$ is an even function defined in the form
\begin{equation}
w(x)=C_p\big[\, e^{\frac{(p-1)x}{2}}\, +\, e^{-\frac{(p-1)x}{2}}\, \big]^{-\frac{2}{p-1}},
\quad \forall\, x\in {\mathbb R}.
\end{equation}
We consider the associated linearized eigenvalue problem,
\begin{equation}
\label{eigenvalue}
h{''}-h+pw^{p-1}h\, =\, \lambda h\quad\mbox{in } {\mathbb R}, \quad h(\pm \infty)\, =\, 0.
\end{equation}
It is well known that this equation possesses a unique positive eigenvalue $\lambda_0$,
with associated eigenfunction $Z$ (even and positive)  in $H^1({\mathbb R})$, which can be normalized so that $\int_{\mathbb R} Z^2=1$.
In fact, a simple computation shows that
\begin{align}
\label{lambda0}
\lambda_0\, =\, \frac 14(p-1)(p+3), \quad Z\, =\, \frac {1}{\sqrt{\int_{\mathbb R} w^{p+1}}\, }w^{\frac {p+1}{2}}.
\end{align}
In this section, the strategy to prove Theorem \ref{theorem 1.1} will be provided step by step.
\subsection{The gluing procedure}
Recall that $\delta_0>0$ is a small constant in (\ref{Fermicoordinates-modified}) and define a cut-off function $\eta_{3\delta}^{\varepsilon}(s)=\eta_{3\delta}(\varepsilon |s|)$ where $\eta_{\delta}(t)$  is also a smooth cut-off function defined as
$$
\eta_{\delta}(t)=1, \, \forall\, 0\leq t\leq \delta
\quad\mbox{and}\quad
\eta_{\delta}(t)=0, \, \forall\, t>2\delta,
$$
for a fixed number $\delta<\delta_0/100$.
 For any given approximate solution ${\mathbf w}$ (to be chosen later, cf. (\ref{globalapproximation})) and a perturbation term
$
\tilde{\phi}({\tilde y})\, =\, \eta_{3\delta}^{\varepsilon}(s){\check\phi}({\tilde y})\, +\, {\check\psi}({\tilde y})
$
on $\Omega_\varepsilon$, the function $u({\tilde y})=\mathbf{w}({\tilde y})+\tilde{\phi}({\tilde y})$ satisfies (\ref{problemafterscaling}) if $({\check\phi}, {\check\psi})$ satisfies the following coupled system:
\begin{align}\label{equivalent system-1}
\eta_{3\delta}^{\varepsilon}\, {\mathbf L}({\check\phi})
\, &=\,
\eta_{\delta}^{\varepsilon}\, \big[-{\mathbf N}(\eta_{3\delta}^{\varepsilon}\, {\check\phi}+{\check\psi})
 \,-\, {\mathbf E}
 \,-\, p\, \mathbf{w}^{p-1}{\check\psi}\big]\quad\mbox{in } \Omega_\varepsilon,
\end{align}
\begin{align}
\label{on the boundary-1}
\eta_{3\delta}^{\varepsilon}\nabla_{{\mathfrak a}(\varepsilon {\tilde y})}  {\check\phi}\cdot \nu_\varepsilon
\, +\, \eta_{\delta}^{\varepsilon}\nabla_{{\mathfrak a}(\varepsilon {\tilde y})} \mathbf{w}\cdot \nu_\varepsilon\, =\, 0
\quad\mbox{on } \partial\, \Omega_\varepsilon,
\end{align}
and
\begin{align}
\label{equivalent system-2}
&\sum_{i=1}^2\partial_{{\tilde y}_i}\big({\mathfrak a}_i(\varepsilon{\tilde y}){\check\psi}_{{\tilde y}_i}\big)
\,-\, V(\varepsilon \tilde y) \, {\check\psi}
\, +\, (1-\eta_{\delta}^{\varepsilon})\, p\, \mathbf{w}^{p-1}\, {\check\psi}
\nonumber\\[2mm]
&\, =\,-\, (1-\eta_{\delta}^{\varepsilon})\, {\mathbf E}
\,-\, \sum_{i=1}^2{\mathfrak a}_i(\varepsilon{\tilde y}){\check\phi}_{{\tilde y}_i}\partial_{{\tilde y}_i}\big(\eta_{3\delta}^{\varepsilon}(s)\big)
\,-\, \sum_{i=1}^2\partial_{{\tilde y}_i}\Big[{\mathfrak a}_i(\varepsilon{\tilde y}){\check\phi}\partial_{{\tilde y}_i}\big( \eta_{3\delta}^{\varepsilon}(s) \big) \Big]
\nonumber\\[2mm]
&\quad\, \,-\, (1-\eta_{\delta}^{\varepsilon})\, {\mathbf N}\, (\eta_{3\delta}^{\varepsilon}\, {\check\phi}+{\check\psi})
\quad\mbox{in } \Omega_\varepsilon,
\end{align}
\begin{align}\label{on the boundary-2}
\nabla_{{\mathfrak a}(\varepsilon {\tilde y})} {\check\psi}\cdot \nu_\varepsilon
\, +\, \big(1-\eta_{\delta}^{\varepsilon}\big)\nabla_{{\mathfrak a}(\varepsilon {\tilde y})} \mathbf{w}\cdot \nu_\varepsilon
\, +\, \varepsilon \nabla_{{\mathfrak a}(\varepsilon {\tilde y})} \eta_{3\delta}^{\varepsilon}\cdot \nu_\varepsilon {\check\phi}\, =\, 0
\quad\mbox{on } \partial\, \Omega_\varepsilon, \qquad
\end{align}
where
\begin{equation}\label{globalerror}
{\mathbf E}\, =\, \sum_{i=1}^2\partial_{{\tilde y}_i}\big({\mathfrak a}_i(\varepsilon{\tilde y}){\mathbf w}_{{\tilde y}_i}\big)
\,-\, V(\varepsilon \tilde y) {\mathbf w} + {\mathbf w}^p,
\end{equation}
\begin{equation}
 {\mathbf L}({\check\phi})\, =\, \sum_{i=1}^2\partial_{{\tilde y}_i}\big({\mathfrak a}_i(\varepsilon{\tilde y}){\check\phi}_{{\tilde y}_i}\big)
\,-\, V(\varepsilon \tilde{y}){\check\phi}
\, +\, p\, \mathbf{w}^{p-1}{\check\phi},
\qquad
{\mathbf N}(\tilde{\phi})\, =\, (\mathbf{w}+\tilde{\phi})^p
\,-\, \mathbf{w}^p
\,-\, p\, \mathbf{w}^{p-1}\tilde{\phi}.
\end{equation}

Assume now that ${\check\phi}$ satisfies the following decay property
\begin{align}
\label{decay property}
\big |\nabla {\check\phi}(\tilde{y})\big|+\big|{\check\phi}(\tilde{y})\big|\, \leq\, Ce^{-\rho/\varepsilon}   \quad\mbox{if } \mathrm{dist}({\tilde y}, \Gamma_\varepsilon)>\delta/\varepsilon,
\end{align}
for a certain constant $\rho>0$, and also that
\begin{align}\label{constraint condition}
\mathbf{w}({\tilde y})
\ \mbox{is exponentially small if } \mathrm{dist}({\tilde y}, \Gamma_\varepsilon)>\delta/\varepsilon.
\end{align}
Since ${\mathbf N}$ is power-like with power greater than one, a direct application of Contraction Mapping Principle yields that (\ref{equivalent system-2})-(\ref{on the boundary-2}) has a unique (small) solution ${\check\psi}={\check\psi}({\check\phi})$ with
\begin{align}
\label{contraction-1}
\|{\check\psi}({\check\phi})\|_{L^{\infty}}\, \leq\, C\varepsilon\big[\, \|{\check\phi}\|_{L^{\infty}(|s|>\delta/\varepsilon)}+
\|\nabla {\check\phi}\|_{L^{\infty}(|s|>\delta/\varepsilon)}+e^{-\delta/\varepsilon}\, \big],
\end{align}
where $|s|>\delta/\varepsilon$ denotes the complement in $\Omega_{\varepsilon}$ of $\delta/\varepsilon$-neighborhood of $\Gamma_\varepsilon$. Moreover, the nonlinear operator ${\check\psi}$ satisfies a Lipschitz condition of the form
\begin{align}
\label{psi-lip}
\|{\check\psi}({\check\phi}_1)-{\check\psi}({\check\phi}_2)\|_{L^{\infty}}\, \leq\, C\varepsilon \big[\, \|{\check\phi}_1-{\check\phi}_2\|_{L^{\infty}(|s|>\delta/\varepsilon)}
\, +\,
\|\nabla {\check\phi}_1-\nabla {\check\phi}_2\|_{L^{\infty}(|s|>\delta/\varepsilon)}\, \big].
\end{align}

\medskip
The key observation is that, after solving (\ref{equivalent system-2})-(\ref{on the boundary-2}),
we can concern (\ref{equivalent system-1})-(\ref{on the boundary-1}) as  a local nonlinear problem  involving ${\check\psi}\, =\, {\check\psi}({\check\phi})$, which can be solved in local coordinates in the sense that we can decompose the interaction among the boundary, the concentration set and the terms ${\mathfrak a}$, $ V$, and then construct a good approximate solution and also derive the resolution theory of the nonlinear problem by delicate analysis.
This whole procedure is called a {\bf gluing technique} in \cite{delPKowWei2007}.

\subsection{Local formulation of the problem}\label{section3.2}
As described in the above, the next step is to consider (\ref{equivalent system-1})-(\ref{on the boundary-1}) in the neighbourhood of $\Gamma_\varepsilon$ so that by the relation ${\tilde y}=y/\varepsilon$ in (\ref{rescaling}) close to $\Gamma_\varepsilon$, the variables $y$ can be represented by the modified Fermi coordinates, say $(t, \theta)$ in (\ref{Fermicoordinates-modified}), which
have been set up in Section \ref{section2}.

\subsubsection{Local forms of the problem}\label{section22}

By recalling the local coordinates $(t, \theta)$ in (\ref{Fermicoordinates-modified}),
we can also define the local rescaling,
\begin{equation}\label{coordinatessz}
(t, \theta)=\varepsilon (s, z),
\end{equation}
and then use the results in \eqref{laplacelocal}, \eqref{boundaryoriginal0} and \eqref{boundaryoriginal1} to give local expressions of the problem.
The equation in (\ref{equivalent system-1})
 can be locally recast in $(s, z)$ coordinate system as follows
\begin{align}
\label{s-z-laplace}
\eta_{3\delta}^{\varepsilon}\, \check{L}({\check\phi})=\eta_{\delta}^{\varepsilon}\, \big[-{\mathbf N}(\eta_{3\delta}^{\varepsilon}\, {\check\phi}+{\check\psi})
-{\mathbf E}-p\, \mathbf{w}^{p-1}{\check\psi}\big],
 \quad \forall\, (s, z)\in (-\delta_0/\varepsilon, \delta_0/\varepsilon)\times(0, 1/\varepsilon).
\end{align}

\noindent $\bullet$
The linear operator is
\begin{align}
\label{tilde L check phi}
\check{L}({\check\phi})\, =\, &h_1(\varepsilon z){\check\phi}_{ss}+h_2(\varepsilon z){\check\phi}_{zz}+\varepsilon h_3(\varepsilon z){\check\phi}_{s}+\varepsilon h_4(\varepsilon z) {\check\phi}_{z}
\nonumber\\[2mm]
&+\hat{B}_{1}({\check\phi})+\hat{B}_0({\check\phi})- V(\varepsilon s, \varepsilon z)\, {\check\phi}+p{\mathbf w}^{p-1}{\check\phi},
\end{align}
where
\begin{equation*}
\hat{B}_1(\cdot)=
\varepsilon^2 sh_5(\varepsilon z)\partial_s
+\varepsilon sh_6(\varepsilon z) \partial^2_{sz}
+\varepsilon^2 s^2h_7(\varepsilon z) \partial^2_{ss}
+\varepsilon sh_8(\varepsilon z)\partial^2_{ss},
\end{equation*}
and
\begin{equation}\label{tilteB00}
\hat{B}_0(\cdot)=
\varepsilon sh_9(\varepsilon s, \varepsilon z)\, \partial^2_{zz}
+\varepsilon^2 s^2h_{10}(\varepsilon s, \varepsilon z)\, \partial^2_{zz}
+\varepsilon^2 s^2h_{11}(\varepsilon s, \varepsilon z)\, \partial^2_{sz}
+\varepsilon^2 sh_{12}(\varepsilon s, \varepsilon z)\partial_z.
\end{equation}

\noindent $\bullet$
The error is then expressed in the form
\begin{align}
{\mathbf E}\, =\, &h_1(\varepsilon z){\mathbf w}_{ss}+h_2(\varepsilon z){\mathbf w}_{zz}+\varepsilon h_3(\varepsilon z){\mathbf w}_{s}+\varepsilon h_4(\varepsilon z) {\mathbf w}_{z}
\nonumber\\[2mm]
&+\hat{B}_{1}({\mathbf w})+\hat{B}_0({\mathbf w})- V(\varepsilon s, \varepsilon z)\, {\mathbf w}+{\mathbf w}^{p}.
\end{align}

On the other hand, the boundary condition in (\ref{on the boundary-1}) can also be expressed precisely in local coordinates.
If $z=0$,
\begin{equation}
\label{boundarycondition11}
\eta_{3\delta}^{\varepsilon}\, \mathbb{D}_0({\check\phi})=-\eta^\varepsilon_\delta\, {\mathbf G}_0\qquad{\rm with}\qquad {\mathbf G}_0=\mathbb{D}_0({\mathbf w}).
\end{equation}
And, at $z=1/\varepsilon$ there holds
\begin{equation}
\label{boundarycondition22}
\eta_{3\delta}^{\varepsilon}\, \mathbb{D}_1({\check\phi})=-\eta^\varepsilon_\delta\, {\mathbf G}_1\qquad{\rm with}\qquad {\mathbf G}_1=\mathbb{D}_1({\mathbf w}).
\end{equation}
The operators on the boundary are
\begin{equation}
\label{boundarycondition1}
\mathbb{D}_0\, =\,
{\mathfrak b}_1\partial_z
+{\mathfrak b}_2\varepsilon\, s\, \partial_s
+{\mathfrak b}_3\varepsilon^2 s^2\, \partial_s
+{\mathfrak b}_4\varepsilon s\partial_z
+{\mathfrak b}_5\varepsilon^2 s^2 \partial_z
+{\hat D}_0^0(\cdot),
\end{equation}
  and
\begin{equation}
\label{boundarycondition2}
\mathbb{D}_1\, =\,
{\mathfrak b}_6 \, \partial_z
+{\mathfrak b}_7\varepsilon s\, \partial_s
+{\mathfrak b}_8\varepsilon^2 s^2\, \partial_s
+{\mathfrak b}_9\varepsilon s \, \partial_z
+{\mathfrak b}_{10}\varepsilon^2 s^2  \, \partial_z+{\hat D}_0^1(\cdot),
\end{equation}
where $\hat{D}^0_0\big(\cdot(s, z)\big)=\varepsilon\, \bar{D}^0_0\big(\cdot(t, \theta)\big)$ and $\hat{D}^1_0\big(\cdot(s, z)\big)=\varepsilon\, \bar{D}^1_0\big(\cdot(t, \theta)\big)$.

\medskip
\subsubsection{Further changing of variables}\label{further change}
A further change of variables in equation (\ref{s-z-laplace}) will be chosen in the forms
\begin{equation}\label{vdefine}
{\check\phi}(s, z)\, =\, \alpha(\varepsilon z)\phi(x, z),
\quad\mbox{with}\quad
x\, =\, \beta(\varepsilon z)s,
\end{equation}
where
\begin{equation}
\label{alpha-beta}
\alpha(\theta)\, =\, V(0, \theta)^{\frac {1}{p-1}},
\qquad
\beta(\theta)\, =\, \sqrt{\frac{ V(0, \theta)}{h_1(\theta)}}
=\sqrt{\frac{ V(0, \theta)\big(a_1(0, \theta)|n_1(\theta)|^2+a_2(0, \theta)|n_2(\theta)|^2\big)}
{|a_1(0, \theta)|^2+|a_2(0, \theta)|^2}}.
\end{equation}
It is also convenient to expand
\begin{equation}
\label{Vexpan}
 V(\varepsilon s, \varepsilon z)\, =\, V(0, \varepsilon z)+ V_t(0, \varepsilon z)\cdot \varepsilon s+\frac 12  V_{tt}(0, \varepsilon z)\cdot \varepsilon^2s^2+a_6(\varepsilon s, \varepsilon z)\, \varepsilon^3\, s^3,
\end{equation}
for a smooth function $a_6(t, \theta)$.
In order to express (\ref{tilde L check phi}) and \eqref{boundarycondition11}-(\ref{boundarycondition22}) in terms of
these new coordinates, the following identities will be prepared
 \begin{align*}
{\check\phi}_s=&\, \alpha\beta \phi_x,
 \qquad
{\check\phi}_{ss}=\alpha\beta^2 \phi_{xx},
 \qquad
 {\check\phi}_z=\, \varepsilon \alpha' \phi
 + \varepsilon \alpha \frac{\beta'}{\beta}x\phi_x+\alpha \phi_z,
 \\[2mm]
{\check\phi}_{sz}\, =\, &\, \varepsilon \alpha'\beta \phi_{x}+ \varepsilon\alpha\beta'\phi_{x}
+\varepsilon\alpha\beta'x\phi_{xx}+\alpha\beta \phi_{xz},
\end{align*}
and
\begin{align*}
{\check\phi}_{zz}
\, =\, &\varepsilon^2 \alpha''\phi+2\varepsilon^2\alpha'\frac{\beta'}{\beta}x\phi_{x}+\varepsilon^2\alpha\frac{\beta''}{\beta}x\phi_{x}
\\[2mm]
&+\varepsilon^2\alpha\Big( \frac{\beta'}{\beta} \Big)^2x^2\phi_{xx}
+2\varepsilon\alpha\frac{\beta'}{\beta}x\phi_{xz}+2\varepsilon\alpha'\phi_{z}+\alpha \phi_{zz}.
\end{align*}

\medskip
We can deduce that
\begin{align}
\frac{1}{\alpha\beta^2}{\check L}({\check\phi})
 =\frac{h_2}{\beta^2}\phi_{zz}+h_1 \phi_{xx}-h_1 \phi+\beta^{-2}p\mathbf{w}^{p-1}\phi+B_2(\phi)+B_3(\phi)\equiv L(\phi),
\end{align}
where $B_3(\phi)$ is a linear differential operator defined by
\begin{align}
B_3(\phi)\, =\, & \varepsilon \frac{h_3}{\beta} \phi_x
+\varepsilon \frac{h_4}{\alpha \beta^2}\big[\, \varepsilon \alpha' \phi
 + \varepsilon \alpha \frac{\beta'}{\beta}x\phi_x+\alpha \phi_z\big]
\nonumber\\[2mm]
&+\frac{h_2}{\beta^2}\left[\, \varepsilon ^2\Big|\frac {\beta'}{\beta}\Big|^2x^2\phi_{xx}+2\varepsilon  \frac {\beta'}{\beta}x \phi_{xz}
+\varepsilon^2\frac {\beta{''}}{\beta}x\phi_x\, \right]
\nonumber\\[2mm]
&+\varepsilon^2 \frac{h_2\alpha''}{\alpha \beta^2}\phi
+2\varepsilon^2 \frac{h_2\alpha'}{\alpha \beta^2}\frac {\beta'}{\beta}x \phi_x
+2\varepsilon \frac{h_2\alpha'}{\alpha \beta^2}\phi_z
+\varepsilon^2\frac{h_5}{\beta^2}x\phi_x
\nonumber\\[2mm]
&+\varepsilon \frac{h_6}{\alpha \beta^2}\frac{x}{\beta}\big[\varepsilon \alpha'\beta \phi_{x}+ \varepsilon\alpha\beta'\phi_{x}+\varepsilon\alpha\beta'x\phi_{xx}+\alpha\beta \phi_{xz}\big]
+\varepsilon^2h_7\Big(\frac{x}{\beta}\Big)^2\phi_{xx}
\nonumber\\[2mm]
&+ \varepsilon h_8\frac{x}{\beta}\phi_{xx}
-\left[\, \varepsilon \frac{V_t(0, \varepsilon z)}{\beta^2}\frac{x}{\beta}
+\frac{\varepsilon^2}{2} V_{tt}(0, \varepsilon z)\beta^{-2}\Big(\frac{x}{\beta}\Big)^2\right]\phi.
\label{B3v}
\end{align}
Here
\begin{equation}
\label{B2v}
B_2(\phi)\, =\, \frac{1}{\alpha\beta^2}\, \hat{B}_0({\check\phi})+\frac{1}{\alpha\beta^2}\, a_6(\varepsilon s, \varepsilon z)\, \varepsilon^3 \, s^3\, \alpha\, \phi,
\end{equation}
and $\hat{B}_0({\check\phi})$ is the operator in (\ref{tilteB00}) where derivatives are expressed in terms of $s$ and $z$ through (\ref{vdefine}), $a_6$ is given by (\ref{Vexpan}), and $s$ is replaced by $\beta^{-1}x$.

\medskip
In the coordinates $(x, z)$, the boundary conditions in (\ref{boundarycondition11})-(\ref{boundarycondition22}) can be recast in the sequel.
For $z=0$,
\begin{equation}
\label{boundary-1-x-z}
\eta_{3\delta}^{\varepsilon}\, \big[D_3^0(\phi)+{\mathfrak b}_1\phi_z+D_2^0(\phi)\big]\, =\, -\frac{1}{\alpha}\, \eta^\varepsilon_\delta\, {\mathbf G}_0,
\end{equation}
where
\begin{align}
\label{D30}
D_3^0(\phi)\, =\, &\varepsilon \Big[{\mathfrak b}_2+{\mathfrak b}_1\frac{\beta'}{\beta} \Big]x\phi_{x}
+\varepsilon {\mathfrak b}_1\frac{\alpha'}{\alpha}\phi
+\varepsilon {\mathfrak b}_4 \frac{x}{\beta}\phi_z
\nonumber\\[2mm]
&+\varepsilon^2\Big[{\mathfrak b}_3\Big(\frac x \beta \Big)^2\beta+{\mathfrak b}_4\Big(\frac{x}{\beta}\Big)\Big(\frac{\beta'}{\beta}x\Big)\Big]\phi_x
+\varepsilon^2{\mathfrak b}_4\frac{\alpha'}{\alpha}\Big(\frac{x}{\beta}\Big)\phi
+\varepsilon^2{\mathfrak b}_5\Big(\frac{x}{\beta}\Big)^2\phi_z,
\end{align}
and
\begin{equation}\label{D20}
D_2^0(\phi)\, =\, \frac{1}{\alpha}{\hat D}_0^0({\check\phi})+\varepsilon^3\, {\mathfrak b}_5\, \Big(\frac{x}{\beta}\Big)^2\frac{\alpha'}{\alpha}\phi
+\varepsilon^3\, {\mathfrak b}_5\, \beta'\, \Big(\frac{x}{\beta}\Big)^3\phi_x.
\end{equation}
Similarly, for $z=1/\varepsilon$,  we have
\begin{equation}
\label{boundary-2-x-z}
\eta_{3\delta}^{\varepsilon}\, \big[D_3^1(\phi)+{\mathfrak b}_6\phi_z+D_2^1(\phi)\big]\, =\, -\frac{1}{\alpha}\, \eta^\varepsilon_\delta\, {\mathbf G}_1,
\end{equation}
where
\begin{align}
\label{D31}
D_3^1(\phi)\, =\, &\varepsilon \Big[{\mathfrak b}_7+{\mathfrak b}_6\frac{\beta'}{\beta} \Big]x\phi_{x}
+\varepsilon {\mathfrak b}_6\frac{\alpha'}{\alpha}\phi
+\varepsilon {\mathfrak b}_9 \frac{x}{\beta}\phi_z
\nonumber\\[2mm]
&+\varepsilon^2\Big[{\mathfrak b}_8\Big(\frac x \beta \Big)^2\beta+{\mathfrak b}_9\Big(\frac{x}{\beta}\Big)\Big(\frac{\beta'}{\beta}x\Big)\Big]\phi_x
+\varepsilon^2{\mathfrak b}_9\frac{\alpha'}{\alpha}\Big(\frac{x}{\beta}\Big)\phi
+\varepsilon^2{\mathfrak b}_{10}\Big(\frac{x}{\beta}\Big)^2\phi_z,
\end{align}
and
\begin{equation}\label{D21}
D_2^1(\phi)\, =\, \frac{1}{\alpha}{\hat D}_0^1({\check\phi})+\varepsilon^3\, {\mathfrak b}_{10}\Big(\frac{x}{\beta}\Big)^2\frac{\alpha'}{\alpha}\phi
+\varepsilon^3\, {\mathfrak b}_{10}\, \beta'\, \Big(\frac{x}{\beta}\Big)^3\, \phi_x.
\end{equation}

\medskip
As a conclusion, in local coordinates $(x, z)$, (\ref{equivalent system-1})-(\ref{on the boundary-1}) become
\begin{equation}\label{localproblem1}
 \eta_{3\delta}^{\varepsilon}\, L(\phi)\, =\, (\alpha \beta^{-2})^{-1}\eta_{\delta}^{\varepsilon}\, \Big[-{\mathbf N}(\eta_{3\delta}^{\varepsilon}\, {\check\phi}+{\check\psi})
-{\mathbf E}-p\, \mathbf{w}^{p-1}{\check\psi}\Big],
\end{equation}
\begin{equation}\label{localproblem2}
\eta_{3\delta}^{\varepsilon}\, \Big[D_3^0(\phi)
+{\mathfrak b}_1\phi_z+D_2^0(\phi)\Big]\, =\, -\frac{1}{\alpha}\, \eta^\varepsilon_\delta\, {\mathbf G}_0,
\end{equation}
\begin{equation}\label{localproblem3}
\eta_{3\delta}^{\varepsilon}\, \Big[D_3^1(\phi)
+{\mathfrak b}_6\phi_z+D_2^1(\phi)\Big]\, =\, -\frac{1}{\alpha}\, \eta^\varepsilon_\delta\, {\mathbf G}_1.
\end{equation}

\subsection{The projected problem}\label{section3.3}
For the convenience of presentation, we pause here to give some notation.

\medskip
\noindent{\textbf{ Notation 3:}}
{\it
 Observe that all functions involved in (\ref{localproblem1})-(\ref{localproblem3}) are expressed in $(x, z)$-variables, and the natural domain for those variables can be extended to the infinite strip
\begin{align}
\begin{aligned}\label{domainS}
{\mathcal S}\, =\, &\Big\{\, (x, z)\, :\, -\infty<x<\infty, \,   0<z<1/\varepsilon\, \Big\},
\\[2mm]
\partial_0{\mathcal S}\, =\, \Big\{\, (x, z)\, :\, -\infty<x<\infty, \, &   z=0\, \Big\},
\qquad
\partial_1{\mathcal S}\, =\, \Big\{\, (x, z)\, :\, -\infty<x<\infty, \, z=1/\varepsilon\, \Big\}.
\end{aligned}
\end{align}
Accordingly, we define
\begin{align}
\begin{aligned}\label{domainS1}
{\mathcal {\hat{S}}}\, =\, &\Big\{\, (x, \tilde{z})\, :\, -\infty<x<\infty, \, 0<\tilde{z}<{\ell}/\varepsilon\, \Big\},
\\[2mm]
\partial_0{\mathcal {\hat{S}}}\, =\, \Big\{\, (x, \tilde{z})\, :\, -\infty<x<\infty, \, & \tilde{z}=0\, \Big\},
\quad
\partial_1{\mathcal {\hat{S}}}\, =\, \Big\{\, (x, \tilde{z})\, :\, -\infty<x<\infty, \, \tilde{z}={\ell}/\varepsilon\, \Big\},
\end{aligned}
\end{align}
where ${\ell}$ is a constant defined as
\begin{equation}
\label{ell}
\ell\equiv \int_0^1 \mathcal{Q}(\theta){\rm d}\theta,
\quad\mbox{with }\,
\mathcal{Q}(\theta)=\sqrt{\frac{V(0, \theta)}{h_2(\theta)}}
=\sqrt{ \frac{V(0, \theta)\Big({\mathfrak a}_1(0, \theta)|n_1(\theta)|^2+{\mathfrak a}_2(0, \theta)|n_2(\theta)|^2\Big)} {{\mathfrak a}_1(0, \theta){\mathfrak a}_2(0, \theta)}}.
\end{equation}

\medskip
In all what follows, we will introduce some parameters $\{{f}_j\}_{j=1}^N$ and $\{e_j\}_{j=1}^N$
and assume the validity of the following constraints, for $j=1, \cdots, N$,
\begin{equation}
\label{constraints of f}
\|f_j\|_{H^2(0, 1)}<C|\ln\varepsilon|^2,
\qquad
\, f_{j+1}(\theta)-f_j(\theta)>2|\ln\varepsilon|-4\ln|\ln\varepsilon|,
\end{equation}
\begin{equation}
\label{enorm}
\|e_j\|_{**}\, \equiv\, \|e_j\|_{L^{\infty}(0, 1)}+\varepsilon\|e_j'\|_{L^2(0, 1)}+\varepsilon^2\|e''_j\|_{L^2(0, 1)}\, \leq\, \varepsilon^{\frac{1}{2}},
\end{equation}
where we have used the convention $f_0=-\infty$ and $f_{N+1}=\infty$.
Set
\begin{align}
\textbf{f}&=(f_1, \cdots, f_N),
\quad
f_j \mbox{'s  satisfy bounds } (\ref{constraints of f}), \nonumber
\\[2mm]
\textbf{e}&=(e_1, \cdots, e_N),
\quad
e_j \mbox{'s satisfy  bounds } (\ref{enorm}), \nonumber
\\[2mm]
\textbf{c}&=(c_1, \cdots, c_N),
\quad
\textbf{d}=(d_1, \cdots, d_N),
\quad
c_j \mbox{'s }\, \text{ and } \, d_j \mbox{'s are in } {L^2{(0, 1)}}, \nonumber
\\[2mm]
\textbf{l}_0&=(l_{0, 1}, \cdots, l_{0, N}),
\quad
\textbf{l}_1=(l_{1, 1}, \cdots, l_{1, N}),
\quad
l_{0, j}\mbox{'s }\, \text{ and } \, l_{1, j}\mbox{'s } \mbox { are constants}, \nonumber
\\[2mm]
\textbf{m}_0&=(m_{0, 1}, \cdots, m_{0, N}),
\quad
\textbf{m}_1=(m_{1, 1}, \cdots, m_{1, N}),
\quad
m_{0, j}\mbox{'s }\, \text{ and } \, m_{1, j}\mbox{'s } \mbox { are constants}, \nonumber
\end{align}
and
\begin{equation}
\label{mathcalF}
\mathcal F\, =\, \big\{\, (\textbf{f}, \textbf{e})\, :\, \{f_j\}_{j=1}^N \mbox{ and } \{e_j\}_{j=1}^N \mbox{ satisfy } (\ref{constraints of f})
\mbox{ and } (\ref{enorm}) \mbox{ respectively}\, \big\}.
\end{equation}
\qed
}

One of the left job is to find the local forms of the approximate solution $\mathbf{w}$ with the constraint (\ref{constraint condition}) and also of the error $\mathbf{E}$. We recall the transformation in (\ref{vdefine})-(\ref{alpha-beta}), and then define the local form of the approximate solution ${\mathbf w}$ by the relation
\begin{equation}
\label{vdefine1}
\eta^{\varepsilon}_{10\delta}(s)\, {\mathbf w}\, =\, \eta^{\varepsilon}_{3\delta}(s)\, \alpha(\varepsilon z)\, v(x, z)
\quad\mbox{with}\quad
x\, =\, \beta(\varepsilon z)s.
\end{equation}
The error $\mathbf{E}$ can be locally recast in $(x, z)$ coordinate system by the relation
 \begin{equation}\label{locallyE}
 \frac{1}{\alpha\beta^2}\, \eta_{\delta}^{\varepsilon}(s)\, \mathbf{E}
 =\eta_{\delta}^{\varepsilon}(s)\, \mathcal{E},
 \end{equation}
 where
 \begin{equation}\label{sv}
 \mathcal{E}\, =\, S(v) \quad {\text{with}}\quad
  S(v)\, =\, \frac{h_2}{\beta^2}v_{zz}+h_1\big[v_{xx}-v+v^{p}\big]+B_2(v)+B_3(v),
 \end{equation}
 with the operator $B_3$ and $B_2$ defined in (\ref{B3v})-(\ref{B2v}).
 In the coordinates $(x, z)$, the boundary errors can be recast as follows. For $z=0$,
\begin{equation}\label{g0}
\frac{1}{\alpha}\eta_\delta^\varepsilon(s)\, {\mathbf G}_0\, =\, \eta_\delta^\varepsilon(s)\, g_0
\quad\mbox{with}\quad
g_0=D^0_3(v)+{\mathfrak b}_1v_z+D_2^0(v),
\end{equation}
and also for $z=1/\varepsilon$,
\begin{align}\label{g1}
\frac{1}{\alpha}\eta_\delta^\varepsilon(s)\, {\mathbf G}_1\, =\, \eta_\delta^\varepsilon(s)\, g_1
\quad\mbox{with}\quad
g_1=D^1_3(v)+{\mathfrak b}_6v_z+D_2^1(v).
\end{align}
It is of importance that (\ref{locallyE}), (\ref{g0}) and (\ref{g1}) hold only in a small neighbourhood of $\Gamma_\varepsilon$. Hence we will consider $v$, $S(v)$ as functions of the variables $x$ and $z$ on $\mathcal S$, and also $g_0, ~g_1$ on $\partial_0{\mathcal S}$ and $\partial_1{\mathcal S}$ in the sequel.
We will find $v=v_4$ in (\ref{basic approximate}) step by step in Section 4, so that $\mathbf{w}$ will be given in (\ref{globalapproximation}) with the property in (\ref{constraint condition}).
In fact, to deal with the resonance, in addition to the parameters $\mathbf{f}$ and $\mathbf{h}$, we shall add one more parameter, say $\mathbf{e}$, in the approximate solution $v_4$.
The exact forms of the error terms $S(v_4)$, $g_0$ and $g_1$ will be given in (\ref{E1-d}) and (\ref{E11}).

\medskip
To make suitable extension of (\ref{localproblem1})-(\ref{localproblem3}), we define an operator on the whole strip $\mathcal S$ in the form
\begin{align}
\label{nonlocal-1}
\begin{split}
\mathcal{L}(\phi)\, \equiv\, &\, \frac{h_2}{\beta^2}\, \phi_{zz}+h_1\big[\phi_{xx}-\phi+pw^{p-1}\phi\big]+\chi(\varepsilon|x|)\, B_3(\phi)
\,-\, h_1(1-\eta_{3\delta}^{\varepsilon})pv_4^{p-1}\, \phi\,
\quad\mbox{in } \mathcal S,
\end{split}
\end{align}
and also the operators
\begin{align}
{\mathcal D}_1(\phi)\, =\, \chi(\varepsilon|x|)\, D_3^1(\phi)+{\mathfrak b}_6\phi_z+\, \chi(\varepsilon|x|)\, D_2^1(\phi)
\quad{\rm on}\ \partial_1\mathcal S,
\label{D1mathcal}
\\[2mm]
{\mathcal D}_0(\phi)\, =\, \chi(\varepsilon|x|)\, D_3^0(\phi)+{\mathfrak b}_1\phi_z+\, \chi(\varepsilon|x|)\, D_2^0(\phi)
\quad{\rm on}\ \partial_0\mathcal S,
\label{D0mathcal}
\end{align}
where $\chi(r)$ is a smooth cut-off function which equals $1$ for $0\leq r<10\delta$ that vanishes identically for $r>20\delta$.

\medskip
Rather than solving problem (\ref{localproblem1})-(\ref{localproblem3}) directly, we deal with the following projected problem: for each pair of parameters $\mathbf{f}$ and $\mathbf{e}$ in $\mathcal F$,
finding functions $\phi\in H^2(\mathcal S), ~\mathbf{c}, ~\mathbf{d} \in L^2(0, 1)$ and constants $\mathbf{l}_0, ~\mathbf{l}_1, $~$~\mathbf{m}_0, ~\mathbf{m}_1$ such that
\begin{equation}
\label{system-1}
\mathcal{L}(\phi)\, =\, \eta_{\delta}^{\varepsilon}(s)\, \big[
-\mathcal{E}-{\mathcal N}({\phi})\big]+\sum_{j=1}^Nc_j(\varepsilon z)\chi(\varepsilon|x|) w_{j, x}+\sum_{j=1}^Nd_j(\varepsilon z)\chi(\varepsilon|x|) Z_j\quad\mbox{in } \mathcal S,
\end{equation}
\begin{equation}
\label{system-2}
{\mathcal D}_1(\phi)\, =\, \eta_{\delta}^{\varepsilon}(s)g_1+\sum_{j=1}^Nl_{1, j}\chi(\varepsilon|x|) w_{j, x}+\sum_{j=1}^Nm_{1, j}\chi(\varepsilon|x|) Z_j\quad\mbox{on } \partial_1 \mathcal S,
\end{equation}
\begin{equation}
\label{system-3}
{\mathcal D}_0(\phi)\, =\, \eta_{\delta}^{\varepsilon}(s) g_0+\sum_{j=1}^Nl_{0, j}\chi(\varepsilon|x|) w_{j, x}+\sum_{j=1}^Nm_{0, j}\chi(\varepsilon|x|) Z_j\quad\mbox{on } \partial_0 \mathcal S,
\end{equation}
\begin{equation}
\label{system-4}
\int_{{\mathbb R}}\phi(x, z)w_{j, x}\, {\rm d}x \, =\, \int_{{\mathbb R}}\phi(x, z)Z_j\, {\rm d}x \, =\, 0, \quad 0<z<\frac {1}{\varepsilon}, \quad \forall\, j=1, \cdots, N,
\end{equation}
where we have denoted
\begin{equation}
{\mathcal N}(\phi)\, =\, \big[v_4+\phi+\psi(\phi)\big]^p-v_4^{p}-pv_4^{p-1}\phi.
\end{equation}
The functions $w_j, \, Z_j, \, j=1, \cdots, N$ are given in (\ref{wjZj}).
This problem has a unique solution $\phi$ such that ${\check\phi}$ satisfies (\ref{decay property}).
 \begin{proposition}
 \label{prop}There is a number ${\tilde\tau}>0$ such that for all $\varepsilon$ small enough  and all parameters $(\mathbf{f}, \mathbf{e})$ in $\mathcal F$,  problem (\ref{system-1})-(\ref{system-4}) has a unique solution $\phi=\phi(\mathbf{f}, \mathbf{e})$ which satisfies
$$
\|\phi\|_{H^2(\mathcal S)}\, \leq\, {\tilde\tau} \varepsilon^{3/2}|\ln\varepsilon|^q,
$$
$$
\big\|\phi\big\|_{L^{\infty}(|x|>\delta/\varepsilon)}+\big\|\nabla \phi\big\|_{L^{\infty}(|x|>\delta/\varepsilon)}\, \leq\, e^{-\rho\delta/\varepsilon}.
$$
 Moreover, $\phi$ depends Lipschitz-continuously on the parameters $\mathbf{f}$ and $\mathbf{e}$ in the sense of the estimate
\begin{align}
\label{characteriztion}
\|\phi(\mathbf{f}_1, \mathbf{e}_1)\,-\, \phi(\mathbf{f}_2, \mathbf{e}_2)\|_{H^2(\mathcal S)}
\, \leq\,
C\varepsilon^{3/2}|\ln\varepsilon|^q\big [\, \|\mathbf{f}_1-\mathbf{f}_2\|_{H^2(0, 1)}\, +\, \|\mathbf{e}_1-\mathbf{e}_2\|_{**}\, \big].
\end{align}
\end{proposition}

\proof
The proof is similar as that for Proposition 5.1 in \cite{delPKowWei2007}, which will be omitted here.
\qed

\medskip
\medskip
We conclude this section by stating the following announcements:
\\[1mm]
\noindent $\bullet$
As we have said in the above, we shall construct the approximate solution in Section \ref{section4}.
\\[1mm]
\noindent $\bullet$
To find a real solution to \eqref{problemafterscaling},
the reduction procedure will be carried out in Sections \ref{section5} and \ref{sectionsolvingreducedequation} to kill the Langrange multipliers in \eqref{system-1}-\eqref{system-4}. This can be done by suitable choice of the parameters $\mathbf{f}=(f_1,\cdots, f_N)$ and $\mathbf{e}=(e_1,\cdots, e_N)$.
We will first derive the equations involving the parameters $f_1,\cdots, f_N$ and $e_1,\cdots, e_N$ in Section \ref{section5},
and then solve the coupled system involving the equations in Section \ref{sectionsolvingreducedequation}.

\section{The local  approximate solutions}
\label{section4}
\setcounter{equation}{0}

The main objective of this section is to construct the local form $v$ of the approximation ${\mathbf w}$ (cf. (\ref{globalapproximation})) and
then evaluate its error $\mathcal{E}$, $g_0$, $g_1$ in the coordinate system $(x, z)$.

\subsection{The first approximate solution}
Recall the notation in Section \ref{section3.3}.
 For a fixed integer $N>1$, we assume that the locations of $N$ concentration layers are characterized by sets
 $$
 \big\{\, (x, z) \, :\, x=\beta(\varepsilon z) f_j(\varepsilon z)+\beta(\varepsilon z)\, h(\varepsilon z)\, \big\}, \qquad j=1, \cdots, N,
 $$
 in the coordinates $(x, z)$.
The function $ h$ satisfies
\begin{equation}\label{assumptionofh}
{\mathfrak b}_1\, h'(0)\,-\, {\mathfrak b}_2 h(0)\, =\, 0,
\qquad
{\mathfrak b}_6\, h'(1)\,-\, {\mathfrak b}_7\, h(1)\, =\, 0.
\end{equation}
In fact, $h$ will be chosen by solving (\ref{equation of h})
and $f_j$'s can be determined in the reduction procedure.

\medskip
By recalling $w$ given in \eqref{wsolution} and $Z$ in (\ref{lambda0}), we set
\begin{equation}\label{wjZj}
w_j(x)\, =\, w(x_j),
\qquad
Z_j(x)\, =\, Z(x_j),
\end{equation}
with
$$
x_j\, =\, x-\beta(\varepsilon z)\, f_j(\varepsilon z)-\beta(\varepsilon z)\, h(\varepsilon z),
$$
and then define the first approximate solution by
\begin{equation}
\label{vvdefine}
{v_1}(x, z)\, =\sum_{j=1}^Nw_j(x).
\end{equation}

For every fixed $n$ with $1\leq n\leq N$, we consider the following set
\begin{equation}
\mathfrak{A}_n
=\Bigg\{(x, z)\in{\mathcal S} \, :\, \frac {\beta f_{n-1}(\varepsilon z)+\beta f_n(\varepsilon z)}{2}\leq x- \beta h\leq \frac {\beta f_n(\varepsilon z)+\beta f_{n+1}(\varepsilon z)}{2}\Bigg\}.
\end{equation}
For $(x, z)\in \mathfrak{A}_n$, we expand $S(v_1)$ by gathering terms of $\varepsilon$ and those of order $\varepsilon^2:$
\begin{align}\label{sv1-gather}
S(v_1)\, =\,
&\sum_{j=1}^N\varepsilon \, \big[h_8(f_j+h)w_{j, xx}- \frac{V_t(0, \varepsilon z)}{\beta^2}\, (f_j+h)w_j\big]
\nonumber\\
&+\sum_{j=1}^N\frac{\varepsilon}{\beta}\big[h_3w_{j, x}- \, \frac{V_t(0, \varepsilon z)}{\beta^2} x_jw_j+h_8x_{j}w_{j, xx}\big]
\nonumber\\
&-\sum_{j=1}^N\varepsilon^2\Bigg[\Big( -\frac {h_5}{\beta}{f_j}
+h_2\frac {f''_j}{\beta}
+h_2\frac {2\beta'}{\beta^2}{f_j'}
+h_2\frac {2\alpha'}{\alpha \beta}{f_j'} \Big) w_{j, x}
\nonumber\\
&\qquad \qquad +h_2\frac {2\beta'}{\beta^2}{f_j'x_jw_{j, xx}}
-\frac{2h_7}{\beta}f_jx_jw_{j, xx}
+\frac { V_{tt}(0, \varepsilon z)}{\beta^3}{f_jx_jw_j}\Bigg]
\nonumber\\
&+\sum_{j=1}^N\varepsilon^2\Big[
\Big( h_2f_j'^2\, +2h_2f_j'h'-h_6\, f_j\, f_j'\, -h_6\, f_j'\, h+h_7f_j^2\Big)w_{j, xx}
 -\frac{1}{2}\, \beta^{-2}\, V_{tt}(0, \varepsilon z)\, f_j^2\, w_j\Big]
\nonumber\\
&+\sum_{j=1}^N\frac {\varepsilon^2}{\alpha\beta^2}\, \Big[h_6\Big(-\alpha \beta f_j' +\alpha \beta'f_j \Big)x_jw_{j, xx}
+ h_6\Big( \alpha'\beta f_j+\alpha\beta'f_j \Big)w_{j, x}
- h_4\alpha\beta f_j' w_{j, x} \Big]
\nonumber\\
&+\sum_{j=1}^N\frac{\varepsilon^2}{\beta^{2}}
\Big[\Big(h_5+h_2\frac {\beta{''}}{\beta}
+h_2\frac {2\alpha'\beta'}{\alpha\beta} \Big){x_jw_{j, x}}
+h_2\frac {|\beta'|^2}{\beta^2}x_j^2w_{j, xx}
+h_2\beta^2\, h'^2\, w_{j, xx}
+h_2\frac {\alpha{''}}{\alpha}w_j
\nonumber\\
&\quad\quad\quad\quad
-\frac {1}{2\beta^2} V_{tt}(0, \varepsilon z)x_j^2w_j
-\frac 12 V_{tt}(0, \varepsilon z)(2f_jh+h^2)w_j
+h_7\beta^2(2f_jh+h^2)+h_7 x^2w_{j, xx}\Big]
\nonumber\\
&-\sum_{j=1}^N\varepsilon^2\Big[\Big( -\frac {h_5}{\beta}{h}
+h_2\frac {h''}{\beta}
+h_2\frac {2\beta'}{\beta^2}{h'}
+h_2\frac {2\alpha'}{\alpha \beta}{h'} \Big) w_{j, x}
\nonumber\\
&\qquad \qquad +h_2\frac {2\beta'}{\beta^2}{h'x_jw_{j, xx}}
-\frac{2h_7}{\beta}hx_jw_{j, xx}
+\frac { V_{tt}(0, \varepsilon z)}{\beta^3}{hx_jw_j}\Big]
\nonumber\\
&+\sum_{j=1}^N\frac {\varepsilon^2}{\alpha \beta^2}
\Big[ h_6\Big( \alpha'+\frac {\alpha \beta'}{\beta} \Big)x_jw_{j, x}
+h_6\frac {\alpha \beta'}{\beta}{x_j^2w_{j, xx}}
- h_6\alpha \beta^2(f_jh'+hh')w_{j, xx}
+h_4 \alpha'w_j\Big]
\nonumber\\
&+\sum_{j=1}^N\, \frac{\varepsilon^2}{\alpha \beta^2}\Big[h_6\Big(-\alpha \beta h' +\alpha \beta'h \Big)x_jw_{j, xx}
+ h_6\Big( \alpha'\beta h+\alpha\beta'h\Big)w_{j, x}
- h_4\alpha\beta h' w_{j, x} \Big]
+B_4(v_1)
\nonumber\\
\, \equiv\, & \varepsilon\sum_{j=1}^N S_{1, j}\, +\, \varepsilon \sum_{j=1}^NS_{2, j}\, +\, \varepsilon^2\sum_{j=1}^NS_{3, j}\, +\, \varepsilon^2\sum_{j=1}^NS_{4, j}\, +\, \varepsilon^2\sum_{j=1}^NS_{5, j}
\, +\, \varepsilon^2\sum_{j=1}^NS_{6, j}\,
\nonumber\\
&+\, \varepsilon^2\sum_{j=1}^NS_{7, j}\, +\, \varepsilon^2\sum_{j=1}^NS_{8, j}\, +\, \varepsilon^2\sum_{j=1}^NS_{9, j}\, +\, B_4(v_1),
\end{align}
where
\begin{align}\label{B2v1}
B_4(v_1)=&\, \frac{1}{\alpha\beta^2} \big[\, \hat{B}_0(v_1)+a_6(\varepsilon s, \varepsilon z)\, \varepsilon^3 \, s^3\, {v_1}\big]
\nonumber\\[2mm]
&+h_1\Big[pw_n^{p-1}(v_1-w_n)-\sum_{j\neq n}{w_j}^p+\frac{1}{2}p(p-1)w_n^{p-2}(v_1-w_n)^2\Big]
+\max_{j\neq n}O(e^{-3 |\beta f_j-x|}).
\end{align}
Here $B_4(v_1)$ turns out to be of size $(\varepsilon^3+\varepsilon^{\frac {3}{2}}|\ln\varepsilon|^q)$. Let us observe that the quantities $S_{1, j}$, $S_{3, j}$, $S_{5, j}$, $S_{7, j}$ and $S_{9, j}$ are odd functions of $x_j$, while $S_{2, j}$, $S_{4, j}$, $S_{6, j}$ and $S_{8, j}$ are even functions of $x_j$.

\medskip
Using  the assumptions of $h$ in (\ref{assumptionofh}), the boundary errors can be formulated as follows.
For $z=0$, the error terms have the expressions
\begin{align}
\label{boundary-1-v_1-x-z}
&\varepsilon \sum_{j=1}^N\beta\Big[{\mathfrak b}_2 f_j-{\mathfrak b}_1f_j'\Big] w_{j, x}
+\varepsilon \Big[{\mathfrak b}_2+{\mathfrak b}_1 \frac{\beta'}{\beta} \Big]\sum_{j=1}^Nx_jw_{j, x}
+\varepsilon {\mathfrak b}_1\frac{\alpha'}{\alpha}\sum_{j=1}^Nw_{j}
\nonumber\\[2mm]
&+\varepsilon^2 {\mathfrak b}_4\sum_{j=1}^N\big( -\beta'\, f_j-\beta\, f_j'-\beta'\, h-\beta\, h'\big)\Big(\frac {x_j} \beta+f_j+h\Big)w_{j, x}
\nonumber\\[2mm]
&+\varepsilon^2\, \sum_{j=1}^N\Big[{\mathfrak b}_3\, \Big(\frac {x_j} \beta+f_j+h\Big)^2\, \beta+{\mathfrak b}_4\, \Big(\frac {x_j} \beta+f_j+h\Big)^2\beta'\Big]\, w_{j, x}
\nonumber\\[2mm]
&+\varepsilon^2{\mathfrak b}_4\frac{\alpha'}{\alpha}\sum_{j=1}^N\Big(\frac {x_j} \beta+f_j+h\Big)w_{j}
\nonumber\\[2mm]
&+\varepsilon^3{\mathfrak b}_5\sum_{j=1}^N\big( -\beta'\, f_j-\beta\, f_j'-\beta'\, h-\beta\, h'\big)\Big(\frac {x_j} \beta+f_j+h\Big)^2w_{j, x}
+D_2^0(v_1).
\end{align}
Similarly, for $z\, =\, 1/\varepsilon$, we have the terms
\begin{align}\label{boundary-2-v_1-x-z}
&\varepsilon \sum_{j=1}^N\, \beta\Big[{\mathfrak b}_7f_j-{\mathfrak b}_6f_j'\Big] w_{j, x}
+\varepsilon \Big[{\mathfrak b}_7+{\mathfrak b}_6\frac{\beta'}{\beta} \Big]\sum_{j=1}^Nx_jw_{j, x}
+\varepsilon {\mathfrak b}_6\frac{\alpha'}{\alpha}\sum_{j=1}^Nw_{j}
\nonumber\\[2mm]
&+\varepsilon^2 {\mathfrak b}_9\sum_{j=1}^N\big( -\beta'\, f_j-\beta\, f_j'-\beta'\, h-\beta\, h'\big)\Big(\frac {x_j} \beta+f_j+h\Big)w_{j, x}
\nonumber\\[2mm]
&+\varepsilon^2\, \sum_{j=1}^N\Big[{\mathfrak b}_8\, \Big(\frac {x_j} \beta+f_j+h\Big)^2\, \beta+{\mathfrak b}_9\, \Big(\frac {x_j} \beta+f_j+h\Big)^2\beta'\Big]\, w_{j, x}
\nonumber\\[2mm]
&+\varepsilon^2{\mathfrak b}_9\frac{\alpha'}{\alpha}\sum_{j=1}^N\Big(\frac {x_j} \beta+f_j+h\Big)w_{j}
\nonumber\\[2mm]
&+\varepsilon^3{\mathfrak b}_{10}\sum_{j=1}^N\big( -\beta'\, f_j-\beta\, f_j'-\beta'\, h-\beta\, h'\big)\Big(\frac {x_j} \beta+f_j+h\Big)^2w_{j, x}
+D_2^1(v_1).
\end{align}

\subsection{ Interior correction layers}

We now want to construct correction terms and establish a further approximation to a real
solution that eliminates the terms of order $\varepsilon$ in the errors.
Inspired by the method in Section 2 of \cite{delPKowWei2007}, for fixed $z$,  we need a solution of
\begin{equation}\label{equationofphi1}
-\phi_{1, xx}+\phi_{1}-p{w}^{p-1}\phi_{1}\, =\, \sum_{j=1}^N S_{1, j}\, +\, \sum_{j=1}^NS_{2, j},
\qquad
\phi_1(\pm \infty)=0.
\end{equation}
As it is well known, this problem is solvable provided that
\begin{equation}
\label{orth-condition}
\int_{\Bbb R}(S_{1, j}+S_{2, j})w_{j, x}{\rm d}x=0.
\end{equation}
Furthermore, the solution is unique under the constrain
\begin{equation}
\label{understrait}
\int_{\Bbb R}\phi_1w_{j, x}{\rm d}x=0.
\end{equation}
Since $S_{1, j}$ is odd in the variable $x_j$, we have
\begin{equation*}
\begin{split}
&\int_{\Bbb R}(S_{1, j}+S_{2, j})w_{j, x}{\rm d}x =\int_{\Bbb R}S_{2, j}w_{j, x}{\rm d}x
\\[2mm]
&= \frac{1}{\beta}\Big[   h_3\int_{\Bbb R}w_{j, x}^2{\rm d}x- \, V_t(0, \varepsilon z) \beta^{-2} \int_{\Bbb R}x_jw_jw_{j, x}{\rm d}x+h_8\int_{\Bbb R}x_{j}w_{j, xx}w_{j, x}{\rm d}x\Big]
\\[2mm]
&= \frac{1}{\beta}\Big[   h_3\, +\, \sigma V_t(0, \varepsilon z) \beta^{-2} -\frac{1}{2} h_8\Big]\int_{\Bbb R}w_{j, x}^2{\rm d}x,
\end{split}
\end{equation*}
where we have used the fact
\begin{equation}\label{relationofwx}
-2\int_{{\mathbb R}}xww_x\, {\rm d}x\,
=\int_{{\mathbb R}}w^2\, {\rm d}x
\, =\, 2\sigma \int_{{\mathbb R}}w_x^2\, {\rm d}x,
\qquad
\int_{{\mathbb R}}w_x^2\, {\rm d}x
 \, =\, -2\int_{{\mathbb R}}xw_xw_{xx}\, {\rm d}x.
\end{equation}
Thanks to the stationary condition \eqref{stationary1}, then we obtain
\begin{equation} \label{relation-1}
h_3\, = \, \frac{1}{2}h_8 -\sigma \frac{V_t(0, \varepsilon z)}{\beta^2}.
\end{equation}
For more details, the reader can refer the computations in Appendix \ref{appendixA}.
Therefore, we have verified the condition \eqref{orth-condition}.

\medskip
Then the solution $\varepsilon\phi_1$ can be written in the form
 \begin{equation}
 \label{phi1}
 \varepsilon\, \phi_1\, =\, \varepsilon\, \sum_{j=1}^N \phi_{1, j}
 =\, \varepsilon\sum_{j=1}^N(\phi_{10, j}+\phi_{11, j}+\phi_{12, j}+\phi_{13, j}),
 \end{equation}
where
\begin{align}
\label{phi10}
\phi_{10, j}(x, z)\, =\, a_{10}(\varepsilon z)\omega_{0, j}(x)=\, a_{10}(\varepsilon z)\omega_{0}(x_j),
\end{align}
\begin{align}
\label{phi11}
\phi_{11, j}(x, z)\, =\, a_{11}(\varepsilon z)\omega_{1, j}(x)=\, a_{11}(\varepsilon z)\omega_{1}(x_j),
\end{align}
\begin{align}\label{phi12}
\phi_{12, j}(x, z)\, =\, \big[f_j(\varepsilon z)+h(\varepsilon z)\big]a_{12}(\varepsilon z)\omega_{2, j}(x)=\, \big[f_j(\varepsilon z)+h(\varepsilon z)\big]a_{12}(\varepsilon z)\omega_{2}(x_j),
\end{align}
\begin{align}\label{phi13}
\phi_{13, j}(x, z)\, =\, \big[f_j(\varepsilon z)+h(\varepsilon z)\big]a_{13}(\varepsilon z)\omega_{3, j}(x)=\, \big[f_j(\varepsilon z)+h(\varepsilon z)\big]a_{13}(\varepsilon z)\omega_{3}(x_j),
\end{align}
with
\begin{align}
\label{a10a11}
a_{10}(\theta)\, =\, \frac{h_3(\theta)}{\beta(\theta) h_1(\theta)},
\qquad
&a_{11}(\theta)\, =\, \frac{h_8(\theta)}{\beta(\theta) h_1(\theta)},
\\[2mm]
\label{a12a13}
a_{12}(\theta)\, =\, - \frac{V_t(0, \theta)}{\beta^2(\theta) h_1(\theta)},
\qquad
&a_{13}(\theta)\, =\, \frac{h_8(\theta)}{h_1(\theta)}.
\end{align}
The functions $\omega_{0}, \omega_{1}$ are respectively the unique odd solutions to
\begin{equation}
\label{w{0}}
-\omega_{0, xx}+\omega_{0}-p{w}^{p-1}\omega_{0}
\, =\, w_{x}+\sigma^{-1} xw, \qquad\int_{{\mathbb R}}\omega_{0}w_{x}\, {\rm d}x
\, =\, 0,
\end{equation}
\begin{equation}
\label{w{1}}
-\omega_{1, xx}+\omega_{1}-p{w}^{p-1}\omega_{1}
\, =\, -\frac{1}{2 \sigma} xw+xw_{xx}, \qquad\int_{{\mathbb R}}\omega_{1}w_{x}\, {\rm d}x
\, =\, 0,
\end{equation}
and $\omega_{2}, \omega_{3}$ are respectively the unique even functions satisfying
\begin{align}
\label{w{2}}
-\omega_{2, xx}+\omega_{2}-p{w}^{p-1}\omega_{2}\, =\, w, \qquad\int_{{\mathbb R}}\omega_{2}w_{x}\, {\rm d}x
\, =\, 0,
\end{align}
\begin{align}
\label{w{3}}
-\omega_{3, xx}+\omega_{3}-p{w}^{p-1}\omega_{3}\, =\, w_{xx}, \qquad\int_{{\mathbb R}}\omega_{3}w_{x}\, {\rm d}x
\, =\, 0.
\end{align}
Observe that $\varepsilon\, \phi_1$ is of size $O(\varepsilon)$.

\subsection{Boundary corrections}\
In the following, we want to cancel the boundary error terms of first order in $\varepsilon$  given in (\ref{boundary-1-v_1-x-z}) and (\ref{boundary-2-v_1-x-z}), i.e.,
\begin{equation*}
\varepsilon \Big[{\mathfrak b}_2+{\mathfrak b}_1\frac{\beta'(0)}{\beta(0)} \Big]x_jw_{j, x}
+\varepsilon\, {\mathfrak b}_1\frac{\alpha'(0)}{\alpha(0)}w_j,
\end{equation*}
and
\begin{equation*}
\varepsilon \Big[{\mathfrak b}_7+{\mathfrak b}_6 \frac{\beta'(1)}{\beta(1)} \Big]x_jw_{j, x}
+\varepsilon\, {\mathfrak b}_6\frac{\alpha'(1)}{\alpha(1)}w_j.
\end{equation*}
On the other hand, the boundary terms
$$
\varepsilon\, \Big[{\mathfrak b}_7f_j(1)\,-\, {\mathfrak b}_6\, f_j'(1)\Big] w_{j, x}
\quad\mbox{ and }\quad
\varepsilon \, \Big[{\mathfrak b}_2f_j(0)\,-\, {\mathfrak b}_1\, f_j'(0)\Big] w_{j, x}
$$
will be dealt with by the standard reduction  procedure in Sections \ref{section5}-\ref{sectionsolvingreducedequation}.

\medskip
This can be done by the methods in Section 2.2 of \cite{JWeiYang2007}.
By defining two constants
\begin{equation*}
{\bf c}_{1}\, =\, \Big[{\mathfrak b}_7+{\mathfrak b}_6\frac{\beta'(1)}{\beta(1)} \Big]\int_{{\mathbb R}}x\, w_x\, Z\, {\rm d}x
\, +\,
{\mathfrak b}_6\frac{\alpha'(1)}{\alpha(1)}\int_{{\mathbb R}}\, w\, Z\, {\rm d}x,
\end{equation*}
\begin{equation*}
{\bf c}_{0}\, =\, \Big[{\mathfrak b}_2+{\mathfrak b}_1\frac{\beta'(0)}{\beta(0)} \Big]\int_{{\mathbb R}}x\, w_x\, Z\, {\rm d}x
\, +\,
{\mathfrak b}_1\frac{\alpha'(0)}{\alpha(0)}\int_{{\mathbb R}}\, w\, Z\, {\rm d}x,
\end{equation*}
and also a function
\begin{align}\label{btheta}
A({\tilde\theta})\, =\, &\, \frac {{\bf{c}_{0}}\cos\big(\sqrt{\lambda_0}\, {\ell}/\varepsilon\big)-\bf{c}_{1}}{\sqrt{\lambda_0}\sin\big(\sqrt{\lambda_0}\, {\ell}/\varepsilon\big)}\cos\big(\varepsilon^{-1}\, \lambda_0\, {\tilde\theta}\big)
\, +\,
\frac {\bf{c}_{0}}{\sqrt{\lambda_0}}\sin\big(\varepsilon^{-1}\, \lambda_0\, {\tilde\theta}\big),
\end{align}
where the constant $\ell$ is given in (\ref{ell}), we choose
\begin{align}
\label{tildephi}
\phi_{21, j}(x, z)\, =\, A\big({\tilde{\bf d}}(\varepsilon z)\big)Z_j=\, A\big({\tilde{\bf d}}(\varepsilon z)\big)Z(x_j),
\quad\mbox{with}\ {\tilde{\bf d}}(\theta)\, =\, \int_0^{\theta}\mathcal{Q}(r)\, {\rm d}r,
\end{align}
where ${\mathcal Q}$ is the function given in \eqref{ell}.
On the other hand, by Corollary 2.4 in \cite{JWeiYang2007}, we then  find a unique solution  $\phi_{*}$ of the following problem:
\begin{align*}
\Delta \phi_{*}-\phi_{*}+p{w}^{p-1}\phi_{*}\, =\, 0\quad&\mbox{in } \hat{\mathcal S},
\\[2mm]
\frac {\partial \phi_{*}}{\partial {\tilde z}}\, =\, \Big[{\mathfrak b}_7+{\mathfrak b}_6\frac{\beta'(1)}{\beta(1) }\Big]xw_{x}+\, {\mathfrak b}_6 \frac{\alpha'(1)}{\alpha(1)}w&-{\bf{c}_{1}}Z
\quad\mbox{on } \partial_1\hat{\mathcal S},
\\[2mm]
\frac {\partial \phi_{*}}{\partial {\tilde z}}\, =\, \Big[{\mathfrak b}_2+{\mathfrak b}_1\frac{\beta'(0)}{\beta(0) }\Big]xw_{x}+\, {\mathfrak b}_1 \frac{\alpha'(0)}{\alpha(0)}w&-{\bf{c}_{0}}Z
\quad\mbox{on } \partial_0\hat{\mathcal S},
\end{align*}
where $\hat{\mathcal S}$, $\partial_0\hat{\mathcal S}$
and $\partial_1\hat{\mathcal S}$
are defined in (\ref{domainS1}). Moreover, $\phi_{*}$ is even in $x$.
By the diffeomophism
\begin{align}
\label{mathfrak-a-function}
{\Upsilon}: [0,  {1}/{\varepsilon}]\rightarrow [0, {{\ell}}/{\varepsilon}],
\qquad
{\Upsilon}(z)\, =\,
\varepsilon^{-1}\int_0^{\varepsilon z}\mathcal{Q}(\theta){\rm d}\theta,
\end{align}
where ${\mathcal Q}$ is the function given in \eqref{ell}, we define
$$
\phi_{22, j}(x, z)\, =\, \phi_{*, j}\, =\, \phi_{*}\big(x_j, {\Upsilon}(z)\big).
$$
Hence, $\phi_{22, j}$ satisfies the following problem:
\begin{align}\label{boundary-1}
\frac{h_2}{\beta^2}\partial_{zz}\phi_{22, j}&+h_1\big[\partial_{xx}\phi_{22, j}-\phi_{22, j}+p{w_j}^{p-1}\phi_{22, j}\big]\, =\, \varepsilon \frac{h_2}{\beta^2}\mathcal{Q}{'}\phi_{{*}, \tilde z}\big(x_j, {\Upsilon}(z)\big)
\quad\mbox{in }   {\mathcal S},
\nonumber\\[2mm]
\frac {\partial \phi_{22, j}}{\partial z}\, =\, &-\frac{\mathcal{Q}(1)}{{\mathfrak b}_6}\Big\{\Big[{\mathfrak b}_7+{\mathfrak b}_6\frac{\beta'(1)}{\beta(1) }\Big]x_jw_{j, x}+\, {\mathfrak b}_6 \frac{\alpha'(1)}{\alpha(1)}w_j-{\bf{c}_{1}}Z_j\, \Big\}
\quad\mbox{on }  \partial_1{\mathcal S},
\\[2mm]
\frac {\partial \phi_{22, j}}{\partial z}\, =\, &- \frac{\mathcal{Q}(0)}{{\mathfrak b}_1}\Big\{\Big[{\mathfrak b}_2+{\mathfrak b}_1\frac{\beta'(0)}{\beta(0) }\Big]x_jw_{j, x}+\, {\mathfrak b}_1 \frac{\alpha'(0)}{\alpha(0)}w_j-{\bf{c}_{0}}Z_j\Big\}
\quad\mbox{on }  \partial_0{\mathcal S},
\nonumber
  \end{align}
where $\mathcal S$, $\partial_0{\mathcal S}$ and $\partial_1{\mathcal S}$ are defined in (\ref{domainS}).
We finally set the boundary correction  term
\begin{align}
 \label{phi2}
\varepsilon\ \phi_2(x, z)
\, =\varepsilon\, \sum_{j=1}^N\phi_{2, j}(x, z)
=\varepsilon\, \sum_{j=1}^N\, \xi(\varepsilon z)\, \Big[\phi_{21, j}(x, z)+\phi_{22, j}(x, z)\Big],
\end{align}
where
\begin{equation*}
\xi(\theta)\, \equiv\, \frac{\chi_0(\theta)}{\mathcal{Q}(0)}
\, +\,
\frac{1-\chi_0(\theta)}{\mathcal{Q}(1)},
\end{equation*}
and the smooth cut-off function $\chi_0$ is defined by
$$
\chi_0(\theta)\, =\, 1  \quad\mbox{if } |\theta|<\frac 18,
\quad {\rm and}\quad
\chi_0(\theta)\, =\, 0  \quad\mbox{if } |\theta|>\frac 38.
$$
Note that $\varepsilon\, \phi_2(x, z)$ is  of size $O(\varepsilon)$ under the gap condition (\ref{gapconditionofve}).

\medskip
Let $v_2(x, z)\, =\, v_1+\varepsilon\, \phi_1+\varepsilon\, \phi_2$
 be the second approximate solution. A careful computation can indicate that
 the new boundary error takes the following form
\begin{align}
&\varepsilon\sum_{j=1}^N\, \Big[{\mathfrak b}_2\beta f_j\,-\, {\mathfrak b}_1\beta\, f_j'\Big]w_{j, x}
\nonumber\\[2mm]
&+\varepsilon^2 \Big[{\mathfrak b}_2+{\mathfrak b}_1 \frac{\beta'}{\beta} \Big]\sum_{j=1}^N\Big(\frac{x_j}{\beta}+f_j+h\Big)\Big[a_{10}\omega_{0, j, x}+a_{11}\omega_{1, j, x}+(f_j+h)(a_{12}\omega_{2, j, x}+a_{13}\omega_{3, j, x})
\nonumber\\[2mm]
&\qquad \qquad \qquad \qquad \qquad \qquad \qquad \qquad \qquad +\mathcal{Q}^{-1}\big(A(0)Z_{j, x}+\phi_{22, j, x}\big)\Big]
\nonumber\\[2mm]
&+\varepsilon^2\sum_{j=1}^N\, {\mathfrak b}_1
\Big[a'_{10}\, \omega_{0, j}
+a'_{11}\, \omega_{1, j}
+(a'_{12}\, \omega_{2, j}+a'_{13}\, \omega_{3, j})(f_j+h)
+(a_{12}\omega_{2, j}+a_{13}\omega_{3, j})(f_j'+h')\Big]
\nonumber\\[2mm]
&-\varepsilon^2\sum_{j=1}^N{\mathfrak b}_1
\big(\beta'f_j+\beta f_j'+\beta'h+\beta h'\big)\Big[a_{10}\omega_{0, j, x}+a_{11}\omega_{1, j, x}+(f_j+h)(a_{12}\omega_{2, j, x}+a_{13}\omega_{3, j, x})
\nonumber\\[2mm]
&\qquad \qquad \qquad \qquad \qquad \qquad \qquad\qquad+\mathcal{Q}^{-1}\big(A(0)Z_{j, x}+\phi_{22, j, x}\big)\Big]
\nonumber\\[2mm]
&-\varepsilon^2 \sum_{j=1}^N{\mathfrak b}_4\big(\beta'f_j+\beta f_j'+\beta'h+\beta h'\big) \Big(\frac{x_j}{\beta}+f_j+h\Big)w_{j, x}
+\varepsilon^2\sum_{j=1}^N {\mathfrak b}_4\frac{\alpha'}{\alpha}\Big(\frac{x_j}{\beta}+f_j+h\Big)w_j
\nonumber\\[2mm]
&+\varepsilon^2\sum_{j=1}^N\Big[{\mathfrak b}_3\Big(\frac{x_j}{\beta}+f_j+h\Big)^2\beta+{\mathfrak b}_4\Big(\frac{x_j}{\beta}+f_j+h\Big)\Big(\frac{\beta'}{\beta}x_j+\beta'f_j+\beta'h\Big)\Big]w_{j, x}
\nonumber\\[2mm]
&+\varepsilon^2 \sum_{j=1}^N {\mathfrak b}_1 \frac{\alpha'}{\alpha} \Big[a_{10}\omega_{0, j}+a_{11}\omega_{1, j}+(f_j+h)(a_{12}\omega_{2, j}+a_{13}\omega_{3, j})+\mathcal{Q}^{-1}(A(0)Z_{j}+\phi_{22, j})\Big]
\nonumber\\[2mm]
&+\varepsilon^2\sum_{j=1}^N {\mathfrak b}_4 \Big(\frac{x_j}{\beta}+f_j+h\Big)\mathcal{Q}^{-1}(\varepsilon A'(0)\mathcal{Q} Z_{j}+\phi_{22, j, z})
\nonumber\\[2mm]
&+D_2^0(v_1+\phi_1+\phi_2)+O(\varepsilon^3)\quad \quad\mbox{on }  \partial_0{\mathcal S},
\label{second-approximated-boundary-new-1}
 \end{align}
where the functions are evaluated at $\theta=0$.
Similar estimate holds on  $\partial_1{\mathcal S}$.

\subsection{The second improvement}
To deal with the resonance phenomena, which were described  in Section 1 of \cite{delPKowWei2007},
and improve the approximation for a solution still keeping the term of $\varepsilon^2$,
we need to introduce  new parameters $\{e_j\}_{j=1}^N$.
In other words, as the methods in \cite{delPKowWei2007} we shall set an improved approximate solution as follows
\begin{equation*}
 \, v_1+\varepsilon\, \phi_1+\varepsilon\, \phi_2+\varepsilon\, \sum_{j=1}^N e_j(\varepsilon z)Z_j(x).
\end{equation*}

To decompose the coupling of the parameters $\{f_j\}_{j=1}^N$ and $\{e_j\}_{j=1}^N$ on the boundary of $\mathcal S$ (in the sense of projection against $Z$ in $L^2$), by Lemma 2.2 in \cite{JWeiYang2007}, we
introduce a new term $\phi^{*}$ (even in $x$) defined by the following problem
\begin{align}\label{equationofphi*}
\begin{split}
\Delta \phi^{*}-\tilde{K}\phi^{*}+p\, {w}^{p-1}\phi^{*}\, =\, &0\quad\mbox{in } \hat{\mathcal S},
\\[2mm]	
\phi^{*}_{{\tilde z}}\, =\, H_{1}(x)\quad\mbox{on }  \partial_0\hat{\mathcal S},
\qquad
\phi^{*}_{{\tilde z}}\, =\, &H_{2}(x)\quad\mbox{on }  \partial_1\hat{\mathcal S},
\end{split}
\end{align}
where $\tilde{K}$ is a large positive constant, and the function
$H_{1}(x)$ is given by the following
\begin{align}\label{H1new}
H_1(x)\, =\,
&\Big[{\mathfrak b}_2+{\mathfrak b}_1 \frac{\beta'(0)}{\beta(0)} \Big]\big(f(0)+h(0)\big)\Big\{\big[a_{10}(0)\, \omega_{0, x}+a_{11}(0)\, \omega_{1, x}\big]
\nonumber\\[2mm]
&\, +\, x\big[a_{12}(0)\, \omega_{2, x}+a_{13}(0)\, \omega_{3, x}\big]
\, +\, x \frac{1}{\beta(0)} \mathcal{Q} \big[A(0)Z_x+\phi_{22, x}\big]\Big\}
\nonumber\\[2mm]
&\, +\, 2{\mathfrak b}_3\big[f(0)+h(0)\big]x\, w_x
\, +\, {\mathfrak b}_4\left[\frac{\beta'(0)}{\beta(0)}\big(f(0)+h(0)\big)+\Big(f'(0)+\frac{\beta'(0)}{\beta(0)}f(0)\Big)\right]x\, w_x
\nonumber\\[2mm]
&\, +\, {\mathfrak b}_1\Big[\beta'(0)f(0)+\beta(0)f'(0)+\beta'(0)h(0)+\beta(0)h'(0)\Big]\, \Big[a_{10}(0)\, \omega_{0, x}+a_{11}(0)\, \omega_{1, x}\Big]
\nonumber\\[2mm]
&\, +\, {\mathfrak b}_1\frac{\alpha'(0)}{\alpha(0)}\Big[(f(0)+h(0))\, (a_{12}(0)\, \omega_2+a_{13}(0)\, \omega_3)
+\frac{1}{\beta(0)}\Big(A(0)Z+\phi_{22}(x, 0)\Big)\Big]
\nonumber\\[2mm]
&+{\mathfrak b}_1\big[f'(0)+h'(0)\big]\Big[a_{12}(0)\omega_2+a_{13}(0)\omega_3\Big]
+{\mathfrak b}_1\big[f(0)+h(0)\big]\Big[a'_{12}(0)\omega_2+a'_{13}(0)\omega_3\Big]
\nonumber\\[2mm]
&+{\mathfrak b}_4\frac{1}{\beta(0)}\big[f(0)+h(0)\big]\big[\varepsilon  A'(0)\beta(0)Z+\phi_{22, z}(x, 0)\big]+{\mathfrak b}_3\big[f(0)+h(0)\big]\frac{\alpha'(0)}{\alpha(0)}w.
\end{align}
The function $H_2(x)$ has a similar expression.
We define a boundary correction term again
\begin{align}
\label{boundary layer}
\varepsilon^2\phi_3(x, z)\, =\, \varepsilon^2\sum_{j=1}^N\phi_{3, j}(x, z)\, =\, \varepsilon^2\sum_{j=1}^N\xi(\varepsilon z)\phi^{*}(x_j, {\Upsilon}(z)).
\end{align}
Note that for $j=1, \cdots, N$, $\varepsilon^2\phi_{3, j}$ is an exponential decaying function which is of order $\varepsilon^2$ and even in the variable $x_j$. Then define the third approximate solution to the problem near $\Gamma_\varepsilon$ as
\begin{equation}
\label{u3define}
v_3(x, z)\, =\, v_1+\varepsilon\, \phi_1+\varepsilon\, \phi_2+\varepsilon\, \sum_{j=1}^N e_j(\varepsilon z)Z_j(x)+\varepsilon^2\phi_3.
\end{equation}

\subsection{The third improvement}
By choosing $h$, we will construct a further approximation that eliminates  the even terms (in $x_j$'s) in the error $ S(v_3) $. This can be fulfilled by adding a term $\Phi=\sum_{j=1}^N\Phi_j$ and then considering the following term
\begin{align}
\label{v_1+phi1+phi2+epsilon sum_{j=1}^N{e_j}Z_j+phi3+Phi}
 S(v_3+\Phi)
\, =\, &S(v_1)+L_0(\varepsilon\, \phi_1)+L_0(\varepsilon\, \phi_2)+L_0 \Big( \varepsilon\sum_{j=1}^N{e_j}Z_j \Big)+L_0(\varepsilon^2\, \phi_3)+L_0(\Phi)
\nonumber\\
&+B_3(\varepsilon\, \phi_1)+B_3(\varepsilon\, \phi_2)+B_3\Big(\varepsilon\sum_{j=1}^N{e_j}Z_j\Big)
+B_3(\varepsilon^2\, \phi_3)+B_3(\Phi)
\nonumber\\
&+N_{0}\Big(\varepsilon\, \phi_1+\varepsilon\, \phi_2+\varepsilon\sum_{j=1}^N{e_j}Z_j+\varepsilon^2\, \phi_3+\Phi\Big),
\end{align}
where
\begin{equation}
L_0(\phi)= \frac{h_2}{\beta^2}\phi_{zz}+h_1\big[\phi_{xx}-\phi+p{v_1}^{p-1}\phi\big],
\qquad
N_0(\phi)= h_1\big[(v_1+\phi)^p- v_1^p-p v_1^{p-1} \phi\big].
\end{equation}
The details will be given in the sequel.

\subsubsection{Rearrangements of the error components}
The first objective of this part is to given the details
to compute the terms in formula (\ref{v_1+phi1+phi2+epsilon sum_{j=1}^N{e_j}Z_j+phi3+Phi}).

\noindent $\bullet$
It is easy to compute that
\begin{align}\label{L0phi1}
L_0(\varepsilon\, \phi_1)&\, =\, \varepsilon\, \Big[\frac{h_2}{\beta^2}\phi_{1, zz}+h_1\big(\phi_{1, xx}-\phi_1+p{v_1}^{p-1}\phi_1\big)\, \Big]
\nonumber\\
&\, =\, -\varepsilon\sum_{j=1}^N(S_{1, j}+S_{2, j})
+\varepsilon\, \frac{h_2}{\beta^2}\phi_{1, zz}
+\varepsilon h_1 \Big( \, p\, {v_1}^{p-1}{\phi_1}-\sum_{j=1}^N\, p\, {w_j}^{p-1}{\phi_{1, j}} \Big).
 \end{align}

\noindent $\bullet$
Recall the expression of $\varepsilon\, \phi_2$ and $A({\tilde\theta})$ defined in (\ref{phi2}) and (\ref{btheta}).
By using the equation of $\phi_{22, j}$ in (\ref{boundary-1}) and the equation of $Z$ in (\ref{eigenvalue}),
we get
\begin{align}
\label{L0phi2}
L_0(\varepsilon\, \phi_2)&\, =\, \varepsilon\, \Big[\frac{h_2}{\beta^2}\, \phi_{2, zz}+h_1\big(\phi_{2, xx}-\phi_2+p\, {v_1}^{p-1}\, \phi_2\big)\Big]
\nonumber\\[2mm]
&=\sum_{j=1}^N M_{11, j}(x, z)+M_{12}(x, z)+\varepsilon\, h_1\Big(p \, {v_1}^{p-1}{\phi_2}-\, \sum_{j=1}^N\, p\, {w_j}^{p-1}{\phi_{2, j}}\Big),
 \end{align}
where
 \begin{align}
 \label{hatL0phi2}
M_{11, j}(x, z)\, \equiv\, &\frac{\varepsilon^{2}}{\beta^{2}}h_2\Big\{2\xi'(\varepsilon z)\big[\varepsilon A'\big({\tilde{\bf d}}(\varepsilon z)\big)\mathcal{Q} Z_j+\phi_{22, j, z}\big]
+\mathcal{Q}'\xi(\varepsilon z)\big[\varepsilon A'\big({\tilde{\bf d}}(\varepsilon z)\big) Z_j+\phi_{*, j, \tilde{z}}\big]\Big\},
 \end{align}
 and
\begin{equation}\label{tildeL0phi2}
M_{12}(x, z)\, \equiv\, \sum_{j=1}^N\frac{\varepsilon^3}{\beta^{2}}\, h_2\, \xi''(\varepsilon z)\big[A'\big({\tilde{\bf d}}(\varepsilon z)\big)Z_j+\phi_{22, j}\big].
\end{equation}

\noindent $\bullet$
According to the equation of $Z$, we obtain
\begin{align*}
L_0 \Big( \varepsilon\sum_{j=1}^N{e_j}Z_j \Big)
&=\varepsilon\sum_{j=1}^N \Big\{ \frac{h_2}{\beta^2} ({e_j}Z_j)_{zz}+h_1\big[({e_j}Z_j)_{xx}-{e_j}Z_j+p{v_1}^{p-1}{e_j}Z_j \big] \Big\}
\\[2mm]
&=\varepsilon\sum_{j=1}^N \Big[\, \varepsilon^2 \frac{h_2}{\beta^2} e_j''Z_j+ h_1\lambda_0{e_j}Z_j\, \Big]
+\varepsilon\, h_1\, p{v_1}^{p-1}\sum_{j=1}^N {e_j}Z_j
\\[2mm]
&\quad  -\varepsilon\, h_1\, \sum_{j=1}^N p{w_j}^{p-1}{e_j}Z_j+\varepsilon^3 a_7(\varepsilon s, \varepsilon z).
\end{align*}

\noindent $\bullet$
Recalling the expression of $\varepsilon^2 \phi_3$ defined in \eqref{boundary layer} and the equation of $\phi^{*}$ in \eqref{equationofphi*},
it follows that
\begin{align}\label{L0phi3}
L_0(\varepsilon^2 \phi_3)
&=\, \varepsilon^2\, \Big[\frac{h_2}{\beta^2}\, \phi_{3, zz}+h_1\big(\phi_{3, xx}-\phi_3+p\, v_1^{p-1}\, \phi_3\big)\Big]
\nonumber\\[2mm]
&=\sum_{j=1}^N\, \varepsilon^2\, h_1(\tilde{K}-1)\phi_{3, j}
+\varepsilon^2\, \Big[\frac{h_2}{\beta^2}\phi_{3, j, zz}+h_1 \big(\phi_{3, j, xx}
-\tilde{K}\phi_{3, j}+p\, {v_1}^{p-1}\, \phi_{3, j} \big) \Big]
\nonumber\\[2mm]
&=\sum_{j=1}^N M_{21, j}(x, z)+M_{22}(x, z)
+ \varepsilon^2\, h_1 \Big( p\, {v_1}^{p-1}\, \phi_3-\sum_{j=1}^N\, p\, {w_j}^{p-1}\, \phi_{3, j}\Big),
\end{align}
where
\begin{equation}
\label{hatL0e^2hatphi}
M_{21, j}(x, z)\, \equiv\, \varepsilon^2\, h_1(\tilde{K}-1)\, \phi_{3, j},
\end{equation}
and
\begin{align}
\label{tildeL0e^2hatphi}
M_{22}(x, z)\, \equiv\, &\sum_{j=1}^N\frac{h_2}{\beta^2}\Big[\varepsilon^4\xi''(\varepsilon z)\phi^{*}(x_j, {\Upsilon}(z))
+2\varepsilon^3\xi'(\varepsilon z)\phi^{*}_{\tilde{z}}(x_j, {\Upsilon}(z))\mathcal{Q}(\varepsilon z)
\nonumber\\
&\qquad \qquad \quad+\varepsilon^3\xi(\varepsilon z)\phi^{*}_{\tilde{z}}(x_j, {\Upsilon}(z)) \mathcal{Q}'(\varepsilon z)\Big].
\end{align}

\noindent $\bullet$
Recall the expression of $B_3$ in (\ref{B3v}), we obtain that
\begin{align}
\label{Bphi1}
B_3(\varepsilon\, \phi_1)
\, =\, &\varepsilon^2 \frac{h_3}{\beta} \, \sum_{j=1}^N \phi_{1, j, x}
+ \varepsilon^2 h_8\, \sum_{j=1}^N \frac{x_j}{\beta}\phi_{1, j, xx}
+ \varepsilon^2 h_8\, \sum_{j=1}^N (f_j+h)\phi_{1, j, xx}
\nonumber\\[2mm]
&-\varepsilon^2\, \frac{ V_t(0, \varepsilon z)}{\beta^{2}}\, \sum_{j=1}^N \frac{x_j}{\beta}\phi_{1, j}
-\varepsilon^2\, \frac{ V_t(0, \varepsilon z)}{\beta^{2}}\, \sum_{j=1}^N (f_j+h)\phi_{1, j}
+\varepsilon^3a_8(\varepsilon s, \varepsilon z)
\nonumber\\[2mm]
\, =\, &\sum_{j=1}^N M_{31, j}(x, z)+M_{32}(x, z),
\end{align}
where
\begin{align}
\label{M31}
M_{31, j}(x, z)
\, \equiv\, &\varepsilon^2 \Bigg\{\frac{h_3}{\beta} \big(a_{10}\omega_{0, j, x} +a_{11}\omega_{1, j, x}+a_{12}h\omega_{2, j, x} +a_{13}h\omega_{3, j, x}\big)
\nonumber\\[2mm]
&\qquad+\frac{h_8}{\beta}\big( a_{10}x_j\omega_{0, j, xx} +a_{11}x_j\omega_{1, j, xx}+a_{12}hx_j\omega_{2, j, xx} +a_{13}hx_j\omega_{3, j, xx}\big)
\nonumber\\[2mm]
&\qquad+ h_8\big( a_{12}f_jh\omega_{2, j, xx} +a_{13}f_jh\omega_{3, j, xx}\big)
\nonumber\\[2mm]
&\qquad+ h_8\big[ ha_{10}\omega_{0, j, xx} +ha_{11}\omega_{1, j, xx}+a_{12}(f_jh+h^2)\omega_{2, j, xx} +a_{13}(f_jh+h^2)\omega_{3, j, xx}\big]
\nonumber\\[2mm]
&\qquad-\frac{ V_t(0, \varepsilon z)}{\beta^{2}}\frac{1}{\beta}\big( a_{10}x_j\omega_{0, j} +a_{11}x_j\omega_{1, j}+a_{12}hx_j\omega_{2, j} +a_{13}hx_j\omega_{3, j}\big)
\nonumber\\[2mm]
&\qquad-\frac{ V_t(0, \varepsilon z)}{\beta^{2}}\big[h a_{10}\omega_{0, j} +ha_{11}\omega_{1, j}+a_{12}(f_jh+h^2)\omega_{2, j} +a_{13}(f_jh+h^2)\omega_{3, j}\big]
\nonumber\\[2mm]
&\qquad-\frac{ V_t(0, \varepsilon z)}{\beta^{2}}\big(a_{12}f_jh\omega_{2, j} +a_{13}f_jh\omega_{3, j}\big) \Bigg\},
\end{align}
and
\begin{align}
\label{M32}
M_{32}(x, z)\, \equiv\, &\sum_{j=1}^N \, \varepsilon^2\Bigg\{ \frac{h_3}{\beta} \, f_j\big(a_{12}\omega_{2, j, x} +a_{13}\omega_{3, j, x}\big)
 +\, \frac{h_8}{\beta}f_j\big( a_{12}x_j\omega_{2, j, xx} +a_{13}x_j\omega_{3, j, xx}\big)
\nonumber\\[2mm]
&\qquad \quad+\, h_8\, \big( f_ja_{10}\omega_{0, j, xx} +f_ja_{11}\omega_{1, j, xx}+a_{12}f_j^2\omega_{2, j, xx} +a_{13}f_j^2\omega_{3, j, xx}\big)
\nonumber\\[2mm]
&\qquad \quad-\frac{V_t(0, \varepsilon z)}{\beta^3}f_j\big( a_{12}x_j\omega_{2, j} +a_{13}x_j\omega_{3, j}\big)
\nonumber\\[2mm]
&\qquad \quad- \frac{V_t(0, \varepsilon z)}{\beta^2}\, \big(f_j a_{10}\omega_{0, j} +f_ja_{11}\omega_{1, j}+a_{12}f_j^2\omega_{2, j} +a_{13}f_j^2\omega_{3, j}\big)\Bigg\}
\,+\,
\varepsilon^3a_8(\varepsilon s, \varepsilon z).
\end{align}

\noindent $\bullet$
Moreover, we can decompose $B_3(\varepsilon\, \phi_2)$ as follows
\begin{align*}
B_3(\varepsilon\, \phi_2)
=&\varepsilon^2 \frac{h_3}{\beta}  \, \sum_{j=1}^N \phi_{2, j, x}
+\varepsilon^2 \frac{ h_4}{\beta^2}\sum_{j=1}^N \phi_{2, j, z}
+2\varepsilon^2 h_2 \frac{\beta'}{\beta^{2}}\sum_{j=1}^N \Big(\frac{x_j}{\beta} + f_j +h\Big) \phi_{2, j, xz}
\\[2mm]
&+2\varepsilon^2 \frac{h_2\alpha'}{\alpha \beta^2}\sum_{j=1}^N \phi_{2, j, z}
+\varepsilon^2 \frac{h_6}{\beta}
\sum_{j=1}^N \Big(\frac{x_j}{\beta} + f_j +h\Big)\phi_{2, j, xz}
\\[2mm]
&+ \varepsilon^2 h_8\, \sum_{j=1}^N \Big( \frac{x_j}{\beta}+ f_j+h\Big) \phi_{2, j, xx}
-\varepsilon^2\, \frac{ V_t(0, \varepsilon z)}{\beta^{2}}\, \sum_{j=1}^N \Big( \frac{x_j}{\beta}+ f_j+h\Big)\phi_{2, j}
+\varepsilon^3a_9(\varepsilon s, \varepsilon z)
\\
\, =\, & \, \sum_{j=1}^N M_{41, j}(x, z)+M_{42}(x, z),
\end{align*}
where
\begin{align}\label{M41}
M_{41, j}(x, z)
\equiv\, & \varepsilon^2 \frac{\xi(\varepsilon z)}{\beta}
\Bigg\{h_3\big[A\big({\tilde{\bf d}}(\varepsilon z)\big)Z_{j, x}+\phi_{22, j, x}\big]+\frac{ h_4}{ \beta}\big[A'\big({\tilde{\bf d}}(\varepsilon z)\big)\varepsilon \mathcal{Q} Z_{j}+\phi_{22, j, z}\big]
\nonumber\\[2mm]
&\quad+2h_2 \frac{\beta'}{\beta}\Big(\frac{x_j}{\beta}+h\Big)\big[ A'\big({\tilde{\bf d}}(\varepsilon z)\big)\varepsilon \mathcal{Q} Z_{j, x}+\phi_{22, j, xz}\big]
+2\frac{h_2\alpha'}{\alpha \beta}\big[ A'\big({\tilde{\bf d}}(\varepsilon z)\big)\varepsilon \mathcal{Q} Z_{j}+\phi_{22, j, z}\big]
\nonumber\\[2mm]
&\quad+ h_6\Big(\frac{x_j}{\beta} +h\Big)\big[ A'\big({\tilde{\bf d}}(\varepsilon z)\big)\varepsilon \mathcal{Q}Z_{j, x}+\phi_{22, j, xz}\big]
+ \frac{h_8}{\beta}\Big(\frac{x_j}{\beta}+h\Big)\big[ A\big({\tilde{\bf d}}(\varepsilon z)\big)Z_{j, xx}+\phi_{22, j, xx}\big]
\nonumber\\[2mm]
&\quad-\, \frac{ V_t(0, \varepsilon z)}{\beta}\Big(\frac{x_j}{\beta}+h\Big)\big[ A\big({\tilde{\bf d}}(\varepsilon z)\big)Z_{j}+\phi_{22, j}\big]\Bigg\},
\end{align}
and
\begin{align}\label{M42}
M_{42}(x, z)\, \equiv\, &\sum_{j=1}^N \varepsilon^2 \xi(\varepsilon z)f_j\frac{1}{\beta}
\Bigg\{2h_2 \frac{\beta'}{\beta}\big[ A'\big({\tilde{\bf d}}(\varepsilon z)\big)\varepsilon \mathcal{Q} Z_{j, x}+\phi_{22, j, xz}\big]
+h_6\big[ A'\big({\tilde{\bf d}}(\varepsilon z)\big)\varepsilon \mathcal{Q}Z_{j, x}
+\phi_{22, j, xz}\big]
\nonumber\\[2mm]
&+\frac{h_8}{\beta}\, \big[ A\big({\tilde{\bf d}}(\varepsilon z)\big)Z_{j, xx}+\phi_{22, j, xx}\big]
-\, \frac{ V_t(0, \varepsilon z)}{\beta}\big[ A\big({\tilde{\bf d}}(\varepsilon z)\big)Z_{j}+\phi_{22, j}\big]\Bigg\}
+\varepsilon^3a_9(\varepsilon s, \varepsilon z).
\end{align}

\noindent $\bullet$
The computations for next term is the following:
\begin{align*}
B_3\Big(\varepsilon\sum_{j=1}^N{e_j}Z_j \Big)
\, =\, &\varepsilon^2 \frac{h_3}{\beta} \, \sum_{j=1}^N {e_j}Z_{j, x}
+ \varepsilon^2 h_8\, \sum_{j=1}^N \Big( \frac{x_j}{\beta}+f_j+h \Big) Z_{j, xx}
\\[2mm]
&-\varepsilon^2\, \frac{ V_t(0, \varepsilon z)}{\beta^{2}}\, \sum_{j=1}^N\Big( \frac{x_j}{\beta}+f_j+h \Big){e_j}Z_{j}
+\varepsilon^3a_{10}(\varepsilon s, \varepsilon z)
  \\[2mm]
\, =\, &\sum_{j=1}^N M_{51, j}(x, z)+M_{52}(x, z),
\end{align*}
where
\begin{equation}
\label{M51}
M_{51, j}(x, z) \, \equiv\, \varepsilon^2 (f_j+h ){e_j} \Big[\, h_8\, Z_{j, xx}-\, \frac{ V_t(0, \varepsilon z)}{\beta^{2}}Z_{j} \, \Big],
\end{equation}
\begin{equation}
\label{M52}
M_{52}(x, z)\, \equiv\, \sum_{j=1}^N\frac{\varepsilon^2}{\beta} {e_j}\Big[\,   h_3  \, Z_{j, x}
+h_8\, x_jZ_{j, xx}-\, \frac{ V_t(0, \varepsilon z)}{\beta^{2}}\, x_jZ_{j}\, \Big] +\varepsilon^3a_{10}(\varepsilon s, \varepsilon z).
\end{equation}

\noindent $\bullet$
The main order of the nonlinear term $N_0\Big(\varepsilon\, \phi_1+\varepsilon\, \phi_2+\varepsilon\sum_{j=1}^N{e_j}Z_j+\varepsilon^2\, \phi_3+\Phi\Big)$
is in the following
$$
\frac {p(p-1)}{2}{v_1}^{p-2}\Big(\varepsilon\, \phi_1+\varepsilon\, \phi_2+\varepsilon\sum_{j=1}^N{e_j}Z_j+\varepsilon^2\, \phi_3+\Phi\Big)^2,
$$
which can be decomposed in the form
\begin{align*}
&\quad\frac {p(p-1)}{2}h_1{v_1}^{p-2}\Big(\varepsilon\, \phi_1+\varepsilon\, \phi_2
+\varepsilon\sum_{j=1}^N{e_j}Z_j+\varepsilon^2\, \phi_3+\Phi\Big)^2
= M_{61}(x, z)\, +\, M_{62}(x, z)\, +\, O(\varepsilon^2),
\end{align*}
where
\begin{align}\label{M61}
& M_{61}(x, z)=\sum_{j=1}^N M_{61, j}(x, z)
\nonumber\\[2mm]
&\equiv\sum_{j=1}^N\varepsilon^2h_1\frac {p(p-1)}{2}{w_j}^{p-2}
\Big\{
a_{10}^2\omega_{0, j}^2+a_{11}^2\omega_{1, j}^2
+2a_{10}a_{11}\omega_{0, j}\omega_{1, j}
\nonumber\\[2mm]
&\qquad\qquad\qquad\qquad\quad\qquad+2h\big[a_{10}a_{12}\omega_{0, j}\omega_{2, j}
+a_{10}a_{13}\omega_{0, j}\omega_{3, j}\big]
\nonumber\\[2mm]
&\qquad\qquad\qquad\qquad\quad\qquad+2h\big[a_{11}a_{12}\omega_{1, j}\omega_{2, j}
+a_{11}a_{13}\omega_{1, j}\omega_{3, j}\big]
\nonumber\\[2mm]
&\qquad\qquad\qquad\qquad\quad\qquad+(2f_jh+h^2)\big[a_{12}^2\omega_{2, j}^2
+a_{13}^2\omega_{3, j}^2
+2a_{12}a_{13}\omega_{2, j}\omega_{3, j}\big]
\nonumber\\[2mm]
&\qquad\qquad\qquad\qquad\quad\qquad+e^2Z_j^2
\, +\,
\xi^2(\varepsilon z)\left[A\big({\tilde{\bf d}}(\varepsilon z)\big)Z_j+\phi_{22, j}\right]^2
\nonumber\\[2mm]
&\qquad\qquad\qquad\qquad\quad\qquad+2(a_{10}\omega_{0, j}+a_{11}\omega_{1, j})\xi(\varepsilon z)\left[A\big({\tilde{\bf d}}(\varepsilon z)\big)Z_j+\phi_{22, j}\right]
\nonumber\\[2mm]
&\qquad\qquad\qquad\qquad\quad\qquad+2 (a_{12}\omega_{2, j} + a_{13}\omega_{3, j})(f_j+h)   \left[ e_jZ_j+ \xi(\varepsilon z)\big(A\big({\tilde{\bf d}}(\varepsilon z)\big)Z_j+\phi_{22, j}\big)\right]
\nonumber\\[2mm]
&\qquad\qquad\qquad\qquad\quad\qquad +2\, \xi(\varepsilon z)\left[A\big({\tilde{\bf d}}(\varepsilon z)\big)Z_j+\phi_{22, j}\right]e_j Z_j
\Big\},
\end{align}
and
\begin{align}\label{M62}
&M_{62}(x, z)=\sum_{j=1}^N M_{62, j}(x, z)
\nonumber\\[2mm]
&\equiv\sum_{j=1}^N\frac {p(p-1)}{2}{w_j}^{p-2}
\Big\{
2\, \varepsilon^2f_j\big(a_{10}a_{12}\omega_{0, j}\omega_{2, j}
+a_{10}a_{13}\omega_{0, j}\omega_{3, j}
+a_{11}a_{12}\omega_{1, j}\omega_{2, j}
+a_{11}a_{13}\omega_{1, j}\omega_{3, j}\big)
\nonumber\\[2mm]
&\qquad\qquad\qquad\qquad\quad\quad+\varepsilon^2a_{12}^2\omega_{2, j}^2f_j^2
+\varepsilon^2a_{13}^2\omega_{3, j}^2f_j^2
+2\, \varepsilon^2\, a_{11}\omega_{1, j} e_jZ_j
\nonumber\\[2mm]
&\qquad\qquad\qquad\qquad\quad\quad+\varepsilon^4\, \phi_3^2(x_j, z)+\Phi^2_j
+2\, \varepsilon^2\, \phi_{1, j}\, \phi_3+2\, \phi_{1, j}\, \Phi_j+2\, \varepsilon^2\, \phi_2\, \phi_3+2\, \phi_{2, j}\, \Phi_j
\nonumber\\[2mm]
&\qquad\qquad\qquad\qquad\quad\quad  +2\, \varepsilon^3\, e_j\, Z_j\, \xi(\varepsilon z)\phi^{*}(x_j, z)+2\, \varepsilon \, e_j\, Z_j\, \Phi_j+2\, \varepsilon^2 \, \phi_3\, \Phi_j
\Big\}.
\end{align}
For the convenience of notation, we also denote
\begin{align}
M_{63}(x, z)=N_0\Big(\varepsilon\, \phi_1+\varepsilon\, \phi_2+\varepsilon\sum_{j=1}^N{e_j}Z_j+\varepsilon^2\, \phi_3+\Phi\Big)-M_{61}(x, z)-M_{62}(x, z).
\end{align}

Whence, according to the above rearrangements, we rewrite (\ref{v_1+phi1+phi2+epsilon sum_{j=1}^N{e_j}Z_j+phi3+Phi}) in terms of
\begin{align}\label{S(v_1+phi1+phi2+epsilon sum_{j=1}^N{e_j}Z_j+phi3+Phi+varphi)}
S(v_3+\Phi)
&=\sum_{j=1}^N
\Big[
\varepsilon^2S_{6, j}+\varepsilon^2S_{7, j}+\varepsilon^2S_{8, j}+\varepsilon^2S_{9, j}+ h_1\big(\Phi_{j, x_jx_j}-\Phi_j+p{w_j}^{p-1}\Phi_j\big)  +M_{11, j}(x, z)
\nonumber\\[2mm]
&\quad\qquad+M_{21, j}(x, z)+M_{31, j}(x, z)+M_{41, j}(x, z)+M_{51, j}(x, z)+M_{61, j}(x, z)
\Big]
\nonumber\\[2mm]
&\quad+\sum_{j=1}^N\varepsilon^2 \Big( S_{3, j}+S_{4, j}+S_{5, j} \Big)
+\sum_{j=1}^N \varepsilon\Big( \varepsilon^2 \frac{h_2}{\beta^2} e_j''Z_j+ h_1\lambda_0{e_j}Z_j\Big)
+B_4(v_1)+\varepsilon\frac{h_2}{\beta^2}\phi_{1, zz}
\nonumber\\
&\quad+\frac{h_2}{\beta^2}\Phi_{zz}+M_{12}(x, z)+M_{22}(x, z)+M_{32}(x, z)+M_{42}(x, z)+M_{52}(x, z)
\nonumber\\[2mm]
&\quad+M_{62}(x, z)+M_{63}(x, z)+B_{3}(\varepsilon^2 \phi_3)+B_3(\Phi)+\varepsilon h_1 \Big( p{v_1}^{p-1}{\phi_1}-\sum_{j=1}^N p {w_j}^{p-1}{\phi_{1, j}} \Big)
\nonumber\\
&\quad+\varepsilon h_1 \Big( p{v_1}^{p-1}\sum_{j=1}^N {e_j}Z_j -\sum_{j=1}^N p{w_j}^{p-1}{e_j}Z_j \Big)
+ \varepsilon h_1\Big( p {v_1}^{p-1}{\phi_2}-\sum_{j=1}^N p {w_j}^{p-1}{\phi_{2, j}} \Big)
\nonumber\\[2mm]
&\quad+\varepsilon^2 h_1 \Big(p{v_1}^{p-1}\phi_3-\sum_{j=1}^Np {w_j}^{p-1} \phi_{3, j} \Big)
+h_1\Big(p {v_1}^{p-1}\Phi-\sum_{j=1}^Np{w_j}^{p-1}\Phi \Big).
\end{align}

\subsubsection{Finding new correction terms and defining the basic approximation}

In order to eliminate the terms between the first brackets in (\ref{S(v_1+phi1+phi2+epsilon sum_{j=1}^N{e_j}Z_j+phi3+Phi+varphi)}),
for fixed $z$, we need a solution to the problem
\begin{align}
\label{equation of Phi}
h_1(-\Phi_{j, xx}+\Phi_j-p{w_j}^{p-1}\Phi_j)=&\, \varepsilon^2S_{6, j}+\varepsilon^2S_{7, j}+\varepsilon^2S_{8, j}+\varepsilon^2S_{9, j}
\nonumber\\[2mm]
&+M_{11, j}(x, z)+M_{21, j}(x, z)+M_{31, j}(x, z)
\nonumber\\[2mm]
&+M_{41, j}(x, z)+M_{51, j}(x, z)+M_{61, j}(x, z),
\quad
\forall x\in {\mathbb R}.
\end{align}
It is well-known that the above problem  is solvable provided that
\begin{align}
\label{solvable of varphi1}
\int_{{\mathbb R}}\big[&\varepsilon^2S_{6, j}+\varepsilon^2S_{7, j}+\varepsilon^2S_{8, j}+\varepsilon^2S_{9, j}+M_{11, j}(x, z)+M_{21, j}(x, z)
\nonumber\\
&+M_{31, j}(x, z)+M_{41, j}(x, z)+M_{51, j}(x, z)+M_{61, j}(x, z)\big]w_{j, x}\, {\rm d}x \, =\, 0.
\end{align}
The computations in Appendix \ref{appendixA} give that the validity of (\ref{solvable of varphi1}) holds
if the following problem
\begin{align}
\label{equation of h}
  \begin{split}
\mathcal{H}_{1}(\varepsilon z)h''
+\mathcal{H}_{1}'(\varepsilon z)h'
+\Big[ \mathcal{H}_{2}'(\varepsilon z)-\mathcal{H}_{3}(\varepsilon z) +\frac{\alpha_1(z)}{\zeta(\varepsilon z)}\Big]h
=\, \frac{1}{\zeta(\varepsilon z)}\big[G_{1}(z)+G_{3}(z)\big],
  \\[2mm]
{\mathfrak b}_1\, h'(0)\,-\, {\mathfrak b}_2 h(0)\, =\, 0,
\qquad
{\mathfrak b}_6\, h'(1)\,-\, {\mathfrak b}_7\, h(1)\, =\, 0,
  \end{split}
\end{align}
 has a solution.
Here the functions $\mathcal{H}_{1}$, $\mathcal{H}_{2}$, $\mathcal{H}_{3}$, $\zeta$, $\alpha_1$,  $G_1$ and $G_2$
are given in (\ref{mathcalH1})-(\ref{mathcalH3}), \eqref{zeta}, (\ref{alpha1}),  (\ref{G1}) and (\ref{G3}).
In fact, for the solvability of problem (\ref{equation of h}) the reader can refer to Lemma 6.1 in \cite{weixuyang}. Moreover, $h$ has the following estimate
$$
 \|h\|_{H^2(0, 1)}\leq C\varepsilon^{\frac{1}{2}}.
$$

\medskip
Now, we can find a function defined by $\Phi=\varepsilon^2\phi_4(x, \varepsilon z)$, such that the terms between the first brackets in (\ref{S(v_1+phi1+phi2+epsilon sum_{j=1}^N{e_j}Z_j+phi3+Phi+varphi)}) disappear. Finally, our basic approximate solution to the problem near the curve $\Gamma_\varepsilon$ is
\begin{align}
\label{basic approximate}
v_4\, =v_1+\varepsilon\phi_1+\varepsilon\phi_2+\varepsilon\sum_{j=1}^N{e_j}Z_j+\varepsilon^2\phi_3+\varepsilon^2\phi_4.
\end{align}

\subsection{The global approximate solution and errors}
Recall the coordinates $(s, z)$ in  (\ref{coordinatessz}), $(t, \theta)$ in (\ref{Fermicoordinates-modified}), and also the local approximate solution $v_{4}(s, z)$ in (\ref{basic approximate}), which is constructed near the curve $\Gamma_\varepsilon$ in the coordinates $(s, z)$.
By the relations in (\ref{vdefine}),
we then make an extension and simply define the approximate solution to (\ref{problemafterscaling}) in the form
\begin{align}
\label{globalapproximation}
\mathbf{w}(\tilde y)\, =\, \eta_{3\delta}^{\varepsilon}(s)\, \alpha(\varepsilon z)\,
v_4(x, z).
\end{align}
Note that, in the coordinates $(\tilde {y}_1, \tilde {y}_2)$ introduced in (\ref{rescaling}), $\mathbf{w}$ is a function defined on $\Omega_{\varepsilon}$ which is extended globally as $0$ beyond the $6\delta/\varepsilon$-neighborhood of $\Gamma_\varepsilon$.
The interior error can be arranged as follows
\begin{align}\label{new error-2}
{\mathcal E}\, \equiv\, S(v_4)
\, =&\, \varepsilon^2
\sum_{j=1}^N\big(\, S_{3, j}+\, S_{4, j}+\, S_{5, j}\big)
+B_4(v_1)
+\sum_{j=1}^N\Big[\, \varepsilon^3 \frac{h_2}{\beta^2} e_j''Z_j+  \varepsilon h_1\lambda_0{e_j}Z_j\, \Big]
\nonumber\\[2mm]
\, &+\varepsilon\frac{h_2}{\beta^2}\phi_{1, zz}
+\varepsilon^2\frac{h_2}{\beta^2}\phi_{4, zz}
+M_{12}(x, z)+M_{22}(x, z)+M_{32}(x, z)
\nonumber\\[2mm]
\, &+M_{42}(x, z)+M_{52}(x, z)+M_{62}(x, z)+M_{63}(x, z)
+B_3(\varepsilon^2\phi_3)+B_3(\varepsilon^2\phi_4)
\nonumber\\
&+\varepsilon h_1\Big(p\, {v_1}^{p-1}{\phi_1}-\sum_{j=1}^Np{w_j}^{p-1}{\phi_{1, j}}\Big)
+\varepsilon h_1 \Big( p\, {v_1}^{p-1}{\phi_2}-\sum_{j=1}^Np{w_j}^{p-1}{\phi_{2, j}} \Big)
\nonumber\\
&+\varepsilon^2 h_1 \Big(p\, {v_1}^{p-1}\phi_3-\varepsilon^2\sum_{j=1}^N\, p\, {w_j}^{p-1}\, \phi_{3, j}\Big)
+\varepsilon^2 h_1 \Big(p\, {v_1}^{p-1}\phi_4-\varepsilon^2\sum_{j=1}^N\, p\, {w_j}^{p-1}\, \phi_{4, j}\Big),
\end{align}
where we have used (\ref{S(v_1+phi1+phi2+epsilon sum_{j=1}^N{e_j}Z_j+phi3+Phi+varphi)}) and the equation of $\phi_4$ in (\ref{equation of Phi}).
The boundary error term $g_0$ has the form
\begin{align}\label{gonew}
g_0(x)\, =&\,
\varepsilon\sum_{j=1}^N\Big[{\mathfrak b}_2\beta f_j\,-\, {\mathfrak b}_1\beta\, f_j'\Big]w_{j, x}
\nonumber\\
&+\varepsilon^2 \Big[{\mathfrak b}_2+{\mathfrak b}_1 \frac{\beta'}{\beta}\Big]\sum_{j=1}^N(f_j+h)
\Big\{
a_{12}\omega_{2, j, x}+a_{13}\omega_{3, j, x}
 +\mathcal{Q}^{-1}\big(A(0)Z_{j, x}+\phi_{22, j, x}\big)
\Big\}
\nonumber\\
&+\varepsilon^2 \Big[{\mathfrak b}_2+{\mathfrak b}_1 \frac{\beta'}{\beta}\Big]\sum_{j=1}^N\Big(\frac{x_j}{\beta}\Big)\big[a_{10}\omega_{0, j, x}+a_{11}\omega_{1, j, x}
\big]
\nonumber\\
&-\varepsilon^2\sum_{j=1}^N{\mathfrak b}_1\big(\beta'f_j+\beta f_j'+\beta'h+\beta h'\big)
\Big\{ (a_{12}\omega_{2, j, x}+a_{13}\omega_{3, j, x})(f_j+h)
\nonumber\\
&\qquad\qquad+\, \frac{1}{\beta}\, \big(A(0)\, Z_{j, x}+\phi_{22, j, x}(x_j, 0)\big)
+ e_j Z_{j, x}\Big\}
\nonumber\\
&-\varepsilon^2\sum_{j=1}^N{\mathfrak b}_1\frac{\alpha'}{\alpha}
\Big\{ (a_{10}\omega_{0, j}+a_{11}\omega_{1, j})+ e_j Z_j\Big\}
+\varepsilon^2{\mathfrak b}_4\frac{\alpha'}{\alpha}\frac{x_j}{\beta}w_j
\nonumber\\
 &-\varepsilon^2\sum_{j=1}^N\Big\{{\mathfrak b}_3\beta\big(\frac{x_j^2}{\beta^2} +f_j^2+h^2+2f_jh\big)
+{\mathfrak b}_4\big[\frac{\beta'}{\beta^2}x_j^2+\beta(f_j+h)^2\big]\Big\}w_{j, x}
\nonumber\\
 &+\varepsilon^2{\mathfrak b}_1\sum_{j=1}^N\Big[(a'_{10}\omega_{0, j}+a'_{11}\omega_{1, j})+e_j'Z_j+{\mathfrak b}_5\frac{x_j}{\beta^2}
\big(\varepsilon A'(0)\beta(0)Z_j+\phi_{22, j, z}(x, 0)\big)\Big]
\nonumber\\
 &- \varepsilon^2\sum_{j=1}^N (\beta f_j' +\beta 'f_j) \Big[{\mathfrak b}_4(f_j+h)- \varepsilon {\mathfrak b}_5\Big(\frac{x_j}{\beta}+f_j+h\Big)^2 \Big] w_{j, x}
+O(\varepsilon^3).
 \end{align}
The term $g_1$ has a similar expression.

\medskip
We decompose
\begin{equation}
\label{E1-d}
{\mathcal E}\, =\, {\mathcal E}_{11}+{\mathcal E}_{12}, \qquad g_0\, =\, g_{01}+g_{02}, \qquad g_1\, =\, g_{11}+g_{12},
\end{equation}
with
\begin{align}
\label{E11}
{\mathcal E}_{11}\, =\, \sum_{j=1}^N {\mathcal E}_{11, j}
\, =\sum_{j=1}^N\big(\varepsilon^3\, \frac{h_2}{\beta^2} e''_j\, Z_j+\varepsilon\, h_1 \lambda_0\, e_j\, Z_j\big),
\quad &{\rm and}\quad
{\mathcal E}_{12}\, =\, {\mathcal E}-{\mathcal E}_{11},
\nonumber\\
g_{01}\, =\, \varepsilon\sum_{j=1}^N\beta(0)\big[{\mathfrak b}_2 f_j\,-\, {\mathfrak b}_1f_j'\big]w_{j, x}+\sum_{j=1}^N\varepsilon^2\, {\mathfrak b}_1\, e_j'\, Z_j,
\quad &{\rm and}\quad
g_{02}\, =\, g_0-g_{01},
\\
g_{11}\, =\, \varepsilon\sum_{j=1}^N\beta(0)\big[{\mathfrak b}_7f_j\,-\, {\mathfrak b}_6\, f_j'\big]w_{j, x}+\sum_{j=1}^N\varepsilon^2\, {\mathfrak b}_6\, e_j'\, Z_j,
\quad &{\rm and}\quad g_{12}\, =\, g_1-g_{11}.\nonumber
\end{align}
For further references, it is useful to estimate the $L^2(\mathcal S)$ norm of ${\mathcal E}$.
From the uniform bound of $e_1, \cdots, e_N$ in (\ref{enorm}), it is easy to see that
\begin{align}
\label{E11norm}
\|{\mathcal E}_{11}\|_{L^2(\mathcal S)}\, \leq\, C\varepsilon^{1/2}|\ln\varepsilon|^q.
\end{align}
Since $\varepsilon\, \phi_1, \varepsilon\, \phi_2$ and $\varepsilon e_jZ_j$ are of size $O(\varepsilon)$,  all terms in ${\mathcal E}_{12}$ carry $\varepsilon^2$ in front. We claim that
\begin{align}
\label{E12norm}
\|{\mathcal E}_{12}\|_{L^2(\mathcal S)}\, \leq\, C\varepsilon^{ 3/2}|\ln\varepsilon|^q.
\end{align}
Similarly, we have the following estimate
\begin{align}
\label{g02norm}
\|g_{02}\|_{L^2({\mathbb R})}+\|g_{12}\|_{L^2({\mathbb R})}\, \leq\, C\varepsilon^{3/2}|\ln\varepsilon|^q.
\end{align}
Moreover, for the Lipschitz dependence of the term of error ${\mathcal E}_{12}$ on the parameters $\mathbf{f}$ and $\mathbf{e}$ for the norm defined in (\ref{constraints of f}) and (\ref{enorm}), we have the validity of the estimate
\begin{align}
\label{E12L}
\|{\mathcal E}_{12}(\mathbf{f}_1, \mathbf{e}_1)-{\mathcal E}_{12}(\mathbf{f}_2, \mathbf{e}_2)\|_{L^2(\mathcal S)}
\, \leq\, C\varepsilon^{ 3/2}|\ln\varepsilon|^q\big[\, \|\mathbf{f}_1-\mathbf{f}_2\|_{H^2(0, 1)}+\|\mathbf{e}_1-\mathbf{e}_2\|_{**}\, \big].
\end{align}
Similarly, we obtain
\begin{align}
\label{g02L}
\begin{split}
&\|g_{02}(\mathbf{f}_1, \mathbf{e}_1)-g_{02}(\mathbf{f}_2, \mathbf{e}_2)\|_{L^2({\mathbb R})}\, +\, \|g_{12}\, (\mathbf{f}_1, \mathbf{e}_1)-g_{12}(\mathbf{f}_2, \mathbf{e}_2)\|_{L^2({\mathbb R})}
\\[2mm]
&\, \leq\, C\varepsilon^{3/2}|\ln\varepsilon|^q\big[\, \|\mathbf{f}_1-\mathbf{f}_2\|_{H^2(0, 1)}+\|\mathbf{e}_1-\mathbf{e}_2\|_{**}\, \big].
\end{split}
\end{align}

 \section{Derivation of the reduced equations: Toda system}\label{section5}
\setcounter{equation}{0}

In this section, we will set up equations for the parameters $\mathbf{f}$ and $\mathbf{e}$ which are equivalent to making $\mathbf{c}(\varepsilon z)$, $\mathbf{d}(\varepsilon z)$, $\mathbf{l}_1$, $\mathbf{l}_0$, $\mathbf{m}_1$ and $\mathbf{m}_0$ are identically zero in the system (\ref{system-1})-(\ref{system-4}).
 The equations
 $$
\mathbf{c}(\varepsilon z)=0, \quad \mathbf{l}_1=0, \quad \mathbf{l}_0=0,
 $$
 are then equivalent to the relations, for $ n=1, \cdots, N$,
 \begin{align}
 \int_{{\mathbb R}}\Big[\, \eta^\varepsilon_{\delta}(s)\, {\mathcal E}+\eta^\varepsilon_{\delta}(s)\, {\mathcal N}(\phi)+\chi(\varepsilon|x|) B_3(\phi)+h_1p\big(\beta^{-2}\chi(\varepsilon |x|){\mathbf w}^{p-1}&-v_1^{p-1}\big)\phi \, \Big]\, w_{n, x}\, {\rm d}x \, =\, 0,
 \label{c=0}
 \\[2mm]
 \int_{{\mathbb R}}\Big[\, \eta^\varepsilon_{\delta}(s) g_1-\chi(\varepsilon|x|)D_3^1(\phi)+\chi(\varepsilon|x|) D_2^1(\phi)\, \Big]\, w_{n, x}\, {\rm d}x\, &=\, 0,
  \quad z=1/\varepsilon, \label{l1=0}
  \\[2mm]
  \int_{{\mathbb R}}\Big[\, \eta^\varepsilon_{\delta}(s) g_0-\chi(\varepsilon|x|) D_3^0(\phi)+\chi(\varepsilon|x|) D_2^0(\phi)\, \Big]\, w_{n, x}\, {\rm d}x\, &=\, 0, \quad z=0. \label{l2=0}
 \end{align}
Similarly,
$$
\mathbf{d}(\varepsilon z)=0,
\quad
\mathbf{m}_1=0, \quad
\mathbf{m}_0=0,
$$
if and only if for $n=1, \cdots, N$,
\begin{align}
   \int_{{\mathbb R}}\Big[\, \eta^\varepsilon_{\delta}(s)\, {\mathcal E}+\eta^\varepsilon_{\delta}(s)\, {\mathcal N}(\phi)+\chi(\varepsilon|x|) B_3(\phi)+h_1p\big(\beta^{-2}\chi(\varepsilon |x|){\mathbf w}^{p-1}&-v_1^{p-1}\big)\phi \, \Big]\, Z_n{\rm d}x\, =\, 0, \label{d=0}
    \\[2mm]
    \int_{{\mathbb R}}\Big[\, \eta^\varepsilon_{\delta}(s) g_1-\chi(\varepsilon|x|)D_3^1(\phi)+\chi(\varepsilon|x|) D_2^1(\phi)\, \Big]\, Z_n\, {\rm d}x\, &=\, 0,
  \quad z=1/\varepsilon, \label{m1=0}
  \\[2mm]
  \int_{{\mathbb R}}\Big[\, \eta^\varepsilon_{\delta}(s) g_0-\chi(\varepsilon|x|) D_3^0(\phi)+\chi(\varepsilon|x|) D_2^0(\phi)\, \Big]\, Z_n\, {\rm d}x\, &=\, 0,
   \quad z=0.\label{m0=0}
\end{align}

\subsection{Estimates for projections of the error}
 For the pair $(\mathbf{f}, \mathbf{e})$ satisfying (\ref{constraints of f}) and (\ref{enorm}),
  we denote  ${\mathbf b}_{1\varepsilon}$ and ${\mathbf b}_{2\varepsilon}$,  generic, uniformly bounded continuous functions of the form
$$
{\mathbf b}_{l\varepsilon j}
    \, =\, {\mathbf b}_{l \varepsilon j}\big(z, {\mathbf f}(\varepsilon z), {\mathbf  e}(\varepsilon z), {\mathbf f}'(\varepsilon z), \varepsilon\, {\mathbf e}'(\varepsilon z)\big), \ l=1, 2,
$$
where ${\mathbf b}_{1\varepsilon j}$ is uniformly Lipschitz in its four last arguments,
and introduce the notation
\begin{align}
{\mathfrak S}=\{x\in{\mathbb R}:(x, z)\in \mathcal S \},
\qquad
{\mathfrak S}_n=\{x\in{\mathbb R}:(x, z)\in {\mathfrak A}_n\}
\quad
n=1, \cdots, N.
\end{align}

\medskip
The computations in Appendix \ref{appendixB} lead to  the estimate, for $n=1, \cdots, N$,
\begin{align}
\int_{\mathfrak S}\, {\mathcal E}w_{n, x}{\rm d}x
=\, &-\varrho_2h_1\Bigg\{\varepsilon^2\varsigma\Big[\, \mathcal{H}_1f''_n+\mathcal{H}_1'f'_n+ \big(\mathcal{H}_2'-\mathcal{H}_3+\alpha_2(z)\big)f_n\, \Big]
\\
&\qquad \qquad+e^{-\beta(f_n-f_{n-1})}-e^{-\beta(f_{n+1}-f_n)} \Bigg\}
\, +\, \mathrm{P}_n(\varepsilon z),
\end{align}
where
\begin{equation}\label{varsigma}
\gamma_1(\theta)=\frac{\varrho_1}{\beta(\theta)},
\qquad
\varsigma(\theta)= \frac{\gamma_1(\theta)}{\varrho_2\, h_1(\theta)},
\qquad
\alpha_2(z)=\frac{\alpha_1(z)}{\zeta(\varepsilon z)},
\end{equation}
\begin{align*}
\mathrm{P}_n(\varepsilon z)=&\, \varepsilon^2\, \gamma_1(\theta)\, \big[\hbar_3(\varepsilon z)\, e_n+\varepsilon^2\, \hbar_4({\varepsilon}z)\, e''_n\big]
+\varepsilon^{\mu_1}\max_{j\neq n}O(e^{-\beta|f_j-f_n|})
\\
&+\varepsilon^3\sum_{j=1}^N\, \big({\mathbf b}_{1\varepsilon j}\, e_j'+{\mathbf b}_{1\varepsilon j}^2\, f''_j+{\mathbf b}_{2\varepsilon j}\big)
+O(\varepsilon^3)\sum_{j=1}^N\, (f_j+f'_j+f''_j+f_j^2),
\end{align*}
where the functions $\mathcal{H}_1$, $\mathcal{H}_2$, $\mathcal{H}_3$ are defined in
\eqref{mathcalH1}-\eqref{mathcalH3}, and $\zeta$ in \eqref{zeta}, $\alpha_1(z)$ in (\ref{alpha1}), $\varrho_1$ and $\varrho_2$ in (\ref{varrho1})-(\ref{varrho2}),
$\hbar_3({\varepsilon}z)$ and $\hbar_4({\varepsilon}z)$ in (\ref{hbar3})-(\ref{hbar4}), $h_1$ in \eqref{h1}.
For further references we observe that
\begin{equation}
\|\mathrm{P}_n\|_{L^2(0, 1)}\leq\, C\varepsilon^{2+\mu_2}, \quad\mbox{for some }\, \mu_2>0, \, n=1, \cdots, N.
\end{equation}

\medskip
On the other hand, the computations in Appendix \ref{appendixC} lead to  the estimate, for $n=1, \cdots, N$,
\begin{align}
\label{E1Z}
\int_{{\mathfrak S}_n}{\mathcal E}Z_n\, {\rm d}x
\, =\, &\varepsilon^3\frac{h_2}{\beta^2}e''_n+\varepsilon\lambda_0h_1e_n+\varepsilon^3\hbar_5e_n'+\varrho_4h_1\big[e^{-\beta(f_n-f_{n-1})}-e^{-\beta(f_{n+1}-f_n)}\big]
\nonumber\\[2mm]
&+\varepsilon^2\varrho_{3}\, (f'_n)^2+\varepsilon^3\rho_1\, (\varepsilon z)+\varepsilon^2\rho_2(\varepsilon z)+\varepsilon^2\rho_3(\varepsilon z)+\mathrm{R}_n(\varepsilon z),
\end{align}
where
\begin{align*}
\mathrm{R}_n(\varepsilon z)=\, &\varepsilon^4 f_n\beta^{-2}\, \hbar_6(\varepsilon z)e''_n
+2\varepsilon^2\varrho_{3}(f'_nh'+f'_nh)
+\varepsilon^3\beta^{-2}\xi''(\varepsilon z)A'\big({\tilde{\bf d}}(\varepsilon z)\big)
\\[2mm]
&+O(\varepsilon^3)\sum_{j=1}^N\big(f_j^2+{f'_j}^2\big)+O(\varepsilon^5)\sum_{j=1}^N\, e''_j
+\varepsilon^{{\hat\tau}_1}\max_{j\neq n}O(e^{-\beta|f_j-f_n|})
\\
&+\varepsilon^3\sum_{j=1}^N\big({\mathbf b}_{1\varepsilon j}f''_j+{\mathbf b}_{2\varepsilon j}e'_j+{\mathbf b}_{2\varepsilon j}\big)+O(\varepsilon^3).
\end{align*}
Here the constants $\lambda_0$, $\varrho_3$, $\varrho_4$ are given \eqref{lambda0}, (\ref{varrho3}) and (\ref{varrho4}),
while the functions $h_1$, $h_2$, $\hbar_5$, $\hbar_6$, $\rho_1$, $\rho_2$  and $\rho_3$ are given
in \eqref{h1}, \eqref{h2}, (\ref{hbar5}), (\ref{hbar6}), (\ref{rho1}), (\ref{rho2}) and \eqref{rho3}.
 For further references we observe that
\begin{align}
\|\mathrm{R}_n\|_{L^2(0, 1)}\leq\, C\varepsilon^{2+\mu_3}, \quad\mbox{for some }\mu_3>0, \, n=1, \cdots, N.
\end{align}

\subsection{Projection of errors on the boundary}
In this section, we compute the projection of errors on the boundary.
Without loss of generality,  only the projections of the error components on $\partial_0{\mathcal S}$ will be given.
According to the expression of $g_0$ as in \eqref{gonew}, the main errors on the boundary integrated against $w_{n, x}$ and $Z_n$ in the variable $x_n$ can be computed as the following:
\begin{align*}
\int_{\mathbb R}& g_0(x)w_{n, x}\, {\rm d}x\,
=\, \varepsilon\sum_{j=1}^N \beta(0)\Big[{\mathfrak b}_2f_j\,-\, {\mathfrak b}_1f_j'\Big]\int_{\mathbb R} w_{j, x}w_{n, x}\, {\rm d}x+\varepsilon^2{\mathcal M}_{0, n}^1(\mathbf{f}, \mathbf{e})
+O(\varepsilon^3).
\end{align*}
Using the following formulas
\begin{equation*}
\int_{\mathbb R} Z^2\, {\rm d}x\, =\, -2\int_{\mathbb R} x\, Z_x\, Z\, {\rm d}x\, =\, 1,
\end{equation*}
we get the following two estimates
\begin{equation*}
\begin{split}
\int_{\mathbb R} g_0(x)Z_n\, {\rm d}x\, =\, &\varepsilon^2 {\mathfrak b}_1\left[e_n'(0) \, +\, \frac{\alpha'(0)}{\alpha(0)}\, e_n(0)  \right]
+O(\varepsilon^3).
\end{split}
\end{equation*}

\medskip
Higher order errors can be proceeded as follows:
\begin{align*}
\int_{\mathbb R}& D_3^0(\phi(x, 0))w_{n, x}\, {\rm d}x\,
\nonumber
\\
=
&\, \sum_{j=1}^N\varepsilon\Big[{\mathfrak b}_2+{\mathfrak b}_1\frac{\beta'}{\beta}\Big]\int_{\mathbb R} x\phi_{x}(x, 0)w_{n, x}\, {\rm d}x
-\sum_{j=1}^N\varepsilon\frac{\alpha'}{\alpha}\int_{\mathbb R} \phi(x, 0)w_{n, x}\, {\rm d}x
 \\
&\, +\, \sum_{j=1}^N\varepsilon {\mathfrak b}_4\int_{\mathbb R}\Big(\frac{x_j}{\beta}+f_j+h\Big)\phi_{z}(x, 0)w_{n, x}\, {\rm d}x
\\
&+\sum_{j=1}^N\varepsilon^2\int_{\mathbb R} \Big[{\mathfrak b}_1\Big(\frac{x_j}{\beta}+f_j+h\Big)^2\beta+{\mathfrak b}_4\Big(\frac{x_j}{\beta}
+f_j+h\Big) \Big(\frac{\beta'}{\beta}x_j+\beta' f_j+\beta' h\Big)\Big]\phi_{x}(x, 0)w_{n, x}\, {\rm d}x
\nonumber
\\
&+\sum_{j=1}^N\varepsilon^2{\mathfrak b}_4\frac{\alpha'}{\alpha}\int_{\mathbb R}\Big(\frac{x_j}{\beta}+f_j+h\Big)\phi(x, 0)w_{n, x}\, {\rm d}x
\, +\, \sum_{j=1}^N\varepsilon^2 {\mathfrak b}_5\int_{\mathbb R}  \Big(\frac{x_j}{\beta}+f_j+h\Big)^2\phi_{z}(x, 0)w_{n, x}\, {\rm d}x
\nonumber
\\
=&\, O(\varepsilon^{2+\mu}),
\end{align*}
and also
\begin{equation*}
 \int_{\mathbb R}  D_3^0(\phi(x, 0))Z_n\, {\rm d}x\, =  O(\varepsilon^{2+\mu})
\end{equation*}
The term $D_2^0(\phi)$ on the boundary integrated against $w_{n, x}$ and $Z_n$ in the variable $x_n$ are of size of order $O(\varepsilon^3)$.

\subsection{The system involving $(\textbf{f}, \textbf{e})$}
\label{The system for fe}

As done in \cite{delPKowWei2007} and \cite{delPKowWei2008}, we can estimate the terms that involve $\phi$ in (\ref{c=0}) and (\ref{d=0}) integrated against the functions $w_{n, x}$ and $Z_n$ in the variable $x_n$ in the similar ways.
As a conclusion, for $n=1, \cdots, N$, there holds the following estimate
\begin{align}
-\varepsilon^2\varsigma\Big[\, \mathcal{H}_1\, f''_n+\mathcal{H}_1'\, f'_n+ \big(\mathcal{H}_2'-\mathcal{H}_3+\alpha_2(z)\big)\, f_n\, \Big] +e^{-\beta(f_n-f_{n-1})}-e^{-\beta(f_{n+1}-f_n)} \, +\, {\mathbb M}_n\, =\, 0.
\end{align}
Moreover, ${\mathbb M}_n$ can be decomposed in the following way
\begin{align}
\begin{aligned}
{\mathbb M}_n={\mathbb M}_{n1}(\theta, \mathbf{f, f', f'', e, e', e''})+{\mathbb M}_{n2}(\theta, \mathbf{f, f', e, e'}),
\end{aligned}
\end{align}
where ${\mathbb M}_{n1}$ and ${\mathbb M}_{n2}$ are continuous of their arguments.
 Functions ${\mathbb M}_{n1}$ and ${\mathbb M}_{n2}$ satisfy the following properties for $n=1, \cdots, N$
$$
\|{\mathbb M}_{n1}\|_{L^2(0, 1)}\, \leq\, C \varepsilon^{2+\mu_0},
\qquad
\|{\mathbb M}_{n2}\|_{L^2(0, 1)}\, \leq\, C \varepsilon^{2+\mu_0},
$$
where $\mu_0$ is a positive constant. For $n=1, \cdots, N$, there also holds the following estimate
\begin{align}
\begin{aligned}
&\varepsilon^3\frac{h_2}{\beta^2}e''_n
+\varepsilon h_1\lambda_0e_n
+\varepsilon^3\hbar_5e_n'
+h_1\varrho_4\big[e^{-\beta(f_n-f_{n-1})}-e^{-\beta(f_{n+1}-f_n)}\big]
\\[2mm]
&+\varepsilon^2\varrho_{3}\, (f'_n)^2
+\varepsilon^3\rho_1
+\varepsilon^2\rho_2
+\varepsilon^2\rho_3
+{\mathbf M}_n\, =\, 0.
\end{aligned}
\end{align}
Moreover, ${\mathbf M}_n$ can be decomposed in the following way
\begin{align}
\begin{aligned}
{\mathbf M}_n={\mathbf M}_{n1}(\theta, \mathbf{f, f', f'', e, e', e''})+{\mathbf M}_{n2}(\theta, \mathbf{f, f', e, e'}),
\end{aligned}
\end{align}
where ${\mathbf M}_{n1}$ and ${\mathbf M}_{n2}$ are continuous of their arguments.
Functions ${\mathbf M}_{n1}$ and ${\mathbf M}_{n2}$ satisfy the following properties for $n=1, \cdots, N$
$$
\|{\mathbf M}_{n1}\|_{L^2(0, 1)}\, \leq\, C \varepsilon^{2+{\hat\tau}_0},
\qquad
\|{\mathbf M}_{n2}\|_{L^2(0, 1)}\, \leq\, C \varepsilon^{2+{\hat\tau}_0},
$$
where ${\hat\tau}_0$ is a positive constant.

\medskip
Therefore, using $\theta=\varepsilon z$  and defining the operators
\begin{align}
{\mathbb L}_{n}(\mathbf{f})
\, \equiv\, &\varepsilon^2\varsigma(\theta)\Big[\, \mathcal{H}_1(\theta)\, f''_n+\mathcal{H}_1'(\theta)\, f'_n+ \big(\mathcal{H}_2'(\theta)-\mathcal{H}_3(\theta)+\alpha_2(\theta/\varepsilon)\big)\, f_n\, \Big]\nonumber
\\[2mm]
&-e^{-\beta(\theta) (f_n-f_{n-1})}
+e^{-\beta(\theta)(f_{n+1}-f_n)},
\quad
n=1, \cdots, N,
\label{mathbbLn}
\end{align}
and
\begin{equation}
{\mathbb L}(e)\, \equiv\,-\, \varepsilon^2\, h_2(\theta)e''\,-\, \varepsilon^2\tilde{\alpha}(\theta)\, e'\,-\, |\beta(\theta)|^2h_1(\theta)\lambda_0 e,
\label{mathbbL}
\end{equation}
we derive the following nonlinear system of differential equations for the parameters $\mathbf{f}$ and  $\mathbf{e}$
\begin{equation}\label{f}
{\mathbb L}_{n}(\mathbf{f})={\mathbb M}_n, \quad n=1, \cdots, N,
\end{equation}
\begin{equation}
{\mathbb L}(e_n)=\, \alpha_5+\alpha_{6, n}+\varepsilon^{-1}\, \beta^2\, {\mathbf M}_n, \quad n=1, \cdots, N,
 \label{e}
\end{equation}
with the boundary conditions, $n=1, \cdots, N$,
\begin{equation}
\label{boundary condition 1}
 {\mathfrak b}_6f'_n(1)\,-\, {\mathfrak b}_7\, f_n(1)+{\mathcal M}_{1, n}^1(\mathbf{f}, \mathbf{e})=\, 0,
 \qquad
 {\mathfrak b}_1f'_n(0)\,-\, {\mathfrak b}_2\, f_n(0)\, +{\mathcal M}_{1, n}^2(\mathbf{f}, \mathbf{e})=\, 0,
\end{equation}
\begin{equation}
\label{boundary condition 3}
 e'_n(1)\, +\, \tilde{b}_6\, e_n(1)\, +{\mathcal M}_{2, n}^1(\mathbf{f}, \mathbf{e})=\, 0,
 \qquad
 e'_n(0)\, +\, \tilde{b}_5\, e_n(0)\, +{\mathcal M}_{2, n}^2(\mathbf{f}, \mathbf{e})=\, 0.
\end{equation}
The constants ${\tilde b}_6$ and ${\tilde b}_5$ are given by
$$\tilde{b}_6\, =\, \frac{\alpha'(1)}{\alpha(1)},
 \qquad
\tilde{b}_5\, =\, \frac{\alpha'(0)}{\alpha(0)},
 $$
where ${\mathcal M}_{j, n}^{i}$'s are some terms of order $O(\varepsilon^{1/2})$.
The functions $\varsigma(\theta)\, =\, \frac{\gamma_1(\theta)}{\varrho_2 h_1(\theta) }>0$,
and $\alpha_2(\theta/\varepsilon)$ are defined in (\ref{varsigma}) by the relation $\theta=\varepsilon z$.
Moreover, we have denoted
\begin{equation}
\tilde{\alpha}(\theta)\, =\, \beta^2(\theta)\, \hbar_5(\theta),
\end{equation}

\begin{equation}
\alpha_5(\theta)\, =\, \beta^2(\theta) \big[ \varepsilon^2\, \rho_1(\theta)\, +\, \varepsilon\rho_2\, (\theta)\, +\, \varepsilon\, \rho_3(\theta) \big],
\label{alpha5}
\end{equation}

\begin{equation}
\alpha_{6, n}(\theta)
\, =\, \varepsilon\, \varrho_{3}\, \beta^2(\theta)\, |f'_n(\theta)|^2
\, +\,
\varepsilon^{-1}\varrho_{4}\, \beta^{2}(\theta)h_1(\theta) \big[e^{-\beta(\theta)(f_n-f_{n-1})}-e^{-\beta(\theta)(f_{n+1}-f_n)}\big],
\end{equation}
where $\hbar_5(\theta)$, $\rho_1(\theta)$, $\rho_2(\theta)$, $\rho_3(\theta)$ are defined in (\ref{hbar5}), (\ref{rho1}), (\ref{rho2}), (\ref{rho3}).

\section{Suitable choosing of parameters}\label{sectionsolvingreducedequation}
\setcounter{equation}{0}

\subsection{Solving the system of reduced equations}
Before solving (\ref{f})-(\ref{boundary condition 3}), some basic facts about the invertibility of corresponding operators will be derived.

\medskip
Firstly, we consider the following problem
\begin{align}\label{equation of e}
\begin{split}
{\mathbb L}(e)=\, \tilde{g}(\theta), \qquad\forall\, 0<\theta<1,
\\
e'(1)\, +\, \tilde{b}_6\, e(1)\, =\, 0,
\qquad\quad
e'(0)\, +\, \tilde{b}_5\, e(0)\, =\, 0.
\end{split}
\end{align}

\begin{proposition}\label{proposition7point1}
If $\tilde{g}\in L^{2}(0, 1)$, then for all small $\varepsilon$ satisfying (\ref{gapconditionofve})
there is a unique solution $e\in H^2(0, 1)$ to problem (\ref{equation of e}), which satisfies
 $$\|e\|_{**}\leq C\, \varepsilon^{-1}\, \|\tilde{g}\|_{L^2(0, 1)}.$$
 Moreover, if $\tilde{g}\in H^{2}(0, 1)$, then
 \begin{align}
\varepsilon^2\, \|e''\|_{L^2(0, 1)}\, +\, \varepsilon\, \|e'\|_{L^2(0, 1)}\, +\, \|e\|_{L^\infty(0, 1)}\, \leq\, C\, \|\tilde{g}\|_{H^2(0, 1)}.
 \end{align}
 \end{proposition}
\begin{proof}
The proof is similar as that for Lemma 8.1 in \cite{delPKowWei2007}.
\end{proof}

Secondly, we consider the following problem
\begin{equation}
\label{ts2a}
{\mathbb L}_{n}({\mathbf f})=\varepsilon^2 \tilde{h}_n,
\end{equation}
\begin{equation}
\label{boundary0000}
{\mathfrak b}_6f'_n(1)-{\mathfrak b}_7\, f_n(1)=0,
\qquad
{\mathfrak b}_1f'_n(0)-{\mathfrak b}_2\, f_n(0)=0,
\end{equation}
for $n=1, \cdots, N$, where $f_0=-\infty, \, f_{N+1}=\infty$.

\begin{proposition}\label{proposition7point2}
For given ${\tilde{\mathbf h}}=(\tilde{h}_1, \cdots, \tilde{h}_N)^{T}\in L^2(0, 1)$,  there exists a sequence $\{\varepsilon_l:l\in{\mathbb N}\}$
from those $\varepsilon$ satisfying the gap condition (\ref{gapconditionofve}) and
approaching $0$ such that problem (\ref{ts2a})-(\ref{boundary0000}) admits a
solution ${\bf f}=(f_1, \cdots, f_N)^{T}$ with the form:
\begin{align}
{\bf f}\, =\,
\frac{1}{\beta}\Bigg\{
\rho_{\varepsilon_l}\Big(1-\frac{N+1}{2},\, 2-\frac{N+1}{2},\, \cdots,\, N-\frac{N+1}{2} \Big)^{T}
\,+\,\bf{\ddot f}
\,+\,\mathbf{P}^T{\hat{\mathfrak u}}
\,+\,\mathbf{P}^T{\bf {\tilde w}}
\,+\,\mathbf{P}^T{\check{\mathfrak u}}
\Bigg\},
\end{align}
where the invertible matrix $\bf P$ is defined in \eqref{mathbfP} and the function $\rho_{\varepsilon_l} (\theta) $ satisfies
\begin{align}\label{def rhoeps}
e^{-\rho_{\varepsilon_l}(\theta)}\, =\, \varepsilon^2 \frac{{\varsigma}(\theta)}{\beta(\theta)}\tau_2 (\theta) \rho_{\varepsilon_l}(\theta),
\end{align}
with $\tau_2$ given in \eqref{tau2}, and in particular
\begin{align*}
\rho_{\varepsilon_l}(\theta)=2|\ln\varepsilon_l|-\ln({2}|\ln\varepsilon_l|)
-\ln\Big(\frac{ \varsigma(\theta) \tau_2(\theta)}{\beta(\theta)}\Big)+O\Big(\frac{\ln\big({2}|\ln\varepsilon_l|\big)}{|\ln\varepsilon_l|}\Big).
	\end{align*}
The vectors  ${\bf{\ddot f}}=(\ddot{f}_{1}, \cdots, \ddot{f}_{N})^T$ defined in Lemma \ref{lemma6point3},
${\hat{\mathfrak u}}=({\hat{\mathfrak u}}_1, \cdots, {\hat{\mathfrak u}}_N)^T$ defined by \eqref{hbar}
and ${\bf {\tilde w}}=({\bf {\tilde w}}_1, \cdots, {\bf {\tilde w}}_N)^{T}$ in (\ref{equationmathbfwn})-(\ref{boundarymathbfw})
do not depend on $\bf{\tilde h}$.
There hold the estimates
$$
\ddot{f}_{j}=O(1),
\qquad
{\hat{\mathfrak u}}_j(\theta) =O\Big(\Big(\frac{1}{\ln|\varepsilon_l|-\ln(\ln|\varepsilon_l|)}\Big)^{\frac{1}{2}}\Big), \qquad j=1, \cdots, N,
$$
\begin{align}
\frac{1}{|\ln\varepsilon_l|}\, \|{\bf {\tilde w}}_j''\|_{L^2(0, 1)}
+ \|{\bf {\tilde w}}_j'\|_{L^2(0, 1)}
+\|{\bf {\tilde w}}_j\|_{L^{2}(0, 1)}
\, \leq\,
\frac{C}{|\ln\varepsilon_l|},
\quad
j=1, \cdots, N-1,
\label{mathbfwn}
\end{align}
\begin{equation}\label{mathbfwN}
\|{\bf {\tilde w}}_N\|_{H^2(0, 1)} \leq C. 
\end{equation}
For the vector ${\check{\mathfrak u}}=({{\check{\mathfrak u}}_1, \cdots, {\check{\mathfrak u}}_N})^T$,  we have
\begin{align*}
\frac{1}{|\ln\varepsilon_l|}\, \|{\check{\mathfrak u}}''\|_{L^2(0, 1)}
+\frac{1}{\sqrt{|\ln\varepsilon_l|}}\, \|{\check{\mathfrak u}}'\|_{L^2(0, 1)}
+\|{\check{\mathfrak u}}\|_{L^{2}(0, 1)}
\leq
C\, \varepsilon^{\mu}\|\tilde{\mathfrak h}\|_{L^2(0, 1)}\, +\, \frac{C}{|\ln\varepsilon_l|}.
\end{align*}
\end{proposition}
\qed

The proof of Proposition \ref{proposition7point2} will be provided in Section \ref{section6.2}.
By accepting this, we here want to finish the proof of Theorem \ref{theorem 1.1}.
\\[1mm]
{\textbf {Proof of Theorem \ref{theorem 1.1}.}}
The profile of the solution given in \eqref{taketheform} can be determined by the approximate solution given in Section \ref{section4}, see \eqref{globalapproximation}.
The properties of the parameters $f_j$'s in \eqref{fproperties1}-\eqref{fproperties2} can be derived from Proposition \ref{proposition7point2}.

\medskip
As we have stated in Section \ref{section3}, we  shall complete the last step of suitable choosing the parameters $\mathbf{f}$ and  $\mathbf{e}$ by solving (\ref{f})-(\ref{boundary condition 3}).
If $\hat{e}$ solves
\begin{align}
 \begin{split}
\mathbb{ L} (\hat{e})=\, \alpha_5(\theta), \qquad\forall\, 0<\theta<1,
\\[2mm]
\hat{e}'(1)\, +\, {\tilde{b}_6}\, \hat{e}(1)\, =\, 0,
\qquad\quad
\hat{e}'(0)\, +\, {\tilde{b}_5}\, \hat{e}(0)\, =\, 0,
\end{split}
\end{align}
from the definition of $\alpha_5(\theta)$ in \eqref{alpha5}, we get
\begin{equation*}
 \|\hat{e}\|_{H^2(0, 1)}\leq C\varepsilon^\frac{1}{2}.
\end{equation*}
Replacing $e_n$ by $\hat{e}+\tilde{e}_n$, the system (\ref{f})-(\ref{boundary condition 3})
keeps the same form except that the term $\alpha_5(\theta)$ disappear. Moreover, let $\tilde{e}_n$ solves
\begin{align}
 \begin{split}
 \mathbb{L}_2(\tilde{e}_n)=\, \alpha_{6, n}(\theta),
  \qquad\forall\, 0<\theta<1,
  \\[2mm]
\tilde{e}_n'(1)\, +\, {\tilde{b}_6}\, \tilde{e}_n(1)\, =\, 0,
    \qquad\quad
\tilde{e}_n'(0)\, +\, {\tilde{b}_5}\, \tilde{e}_n(0)\, =\, 0,
  \end{split}
\end{align}
then it derives
\begin{equation*}
  \|\tilde{e}_n\|_{H^2(0, 1)}\leq C\varepsilon^\mu.
\end{equation*}

Define the set
$$
\mathcal{D}=\Big\{\, \mathbf{f}, \mathbf{e}\in H^2(\mathcal{S}): \|\mathbf{f}\|_{H^2(0, 1)}\leq D|\ln\varepsilon|^2,
\quad \|\mathbf{e}\|_{**}\leq C\varepsilon^\mu\, \Big\}.
 $$
For $(\bar{\mathbf{f}}, \bar{\mathbf{e}})\in \mathcal{D}$, we can set for $n=1, \cdots, N$
\begin{equation*}
\tilde{h}_n(\mathbf{f}, \mathbf{e})\equiv \varepsilon^{-2}{\mathbb M}_{n1}(\mathbf{f}, \mathbf{f}', \mathbf{f}'', \mathbf{e}, \mathbf{e}', \mathbf{e}'')
\, +\, \varepsilon^{-2}
{\mathbb M}_{n2}(\bar{\mathbf{f}}, \bar{\mathbf{f}}', \bar{\mathbf{e}}, \bar{\mathbf{e}}'),
\end{equation*}
\begin{equation*}
 \tilde{g}_n(\mathbf{f}, \mathbf{e})\, \equiv\, \varepsilon^{-1}\, \beta^2 {\mathbf M}_{n1}(\mathbf{f}, \mathbf{f}', \mathbf{f}'', \mathbf{e}, \mathbf{e}', \mathbf{e}'')
\, +\,
\varepsilon^{-1}\, \beta^2 {\mathbf M}_{n2}(\bar{\mathbf{f}}, \bar{\mathbf{f}}', \bar{\mathbf{e}}, \bar{\mathbf{e}}').
\end{equation*}
We now use Contraction Mapping Principle and Schauder Fixed Point Theorem to solve (\ref{f})-(\ref{boundary condition 3}) with the right hand replacing by $\tilde{h}_n$ and $\tilde{g}_n$.
Whence, by the fact that ${\mathbb M}_{n1}$, ${\mathbf M}_{n1}$ are contractions on $\mathcal{D}$,
making use of the argument developed in Propositions \ref{proposition7point1} and  \ref{proposition7point2}, and the Contraction Mapping Principle, we find $\mathbf{f}$ and $\mathbf{e}$ for a fixed $\bar{\mathbf{f}}$ and $\bar{\mathbf{e}}$.
In this way, we define a mapping $\mathcal{Z}(\bar{\mathbf{f}}, \bar{\mathbf{e}})=(\mathbf{f}, \mathbf{e})$ and the solution of our problem is simply a fixed point of $\mathcal{Z}$. Continuity of ${\mathbb M}_{n2}$ and ${\mathbf M}_{n2}$, $n=1, \cdots, N$, with respect to its parameters and a standard regularity argument allows us to conclude that $\mathcal{Z}$ is compact as mapping from $H^2(0, 1)$ into itself. The Schauder Fixed Point Theorem applies to yield the existence of a fixed point of $\mathcal{Z}$ as required. This ends the proof of Theorem \ref{theorem 1.1}.
\qed

\subsection{Proof of Proposition \ref{proposition7point2}}\label{section6.2}

Note that (\ref{ts2a})-(\ref{boundary0000}) can be concerned as a small perturbation
of a simpler problem in the form, for $n=1, \cdots, N$,
\begin{align}\label{simpleproblem0}
\varepsilon ^2\varsigma\Big[\, \mathcal{H}_1\, f''_n+\mathcal{H}_1'\, f'_n+ \big(\mathcal{H}_2'-\mathcal{H}_3\big)\, f_n\, \Big]
-e^{-\beta(f_n-f_{n-1})}
+e^{-\beta(f_{n+1}-f_n)}
=\varepsilon^{2+\mu} \tilde{h}_n,
\end{align}
\begin{align}
\label{boundary10}
{\mathfrak b}_6f'_n(1)-{\mathfrak b}_7\, f_n(1)=0,
\qquad
{\mathfrak b}_1f'_n(0)-{\mathfrak b}_2\, f_n(0)=0,
\end{align}
where $\mathcal{H}_1$, $\mathcal{H}_2$ and $\mathcal{H}_3$ are functions defined in (\ref{mathcalH1})-(\ref{mathcalH3}).
By similar arguments as done in Section 6 of \cite{weixuyang},
we can finish the proof of Proposition \ref{proposition7point2} if we can solve (\ref{simpleproblem0})-(\ref{boundary10}).

\medskip
We now focus on the resolution theory for (\ref{simpleproblem0})-(\ref{boundary10}),
whose proof basically follows the methods in \cite{delPKowWeiYang} and \cite{YangYang2013}.
However, in this paper, the homogeneous boundary conditions in \eqref{boundary10} make the procedure much more complicated,
which will be divided into three steps.
In the first step, we will find an approximate solution by solving an algebraic system and then derive the improved equivalent nonlinear system of (\ref{simpleproblem0})-(\ref{boundary10}), see \eqref{equationpftildee}-\eqref{boundaryconditionoftildee}.
In step 2, by the decomposition  method, the problem can be further transformed into \eqref{equationofui}-\eqref{boundaryconditionofui}.
To cancel the boundary error terms $\tilde{\mathbb G}_{1}$ and $\tilde{\mathbb G}_{2}$ (see \eqref{boundaryconditionofui}),
we need to find more boundary correction terms ${\hat{\mathfrak u}}_n, \, n=1,\cdots, N$ (see \eqref{hbar}) in the expansions of $\mathfrak{u}_n$'s,
which directly leads to the system \eqref{equationofddotu}-\eqref{boundaryconditionofddotu}.
Finally, after giving the linear resolution theory in Lemma \ref{lemma6point4},
the proof  can be finished by the Contraction Mapping Principle in Step 3.

\medskip
\noindent{\bf Step 1:}
By setting
\begin{equation}\label{transformation}
{\check f}_n(\theta)=\beta(\theta)f_n(\theta),
\end{equation}
we get
\begin{equation}\label{simpleproblem}
\varepsilon ^2\frac{\varsigma}{\beta}
\big[\mathcal{H}_1{\check f}''_n+\tau_1{\check f}'_n+\tau_2{\check f}_n\, \big]
\,-\, e^{-({\check f}_n-{\check f}_{n-1})}
\, +\, e^{-({\check f}_{n+1}-{\check f}_n)}
=\varepsilon^{2+\mu} \tilde{h}_n,
\end{equation}
\begin{align}
\label{boundary1}
{\check f}'_n(1)+K_2\, {\check f}_n(1)=0,
\qquad
{\check f}'_n(0)+K_1\, {\check f}_n(0)=0,
\end{align}
where ${\check f}_0=-\infty, \, {\check f}_{N+1}=\infty$.
Here we have denoted
\begin{equation}\label{tau1}
\tau_1(\theta)=\mathcal{H}_1'(\theta)-2\frac{\beta'(\theta)}{\beta(\theta)}\mathcal{H}_1(\theta), 	
\end{equation}
\begin{equation}\label{tau2}
\tau_2(\theta)=\mathcal{H}_2'(\theta)-\mathcal{H}_3(\theta)
+2\frac{|\beta'(\theta)|^2}{\beta^2(\theta)}\mathcal{H}_1(\theta)
-\frac{\beta''(\theta)}{\beta(\theta)}\mathcal{H}_1(\theta)
-\frac{\beta'(\theta)}{\beta(\theta)}\mathcal{H}_1'(\theta),
\end{equation}
\begin{equation}\label{K1K2}
K_1=\frac{\beta'(0)}{\beta(0)}+\frac{{\mathfrak b}_2}{{\mathfrak b}_1},
\qquad
K_2=\frac{\beta'(1)}{\beta(1)}+\frac{{\mathfrak b}_7}{{\mathfrak b}_6}, 	
\end{equation}
so that
$$
\tau_2(\theta)>0,
\qquad
K_1\, =\, K_2\, =\, 0,
$$
due to the assumptions in (\ref{taupositivity})-(\ref{boundaryadmissibility}).

\medskip
Recall that the assumption \eqref{taupositivity} implies that  ${\varsigma}(\theta)\tau_2(\theta)/\beta(\theta)>0$.
Let us define two positive functions $\rho_\varepsilon(\theta)$ and $\delta(\theta)$ by
\begin{align}
e^{-\rho_\varepsilon(\theta)}
\, =\,
\varepsilon^2\frac{{\varsigma}(\theta)}{\beta(\theta)} \tau_2 (\theta)\rho_\varepsilon(\theta),
\qquad
\frac{1}{\delta^2(\theta)}\, =\, \tau_2(\theta)\rho_\varepsilon(\theta).
\end{align}
We can easily obtain that
\begin{align}
\rho_\varepsilon(\theta)=2|\ln\varepsilon|-\ln\big(2|\ln\varepsilon|\big)
-\ln\Big(\frac{\varsigma(\theta) \tau_2(\theta)}{\beta(\theta)}\Big)+O\Big(\frac{\ln\big(2|\ln\varepsilon|\big)}{|\ln\varepsilon|}\Big),
\end{align}
\begin{align}\label{del1}
\frac{1}{\delta^2(\theta)}\, =\,
\tau_2(\theta)\, \left[\, 2|\ln\varepsilon|
-\ln\big(2|\ln\varepsilon|\big)
-\ln\Big(\frac{\varsigma(\theta)\tau_2(\theta)}{\beta(\theta)}\Big)
+\, O\Big(\frac{\ln\big(2|\ln\varepsilon|\big)}{|\ln\varepsilon|}\Big)\, \right].
\end{align}
Then  multiplying equation (\ref{simpleproblem}) by $\varepsilon^{-2}\delta^2(\theta)$ and setting
\begin{equation*}
	{\check f}_n(\theta)\, =\, 	\Big(n-\frac{N}{2}-\frac{1}{2}\Big)\rho_{\varepsilon}(\theta) +{\hat f}_n(\theta), \quad  n=1, \cdots, N,
\end{equation*}
we get an  equivalent system,  for $n=1, \cdots, N$,
\begin{align}\label{ts3a}
 &\delta^2\big[\mathcal{H}_1{\hat f}''_n+\tau_1{\hat f}'_n+\tau_2{\hat f}_n\big]
 \,-\, e^{-({\hat f}_n-{\hat f}_{n-1})}
\, +\, e^{-({\hat f}_{n+1}-{\hat f}_n)}
\nonumber\\[2mm]
&=\varepsilon^\mu\delta^2\frac{\beta}{\varsigma} \tilde{h}_n-	
\delta^2\Big(n-\frac{N}{2}-\frac{1}{2}\Big)\rho_{\varepsilon}''
-\delta^2\tau_1\Big(n-\frac{N}{2}-\frac{1}{2}\Big)\rho_{\varepsilon}'
-\Big(n-\frac{N}{2}-\frac{1}{2}\Big)
,
\end{align}
where ${\hat f}_0=-\infty, \, {\hat f}_{N+1}=\infty$.
The boundary conditions become
\begin{align}\label{boundarysim}
{\hat f}_n'(1)=-\Big(n-\frac{N}{2}-\frac{1}{2}\Big)\rho_{\varepsilon}'(1),
\qquad
{\hat f}_n'(0)=-\Big(n-\frac{N}{2}-\frac{1}{2}\Big)\rho_{\varepsilon}'(0).
\end{align}

\begin{remark}\label{remark61}
Note that the terms of order $O(|\ln\varepsilon|)$ in the right hand sides of the equations in (\ref{boundarysim}) disappear
so that they are of $O(1)$ due to the assumptions $K_1=K_2=0$ which are exactly given in (\ref{boundaryadmissibility}).
\qed
\end{remark}

\medskip
First, we want to cancel the terms of $O(1)$ in right hand side of (\ref{ts3a}). To this end, we will introduce the following lemma

\begin{lemma}\label{lemma6point3}
There exists a solution $\ddot{\bf f}=(\ddot{f}_1, \cdots, \ddot{f}_N)^{T}$ to the following nonlinear algebraic system
\begin{align}\label{equationofe0}
 -e^{- (\ddot{f}_{n}-\ddot{f}_{n-1})}+ e^{-(\ddot{f}_{n+1}-\ddot{f}_{n})}=-\Big(n-\frac{N}{2}-\frac{1}{2}\Big)
\end{align}
with $n$ running from $1$ to $N$, where $\ddot{f}_0=-\infty, \, \ddot{f}_{N+1}=\infty$.
\end{lemma}

\begin{proof}
By setting
\begin{align}\label{deinfinitionofa_n0}
a_0\, =\, a_{N}\, =\, 0, \quad a_n\, =\, e^{-(\ddot{f}_{n+1}-\ddot{f}_{n})}, \quad
n=1, \cdots, N-1,
\end{align}
the proof can be found in the solving method for equation (7.10) in \cite{YangYang2013}.
\end{proof}

\medskip
We set ${\hat f}_n\, =\, \ddot{f}_n + {\tilde f}_n, \, n=1, \cdots, N$,  where $\ddot{f}_n=O(1)$ satisfies the system \eqref{equationofe0}. It is obvious that system (\ref{ts3a})-(\ref{boundarysim}) is equivalent to the following nonlinear system of equations,
\begin{align}\label{equationpftildee}
\begin{split}
\delta^2\Bigl[\mathcal{H}_1\, {\tilde f}_n''
+\tau_1{{\tilde f}_n'}+ \tau_2{{\tilde f}_n}\Bigr]
&\, +\, a_{n-1}\bigl(\, {\tilde f}_n-{\tilde f}_{n-1} \bigr)
\\
&\,-\, a_{n} \bigl(\, {\tilde f}_{n+1}-{\tilde f}_{n}  \bigr)
\, =\, \delta^2\, \varepsilon^{\mu}\frac{\beta}{\varsigma} \tilde{h}_{n}\, +\, \delta^2\tilde{g}_n+{\mathfrak N}_n({\bf{\tilde f}}),
\end{split}
\end{align}
with boundary conditions
\begin{align}\label{boundaryconditionoftildee}
{\tilde f}_n'(0)\, =\, G_{1, n}, \qquad {\tilde f}_n'(1)\, =\, G_{2, n},
\end{align}
where we have denoted
\begin{align}
{\tilde g}_n\, =
\,-\, \Big(n-\frac{N}{2}-\frac{1}{2}\Big)\rho''_{\varepsilon}
\,-\, \tau_1\Big(n-\frac{N}{2}-\frac{1}{2}\Big)\rho'_{\varepsilon}
\,-\, \tau_2\ddot{f}_n,
\end{align}
\begin{align}\label{G1G2}
\, G_{1, n}=-\Big(n-\frac{N}{2}-\frac{1}{2}\Big)\rho'_{\varepsilon}(0), \qquad
\, G_{2, n}=-\Big(n-\frac{N}{2}-\frac{1}{2}\Big)\rho'_{\varepsilon}(1).
\end{align}
Moreover, the nonlinear terms ${\mathfrak N}_n, \, n=1, \cdots, N$, are given by
\begin{align}
\begin{split}
{\mathfrak N}_n({\bf{\tilde f}}) =&\, a_{n-1} \big[e^{-(\, {\tilde f}_n    - {\tilde f}_{n-1} )} -1+ \bigl(\, {\tilde f}_n    - {\tilde f}_{n-1} \bigr)\big]
\,-\, a_n\big[e^{-(\, {\tilde f}_{n+1} - {\tilde f}_{n} )}-1+ \bigl(\, {\tilde f}_{n+1} - {\tilde f}_{n}  \bigr)\big],
\end{split}
\end{align}
and ${\tilde f}_0=-\infty, \, {\tilde f}_{N+1}=\infty$.

\medskip
\medskip
\noindent{\bf Step 2:} The first try is to decompose the above system. We will denote:
\begin{align*}
{\bf {\tilde f}}=({\tilde f}_1, \cdots, {\tilde f}_N)^{T},
\qquad
{\bf {\tilde h}}=({\tilde h}_1, \cdots, {\tilde h}_N)^{T},
\qquad
{\bf {\tilde g}}=({\tilde g}_1, \cdots, {\tilde g}_N)^{T},
\end{align*}
\begin{align*}
{\mathfrak N}({\bf{\tilde f}})=\big({\mathfrak N}_1({\bf{\tilde f}}), \cdots, {\mathfrak N}_N({\bf{\tilde f}})\big)^{T},
\qquad
{\mathbb G}_1=(G_{1, 1}, \cdots, G_{1, N})^{T},
\qquad
{\mathbb G}_2=(G_{2, 1}, \cdots, G_{2, N})^{T}.
\end{align*}
Then system (\ref{equationpftildee})  becomes:
\begin{align}\label{linear1}
\begin{split}
\delta^{2}{\bf I}\Bigl[\, \mathcal{H}_1\, \frac{\mathrm{d}^2}{\mathrm{d}\theta^2}
\, +\, \tau_1\frac{\mathrm{d}}{\mathrm{d}\theta}
\, +\, \tau_2 \, \Bigr]{\bf {\tilde f}}
\, +\, {\mathbb A}{\bf {\tilde f}}
\, =\, \delta^2\, \varepsilon^{\mu} \frac{\beta}{\varsigma}{\bf{\tilde h}}\, +\, \delta^2 {\bf{\tilde g}}\, +\, {\mathfrak N}({\bf{\tilde f}}),
\end{split}
\end{align}
where $\bf I$ is a $ N\times N$ unit Matrix and the Matrix ${\mathbb A}$ defined as
\begin{align}\label{definitionofbfB}
{\mathbb A}\, =\,  \left(
\begin{array}{ccccccccc}
a_1 & -a_1&0&  0&\cdots &0&0&0&0
\\
-a_1 &(a_1+a_2)&-a_2& 0&\cdots &0&0&0&0
\\
\vdots& \vdots&\vdots & \vdots&\ddots& \vdots& \vdots& \vdots&\vdots
\\
0 & 0&0&0&\cdots&  0&-a_{N-2}&(a_{N-2}+a_{N-1})&-a_{N-1}
\\
0 & 0&0&0&\cdots&0&  0&-a_{N-1}&a_{N-1}
\end{array}
\right).
\end{align}
For the symmetric matrix ${\mathbb A}$,  using elementary matrix operations it is easy to prove that there exists an invertible matrix ${\bf Q}$ such that
$$
{\bf Q}{\mathbb A}{\bf Q}^T=\mbox{diag}(a_1, \cdots, a_{N-1}, 0).
$$
Since $a_1, \cdots, a_{N-1}$ are positive constants defined in (\ref{deinfinitionofa_n0}),
then all eigenvalues of the matrix ${\mathbb A}$ are
\begin{align*}
 \lambda_1\geq\lambda_2\geq \cdots\geq\lambda_{N-1}>\lambda_N=0.
\end{align*}
Moreover, since ${\mathbb A}$ is a symmetric matrix, there exists another invertible matrix $\bf P$ independent of $\theta$ with the form
\begin{align}\label{mathbfP}
\mathbf{P}\mathbf{P}^T={\bf I},
\qquad
\mathbf{P}=
\left(
  \begin{array}{cccc}
    p_{11} & \cdots & p_{1N-1} & \frac{1}{\sqrt N} \\
    p_{21} & \cdots & p_{2N-1} & \frac{1}{\sqrt N} \\
    \vdots & \vdots & \vdots & \vdots \\
    p_{N1} & \cdots &p_{NN-1}  & \frac{1}{\sqrt N} \\
  \end{array}
\right),
\end{align}
in  such a way that
\begin{align}\label{P}
{\bf P}^T{\mathbb A}{\bf P}=\mbox{diag}(\lambda_1, \lambda_2, \cdots, \lambda_{N-1}, \lambda_{N}).
\end{align}

\medskip
We denote
\begin{align}\label{alphaanddelta}
\kappa=\frac{1}{ 2|\ln\varepsilon|-\ln(2|\ln\varepsilon|)},
\qquad
\frac{\delta^{-2}(\theta)}
{\,\kappa\, }
\, =\, \tau_2(\theta)\, +\, {\tilde\sigma}(\theta).
\end{align}
By (\ref{del1}) we have
\begin{align}
{\tilde\sigma}(\theta)\, =\, O\Big(\frac{1}{|\ln\varepsilon|}\Big).
\label{sigma ep}
\end{align}
Multiplying (\ref{linear1}) by $\tau_2(\theta)\, +\, {\tilde\sigma}(\theta)$,  we get the following system
\begin{align}\label{linear}
\begin{split}
\kappa {\bf I}\Bigl[\, \mathcal{H}_1\, \frac{\mathrm{d}^2}{\mathrm{d}\theta^2}
+\tau_1\frac{\mathrm{d}}{\mathrm{d}\theta}
\, +\, \tau_2 \Bigr]{\bf {\tilde f}}
\, +\,
\bigl(\tau_2+{\tilde\sigma}\bigr){\mathbb A}{\bf {\tilde f}}
\, =\,
\kappa\varepsilon^{\mu}\frac{\beta}{\varsigma} {\bf{\tilde h}}
\, +\, \kappa {\bf{\tilde g}}
\, +\, \big(\tau_2+{\tilde\sigma}\big){\mathfrak N}({\bf{\tilde f}}).
\end{split}
\end{align}

\medskip
Now,  define six new vectors
\begin{align*}
\mathfrak{u}=(\mathfrak{u}_1, \cdots, \mathfrak{u}_{N})^T={\bf P}^T{\bf {\tilde f}},
\qquad
\tilde{\mathfrak{h}}=(\tilde{\mathfrak{h}}_1, \cdots, \tilde{\mathfrak{h}}_{N})^T= \frac{\beta}{\varsigma}{\bf P}^T {\bf {\tilde h}},
\qquad
\tilde{\mathfrak{g}}=(\tilde{\mathfrak{g}}_1, \cdots, \tilde{\mathfrak{g}}_{N})^T= {\bf P}^T {\bf{\tilde g}},
\end{align*}
\begin{align*}
{\tilde{\mathfrak N}}(\mathfrak{u})=\big({\tilde{\mathfrak N}}_1(\mathfrak{u}), \cdots, {\tilde{\mathfrak N}}_N(\mathfrak{u})\big)^{T}
=\big(\tau_2+{\tilde\sigma}\big){\bf P}^T{\mathfrak N}({\bf{\tilde f}})
=\big(\tau_2+{\tilde\sigma}\big){\bf P}^T{\mathfrak N}({\bf P}\mathfrak{u}),
\end{align*}
\begin{align*}
\tilde{{\mathbb G}}_1\, =\, (\tilde{G}_{1, 1}, \cdots, \tilde{G}_{1, N})^{T}={\bf P}^T {\mathbb G}_1,
\qquad
\tilde{{\mathbb G}}_2\, =\, (\tilde{G}_{2, 1}, \cdots, \tilde{G}_{2, N})^{T}={\bf P}^T {\mathbb G}_2.
\end{align*}
Note that the form of ${\mathbf P}$ in (\ref{mathbfP}) and the expressions of ${\mathbb G}_1$ and ${\mathbb G}_2$ in (\ref{G1G2}) imply that
$$
\tilde{G}_{1, N}=\tilde{G}_{2, N}=0.
$$
Therefore \eqref{equationpftildee}-(\ref{boundaryconditionoftildee}) become
\begin{align}\label{equationofui}
\kappa\, \Big[\mathcal{H}_1\, \mathfrak{u}''+\tau_1\mathfrak{u}'+ \tau_2\mathfrak{u}\Big]
\, +\, \mbox{diag}(\lambda_1, \cdots, \lambda_N)\, (\tau_2+{\tilde\sigma})\mathfrak{u}
=\, \kappa\, \varepsilon^{\mu}\tilde{\mathfrak h}
\, +\, \kappa \tilde{\mathfrak{g}}
\, +\, {\tilde{\mathfrak N}}(\mathfrak{u}),
\end{align}
with boundary conditions
\begin{equation}\label{boundaryconditionofui}
\mathfrak{u}'(0)=\tilde{{\mathbb G}}_1,
\quad
\mathfrak{u}'(1)=\tilde{{\mathbb G}}_2.
\end{equation}

For the convenience of notation, we denote
\begin{align*}
 \ell_{n}\, =\, \Big( \frac{\kappa}{\lambda_n}\Big)^{\frac{1}{2}},
 \qquad
 \Pi(s)\, =\, \tau_2(s)+{\tilde\sigma}(s),
\end{align*}
\begin{equation*}
\vartheta(\theta)=\int_{0}^{\theta} \Pi(s)^{\frac{1}{2}}\, {\mathrm{d}}s,
\qquad
l_0=\int_{0}^{1} \Pi(s)^{\frac{1}{2}}\, {\mathrm{d}}s.
\end{equation*}
In order to cancel the error terms on the boundary in \eqref{boundaryconditionofui}, we introduce the following functions
\begin{align}\label{hbar}
{\hat{\mathfrak u}}_n(\theta)\, =&\chi(\theta)\frac{\tilde{G}_{1, n}\ell_{n}}{\sqrt{\Pi(0)}}\, \sin\Big(\frac{\vartheta(\theta)}{\ell_{n}}\, \Big)
\,-\,
\big(1-\chi(\theta)\big)\frac{\tilde{G}_{2, n}\ell_{n}}{\sqrt{\Pi(1)}}
\sin\Big(\frac{l_0-\vartheta(\theta)}{\ell_{n}}\, \Big), \quad n=1, \cdots, N-1,
\nonumber
\end{align}
and
$$
{\hat{\mathfrak u}}_N(\theta)=0.
$$
In the above, $\chi$ is a smooth cut-off function with the properties
$$
\chi(\theta)=1\quad \mbox{if } |\theta|<1/8
\qquad\mbox{and}\qquad
\chi(\theta)=0\quad \mbox{if } |\theta|>2/8.
$$
It is easy to show
\begin{align}
\|{\hat{\mathfrak u}}_n\|_{L^2(0, 1)}\, \leq\, \frac{C}{\sqrt{|\ln\varepsilon|}}.
\end{align}
For later use, we compute, for $n=1, \cdots, N-1$,
\begin{align*}
{\hat{\mathfrak u}}_n'(\theta)
\, =&\chi(\theta)\frac{\tilde{G}_{1, n}}{\sqrt{\Pi(0)}}\, \sqrt{\Pi(\theta)}\cos\Big(\frac{\vartheta(\theta)}{\ell_{n}}\, \Big)
\, +\,
\big(1-\chi(\theta)\big)\frac{\tilde{G}_{2, n}}{\sqrt{\Pi(1)}}\, \sqrt{\Pi(\theta)}\cos\Big(\frac{l_0-\vartheta(\theta)}{\ell_{n}}\, \Big)
\\[2mm]
&\, +\, \frac{\tilde{G}_{1, n}\ell_{n}}{\sqrt{\Pi(0)}}\, \chi'(\theta)\sin\Big(\frac{\vartheta(\theta)}{\ell_{n}}\, \Big)
\, +\,
\frac{\tilde{G}_{2, n}\ell_{n}}{\sqrt{\Pi(1)}}
\chi'(\theta)\sin\Big(\frac{l_0-\vartheta(\theta)}{\ell_{n}}\, \Big).
\end{align*}
This implies that $\hat{\mathfrak u}_n$ satisfies the following boundary conditions
\begin{align}
{\hat{\mathfrak u}}'_n(0)\, =\, \tilde{G}_{1, n},
\qquad
{\hat{\mathfrak u}}'_n(1)\, =\, \tilde{G}_{2, n}, \quad n=1, \cdots, N.
\end{align}	
For $n=1, \cdots, N-1$, there holds
\begin{align*}
{\hat{\mathfrak u}}_n''(\theta)\, =&
-\chi(\theta)\frac{\tilde{G}_{1, n}}{\sqrt{\Pi(0)}}\,
\, \frac{\Pi(\theta)}{\ell_{n}}\sin\Big(\frac{\vartheta(\theta)}{\ell_{n}}\, \Big)
\,+\,
\big(1-\chi(\theta)\big)\frac{\tilde{G}_{2, n}}{\sqrt{\Pi(1)}}\, \, \frac{\Pi(\theta)}{\ell_{n}}\sin\Big(\frac{l_0-\vartheta(\theta)}{\ell_{n}}\, \Big)
\\[2mm]
&\,+\,\chi(\theta)\frac{\tilde{G}_{1, n}}{\sqrt{\Pi(0)}}\, \frac{\Pi'(\theta)}{\sqrt{\Pi(\theta)}}\, \cos\Big(\frac{\vartheta(\theta)}{\ell_{n}}\, \Big)
\,+\,
\big(1-\chi(\theta)\big)\frac{\tilde{G}_{2, n}}{\sqrt{\Pi(1)}}\, \frac{\Pi'(\theta)}{\sqrt{\Pi(\theta)}}\,
\cos\Big(\frac{l_0-\vartheta(\theta)}{\ell_{n}}\, \Big)
\\[2mm]
&\,+\,\chi'(\theta)\frac{\tilde{G}_{1, n}}{\sqrt{\Pi(0)}}\, \Pi(\theta)\cos\Big(\frac{\vartheta(\theta)}{\ell_{n}}\, \Big)
\,-\,
\chi'(\theta)\frac{\tilde{G}_{2, n}}{\sqrt{\Pi(1)}}\, \Pi(\theta)\cos\Big(\frac{l_0-\vartheta(\theta)}{\ell_{n}}\, \Big)
\\[2mm]
&\,+\,\chi''(\theta)\frac{\tilde{G}_{1, n}\ell_{n}}{\sqrt{\Pi(0)}}\, \sin\Big(\frac{\vartheta(\theta)}{\ell_{n}}\, \Big)
\,+\,
\chi''(\theta)\frac{\tilde{G}_{2, n}\ell_{n}}{\sqrt{\Pi(1)}}\, \sin\Big(\frac{l_0-\vartheta(\theta)}{\ell_{n}}\, \Big).
\end{align*}
Whence, we obtain, for $n=1, \cdots, N$,
\begin{align}
\Big\|\kappa\Big[\mathcal{H}_1\, {\hat{\mathfrak u}}_n''+\tau_1 {\hat{\mathfrak u}}_n'+\tau_2{\hat{\mathfrak u}}_n\Big]
+\lambda_n(\tau_2+{\tilde\sigma}){\hat{\mathfrak u}}_n\Big\|_{L^2(0, 1)}
\, \leq\,
\frac{C}{|\ln\varepsilon|}.
\end{align}

\medskip
Letting $\mathfrak{u}\, =\, {\tilde{\mathfrak u}}\, +\, {\hat{\mathfrak u}}$ with ${\hat{\mathfrak u}}=({\hat{\mathfrak u}}_1, \cdots, {\hat{\mathfrak u}}_N)^T$, the system \eqref{equationofui}-\eqref{boundaryconditionofui}
is equivalent to the following system, for $n=1, \cdots, N$,
\begin{align}\label{equationofddotu}
\begin{split}
\kappa\Big[
\mathcal{H}_1\, {\tilde{\mathfrak u}}_n''
+\tau_1 {\tilde{\mathfrak u}}_n'
+\tau_2{\tilde{\mathfrak u}}_n
\Big]
+\lambda_n(\tau_2+{\tilde\sigma}){\tilde{\mathfrak u}}_n
\, =\, \kappa\varepsilon^{\mu}{\mathfrak h}_n
\, +\, \kappa \tilde {\mathfrak{g}}_n
\, +\, \hat{\mathfrak{g}}_n
\, +\, {\tilde{\mathfrak N}}_n({\tilde{\mathfrak u}+{\hat{\mathfrak u}}}),
\end{split}
\end{align}
with boundary conditions
\begin{equation}\label{boundaryconditionofddotu}
{\tilde{\mathfrak u}}_n'(0)\, =\, 0,
\qquad
{\tilde{\mathfrak u}}_n'(1)\, =\, 0,
\end{equation}
where
\begin{equation*}
 \hat{\mathfrak{g}}_n
\, =\, -
\kappa\Big[\mathcal{H}_1\, {\hat{\mathfrak u}}_n''+\tau_1{\hat{\mathfrak u}}_n'
\, +\, \tau_2{\hat{\mathfrak u}}_n\Big]
\,-\, \lambda_n(\tau_2+{\tilde\sigma}){\hat{\mathfrak u}}_n, \quad
n=1, \cdots, N-1,
\end{equation*}
and also
\begin{equation*}
  \tilde{\mathfrak{g}}_N\, =\, 0.
\end{equation*}
For later use, we will estimate the terms in the right hand of \eqref{equationofddotu}.
For $n=1, \cdots, N-1$, there hold
\begin{align*}
{\tilde{\mathfrak N}}_n({\tilde{\mathfrak u}+{\hat{\mathfrak u}}})
&=\big(\tau_2+{\tilde\sigma}\big)\big({\bf P}^T{\mathfrak N}({\bf P}(\tilde{\mathfrak{u}}+{\hat{\mathfrak u}}))\big)_n
=\big(\tau_2+{\tilde\sigma}\big)\sum_{i=1}^{N}p_{in}{\mathfrak N}_i({\bf P}(\tilde{\mathfrak{u}}+{\hat{\mathfrak u}})),
\end{align*}
and also
\begin{equation*}
 \tilde{{\mathfrak h}}_N\, =\, \frac{1}{\sqrt{N}}\, \sum_{i=1}^{N}\, \tilde{h}_i,
 \qquad
\tilde{{\mathfrak g}}_N\, =\, \frac{1}{\sqrt{N}}\, \sum_{i=1}^{N}\, \tilde{g}_i,
\end{equation*}
\begin{equation*}
 {\tilde{\mathfrak N}}_N({\tilde{\mathfrak u}+{\hat{\mathfrak u}}})
 =\big(\tau_2+{\tilde\sigma}\big)\big({\bf P}^T{\mathfrak N}({\bf P}(\tilde{\mathfrak{u}}+{\hat{\mathfrak u}}))\big)_N
 =\big(\tau_2+{\tilde\sigma}\big)\frac{1}{\sqrt{N}}\sum_{i=1}^{N}{\mathfrak N}_i({\bf P}(\tilde{\mathfrak{u}}+{\hat{\mathfrak u}}))=0.
\end{equation*}
According to the definitions of $\tilde{\mathfrak h}_n$'s and $\tilde{\mathfrak g}_n$'s, we can easily get
\begin{align}\label{evaluefrakg}
 \|\tilde{\mathfrak h}_n\|_{L^2(0, 1)}\, \leq\, C, \quad
\|\tilde{\mathfrak g}_n\|_{L^2(0, 1)}\, \leq\, C, \quad
n=1, \cdots, N.
\end{align}

\medskip
\noindent{\bf Step 3:}
For the purpose of using a fixed point argument to solve  \eqref{equationofddotu}-\eqref{boundaryconditionofddotu},
we concern the following resolution theory for the linear differential equations.

\begin{lemma}\label{lemma6point4}\
{\textbf{(1).}}
Assume that the non-degeneracy condition (\ref{nondegeneracy}) holds. For any small $\varepsilon$, there exists a unique solution $v$
to the equation
\begin{align}\label{equationofvvv}
\Big[\mathcal{H}_1\frac{\mathrm{d}^2}{\mathrm{d}\theta^2}+\tau_1\, \frac{\mathrm{d}}{\mathrm{d}\theta} +\tau_2\Big]v
\, =\, h,
\qquad
v'(0)\, =\, 0,
\qquad
v'(1)\, =\, 0,
\end{align}
with the estimate
\begin{align}
\|v\|_{H^2(0, 1)} \leq C\, {\|h\|_{L^2(0, 1)}}.
\label{l2estimatev}
\end{align}

\noindent{\textbf{(2).}}
Consider the following system, for $n=1, \cdots, N-1$,
\begin{equation}\label{equationofvn}
\begin{split}
\kappa\, \Big[\mathcal{H}_1\frac{\mathrm{d}^2}{\mathrm{d}\theta^2}
+\tau_1\, \frac{\mathrm{d}}{\mathrm{d}\theta}
+\tau_2\Big]v_n
&\, +\,
\lambda_n\Bigl(\tau_2+{\tilde\sigma}\Bigr)v_n
\, =\, \mathfrak{p}_n,
\\
\qquad v_n'(0)\, =\, & 0,
\qquad v_n'(1)\, =\, 0.
\end{split}
\end{equation}
There exists a sequence $\{\varepsilon_l, \, l\in {\mathbb N}\}$ approaching $0$ and satisfying the gap condition (\ref{gapconditionofve}) such that problem (\ref{equationofvn}) has a unique solution $\bf{v}\, =\,  \bf{v}(\bf{\mathfrak{p}})$ and
\begin{align}\label{h2estimate1}
\frac{1}{|\ln\varepsilon_l|}\, \|{\bf v}''\|_{L^2(0, 1)}
+\frac{1}{\sqrt{|\ln\varepsilon_l|}}\, \|{\bf v}'\|_{L^2(0, 1)}
+\|{{\bf v}}\|_{L^{2}(0, 1)}\leq C\, \sqrt{|\ln\varepsilon_l|}\, \|{\bf \mathfrak{p}}\|_{L^2(0, 1)},
\end{align}
where ${\bf v}=(v_1, \cdots, v_{N-1})^{T}$ and ${\bf{\mathfrak{p}}}=(\mathfrak{p}_1, \cdots, \mathfrak{p}_{N-1})^T$.
Moreover,  if ${\bf{ \mathfrak{p}}}\in H^2(0, 1)$ then
\begin{align}
\frac{1}{|\ln\varepsilon_l|}\, \|{\bf v}''\|_{L^2(0, 1)}
+ \|{\bf v}'\|_{L^2(0, 1)}
+\|{{\bf v}}\|_{L^{2}(0, 1)}\leq C\, \|{\bf \mathfrak{p}}\|_{H^2(0, 1)}.
\label{h2estimate2}
\end{align}
\end{lemma}

\begin{proof}
We can use the inverse of the transformation in (\ref{transformation}), i.e., $v\, =\, \beta(\theta)\, \tilde{v}$, and then obtain
\begin{equation}
\frac{\mathrm{d}}{\mathrm{d}\theta}\, v\, =\, \beta'\, \tilde{v}\, +\, \beta\, \frac{\mathrm{d} }{\mathrm{d}\theta}\, \tilde{v}, \qquad
\frac{\mathrm{d}^2}{\mathrm{d}\theta^2}\, v\, =\, \beta''\, \tilde{v}\, +\, 2\, \beta'\, \frac{\mathrm{d} }{\mathrm{d}\theta}\, \tilde{v}\, +\, \beta\, \frac{\mathrm{d}^2}{\mathrm{d}\theta^2}\, \tilde{v}.
\end{equation}
Therefore, we get an equivalent problem of (\ref{equationofvvv})
\begin{equation}\label{equationoftidlev}
\begin{split}
 \mathcal{H}_1\tilde{v}_{\theta\theta}
 \, +\, \mathcal{H}_1'\, \tilde{v}_{\theta}
 \, +\, \big( \mathcal{H}_2'-\mathcal{H}_3\big) \, \tilde{v} & \, =\, \frac{1}{\beta} h,
 \\[2mm]
 {\mathfrak b}_1\tilde{v}'(0) \,-\, {\mathfrak b}_2\, \tilde{v}(0)\, =\, 0,
\quad
 {\mathfrak b}_6\tilde{v}'(1)& \,-\, {\mathfrak b}_7\, \tilde{v}(1)\, =\, 0.
\end{split}
\end{equation}
Recalling the definition of $\mathcal{H}_1$, $\mathcal{H}_2$ and $\mathcal{H}_3$ in \eqref{mathcalH1}-\eqref{mathcalH3} and applying the non-degeneracy condition (\ref{nondegeneracy}), we can solve \eqref{equationoftidlev} directly.
The proof of the second part is similar as that for Claim 1 and Claim 2 in \cite{YangYang2013}.
The details are omitted here.
\end{proof}

\medskip
In order to solve (\ref{equationofddotu})-(\ref{boundaryconditionofddotu}),  we first concern the system, for $n=1, \cdots, N-1$
\begin{align}\label{equationmathbfwn}
&\kappa \Big[\mathcal{H}_1{\bf {\tilde w}}_n''
\, +\, \tau_1\, {\bf {\tilde w}}_n'
\, +\, \tau_2\, {\bf {\tilde w}}_n\Big]
\, +\, \lambda_n\big(\tau_2\, +\, {\tilde\sigma}\big){\bf {\tilde w}}_n
\, =\,
\kappa\tilde{\mathfrak g}_n
\, +\,
\hat{\mathfrak{g}}_n,
\end{align}
\begin{equation}\label{equationmathbfwN}
\mathcal{H}_1{\bf {\tilde w}}_N''
\, +\, \tau_1\, {\bf {\tilde w}}_N'
\, +\, \tau_2\, {\bf {\tilde w}}_N
\, =\, \tilde{\mathfrak g}_N,
\end{equation}
with boundary conditions
\begin{equation}\label{boundarymathbfw}
{\bf {\tilde w}}_n'(0)\, =\, 0,
\qquad
{\bf {\tilde w}}_n'(1)\, =\, 0,
\qquad n=1, \cdots, N.
\end{equation}
Using Lemma \ref{lemma6point4}, we can solve the above system and get the estimates as in (\ref{mathbfwn})-(\ref{mathbfwN}).
The substituting
$$
{\tilde{\mathfrak u}}_n={\bf {\tilde w}}_n+{\check{\mathfrak u}}_n,\quad n=1, \cdots, N,
$$
 will then imply that
the nonlinear problem (\ref{equationofddotu})-(\ref{boundaryconditionofddotu}) can be transformed into the following system
for $n=1, \cdots, N$,
\begin{equation}\label{equationofddotu1}
\kappa\Big[\mathcal{H}_1\, {\check{\mathfrak u}}_n''
\, +\, \tau_1\, {\check{\mathfrak u}}_n'
\, +\, \tau_2\, {\check{\mathfrak u}}_n\Big]
\, +\, \lambda_n\, \big(\tau_2\, +\, \tilde{\sigma}\big)
{\check{\mathfrak u}}_n
\, =\,
\kappa\varepsilon^{\mu}{\mathfrak h}_n
\, +\, \tilde{\mathfrak N}_n({{\check{\mathfrak u}}+{\bf {\tilde w}}+\hat{\mathfrak u}}),
\end{equation}
with boundary conditions
\begin{equation}\label{boundaryconditionofddotu1}
{\check{\mathfrak u}}_n'(0)\, =\, 0,
\qquad
{\check{\mathfrak u}}_n'(1)\, =\, 0
\qquad n=1, \cdots, N.
\end{equation}

\medskip
Finally, we claim that problem (\ref{equationofddotu1})-(\ref{boundaryconditionofddotu1}) can be solved by using Lemma \ref{lemma6point4} and a Contraction Mapping Principle in the set
\begin{align*}
{\mathcal{X}}\, =\, \left\{\, {\check{\mathfrak u}}\in H^2(0, 1)
\, :\,
\frac{1}{|\ln\varepsilon|}\, \|{\check{\mathfrak u}}''\|_{L^2(0, 1)}
\, +\, \frac{1}{\sqrt{|\ln\varepsilon|}}\, \|{\check{\mathfrak u}}'\|_{L^2(0, 1)}
\, +\, \|{\check{\mathfrak u}}\|_{L^2(0, 1)}
\, \leq\,
\frac{C}{\sqrt{|\ln\varepsilon|}}\, \right\}.
\end{align*}
In fact, this can be done in the following way.
According to the definition of ${\tilde{\mathfrak N}} $, we obtain, for any $ \, {\check{\mathfrak u}}^0 \in 	{\mathcal{X}}$
\begin{align*}
{\tilde{\mathfrak N}}({\check{\mathfrak u}}^0+{\bf{\tilde w}}+{\hat{\mathfrak u}})
=\big(\tau_2+{\tilde\sigma}\big){\bf P}^T
\left(\begin{array}{c}
{\mathfrak N}_1\big({\bf P}(\check{\mathfrak u}^0+{\bf {\tilde w}}+{\hat{\mathfrak u}})\big)
\\[2mm]
\vdots
\\[2mm]
{\mathfrak N}_N\big({\bf P}(\check{\mathfrak u}^0+{\bf {\tilde w}}+{\hat{\mathfrak u}})\big)
\end{array}
\right).
\end{align*}
This implies that
\begin{equation*}
\big\|	{\tilde{\mathfrak N}}({\check{\mathfrak u}}^0+{\bf{\tilde w}}+{\hat{\mathfrak u}})\big\|_{L^2(0, 1)}
\leq C\sum_{n=1}^{N} \big\|{\mathfrak N}_{n}\big({\bf P}(\check{\mathfrak u}^0+{\bf {\tilde w}}+{\hat{\mathfrak u}})\big)\big\|_{L^2(0, 1)},
	\end{equation*}
where the expression of ${\mathfrak N}_{n} $ is
\begin{align*}
{\mathfrak N}_{n}\big({\bf P}(\check{\mathfrak u}^0+{\bf {\tilde w}}+{\hat{\mathfrak u}})\big)
=&\, a_{n-1} \Bigg\{e^{-\left[\, ( {\bf P}({\check{\mathfrak u}}^0+{\bf {\tilde w}}+{\hat{\mathfrak u}}))_n
-\big({\bf P}(\check{\mathfrak u}^0+{\bf {\tilde w}}+{\hat{\mathfrak u}})\big)_{n-1} \right]}
-1
\\
& \qquad \quad
+\Big[( {\bf P}({\check{\mathfrak u}}^0+{\bf {\tilde w}}+{\hat{\mathfrak u}}))_n
- ( {\bf P}({\check{\mathfrak u}}^0+{\bf {\tilde w}}+{\hat{\mathfrak u}}))_{n-1} \Big]\Bigg\}
\\[2mm]
&
\,-\, a_n\Bigg\{e^{-\left[\, ( {\bf P}({\check{\mathfrak u}}^0+{\bf {\tilde w}}+{\hat{\mathfrak u}}))_{n+1} - ( {\bf P}({\check{\mathfrak u}}^0+{\bf {\tilde w}}+{\hat{\mathfrak u}}))_{n}  \right]}
-1
\\
& \qquad \quad+\Big[\, ( {\bf P}({\check{\mathfrak u}}^0+{\bf {\tilde w}}+{\hat{\mathfrak u}}))_{n+1} - ( {\bf P}({\check{\mathfrak u}}^0+{\bf {\tilde w}}+{\hat{\mathfrak u}}))_{n}  \Big]\Bigg\}
\\[2mm]
=&O\Big(\Big|\big( {\bf P}({\check{\mathfrak u}}^0+{\bf {\tilde w}}+{\hat{\mathfrak u}}))_n    - ( {\bf P}({\check{\mathfrak u}}^0+{\bf {\tilde w}}+{\hat{\mathfrak u}})\big)_{n-1} \Big|^2\Big)
\\[2mm]
&+ O\Big(\Big| \big( {\bf P}({\check{\mathfrak u}}^0+{\mathbf w}+{\hat{\mathfrak u}}))_{n+1} - ( {\bf P}({\check{\mathfrak u}}^0+{\bf {\tilde w}}+{\hat{\mathfrak u}})\big)_{n}  \Big|^2\Big).
\end{align*}
The definitions of ${\hat{\mathfrak u}}$ in \eqref{hbar} and ${\bf P}$ in (\ref{mathbfP}) will imply that
\begin{equation*}
\big\| |({\bf P}{\hat{\mathfrak u}})_n|^2\big\|_{L^2(0, 1)}\leq \frac{C}{|\ln\varepsilon|},
\qquad
 \big\| |({\bf P}{\check{\mathfrak u}}^0)_n|^2\big\|_{L^2(0, 1)}\leq \frac{C}{|\ln\varepsilon|},
 \quad n=1, \cdots, N.
\end{equation*}
Gathering the above estimates, we get the following estimate
\begin{equation*}
\big\|{\tilde{\mathfrak N}}({\check{\mathfrak u}}^0+{\bf{\tilde w}}+{\hat{\mathfrak u}})\big\|_{L^2(0, 1)}\leq \frac{C}{|\ln\varepsilon|}.
\end{equation*}
Therefore, for any $n=1, \cdots, N-1$, 	using (\ref{evaluefrakg}) and Lemma \ref{lemma6point4}, we can get a solution ${\check{\mathfrak u}}_n$ to
\begin{align*}
\kappa\Big[\mathcal{H}_1\, {\check{\mathfrak u}}_n''
\, +\, \tau_1\, {\check{\mathfrak u}}_n'
\, +\, \tau_2\, {\check{\mathfrak u}}_n\Big]
\, +\, &\lambda_n\, \big(\tau_2\, +\, \tilde{\sigma}\big){\check{\mathfrak u}}_n
\, =\, \kappa\varepsilon^{\mu}\tilde{\mathfrak h}_n\, +\,
{\tilde{\mathfrak N}}_n({\check{\mathfrak u}}^0+{\bf{\tilde w}}+{\hat{\mathfrak u}}),
\\[2mm]
{\check{\mathfrak u}}_n'(0)\, =\, &0,
\qquad \qquad
{\check{\mathfrak u}}_n'(1)\, =\, 0,
\end{align*}
with the following estimate
\begin{align*}
 &\frac{1}{|\ln\varepsilon_l|}\, \|{\check{\mathfrak u}}_n''\|_{L^2(0, 1)}
+\frac{1}{\sqrt{|\ln\varepsilon_l|}}\, \|{\check{\mathfrak u}}_n'\|_{L^2(0, 1)}
+\|{\check{\mathfrak u}}_n\|_{L^{2}(0, 1)}
\\[2mm]
&\leq C\, \sqrt{|\ln\varepsilon_l|}\,
\bigg\|\kappa\varepsilon^{\mu}\tilde{\mathfrak h}_n \, +\,
{\tilde{\mathfrak N}}_n({\check{\mathfrak u}}^0+{\bf{\tilde w}}+{\hat{\mathfrak u}})\bigg\|_{L^2(0, 1)}
\\
&\, \leq\, \varepsilon^{\mu}\|\tilde{\mathfrak h}_n\|_{L^2(0, 1)}
\, +\, \frac{C}{\sqrt{|\ln\varepsilon_l|}} \, \leq \, \frac{C}{\sqrt{|\ln\varepsilon_l|}}.
\end{align*}
Concerning the $N$-th equation in \eqref{equationofddotu}-(\ref{boundaryconditionofddotu}), i.e.,
\begin{equation*}
\mathcal{H}_1\, {\check{\mathfrak u}}_N''
\, +\, \tau_1\, {\check{\mathfrak u}}_N'	
\, +\, \tau_2\, {\check{\mathfrak u}}_N
\, =\, \varepsilon^{\mu}\tilde{\mathfrak h}_N,
\qquad
{\check{\mathfrak u}}_N'(0)\, =\, 0,
\qquad
{\check{\mathfrak u}}_N'(1)\, =\, 0,
\end{equation*}
using (\ref{evaluefrakg}) and Lemma \ref{lemma6point4}, we can also find a solution satisfying
\begin{equation*}
\begin{split}
\|{\check{\mathfrak u}}_N\|_{H^2(0, 1)}
\leq \, \varepsilon^{\mu}\|\tilde{\mathfrak h}_N
\|_{L^2(0, 1)}
\leq \frac{C}{\sqrt{|\ln\varepsilon_l|}}.
\end{split}
\end{equation*}
Now,  the result follows by a straightforward application of Contraction Mapping Principle and Lemma \ref{lemma6point4}.
The proof of Proposition \ref{proposition7point2} is complete.
\qed

\bigskip
\medskip
\medskip
\noindent {\bf Acknowledgements:}
S. Wei was supported by NSFC (No. 12001203) and Guangdong Basic and Applied Basic Research Foundation (No. 2020A1515110622);
J. Yang was supported by NSFC (No. 11771167 and No. 11831009).
\qed

\begin{appendices}

\section{The derivation of the equation for $h$}\label{appendixA}
\setcounter{equation}{0}

The computations of (\ref{solvable of varphi1}) can be showed as follows.
Since $S_{6, j}$, $S_{8, j}$, $M_{11, j}(x, z)$, $M_{21, j}(x, z)$, $M_{51, j}(x, z)$ are even functions of $x_j$, then integration against $ w_{j, x}$ therefore just vanish.
This gives that
\begin{equation*}
\begin{split}
\text{LHS of \eqref{solvable of varphi1}}&=\int_{{\mathbb R}}\big[\varepsilon^2S_{7, j}+\varepsilon^2S_{9, j}+M_{31, j}(x, z)+M_{41, j}(x, z)+M_{61, j}(x, z)\big]w_{j, x}{\rm d}x
\\[2mm]
&\equiv J_1\, +\, J_2\, +\, J_3\, +\, J_4\, +\, J_5.
\end{split}
\end{equation*}
These terms can be computed in the sequel.

\noindent $\bullet$
Recalling the expression of $S_{7, j}$ in (\ref{sv1-gather}), direct computation leads to
\begin{align}
\label{J1}
J_{1}\, =\, &-\varepsilon^2\Big( -\frac {h_5}{\beta}{h}
+h_2\frac {h''}{\beta}
+h_2\frac {2\beta'}{\beta^2}{h'}
+h_2\frac {2\alpha'}{\alpha \beta}{h'} \Big)\int_{{\mathbb R}}w_{j, x}^{2}\, {\rm d}x
\nonumber\\[2mm]
&-\varepsilon^2h_2\frac {2\beta'}{\beta^2}{h'}\int_{{\mathbb R}}x_j w_{j, x}w_{j, xx}\, {\rm d}x
+\varepsilon^2\frac{2h_7}{\beta}h\int_{{\mathbb R}}x_j w_{j, x}w_{j, xx}\, {\rm d}x
-\varepsilon^2\frac { V_{tt}(0, \varepsilon z)}{\beta^3}{h}\int_{{\mathbb R}}x_j w_{j}w_{j, x}\, {\rm d}x
\nonumber\\[2mm]
\, =\, &\varepsilon^2\, \frac{\varrho_1}{\beta}\Big[\, -h_2\, h''
-h_2(\frac{\beta'}{\beta}+2\frac{\alpha'}{\alpha}){h'}+ (h_5-h_7+\sigma V_{tt}(0, \varepsilon z)\beta^{-2})h\, \Big],
\end{align}
where we have  used  the relations \eqref{relationofwx} and
\begin{equation}
\label{varrho1}
 \varrho_{1}\, \equiv\, \int_{{\mathbb R}}w_x^2\, {\rm d}x
 \, =\, -2\int_{{\mathbb R}}xw_xw_{xx}\, {\rm d}x.
\end{equation}

\noindent $\bullet$
According to the definition of $S_{9, j}$ in (\ref{sv1-gather}) and \eqref{varrho1}, it follows that
\begin{align}
\label{J2}
J_{2}\, =\, &
\frac{\varepsilon^2}{\alpha \beta^2}\int_{{\mathbb R}}\big[h_6\Big(-\alpha \beta h' +\alpha \beta'h \Big)x_jw_{j, xx}
+ h_6\Big( \alpha'\beta h+\alpha\beta'h\Big)w_{j, x}
- h_4\alpha\beta h' w_{j, x} \big]w_{j, x}\, {\rm d}x
\nonumber\\[2mm]
\, =\, &\varepsilon^2\, \frac{\varrho_1}{\beta}\Big[
h_6\Big(\frac{1}{2} \frac{\beta'}{\beta}+\frac{\alpha'}{\alpha}\Big)h
+\frac{1}{2} h_6 h'-h_4 h' \Big].
\end{align}

\noindent $\bullet$
By the definition of $M_{31, j}(x, z)$ in \eqref{M31},  the facts $\omega_{0, j}$, $\omega_{1, j}$, $w_{j, x}$ are odd functions of $x_j$ and $\omega_{2, j}, \omega_{3, j}$ are even functions of $x_j$, we obtain
\begin{align}\label{J3}
J_3
\, =\, &\varepsilon^2 h\int_{{\mathbb R}} \Big\{\frac{h_3}{\beta}  \big(a_{12}\omega_{2, j, x} +a_{13}\omega_{3, j, x}\big)
+ \frac{h_8}{\beta} \big( a_{12}x_j\omega_{2, j, xx} +a_{13}x_j\omega_{3, j, xx}\big)
\nonumber\\[2mm]
&\qquad \qquad -\, \frac{V_t(0, \varepsilon z)}{\beta^3}\big(a_{12}x_j\omega_{2, j} +a_{13}x_j\omega_{3, j}\big)
\nonumber\\[2mm]
&\qquad \qquad+h_8 \big( a_{10}\omega_{0, j, xx} +a_{11}\omega_{1, j, xx}\big)
-\, \frac{V_t(0, \varepsilon z)}{\beta^2}\big( a_{10}\omega_{0, j} +a_{11}\omega_{1, j}\big)\Big\}w_{j, x}\, {\rm d}x
\nonumber\\
\, =\, &\varepsilon^2 h\int_{{\mathbb R}} \Big\{h_1a_{10} \big(a_{12}\omega_{2, j, x} +a_{13}\omega_{3, j, x}\big)
+ h_1a_{11}\big( a_{12}x_j\omega_{2, j, xx} +a_{13}x_j\omega_{3, j, xx}\big)
\nonumber\\
&\qquad \qquad \qquad+h_1( a_{10}-\frac{1}{2}a_{11})\big(a_{12}\sigma^{-1}x_j\omega_{2, j} +a_{13}\sigma^{-1}x_j\omega_{3, j}\big)\Big\}w_{j, x}\, {\rm d}x
\nonumber\\
&+\varepsilon^2 h\int_{{\mathbb R}}\Big\{ h_1a_{13}\big( a_{10}\omega_{0, j, xx} +a_{11}\omega_{1, j, xx}\big)
+\, h_1a_{12}\big( a_{10}\omega_{0, j} +a_{11}\omega_{1, j}\big)\Big\}w_{j, x}\, {\rm d}x.
\end{align}
Here, we have used the  \eqref{relation-1} and the definitions of $a_{10}, a_{11}, a_{12}, a_{13}$ as in \eqref{a10a11}-\eqref{a12a13}.

\noindent $\bullet$
By the definition of $M_{41, j}(x, z)$ in \eqref{M41}, we get
\begin{align}\label{J4}
J_4
\, =\, &\varepsilon^2 \frac{\xi(\varepsilon z)}{\beta} h\int_{{\mathbb R}}
 \big( 2h_2 \frac{\beta'}{\beta}+ h_6 \big) \, \Big[ A'\big({\tilde{\bf d}}(\varepsilon z)\big)\varepsilon \mathcal{Q} Z_{j, x}+\phi_{22, j, xz}\Big]\, w_{j, x}\, {\rm d}x
\nonumber\\
&+\varepsilon^2 \frac{\xi(\varepsilon z)}{\beta} \int_{{\mathbb R}}
\Big\{ h_3\Big[ A\big({\tilde{\bf d}}(\varepsilon z)\big)Z_{j, x}+\phi_{22, j, x}\Big]
+\frac{h_8}{\beta}\Big(\frac{x_j}{\beta}\Big)\, \Big[ A\big({\tilde{\bf d}}(\varepsilon z)\big)Z_{j, xx}+\phi_{22, j, xx}\Big]
\nonumber\\[2mm]
&\qquad \qquad \qquad \qquad -\, \frac{ V_t(0, \varepsilon z)}{\beta}\, \Big(\frac{x_j}{\beta}\Big)\, \Big[ A\big({\tilde{\bf d}}(\varepsilon z)\big)Z_{j}+\phi_{22, j}\Big]\Big\}w_{j, x}\, {\rm d}x
\nonumber\\
\, =\, &\, \varepsilon^2\, \frac{\varrho_1}{\beta}\, \alpha_1(z)\, h
\, +\, \varepsilon^2\, \frac{\varrho_1}{\beta}\, G_{1}(z),
\end{align}
 where
\begin{align}
\label{alpha1}
\alpha_1(z)\, =\, 2 \frac{\xi(\varepsilon z)}{\varrho_1}
\, \int_{{\mathbb R}}\Big[\varepsilon\, A'\big({\tilde{\bf d}}(\varepsilon z)\big)\, \beta (\varepsilon z)\, Z_{j, x}+\phi_{22, j, xz}(x, z)\Big] w_{j, x}\, {\rm d}x,
\end{align}
\begin{align}
\label{G1}
 G_{1}(z)\, =\, &\varepsilon^2 \frac{\xi(\varepsilon z)}{\varrho_1}\int_{{\mathbb R}}
\Big\{ h_3\Big[ A\big({\tilde{\bf d}}(\varepsilon z)\big)Z_{j, x}+\phi_{22, j, x}\Big]
+\frac{h_8}{\beta}\Big(\frac{x_j}{\beta}\Big)\, \Big[ A\big({\tilde{\bf d}}(\varepsilon z)\big)Z_{j, xx}+\phi_{22, j, xx}\Big]
\nonumber\\[2mm]
&\qquad \qquad \qquad \qquad -\, \frac{ V_t(0, \varepsilon z)}{\beta}\, \Big(\frac{x_j}{\beta}\Big)\, \Big[ A\big({\tilde{\bf d}}(\varepsilon z)\big)Z_{j}+\phi_{22, j}\Big]\Big\}w_{j, x}\, {\rm d}x.
\end{align}

\noindent $\bullet$
The definition of $M_{61, j}$ is given in \eqref{M61}.
Since $\omega_{0, j}$, $\omega_{1, j}$ and $w_{j, x}$ are odd functions of $x_j$,
while $\omega_{2, j}$, $\omega_{3, j}$, $\phi_{21, j}$ and $\phi_{22, j}$ are even functions of $x_j$,
we then obtain that
\begin{align}\label{J5}
J_5
\, =\, &\varepsilon^2h_1\int_{{\mathbb R}}\frac {p(p-1)}{2}{w_j}^{p-2}\Big[2h\big(a_{10}a_{12}\omega_{0, j}\omega_{2, j}
+a_{10}a_{13}\omega_{0, j}\omega_{3, j}
+a_{11}a_{12}\omega_{1, j}\omega_{2, j}
+a_{11}a_{13}\omega_{1, j}\omega_{3, j}\big)
\nonumber\\[2mm]
&\qquad\qquad\qquad\qquad\quad\quad+2\big(a_{10}\omega_{0, j}+a_{11}\omega_{1, j}\big)\xi(\varepsilon z)\big(A\big({\tilde{\bf d}}(\varepsilon z)\big)Z_j+\phi_{22, j}\big)\Big]w_{j, x}\, {\rm d}x
\nonumber\\
\, \equiv\, &\varepsilon^2h_1h\int_{{\mathbb R}} \, p(p-1)\, {w_j}^{p-2} \Big[a_{10}a_{12}\omega_{0, j}\omega_{2, j}
+a_{10}a_{13}\omega_{0, j}\omega_{3, j}
\nonumber\\
&\qquad \qquad \qquad \qquad \qquad \quad+a_{11}a_{12}\omega_{1, j}\omega_{2, j}
+a_{11}a_{13}\omega_{1, j}\omega_{3, j}\Big]w_{j, x}\, {\rm d}x
+\varepsilon^2 \frac{\varrho_1}{\beta}G_{2}(z),
\end{align}
where
\begin{align}
\label{G3}
 G_{2}(z)
 \, =\, \frac{\beta}{\varrho_1}\, \xi(\varepsilon z) h_1
 \, p(p-1)\int_{{\mathbb R}}{w_j}^{p-2}\, \big(a_{10}\omega_{0, j}+a_{11}\omega_{1, j}  \big)\Big[A\big({\tilde{\bf d}}(\varepsilon z)\big)Z_j+\, \phi_{22, j}(x, z)\Big]\, w_{j, x}\, {\rm d}x.
\end{align}

\medskip
By  differentiating the equation (\ref{w{0}}) and using equations (\ref{w{2}}), \eqref{w{3}}, we obtain
\begin{equation}\label{relation1}
\int_{{\mathbb R}}p(p-1){w_j}^{p-2}w_{j, x} \omega_{0, j}\omega_{2, j}{\rm d}x
=-\int_{{\mathbb R}} w_{j, x}\omega_{0, j}{\rm d}x
+\int_{{\mathbb R}}\Big[w_{j, x}+\frac{1}{\sigma}x_j w_{j}\Big]\omega_{2, j, x}{\rm d}x, 	
\end{equation}
\begin{equation}\label{relation2}
 \int_{{\mathbb R}}p(p-1){w_j}^{p-2}w_{j, x} \omega_{0, j}\omega_{3, j}{\rm d}x
=-\int_{{\mathbb R}} w_{j, xxx}\omega_{0, j}{\rm d}x
+\int_{{\mathbb R}}\Big[w_{j, x}+\frac{1}{\sigma}x_j w_{j}\Big]\omega_{3, j, x}{\rm d}x.
\end{equation}
Similarly,  by  differentiating the equation (\ref{w{1}}) and using equations (\ref{w{2}}), \eqref{w{3}}, we obtain
\begin{equation}\label{relation3}
 \int_{{\mathbb R}}p(p-1){w_j}^{p-2}w_{j, x} \omega_{1, j}\omega_{2, j}{\rm d}x
=-\int_{{\mathbb R}} w_{j, x}\omega_{1, j}{\rm d}x
+\int_{{\mathbb R}}\Big[- \frac{1}{2 \sigma }x_j w_{j}+x_j w_{j, xx}\Big]\omega_{2, j, x}{\rm d}x,
\end{equation}
\begin{equation}\label{relation4}
 \int_{{\mathbb R}}p(p-1){w_j}^{p-2}w_{j, x} \omega_{1, j}\omega_{3, j}{\rm d}x
=-\int_{{\mathbb R}} w_{j, xxx}\omega_{1, j}{\rm d}x
+\int_{{\mathbb R}}\Big[- \frac{1}{2 \sigma }x_j w_{j}+x_j w_{j, xx}\Big]\omega_{3, j, x}{\rm d}x.
\end{equation}
%
%
%
Adding (\ref{J3}), (\ref{J5}) and using (\ref{relation1})-\eqref{relation4}, we have
\begin{align}
\label{J3+J5}
J_3+J_5\, =\, &
\varepsilon^2 h_1h\Bigg\{a_{10}a_{12}2\int_{{\mathbb R}} w_{j, x}\omega_{2, j, x}{\rm d}x
-a_{10}a_{12}\frac{1}{\sigma}\int_{{\mathbb R}} w_{j}\omega_{2, j}{\rm d}x
+a_{11}a_{12}\frac{1}{2 \sigma }\int_{{\mathbb R}}w_{j}\omega_{2, j}{\rm d}x
\nonumber\\[2mm]
&\qquad \quad
+a_{10}a_{13}2\int_{{\mathbb R}}w_{j, x}\omega_{3, j, x}{\rm d}x
-a_{10}a_{13}\frac{1}{\sigma}\int_{{\mathbb R}}w_{j}\omega_{3, j}{\rm d}x
+a_{11}a_{13}\frac{1}{2 \sigma }\int_{{\mathbb R}}w_{j}\omega_{3, j}{\rm d}x
\nonumber\\[2mm]
&\qquad \quad
-a_{11}a_{12}\int_{{\mathbb R}}w_{j, x}\omega_{2, j, x}{\rm d}x
-a_{11}a_{13}\int_{{\mathbb R}}w_{j, x}\omega_{3, j, x}{\rm d}x\Bigg\}
\nonumber\\[2mm]
=&\, \varepsilon^2 h_1h\Big\{-a_{10}a_{12}\Big(\frac {2}{p-1}+\frac 12\Big)
-a_{10}a_{12}\Big(\frac 12-\frac {2}{p-1}\Big)
+a_{11}a_{12}\frac{1}{2  }\Big(\frac 12-\frac {2}{p-1}\Big)
\nonumber\\[2mm]
&\qquad \quad
- a_{10}a_{13}
+\frac{1}{2}a_{11}a_{13}
+a_{11}a_{12}\frac{1}{2}\Big(\frac {2}{p-1}+\frac 12\Big)
 \Big\}\int_{{\mathbb R}}w_{x}^2{\rm d}x
+\varepsilon^2 \frac{\varrho_1}{\beta}G_{2}(z)
\nonumber\\[2mm]
\, =&\, \varepsilon^2 \, \rho_1\, h\, h_1 \Big\{-a_{10}a_{12}
+\frac{1}{2}a_{11}a_{12}
-a_{10}a_{13}
+\frac{1}{2}a_{11}a_{13}
 \Big\}
+\varepsilon^2 \frac{\varrho_1}{\beta}G_{2}(z)
\nonumber\\[2mm]
\, =\, &\varepsilon^2 \, \rho_1\, h\Big[- \sigma  \frac{V_t(0, \theta)}{\beta^3} \Big(\frac{ V_t(0, \theta)}{ V(0, \theta)} -\frac{h_8(\theta)}{h_1(\theta)}  \Big)\Big]
+\varepsilon^2 \frac{\varrho_1}{\beta}G_{2}(z),
\end{align}
where we have used (\ref{stationary}) and the following integral identities
\begin{align}
2\int_{{\mathbb R}}\omega_{2, j, x}w_{j, x}\, {\rm d}x
\, =\, -\sigma\int_{{\mathbb R}}(w_{j, x})^2\, {\rm d}x
\, =\, -\sigma\int_{{\mathbb R}}w_{x}^2{\rm d}x,
\end{align}
\begin{align}
\sigma^{-1}\int_{{\mathbb R}}\omega_{2, j}w_j\, {\rm d}x
\, =\, \Big(\frac 12-\frac {2}{p-1}\Big)\int_{{\mathbb R}}(w_{j, x})^2\, {\rm d}x
\, =\, \Big(\frac 12-\frac {2}{p-1}\Big)\int_{{\mathbb R}}w_{x}^2{\rm d}x,
\end{align}
\begin{equation}
2 \int_{{\mathbb R}}\omega_{3, j, x} w_{j, x}{\rm d}x\, =\, -\frac{1}{2} \int_{{\mathbb R}} w_{j, x}^2{\rm d}x\, =\, -\frac{1}{2} \int_{{\mathbb R}} w_{x}^2{\rm d}x,
\end{equation}
\begin{align}
\int_{{\mathbb R}}\omega_{3, j} w_{j}{\rm d}x\, =\, \Big( \frac{1}{p-1}+\frac{1}{4} \Big)\int_{{\mathbb R}} w_{j, x}^2{\rm d}x\, =\, \frac{\sigma}{2}\int_{{\mathbb R}} w_{x}^2{\rm d}x.
\end{align}

\medskip
Finally, denote
\begin{align}\label{zeta}
\zeta(\theta)=\frac{1}{\alpha^2 \beta \sqrt{{\mathfrak h}_1}},
\end{align}
\begin{align}
\label{hbar1}
\hbar_1(\theta)\, =\, h_2\Big[\, \frac {\beta'}{\beta}
+\frac {2\alpha'}{\alpha}\, \Big]
+h_4 -\frac{1}{2} h_6 \, =\, \zeta(\theta)\mathcal{H}_{1}'(\theta),
\end{align}
and
\begin{align}
\label{hbar2}
\hbar_2(\theta)\, =\, &-\Big[\, h_5
-h_7+\sigma\frac { V_{tt}(0, \theta)}{\beta^2}\, \Big]
-h_6\Big(\frac{1}{2} \frac{\beta'}{\beta}+\frac{\alpha'}{\alpha}\Big)
+\sigma  \frac{V_t(0, \theta)}{\beta^2} \Big(\frac{ V_t(0, \theta)}{ V(0, \theta)} -\frac{h_8}{h_1}  \Big)
\nonumber\\[2mm]
\, =\, &\zeta(\theta) \big[\, \mathcal{H}_{2}'(\theta)-\mathcal{H}_{3} (\theta)\, \big],
\end{align}
where the last equalities in \eqref{hbar1} and \eqref{hbar2} will be verified in Appendix \ref{appendixD}.
We infer that equation (\ref{solvable of varphi1}) becomes
\begin{align*}
\text{LHS of \eqref{solvable of varphi1}}
=&-\varepsilon^2 \frac{\varrho_1}{\beta}\Big[h_2h''
+\hbar_1(\varepsilon z)h'+(\hbar_2(\varepsilon z)+\alpha_1(z))h-G_{1}(z)-G_{2}(z)\Big]
\\
=& -\varepsilon^2 \frac{\varrho_1}{\beta}\zeta(\varepsilon z)
\Bigg\{\mathcal{H}_{1}(\varepsilon z)h''
+\mathcal{H}_{1}'(\varepsilon z)h'
+\Big[\big( \mathcal{H}_{2}'(\varepsilon z)-\mathcal{H}_{3}(\varepsilon z) \big)
+\frac{\alpha_1(z)}{\zeta(\varepsilon z)}\Big]h
-\frac{ G_{1}(z)+G_{2}(z)}{\zeta(\varepsilon z)}\Bigg\}
\\
=\, &0.
\end{align*}

\section{The first projection of error}\label{appendixB}
\setcounter{equation}{0}

We do estimates for the term $ \int_{\mathfrak S}{\mathcal E}w_{n, x}\, {\rm d}x$ given in Section \ref{section5},
where ${\mathcal E}$ is defined in (\ref{new error-2}) and $w_{n, x}$ is an odd function of $x_n$. Integration against all even terms of $x_n$, say ${\mathcal E}_{11, n}$ and $S_{4, n}, M_{12, n}, M_{22, n}$  in ${\mathcal E}_{12}$, therefore just vanish. We have
\begin{equation}
\int_{\mathfrak S}\, {\mathcal E}\, w_{n, x}\, {\rm d}x \, =\Big\{\int_{\mathfrak S_n}+\int_{\mathfrak S\setminus\mathfrak S_n}\Big\}\, {\mathcal E}\, w_{n, x}\, {\rm d}x.
\end{equation}
We begin with
\begin{align}
\int_{{\mathfrak S}_n}{\mathcal E}_{12}w_{n, x}{\rm d}x
\, =\, &\,
\sum_{j=1}^N\int_{{\mathfrak S}_n}\varepsilon^2S_{3, j}w_{n, x}{\rm d}x
 +\sum_{j=1}^N\int_{{\mathfrak S}_n}\varepsilon^2S_{5, j}w_{n, x}{\rm d}x
 +\int_{{\mathfrak S}_n}B_4(v_1)w_{n, x}{\rm d}x
\nonumber
\\
&
+\int_{{\mathfrak S}_n}\varepsilon\frac{h_2}{\beta^2}\phi_{1, zz}w_{n, x}{\rm d}x
+\int_{{\mathfrak S}_n}\varepsilon^2\, \frac{h_2}{\beta^2}\phi_{4, zz}\, w_{n, x}\, {\rm d}x
+\int_{{\mathfrak S}_n}M_{32}(x, z)\, w_{n, x}\, {\rm d}x
\nonumber
\\
&+\int_{{\mathfrak S}_n}M_{42}(x, z)\, w_{n, x}\, {\rm d}x
+\int_{{\mathfrak S}_n}M_{52}(x, z)\, w_{n, x}\, {\rm d}x
+\int_{{\mathfrak S}_n}M_{62}(x, z)\, w_{n, x}\, {\rm d}x
\nonumber
\\
&+\int_{{\mathfrak S}_n}M_{63}(x, z)\, w_{n, x}\, {\rm d}x
+\int_{{\mathfrak S}_n}\big[B_3(\varepsilon^2\, \phi_3)
+B_3(\varepsilon^2\, \phi_4)\big]\, w_{n, x}\, {\rm d}x\, +\, O(\varepsilon^3)
\nonumber
\\
\, \equiv\, &\, \textrm{I}_1
\, +\, \textrm{I}_2
\, +\, \textrm{I}_3
\, +\, \textrm{I}_4
\, +\, \textrm{I}_5
\, +\, \textrm{I}_6
\, +\, \textrm{I}_7
\, +\, \textrm{I}_8
\, +\, \textrm{I}_9
\, +\, \textrm{I}_{10}
\, +\, \textrm{I}_{11}\, +\, O(\varepsilon^3).
\label{intE1wx}
\end{align}
These terms will be estimated as follows.

\noindent $\bullet$
By repeating the same computation used in (\ref{J1}) and (\ref{J2}), we get
\begin{equation}
\label{I1}
\textrm{I}_{1}
\, =\, \varepsilon^2\, \frac{\varrho_1}{\beta}\left\{-h_2\, f_{n}''
-h_2\Big[\, \frac {\beta'}{\beta} +\frac {2\alpha'}{\alpha}\, \Big]{f_{n}'}
+\Big[\, h_5-h_7+\sigma\frac { V_{tt}(0, \varepsilon z)}{\beta^2}\, \Big]{f_n} \right\}+O(\varepsilon^3)\sum_{j=1}^N\, (f_j+f'_j+f''_j),
\end{equation}
where  $\varrho_{1}$ is a positive constant defined in \eqref{varrho1}.

\noindent $\bullet$
There also holds
\begin{equation}\label{I2}
\textrm{I}_{2}
 \, =\, \varepsilon^2\frac{\varrho_1}{\beta}\, \Big[
h_6\Big(\frac{1}{2} \frac{\beta'}{\beta}+\frac{\alpha'}{\alpha}\Big)f_n
+\frac{1}{2} h_6f_n'-h_4 f_n' \Big]
 +O(\varepsilon^3)\sum_{j=1}^N\, (f_j+f'_j).
\end{equation}

\noindent $\bullet$
Recall the expression of $B_4(v_1)$ in (\ref{B2v1}), then
\begin{align}\label{I3}
\textrm{I}_{3}
\, =\, &
\int_{{\mathfrak S}_n}\Big\{\, \frac{1}{\alpha\beta^2}\big[\hat{B}_0(v_1)+\, a_6(\varepsilon s, \varepsilon z)\, \varepsilon^3 \, s^3\, {v_1}\big]
\, +\, h_1p(w_n)^{p-1}(v_1-w_n)
\nonumber\\[2mm]
\, &\qquad
\,-\, h_1\sum_{j\neq n}{w_j}^p\, +\, h_1\frac{1}{2}p(p-1)(w_n)^{p-2}(v_1-w_n)^2\, +\, \max_{j\neq n}O(e^{-3|\beta f_j-x|})\, \Big\}\, w_{n, x}\, {\rm d}x
\nonumber\\[2mm]
=\, &\varrho_2h_1\Big[e^{-\beta(f_n-f_{n-1})}-e^{-\beta(f_{n+1}-f_n)}\Big]+\varepsilon^{\mu_1}\max_{j\neq n}O(e^{-\beta|f_j-f_n|})+\varepsilon^3 \sum_{j=1}^N\big({\mathbf b}_{1\varepsilon j}\, f''_j+{\mathbf b}_{2\varepsilon j}\big),
\end{align}
where $\mu_1$ is a small positive constant, and $\varrho_{2}$ is a positive constant given by
\begin{equation}\label{varrho2}
\varrho_2\, =\, p\, C_p
\int_{0}^{\infty}\, w^{p-1}\, w_x(e^{-x}-e^{x})\, {\rm d}x.
\end{equation}

\noindent $\bullet$
Recall the expression of $\varepsilon\, \phi_1$ in \eqref{phi1}, then
\begin{equation}\label{I4}
\textrm{I}_{4}
\, =\, \int_{{\mathfrak S}_n}\varepsilon\, \frac{h_2}{\beta^2}\, \phi_{1, zz}\, w_{n, x}\, {\rm d}x
\, =\, \varepsilon^3\sum_{j=1}^N\big({\mathbf b}_{1\varepsilon j}f''_j+{\mathbf b}_{2\varepsilon j}\big).
\end{equation}

\noindent $\bullet$
It can be derived  that
\begin{equation}
\label{I5}
\textrm{I}_5\, =\, \int_{{\mathfrak S}_n}\varepsilon^2\, \frac{h_2}{\beta^2}\, \phi_{4, zz}\, w_{n, x}\, {\rm d}x\, =\, O(\varepsilon^4).
\end{equation}

\noindent $\bullet$
From the definition of $M_{42}(x, z)$ in \eqref{M42}, we can estimate the term $\textrm{I}_7$ as the following
\begin{align}
\label{I7}
\textrm{I}_7
\, =\, &\varepsilon^2 \xi(\varepsilon z) \frac{1}{\beta}f_n\, \int_{{\mathfrak S}_n}\, \big(2h_2 \frac{\beta'}{\beta}\, + h_6\big)  \Big[ A'\big({\tilde{\bf d}}(\varepsilon z)\big)\varepsilon \mathcal{Q} Z_{n, x}+\phi_{22, n, xz}\Big]\, w_{n, x}\, {\rm d}x
\nonumber\\[2mm]
&+\varepsilon^3\sum_{j=1}^N\big( {\mathbf b}_{1\varepsilon j}f''_j+{\mathbf b}_{2\varepsilon j}\big)
\nonumber\\
\equiv& -\varepsilon^2 \frac{\varrho_1}{\beta}\, \alpha_1(z)f_n  +\varepsilon^3\sum_{j=1}^N \big( {\mathbf b}_{1\varepsilon j}f''_j+{\mathbf b}_{2\varepsilon j} \big) +O(\varepsilon^3)\sum_{j=1}^N\big(f_j+f'_j\big),
\end{align}
where $\alpha_1(z)$ and $\varrho_{1}$ are defined in (\ref{alpha1}) and (\ref{varrho1}).

\noindent $\bullet$
From the definition of $M_{52}(x, z)$ in \eqref{M52}, we need only consider the odd terms and the higher order terms involving $e_j'$ and $e''_j$, so we get
\begin{align}
\label{I8}
 \textrm{I}_8=&\sum_{j=1}^N\frac{\varepsilon^2}{\beta} {e_j}\int_{{\mathfrak S}_n}\Big[ h_3  \, Z_{j, x}+h_8\, x_jZ_{j, xx}-\, \frac{V_t(0, \varepsilon z)}{\beta^2}\, x_jZ_{j}\Big] \, w_{n, x}\, {\rm d}x
\nonumber\\
 &+\sum_{j=1}^N\, \varepsilon^4\, \frac{1}{\alpha\beta^3} h_9\, e''_j(\varepsilon z)\int_{{\mathfrak S}_n} x_j\, w_{n, x}\, Z_j\, {\rm d}x
 +\varepsilon^3\sum_{j=1}^N\, \big({\mathbf b}_{1\varepsilon j}\, e'_j+ {\mathbf b}_{1\varepsilon j}\, f''_j+{\mathbf b}_{2\varepsilon j}\big)
\nonumber\\
\, \equiv\, &\varepsilon^2 \frac{\varrho_1}{\beta}\Big[\hbar_3({\varepsilon}z)\, e_n+\varepsilon^2\hbar_4({\varepsilon}z)\, e''_n\Big]
+\varepsilon^3\sum_{j=1}^N\big({\mathbf b}_{1\varepsilon j}\, e'_j+{\mathbf b}_{1\varepsilon j}\, f''_j+{\mathbf b}_{2\varepsilon j}\big),
\end{align}
where $\hbar_3(\varepsilon z)$ and $\hbar_4(\varepsilon z)$ are defined like the following
\begin{equation}\label{hbar3}
\hbar_3(\varepsilon z)\, =\, -\varrho_{1}
\int_{{\mathbb R}}\Big[h_3\, Z_x+h_8xZ_{xx}- V_t \beta^{-2}\, x\, Z\Big]\, w_x\, {\rm d}x,
\end{equation}
\begin{equation}\label{hbar4}
\hbar_4(\varepsilon z)\, =\, -\frac{1}{\varrho_1\alpha\beta^2}\, h_9
\int_{{\mathbb R}}\, x\, w_x\, Z\, {\rm d}x.
\end{equation}

\noindent $\bullet$
Recalling the definitions of $M_{32}(x, z)$ and $M_{62}(x, z)$ as in \eqref{M32}, \eqref{M62}, we can get
\begin{equation*}
\textrm{I}_6+ \textrm{I}_9
=\, \varepsilon^2 \, \rho_1\, \Big[- \sigma  \frac{V_t(0, \varepsilon z)}{\beta^3} \Big(\frac{ V_t(0, \varepsilon z)}{ V(0, \varepsilon z)} -\frac{h_8}{h_1}  \Big)\Big]f_n\, +\, \varepsilon^3\sum_{j=1}^N\big({\mathbf b}_{1\varepsilon j}\, f''_j+\, {\mathbf b}_{2\varepsilon j}\big).
\end{equation*}

\noindent $\bullet$
According to the fact that the terms in $B_3(\varepsilon^2\, \phi_3)$ and $B_3(\varepsilon^2\, \phi_4)$ are of order $O(\varepsilon^3)$,
it follows that
\begin{equation}
\label{I10}
\textrm{I}_{10}+\textrm{I}_{11}
\, =\, \varepsilon^3\sum_{j=1}^N\big({\mathbf b}_{1\varepsilon j}\, f''_j+\, {\mathbf b}_{2\varepsilon j}\big).
\end{equation}

The above computations lead to  the estimate
\begin{align}
 \int_{{\mathfrak S}_n}{\mathcal E}w_{n, x}{\rm d}x
 \, =\, &-\varepsilon^2\, \frac{\varrho_1}{\beta}\zeta(\varepsilon z)\, \Big\{\mathcal{H}_1(\varepsilon z)\, f''_n+\mathcal{H}_1'(\varepsilon z)\, f'_n+ \Big[\mathcal{H}_2'(\varepsilon z)-\mathcal{H}_3(\varepsilon z)+\frac{\alpha_1(z)}{\zeta(\varepsilon z)}\Big]\, f_n\Big\}
\nonumber\\[2mm]
&+\varepsilon^2 \frac{\varrho_1}{\beta}\big[\hbar_3(\varepsilon z)e_n+\varepsilon^2\hbar_4({\varepsilon}z)e''_n\big]+h_1\, \varrho_2\big[e^{-\beta(f_n-f_{n-1})}-e^{-\beta(f_{n+1}-f_n)}\big]
\nonumber\\[2mm]
&+\varepsilon^{\mu_1}\max_{j\neq n}O(e^{-\beta|f_j-f_n|})+\varepsilon^3\sum_{j=1}^N\, \big({\mathbf b}_{1\varepsilon j}\, e_j'+{\mathbf b}_{1\varepsilon j}^2\, f''_j+{\mathbf b}_{2\varepsilon j}\big)
 \nonumber\\
&+O(\varepsilon^3)\sum_{j=1}^N\big(\, f_j+\, f'_j+\, f''_j+\, f_j^2\big).
\end{align}

On the other hand, to compute $\int_{\mathfrak S\setminus\mathfrak S_n}\, {\mathcal E}\, w_{n, x}\, {\rm d}x$ for fixed $n=1, \cdots, N$,
we notice that for $(x, z)\in {\mathcal S}_{\delta/\varepsilon}\backslash {\mathfrak A}_n$ with
\begin{equation*}
{\mathcal S}_{\delta/\varepsilon}\, =\, \big\{\, -\delta/\varepsilon\, <\, x\, <\, \delta/\varepsilon, \ 0\, <\, z\, <\, 1/\varepsilon\, \big\},
\end{equation*}
there holds
\begin{equation*}
w_{n, x}\, =\, \max_{j\neq n}O(e^{\frac{1}{2}\beta|f_j-f_n|}).
\end{equation*}
Thus we can estimate
\begin{align}
\int_{\mathfrak S\setminus\mathfrak S_n}\, {\mathcal E}\, w_{n, x}\, {\rm d}x=\varepsilon^\frac{1}{2}\max_{j\neq n}O(e^{-\beta|f_j-f_n|})+O(\varepsilon^{\frac{1}{2}})\sum_{i=1}^{11}\textrm{I}_i.
\end{align}

\section{The second projection of error}\label{appendixC}
We estimate  the term $ \int_{\mathfrak S}{\mathcal E}Z_n\, {\rm d}x$ in Section \ref{section5},
where ${\mathcal E}$ and its decomposition are defined in (\ref{new error-2}) and \eqref{E1-d}, and $Z_n$ is an odd function of $x_n$.  We have
$$\int_{{\mathfrak S}}{\mathcal E}\, Z_n\, {\rm d}x
\, =\, \int_{\mathfrak S}{\mathcal E}_{11}\, Z_n\, {\rm d}x +\int_{\mathfrak S}{\mathcal E}_{12}\, Z_n\, {\rm d}x, $$
where
\begin{equation*}
 \int_{\mathfrak S}{\mathcal E}_{11}\, Z_n\, {\rm d}x
 \, =\, \varepsilon^3\, \frac{h_2}{\beta^2}\, e''_n+\varepsilon h_1\lambda_0\, e_n+O(\varepsilon),
\end{equation*}
and
\begin{equation}
\int_{\mathfrak S}\, {\mathcal E}_{12}\, Z_n\, {\rm d}x \, =\Big\{\int_{\mathfrak S_n}+\int_{\mathfrak S\setminus\mathfrak S_n}\Big\}\, {\mathcal E}_{12}
\, Z_n\, {\rm d}x.
\end{equation}
Since $S_{6, j}$, $S_{8, j}$, $M_{11, j}(x, z)$, $M_{21, j}(x, z)$, $M_{51, j}(x, z)$ are even functions of $x_j$, then integration against $ w_{j, x}$ just vanish.
This gives that
\setcounter{equation}{0}
\begin{align}
 \int_{\mathfrak S_n}{\mathcal E}_{12}\, Z_n\, {\rm d}x
\, =\, &
\sum_{j=1}^N\int_{\mathfrak S_n}\varepsilon^2\, S_{4, j}\, Z_n\, {\rm d}x
+\int_{\mathfrak S_n}B_4(v_1)\, Z_n\, {\rm d}x
+\int_{\mathfrak S_n}\varepsilon\, \frac{h_2}{\beta^2}\phi_{1, zz}\, Z_n\, {\rm d}x \nonumber
\\[2mm]
&+\int_{\mathfrak S_n}\varepsilon^2\, \frac{h_2}{\beta^2}\phi_{4, zz}\, Z_n\, {\rm d}x
+\int_{\mathfrak S_n}M_{12}(x, z)\, Z_n\, {\rm d}x
+\int_{\mathfrak S_n}M_{22}(x, z)\, Z_n\, {\rm d}x \nonumber
\\[2mm]
&+\int_{\mathfrak S_n}M_{32}(x, z)\, Z_n\, {\rm d}x
+\int_{\mathfrak S_n}M_{42}(x, z)\, Z_n\, {\rm d}x
+\int_{\mathfrak S_n}M_{52}(x, z)\, Z_n\, {\rm d}x\nonumber
\\[2mm]
&+\int_{\mathfrak S_n}\big(M_{62}(x, z)+M_{63}(x, z)\big)Z_n\, {\rm d}x
+\int_{\mathfrak S_n}B_3(\varepsilon^2\, \phi_3)Z_n\, {\rm d}x
+\int_{\mathfrak S_n}B_3(\varepsilon^2\, \phi_4)Z_n\, {\rm d}x\nonumber
\\[2mm]
\, \equiv\, &
\textrm{II}_1
+\textrm{II}_2
+\textrm{II}_3
+\textrm{II}_4
+\textrm{II}_5
+\textrm{II}_6
+\textrm{II}_7
+\textrm{II}_8
+\textrm{II}_9+\textrm{II}_{10}
+\textrm{II}_{11}
+\textrm{II}_{12}.
\label{E12Z}
\end{align}
Here are the details of computations.

\noindent $\bullet$
According to the expression of $S_{4, j}$ in (\ref{sv1-gather}) and the assumption of $f_j$, it follows that
\begin{align}
\textrm{II}_1
=&\sum_{j=1}^N\varepsilon^2\int_{\mathfrak S_n}\Big[
\Big( h_2f_j'^2\, +2h_2f_j'h'-h_6\, f_j\, f_j'\, -h_6\, f_j'\, h+h_7f_j^2\Big)w_{n, xx}
 -\frac{1}{2}\, \beta^{-2}\, V_{tt}(0, \varepsilon z)\, f_n^2\, w_n\Big]Z_n\, {\rm d}x
\nonumber
\\
\, =\, &\varepsilon^2\, \varrho_{3}\big( h_2f_n'^2\, +2h_2f_n'h'-h_6\, f_n\, f_n'\, -h_6\, f_n'\, h+h_7f_n^2\big)
+O(\varepsilon^3)\sum_{j=1}^N\, \big(f_j^2+\, {f'_j}^2+f'_jf_j \big),
\label{II1}
\end{align}
where
\begin{equation}\label{varrho3}
 \varrho_{3}\, =\, 2\int_{{\mathbb R}}\, w_{xx}\, Z\, {\rm d}x.
\end{equation}

\noindent $\bullet$
Recall the expression of $B_4(v_1)$ in (\ref{B2v1}), then
\begin{align}\label{II2}
\textrm{II}_{2}
\, =\, &\int_{{\mathfrak S}_n}\Big\{\, \frac{1}{\alpha\beta^2}\big[\hat{B}_0(v_1)+\, a_6(\varepsilon s, \varepsilon z)\, \varepsilon^3 \, s^3\, {v_1}\big]
\, +\, h_1p(w_n)^{p-1}(v_1-w_n)
\,-\, h_1\sum_{j\neq n}{w_j}^p
\nonumber\\
\, &\quad\qquad\, +\, h_1\frac{1}{2}p(p-1)(w_n)^{p-2}(v_1-w_n)^2
\, +\, \max_{j\neq n}O(e^{-3|\beta f_j-x|})\, \Big\}\, Z_n\, {\rm d}x
\nonumber\\
=\, &
\varrho_4\, h_1\Big[e^{-\beta(f_n-f_{n-1})}-e^{-\beta(f_{n+1}-f_n)}\Big]+\varepsilon^{{\hat\tau}_1}\max_{j\neq n}O(e^{-\beta|f_j-f_n|})
+\varepsilon^3\sum_{j=1}^N \big({\mathbf b}_{1\varepsilon j}\, f''_j+ {\mathbf b}_{2\varepsilon j} \big),
\end{align}
where ${\hat\tau}_1$ is a small positive constant with $\frac{1}{2}<{\hat\tau}_1<1$ and $\varrho_4$ is positive constant given by
\begin{align}
\varrho_4\, =\, p\, C_p\int_{0}^{\infty}\, w^{p-1}(x)Z(x)[e^{-x}-e^{x}]\, {\rm d}x.
\label{varrho4}
\end{align}

\noindent $\bullet$
According to the definition of $\varepsilon\, \phi_1(x, z)$ as in \eqref{phi1},
we obtain
\begin{align}\label{II3}
 \textrm{II}_{3}
&\, =\, \sum_{j=1}^N\Big[\varepsilon^3 \, \beta^{-2}\, a_{12}(\varepsilon z)\, f''_j\int_{{\mathfrak S}_n} \omega_{2, j}\, Z_n\, {\rm d}x
+\varepsilon^3 \, \beta^{-2}\, a_{13}(\varepsilon z)\, f''_j\int_{{\mathfrak S}_n} \omega_{3, j}\, Z_n\, {\rm d}x
+\varepsilon^3 \, {\mathbf b}_{2\varepsilon j}\Big]
\nonumber\\
&\, =\varepsilon^3 \, \beta^{-2} f''_n\, \int_{{\mathfrak S}_n}\big[\, a_{12}(\varepsilon z) \omega_{2, n}+a_{13}(\varepsilon z)\omega_{3, n}\, \big]\, Z_n\, {\rm d}x
+\varepsilon^3\sum_{j=1}^N \big({\mathbf b}_{1\varepsilon j}\, f''_j+ {\mathbf b}_{2\varepsilon j} \big)
\nonumber\\
&\, =\varepsilon^3\, \rho_1\, (\varepsilon z)+\varepsilon^3\sum_{j=1}^N \big({\mathbf b}_{1\varepsilon j}\, f''_j+ {\mathbf b}_{2\varepsilon j} \big),
\end{align}
where
\begin{equation}\label{rho1}
\rho_1\, (\varepsilon z)= \, \beta^{-2} f''_n\, \int_{{\mathbb R}}\Big[a_{12}(\varepsilon z) \omega_{2, n}+a_{13}(\varepsilon z)\omega_{3, n}\Big]\, Z_n\, {\rm d}x.
\end{equation}

\noindent $\bullet$
It is easy to prove that
\begin{align}
\label{II4}
\textrm{II}_{4}\, =\, \varepsilon^2\beta^{-2} h_2\, \int_{{\mathfrak S}_n}\, \phi_{4, zz}\, Z_n\, {\rm d}x\, =\, O(\varepsilon^4).
\end{align}

\noindent $\bullet$
According to the expression of $M_{12}(x, z)$ and $M_{22}(x, z)$ as in \eqref{tildeL0phi2} and \eqref{tildeL0e^2hatphi}, hence, it is easy to obtain that
\begin{align}
\textrm{II}_{5}
&=\sum_{j=1}^N\frac{\varepsilon^3}{\beta^{2}}\, h_2\, \xi''(\varepsilon z)\int_{{\mathfrak S}_n}\big[A'\big({\tilde{\bf d}}(\varepsilon z)\big)Z_j+\phi_{22, j}\big]\, Z_n\, {\rm d}x, \label{II5}
\nonumber\\[2mm]
&=\varepsilon^3\frac{h_2}{\beta^2}\, \xi''(\varepsilon z)A'\big({\tilde{\bf d}}(\varepsilon z)\big)+O(\varepsilon^3).
\end{align}

\noindent $\bullet$
Similarly, there holds
\begin{equation}\label{II6}
\textrm{II}_{6}=\int_{{\mathfrak S}_n}M_{22}(x, z)\, Z_n\, {\rm d}x=O(\varepsilon^3).
\end{equation}

\noindent $\bullet$
The estimate of $\textrm{II}_7$ can be proved by the same way, i.e.,
\begin{align}\label{II7}
\textrm{II}_7
\, =\, &\varepsilon^2f_n^2\int_{{\mathfrak S}_n}\Big[\, h_8\big(a_{12}\omega_{2, n, xx} +a_{13}\omega_{3, n, xx}\big)
\nonumber\\
&\qquad \qquad -\, \frac{V_t(0, \varepsilon z)}{\beta^2}\, \big(a_{12}\omega_{2, n}+a_{13}\omega_{3, n}\big)\, \Big]Z_n\, {\rm d}x
+\varepsilon^3\sum_{j=1}^N \big({\mathbf b}_{1\varepsilon j}\, f''_j+ {\mathbf b}_{2\varepsilon j} \big)
\nonumber\\
\, =\, &\, \varepsilon^2\, \rho_2\, (\varepsilon z)+\varepsilon^3\sum_{j=1}^N \big({\mathbf b}_{1\varepsilon j}\, f''_j+ {\mathbf b}_{2\varepsilon j} \big),
\end{align}
where
\begin{equation}\label{rho2}
\rho_2\, (\varepsilon z)=f_n^2\int_{{\mathbb R}}\Big[h_8\big(a_{12}\omega_{2, xx} +a_{13}\omega_{3, xx}\big)
-\, \frac{V_t(0, \varepsilon z)}{\beta^2}\, \big(a_{12}\omega_{2}+a_{13}\omega_{3}\big)\Big]Z\, {\rm d}x.
\end{equation}

\noindent $\bullet$
On the other hand, from the definition of $M_{42}(x, z)$ in \eqref{M42}, we can compute that
\begin{align}\label{II8}
\textrm{II}_{8}
&=\, \sum_{j=1}^N \varepsilon^2 \frac{1}{\beta}f_j \xi(\varepsilon z)\, \int_{{\mathfrak S}_n} \Big\{\frac{h_8}{\beta}\, \big[ A\big({\tilde{\bf d}}(\varepsilon z)\big)Z_{j, xx}+\phi_{22, j, xx}\big]
\nonumber\\
&\qquad \qquad \qquad \qquad \qquad \qquad-\, \frac{V_t(0, \varepsilon z)}{\beta}\, \big[ A\big({\tilde{\bf d}}(\varepsilon z)\big)Z_{j}+\phi_{22, j}\big]\Big\}\, Z_n\, {\rm d}x+O(\varepsilon^3)
\nonumber\\[2mm]
&=\, \varepsilon^2\rho_3\, (\varepsilon z)+O(\varepsilon^3),
\end{align}
where
\begin{align}\label{rho3}
\rho_3\, (\varepsilon z)\, =\, &\frac{1}{\beta}f_n \xi(\varepsilon z)\, \int_{\mathbb R} \Big\{\frac{h_8}{\beta}\, \big[ A\big({\tilde{\bf d}}(\varepsilon z)\big)Z_{xx}+\phi_{22, xx}\big]
\nonumber\\
&\qquad \qquad \qquad \quad -\, \frac{V_t(0, \varepsilon z)}{\beta}\, \big[ A\big({\tilde{\bf d}}(\varepsilon z)\big)Z+\phi_{22}\big]\Big\}\, Z\, {\rm d}x.
\end{align}

\noindent $\bullet$
We need only to compute those parts of $M_{52}(x, z)$ in \eqref{M52}, which are even in $x_n$. It is easy to check that
\begin{align}
\begin{split}
\label{II9}
 \textrm{II}_{9}
 &=\sum_{j=1}^N\frac{\varepsilon^2}{\beta} {e_j}\int_{{\mathfrak S}_n}\Big[\, h_3  \, Z_{j, x}
+h_8\, x_jZ_{j, xx}-\, \frac{V_t(0, \varepsilon z)}{\beta^2}\, x_jZ_{j}\, \Big]\, Z_n\, {\rm d}x+O(\varepsilon^3)
=\, O(\varepsilon^3).
\end{split}
\end{align}
Additionally, we also need to consider some higher order terms in $\textrm{II}_{9}$. The ones involving first derivative of $e_j$ are
\begin{align}
\label{M121important}
&\varepsilon^3\sum_{j=1}^N\, e_j'\, \Big(\frac{h_2 \beta'}{\beta^{3}}+\frac{h_6}{\beta^{2}}\, \, \Big)\int_{{\mathfrak S}_n}x_jZ_{j, x}Z_n\, {\rm d}x+\varepsilon^3\sum_{j=1}^N\, e_j'\, \Big(\frac{2\alpha'}{\alpha\beta^{2}}+\frac{h_4}{\beta^{2}}\, \Big)\int_{{\mathfrak S}_n} Z_j Z_n\, {\rm d}x
\nonumber\\
&\, =\, \varepsilon^3\, \Big[\, \Big(\frac{2\alpha'}{\alpha\beta^{2}}+\frac{h_4}{\beta^{2}}\, \Big)- \frac{1}{2}\Big(\frac{h_2 \beta'}{\beta^{3}}+\frac{h_6}{\beta^{2}}\, \Big)\, \Big]e_n'+\varepsilon^3\sum_{j=1}^N \, {\mathbf b}_{2\varepsilon j}e'_j(\varepsilon z)
\nonumber\\
&\, \equiv\, \varepsilon^3\, \hbar_5({\varepsilon}z)\, e_n'+\varepsilon^3\sum_{j=1}^N \, {\mathbf b}_{2\varepsilon j}e'_j(\varepsilon z),
\end{align}
\noindent where
\begin{equation}
\label{hbar5}
\hbar_5(\varepsilon z)\, =\, \Big(\frac{2\alpha'}{\alpha\beta^{2}}+\frac{h_4}{\beta^{2}}\, \Big)- \frac{1}{2}\Big(\frac{h_2 \beta'}{\beta^{3}}+\frac{h_6}{\beta^{2}}\, \, \Big).
\end{equation}
Moreover, the ones involving second derivative of $e_j$ in $\textrm{II}_{9}$ are
\begin{equation}
\varepsilon^4\, f_n\, \beta^{-2}\, \hbar_6(\varepsilon z)e''_n(\varepsilon z)+O(\varepsilon^5)\sum_{j=1}^Ne''_j(\varepsilon z),
\label{hbar6}
\end{equation}
with $O(\varepsilon^2)$ uniform in $\varepsilon$ and $\hbar_6({\varepsilon}z)$ is a smooth functions of their argument.

\noindent $\bullet$
In the terms of $\textrm{II}_{10}$ and $\textrm{II}_{12}$, we need only to consider those parts which are even in $x$.
It is good that the even (in $x$) terms in $\textrm{II}_{10}$ and $\textrm{II}_{12}$ are of  order $O(\varepsilon^3)$.
Moreover, the terms in $B_3(\varepsilon^2\phi_3)$ are of order $O(\varepsilon^3)$. Consequently, we deduce that
\begin{equation}
\begin{split}
  &\textrm{II}_{10}+\textrm{II}_{11}+\textrm{II}_{12}
 =O(\varepsilon^3).
\end{split}
\end{equation}

\medskip
To compute $\int_{\mathfrak S\setminus\mathfrak S_n}\, {\mathcal E}_{12}\, Z_n\, {\rm d}x$, we notice that for $(x, z)\in {\mathcal S}_{\delta/\varepsilon}\backslash {\mathfrak A}_n$,
\begin{equation*}
Z(x_n)=\max_{j\neq n}O(e^{-\frac{p}{2}\beta |f_j-f_n|}),
\end{equation*}
and thus we can estimate
\begin{equation}
\int_{\mathfrak S\setminus\mathfrak S_n}\, {\mathcal E}_{12}\, Z_n\, {\rm d}x=\varepsilon^\frac{1}{2}\max_{j\neq n}O(e^{-\beta|f_j-f_n|})+O(\varepsilon^{\frac{1}{2}})\sum_{i=1}^{12}\textrm{II}_i.
\end{equation}

\section{The computations of \eqref{relation-1}, \eqref{hbar1} and \eqref{hbar2}}\label{appendixD}

We first show the validity of  \eqref{relation-1} under the assumption of stationary condition for $\Gamma$ in {\bf (A3)} of Section \ref{section1}.
In fact,
the stationary assumption means that (c.f. \eqref{stationary1}), i.e.,
\begin{align}\label{relationoff1f0}
\frac{1}{2}{\mathfrak f}_1\, =\, -\sigma\frac{ V_t(0, \theta)}{ V(0, \theta)}{\mathfrak f}_0.
\end{align}
This gives that
\begin{align*}
 \frac{1}{{\mathfrak h}_1}{\mathfrak f}_1\beta
-\frac{1}{2}\frac{{\mathfrak h}_2}{{\mathfrak h}_1^{2}}{\mathfrak f}_0 \beta
\, =\, \frac{1}{2}\frac{1}{{\mathfrak h}_1}{\mathfrak f}_1\beta
-\frac{1}{2}\frac{{\mathfrak h}_2}{{\mathfrak h}_1^2}{\mathfrak f}_0\beta
-\sigma \frac{V_t(0, \theta)}{\beta}.
\end{align*}
According to the expressions of $h_3, h_8$ as in \eqref{h3h4} and \eqref{h8}, we can obtain that
\begin{align*}
h_3\beta\, = \, \frac{1}{2}h_8\beta -\sigma \frac{V_t(0, \theta)}{\beta},
\end{align*}
which is exactly the formula \eqref{relation-1}.

\medskip
Recalling the expression as in \eqref{h1}-\eqref{h2}, we then have
\begin{align}
{\mathcal H}_1(\theta)= V^{\frac{2}{p-1}}(0, \theta)\frac{\sqrt{ V\big(0, \theta\big)}}{\sqrt{{\mathfrak f}_0}}\, {\mathfrak w}_0\, =\, \alpha^2 \beta \sqrt{{\mathfrak h}_1}\, h_2.
\end{align}
This gives that
\begin{align*}
{\mathcal H}_1'(\theta)
 =&\alpha^2 \beta\sqrt{{\mathfrak h}_1}\Big[\, h_2 \frac{\beta'}{\beta}
 +\frac{1}{\sqrt{{\mathfrak h}_1}} \Big( \frac{1}{\sqrt{{\mathfrak h}_1}}{\mathfrak w}_0\Big)'
 +2 \frac{\alpha'}{\alpha} h_2\, \Big].
\end{align*}
Recalling the definitions of $h_2, h_4, h_5$ as in \eqref{h1}-\eqref{h2}, \eqref{h3h4} and \eqref{h5}, the coefficient of $f'$ is
\begin{align*}
h_2\Big[\, \frac {\beta'}{\beta} +\frac {2\alpha'}{\alpha}\, \Big]+h_4
-\frac{1}{2} h_6
\, =\, &h_2\Big[\, \frac {\beta'}{\beta} +\frac {2\alpha'}{\alpha}\, \Big]
+\frac{1}{\sqrt{{\mathfrak h}_1}} \Big( \frac{1}{\sqrt{{\mathfrak h}_1}}{\mathfrak w}_0\Big)'
\\[2mm]
\, =\, &\frac{1}{\alpha^2 \beta \sqrt{{\mathfrak h}_1}} {\mathcal H}_1',
\end{align*}
which is exactly the formula \eqref{hbar1}.

\medskip
According to the definition of $\beta$ in \eqref{alpha-beta} and $\sigma= \frac{p+1}{p-1}-\frac{1}{2}$, then there holds
\begin{align*}
{\mathcal H}_2(\theta)
=\frac{ V^{\sigma}\big(0, \theta\big)}{\sqrt{{\mathfrak f}_0}}{\mathfrak l}_1
=\alpha^2 \beta \frac{1}{\sqrt{{\mathfrak h}_1}}{\mathfrak l}_1.
\end{align*}
This implies that
\begin{align*}
{\mathcal H}_2'(\theta)
=&\alpha^2\Big[\beta'\frac{1}{\sqrt{{\mathfrak h}_1}}{\mathfrak l}_1
\, +\, \beta\frac{1}{\sqrt{{\mathfrak h}_1}}\partial_{\theta} {\mathfrak l}_1
\,-\, \frac{1}{2}\frac{1}{(\sqrt{{\mathfrak h}_1})^3}\partial_{\theta}{\mathfrak h}_1\, \beta{\mathfrak l}_1\Big]
+2\alpha \alpha' \beta \frac{1}{\sqrt{{\mathfrak h}_1}}{\mathfrak l}_1.
\end{align*}
Using the fact \eqref{relationoff1f0} and \eqref{h1}-\eqref{h2}, then
\begin{align*}
{\mathcal H}_2'(\theta)-{\mathcal H}_3(\theta)
\, =\, &\alpha^2\Big[\beta'\frac{1}{\sqrt{{\mathfrak h}_1}}{\mathfrak l}_1
\, +\, \beta\frac{1}{\sqrt{{\mathfrak h}_1}}\partial_{\theta} {\mathfrak l}_1
\,-\, \frac{1}{2}\frac{1}{(\sqrt{{\mathfrak h}_1})^3}\partial_{\theta}{\mathfrak h}_1\, \beta{\mathfrak l}_1\Big]
+2\alpha \alpha' \beta \frac{1}{\sqrt{{\mathfrak h}_1}}{\mathfrak l}_1
\\[2mm]
&-\alpha^2\frac{ V^{\frac{1}{2}}\big(0, \theta\big)}{\sqrt{{\mathfrak f}_0}}{\mathfrak f}_2
-\alpha^2 \sigma \frac{ V_{tt}(0, \theta)\, }{ V^{\frac{1}{2}}\big(0, \theta\big)}\, \sqrt{{\mathfrak f}_0}
-\alpha^2 \sigma(\sigma-1)\frac{\big| V_t(0, \theta)\big|^2\, }{ V^{\frac{3}{2}}(0, \theta)}
\, \sqrt{{\mathfrak f}_0}
\nonumber \\[2mm]
&- \alpha^2 \sigma\frac{ V_t(0, \theta) \, }{ V^{\frac{1}{2}}\big(0, \theta\big)}\, \frac{1}{\sqrt{{\mathfrak f}_0}}{\mathfrak f}_1
\, +\, \alpha^2\frac{1}{4}\frac{ V^{\frac{1}{2}}\big(0, \theta\big)}{\ \big(\sqrt{{\mathfrak f}_0}\, \big)^3\ }{\mathfrak f}_1^2.
\\[2mm]
\, =\, &\alpha^2 \beta \sqrt{{\mathfrak h}_1}
\Bigg\{ \frac{\beta'}{\beta}\frac{1}{{\mathfrak h}_1}{\mathfrak l}_1
\, +\, \frac{1}{{\mathfrak h}_1}\partial_{\theta} {\mathfrak l}_1
\,-\, \frac{1}{2}\frac{1}{{\mathfrak h}_1^{2}}\partial_{\theta}{\mathfrak h}_1{\mathfrak l}_1
\, +\, 2\frac{\alpha'}{\alpha}\frac{1}{{\mathfrak h}_1}{\mathfrak l}_1
\\[2mm]
&\qquad \qquad
\,-\, \frac{1}{{\mathfrak h}_1}{\mathfrak f}_2
\,-\, \sigma \frac{ V_{tt}(0, \theta)}{\beta^{2}}
\,-\, \sigma \frac{|V_t(0, \theta)|^2}{V \beta}
\,-\, \sigma \frac{V_t(0, \theta)}{\beta^2}\frac{{\mathfrak f}_1}{{\mathfrak f}_0} \Bigg\}.
\end{align*}
On the other hand, according to the definitions of $h_1, h_5, h_6, h_7, h_8$ as in \eqref{h1}-\eqref{h2}, \eqref{h5}, \eqref{h6h7}, and \eqref{h8}, the coefficient of $f$ is
\begin{align}\label{defoff}
&-\Big[\, h_5-h_7+\sigma\frac { V_{tt}(0, \theta)}{\beta^2}\, \Big]
-h_6\Big(\frac{1}{2} \frac{\beta'}{\beta}+\frac{\alpha'}{\alpha}\Big)
+\sigma  \frac{V_t(0, \theta)}{\beta^2} \Big(\frac{ V_t(0, \theta)}{ V(0, \theta)} -\frac{h_8(\theta)}{h_2(\theta)}  \Big)
\nonumber\\[2mm]
\, =\, &
-2\frac{1}{{\mathfrak h}_1}{\mathfrak f}_2
+\frac{3}{2}\frac{{\mathfrak h}_2}{{\mathfrak h}_1^{2}}{\mathfrak f}_1
-\frac{{\mathfrak h}_2^2}{{\mathfrak h}_1^{3}}{\mathfrak f}_0
+\frac{{\mathfrak h}_3}{{\mathfrak h}_1^{2}}{\mathfrak f}_0
+\frac{1}{\sqrt{{\mathfrak h}_1}} \Big[\frac{1}{\sqrt{{\mathfrak h}_1}}{\mathfrak l}_1\Big]'
+\Big[\, \frac{1}{{\mathfrak h}_1}{\mathfrak f}_2
+\frac{{\mathfrak h}_2}{{\mathfrak h}_1^2}{\mathfrak f}_1
-\Big(\frac{{\mathfrak h}_2^2}{{\mathfrak h}_1^3}
+\frac{{\mathfrak h}_3}{{\mathfrak h}_1^2}\Big){\mathfrak f}_0\Big]
\nonumber\\[2mm]
&-\sigma \frac{V_{tt}(0, \theta)}{\beta} \frac {1}{\beta}
+2\frac{1}{{\mathfrak h}_1}{\mathfrak l}_1\, \Big(\frac{1}{2} \frac{\beta'}{\beta}+\frac{\alpha'}{\alpha}\Big)
-\sigma  \frac{V_t(0, \theta)}{\beta^2} \Big[\frac{ V_t(0, \theta)}{ V(0, \theta)} -\frac{{\mathfrak f}_1}{{\mathfrak f}_0}
+\frac{{\mathfrak h}_2}{{\mathfrak h}_1} \Big]
\nonumber\\[2mm]
=&-\frac{1}{{\mathfrak h}_1}{\mathfrak f}_2
-\frac{1}{2}\frac{1}{{\mathfrak h}_1^2}\partial_{\theta}{\mathfrak h}_1{\mathfrak l}_1
+\frac{1}{{\mathfrak h}_1}\partial_{\theta}{\mathfrak l}_1
+\frac{\beta'}{\beta}\frac{1}{{\mathfrak h}_1}{\mathfrak l}_1
+2\frac{\alpha'}{\alpha}\frac{1}{{\mathfrak h}_1}{\mathfrak l}_1
-\sigma \frac{ V_{tt}(0, \theta)}{\beta^2}
+\sigma\frac{| V_t(0, \theta)|^2}{ V\beta^2}
-\sigma \frac{ V_t(0, \theta)}{\beta^2}\frac{{\mathfrak f}_1}{{\mathfrak f}_0}.
\end{align}
Therefore, we obtain that
\begin{equation*}
\Big[\, h_5-h_7+\sigma\frac { V_{tt}(0, \theta)}{\beta^2}\, \Big]
+h_6\Big(\frac{1}{2} \frac{\beta'}{\beta}+\frac{\alpha'}{\alpha}\Big)
- \sigma  \frac{V_t(0, \theta)}{\beta^2} \Big[\frac{ V_t(0, \theta)}{ V(0, \theta)} -\frac{h_8}{h_2}  \Big]
\, =\, \frac{{\mathcal H}_2'-{\mathcal H}_3}{\alpha^2 \beta \sqrt{{\mathfrak h}_1}},
\end{equation*}
which is exactly the formula \eqref{hbar2}.

\section{The computations of \eqref{boundaryoffunctional1} and \eqref{boundaryoffunctional2}}\label{appendixE}
Due to the assumptions in \eqref{a1=a2}, we obtain
\begin{equation*}
{\mathfrak a}_1(0, 0)={\mathfrak a}_2(0, 0),
\qquad
{\mathfrak a}_1(0, 1)={\mathfrak a}_2(0, 1),
\end{equation*}
\begin{equation*}
{ \tilde{\mathfrak a}}_1(0)={\tilde{\mathfrak a}}_2(0)=
{\tilde{\mathfrak a}}_1(1)={\tilde{\mathfrak a}}_2(1)=\frac{1}{\sqrt{2}},
\end{equation*}
and
\begin{equation*}
-{\mathfrak a}_2(0, 0)|{\tilde{\mathfrak a}}_1'(0)|^2
+ {\mathfrak a}_1(0, 0)|{\tilde{\mathfrak a}}_2'(0)|^2=0,
\qquad
-{\mathfrak a}_2(0, 1)|{\tilde{\mathfrak a}}_1'(1)|^2
+ {\mathfrak a}_1(0, 1)|{\tilde{\mathfrak a}}_2'(1)|^2=0.
\end{equation*}
Using the above facts and expressions of ${\mathfrak w}_0$, ${\mathfrak l}_1$ as in  \eqref{m11}, \eqref{m12}, we can derive that
\begin{equation}
{\mathfrak w}_0(0)
 \, =\, {\mathfrak a}_1(0, 0){\tilde{\mathfrak a}}_2^2(0) |n_2(0)|^2+{\mathfrak a}_2(0, 0){\tilde{\mathfrak a}}_1^2(0) |n_1(0)|^2
\, =\, \frac{{\mathfrak a}_1(0, 0)}{2},
\label{E24}
\end{equation}
and
\begin{align}
{\mathfrak l}_1(0)
 \, =\, &\big[{\mathfrak a}_2(0, 0)|n_2(0)|^2+{\mathfrak a}_1(0, 0)|n_1(0)|^2\big]\Theta_{tt}(0, 0)
 +\big[{\mathfrak a}_2(0, 0){\tilde{\mathfrak a}}_1(0){\tilde{\mathfrak a}}_1'(0)|n_1(0)|^2
 +{\mathfrak a}_1(0, 0){\tilde{\mathfrak a}}_2(0){\tilde{\mathfrak a}}_2'(0)|n_2(0)|^2\big]
 \nonumber\\[2mm]
 &-k\big[-{\mathfrak a}_2(0, 0)|{\tilde{\mathfrak a}}_1'(0)|^2
 + {\mathfrak a}_1(0, 0)|{\tilde{\mathfrak a}}_2'(0)|^2\big]n_1(0)n_2(0)
 +\big[-\partial_t{\mathfrak a}_2(0, 0){\tilde{\mathfrak a}}_1(0)
 +\partial_t{\mathfrak a}_1(0, 0){\tilde{\mathfrak a}}_2(0)\big]n_1(0)n_2(0)
\nonumber\\[2mm]
\, =\, &{\mathfrak a}_1\big(0, 0\big)\frac{k_1}{2}
\, +\, \frac{{\mathfrak a}_1(0, 0)}{\sqrt{2}}\bigg[{\tilde{\mathfrak a}}_1'(0)|n_1(0)|^2+{\tilde{\mathfrak a}}_2'(0)|n_2(0)|^2\bigg]
\nonumber\\[2mm]
&\, +\, \frac{1}{\sqrt{2}}\bigg[-\partial_t{\mathfrak a}_2(0, 0)
+\partial_t{\mathfrak a}_1(0, 0)\bigg]n_1(0)n_2(0),
\label{E25}
\end{align}
where we have used \eqref{Thetaderivative1}-\eqref{Thetaderivative2} and \eqref{tildek1}.

\medskip
The formulas \eqref{mathfrakh1}, \eqref{g12inverse} and \eqref{g22inverse} imply that
 \begin{equation*}
 {\mathfrak h}_1(0)\, =\, \big[{\tilde{\mathfrak a}}_1(0)|n_1(0)|^2+{\tilde{\mathfrak a}}_2(0)|n_2(0)|^2\big]^2 \, =\, {\tilde{\mathfrak a}}_1^2(0)=\frac{1}{2},
 \end{equation*}
$$
{\mathfrak g}_3(0)=\frac{1}{{\mathfrak h}_1(0)}\big[{\tilde{\mathfrak a}}_1(0)-{\tilde{\mathfrak a}}_2(0)\big]n_1(0)n_2(0)=0,
$$
$$
{\mathfrak g}_6(0)\, =\, \frac{1}{{\mathfrak h}_1(0)}\big[{\tilde{\mathfrak a}}_1^2(0) |n_1(0)|^2 +{\tilde{\mathfrak a}}_2^2(0) |n_2(0)|^2 \big]\, =\, 1,
$$
\begin{align*}
{\mathfrak g}_4(0)=&-\frac{1}{{\mathfrak h}_1(0)}\big[{\tilde{\mathfrak a}}_1(0){\tilde{\mathfrak a}'}_1(0)|n_1(0)|^2
+{\tilde{\mathfrak a}}_2(0){\tilde{\mathfrak a}'}_2(0) |n_2(0)|^2\big]
\, +\,
k\frac{1}{{\mathfrak h}_1(0)}\big[-{\tilde{\mathfrak a}}_1^2(0)+{\tilde{\mathfrak a}}_2^2(0)\big]n_1(0)n_2(0)
\nonumber\\[2mm]
&\qquad-\frac{1}{{\mathfrak h}_1(0)}\big[-q_1(0)n_2(0)+q_2(0)n_1(0)\big]
+{\mathfrak g}_1\big[{\tilde{\mathfrak a}}_1(0)-{\tilde{\mathfrak a}}_2(0)\big]n_1(0)n_2(0)
\nonumber\\[2mm]
\, =\, &-\sqrt{2}\Big({\tilde{\mathfrak a}'}_1(0)|n_1(0)|^2+{\tilde{\mathfrak a}'}_2(0) |n_2(0)|^2\Big) -k_1,
\end{align*}
where we have used \eqref{2.22}, \eqref{Thetaderivative1}-\eqref{Thetaderivative2} and \eqref{tildek1}.
Recalling the definition of ${\mathfrak y}_4$, $\mathfrak{p}_8$ as in  \eqref{sigma2}, \eqref{p8p9p10}, it is easy to obtain

\begin{equation*}
{\mathfrak y}_4(0)\, =\, {\sqrt{{\mathfrak g}_6(0)}}\, =\, 1,
\end{equation*}
\begin{equation*}
 \mathfrak{p}_8(0)\, =\, {\mathfrak a}_2(0, 0){\tilde{\mathfrak a}}_1(0)|n_1(0)|^2+{\mathfrak a}_1(0, 0){\tilde{\mathfrak a}}_2(0)|n_2(0)|^2
 = {\mathfrak a}_1(0, 0) {\tilde{\mathfrak a}}_1(0)=\frac{{\mathfrak a}_1(0, 0)}{\sqrt{2}}.
\end{equation*}
Then by using \eqref{E24}, we get
\begin{align}
 {\mathfrak b}_1=&\frac{1}{\sqrt{{\mathfrak h}_1(0)}}{\mathfrak y}_4(0)\mathfrak{p}_8(0)
 ={\mathfrak a}_1(0, 0) =2{\mathfrak w}_0(0).
\label{b1=mathfrakw0}
 \end{align}
On the other hand, recalling the definition of $\mathfrak{p}_1$, $\mathfrak{p}_4$, as in \eqref{p1p2}, \eqref{p3p4},
it is easy to obtain
\begin{align*}
  \mathfrak{p}_1(0)\, =\, {\mathfrak a}_1(0, 0){\tilde{\mathfrak a}}_1(0)|n_1|^2\, +\, {\mathfrak a}_2(0, 0){\tilde{\mathfrak a}}_2(0)|n_2|^2
  \, =\, {\mathfrak a}_1(0, 0){\tilde{\mathfrak a}}_1(0)
  \, =\, \frac{{\mathfrak a}_1(0, 0)}{\sqrt{2}},
\end{align*}
\begin{align*}
\mathfrak{p}_4(0)\, =\, -\big[\partial_t{\mathfrak a}_1(0, 0)- \partial_t{\mathfrak a}_2(0, 0) \big]n_1(0)n_2(0).
\end{align*}
The term ${\mathfrak y}_2$ given in \eqref{sigma1} will be evaluated as the following
\begin{align*}
{\mathfrak y}_2(0)\, =\, \frac{{\mathfrak g}_4}{{\mathfrak g}_6}-\frac{{\mathfrak g}_7{\mathfrak g}_3}{{\mathfrak g}_6^2}
\, =\, -\sqrt{2}\Big({\tilde{\mathfrak a}'}_1(0)|n_1(0)|^2+{\tilde{\mathfrak a}'}_2(0) |n_2(0)|^2\Big) -k_1.
\end{align*}
Then by using \eqref{E25}, we have
 \begin{align}
{\mathfrak b}_2=&\frac{1}{\sqrt{{\mathfrak h}_1(0)}}{\mathfrak y}_4(0)\mathfrak{p}_4(0)
\, +\,
\frac{1}{\sqrt{{\mathfrak h}_1(0)}}{\mathfrak y}_2(0)\mathfrak{p}_1(0)
\nonumber\\[2mm]
=&-\sqrt{2}\big[\partial_t{\mathfrak a}_1(0, 0)- \partial_t{\mathfrak a}_2(0, 0) \big]n_1(0)n_2(0)
\,-\,
\sqrt{2}\Big({\tilde{\mathfrak a}'}_1(0)|n_1(0)|^2+{\tilde{\mathfrak a}'}_2(0) |n_2(0)|^2\Big){\mathfrak a}_1(0, 0)
\,-\,
k_1 {\mathfrak a}_1(0, 0)
\nonumber\\[2mm]
\, =\, &-2{\mathfrak l}_1(0).
\label{b2=mathfrakw1}
 \end{align}
The formulas \eqref{b1=mathfrakw0} and \eqref{b2=mathfrakw1} will lead to
 \begin{align}
\frac{ V^{\sigma}\big(0, \theta\big)}{\sqrt{{\mathfrak f}_0}}
\Big[{\mathfrak w}_0(0)f'(0)
\, +\,
{\mathfrak l}_1(0)f(0)
\Big]
\, =\,
\frac{ V^{\sigma}\big(0, \theta\big)}{2\sqrt{{\mathfrak f}_0}}
\big[
{\mathfrak b}_1f'(0)
\,-\,
{\mathfrak b}_2f(0)
\big].
 \end{align}
This is the formula \eqref{boundaryoffunctional2}.
Formula \eqref{boundaryoffunctional1} can be verified in the identically same way.

\end{appendices}

\end{document}